\documentclass[8 pt]{amsproc}
\usepackage{amsmath}
\usepackage{graphicx}
\usepackage{amssymb}
\oddsidemargin=-1cm\evensidemargin=-1cm
\begin{document}
\numberwithin{equation}{section}

\def\1#1{\overline{#1}}
\def\2#1{\widetilde{#1}}
\def\3#1{\widehat{#1}}
\def\4#1{\mathbb{#1}}
\def\5#1{\frak{#1}}
\def\6#1{{\mathcal{#1}}}

\def\C{{\4C}}
\def\R{{\4R}}
\def\n{{\4n}}
\def\Z{{\4Z}}
\author{Valentin Burcea  }
\title{CR Singularities and Generalizations of Moser's Theorem II}

 \begin{abstract}  Let 
  a Real-Analytic  Submanifold $M$     formally (holomorphically) equivalent to the following model
 \begin{equation*}w=z_{1}\overline{z}_{1}+z_{2}\overline{z}_{2}+\dots+z_{N}\overline{z}_{N}+\lambda_{1}\left(z_{1}^{2}+\overline{z}_{1}^{2}\right)+\lambda_{2}\left(z_{2}^{2}+\overline{z}_{2}^{2}\right)+\dots+\lambda_{N}\left(z_{N}^{2}+\overline{z}_{N}^{2}\right),\quad \mbox{where $  \lambda_{1},\lambda_{2},\dots, \lambda_{N}\in \left[0,\frac{1}{2}\right)$.}\end{equation*}
 
 It is proven that $M$ is holomorphically equivalent to this model by considering   Spaces of Normalizations, which are also used also in proving other analogues of The Theorem of Moser in certain non-equidimensional situations.  Also, these methods allow   to develop   a  (formal) normal form for a    class  of C.-R. Singular $2$-Codimensional  Real-Formal Submanifolds in Complex Spaces  Moreover,  the defined Spaces of Normalizations are used in order to provide a classification of  the so-called infinitesimal holomorphic   automorphisms in the above studied case. 

Also in the above considered situation, we consider Formal Holomorphic Segre Preserving Equivalences by defining Double-Spaces of Normalizations. Then, there   are proven   other analogues  of The Theorem of Moser using The Approximation Theorem of Artin in the equidimensional case, and then in the non-equidimensional case by reformulating  simple computations  from the standard case.

\end{abstract}

\dedicatory{To the Memory of Professor Donal O'Donovan, who was My Friend}
\address{V. Burcea: INDEPENDENT}
\email{valentin@maths.tcd.ie}
\thanks{ THANKS to  Science Foundation Ireland Grant 10/RFP/MTH 2878 for ENORMOUS FUNDING while I was working in Trinity College Dublin, Ireland}
\thanks{ Special thanks to CAPES for generous funding while I was working in The Federal University of Santa Catarina, Brazil} 
\thanks{ The reference \cite{bu1} was fully supported by  Science Foundation Ireland Grant 06/RFP/MAT 018.}

  \maketitle

\def\Label#1{\label{#1}{\bf (#1)}~}


\def\cn{{\C^n}}
\def\cnn{{\C^{n'}}}
\def\ocn{\2{\C^n}}
\def\ocnn{\2{\C^{n'}}}


\def\dist{{\rm dist}}
\def\const{{\rm const}}
\def\rk{{\rm rank\,}}
\def\id{{\sf id}}
\def\tr{{\bf tr\,}}
\def\aut{{\sf aut}}
\def\Aut{{\sf Aut}}
\def\CR{{\rm CR}}
\def\GL{{\sf GL}}
\def\Re{{\sf Re}\,}
\def\Im{{\sf Im}\,}
\def\span{\text{\rm span}}
\def\Diff{{\sf Diff}}

\def\codim{{\rm codim}}
\def\crd{\dim_{{\rm CR}}}
\def\crc{{\rm codim_{CR}}}

\def\phi{\varphi}
\def\eps{\varepsilon}
\def\d{\partial}
\def\a{\alpha}
\def\b{\beta}
\def\g{\gamma}
\def\G{\Gamma}
\def\D{\Delta}
\def\Om{\Omega}
\def\k{\kappa}
\def\l{\lambda}
\def\L{\Lambda}
\def\z{{\bar z}}
\def\w{{\bar w}}
\def\Z{{\1Z}}
\def\t{\tau}
\def\th{\theta}

\emergencystretch15pt \frenchspacing

\newtheorem{Thm}{Theorem}[section]
\newtheorem{Cor}[Thm]{Corollary}
\newtheorem{Pro}[Thm]{Proposition}
\newtheorem{Lem}[Thm]{Lemma}

\theoremstyle{definition}\newtheorem{Def}[Thm]{Definition}

\theoremstyle{remark}
\newtheorem{Rem}[Thm]{Remark}
\newtheorem{Exa}[Thm]{Example}
\newtheorem{Exs}[Thm]{Examples}

\def\bl{\begin{Lem}}
\def\el{\end{Lem}}
\def\bp{\begin{Pro}}
\def\ep{\end{Pro}}
\def\bt{\begin{Thm}}
\def\et{\end{Thm}}
\def\bc{\begin{Cor}}
\def\ec{\end{Cor}}
\def\bd{\begin{Def}}
\def\ed{\end{Def}}
\def\br{\begin{Rem}}
\def\er{\end{Rem}}
\def\be{\begin{Exa}}
\def\ee{\end{Exa}}
\def\bpf{\begin{proof}}
\def\epf{\end{proof}}
\def\ben{\begin{enumerate}}
\def\een{\end{enumerate}}
\def\beq{\begin{equation}}
\def\eeq{\end{equation}}

\section{Introduction and Main Results}
 This paper  continues the author's research\cite{bu1},\cite{bu2},\cite{bu3},\cite{bu3se},\cite{bu4} regarding  Real Submanifolds in $\mathbb{C}^{N+1}$ defined near C.-R.  Singularities\cite{bi} and the local equivalence problem\cite{ko2},\cite{ko3} in Complex Analysis, which generally asks when two Real Submanifolds are actually equivalent.  Its non-triviality was detected by Moser-Webster\cite{mowe}, Gong\cite{go2} and Kossovskiy-Shafikov\cite{kosh}. They studied the existence and construction of two Real-Analytic Submanifolds in Complex Spaces, which are only formally (holomorphically) equivalent. 
 
 We recall that if $T^{c}_{q}M$ is   the complex tangent space of $M$ at $q$, the point $p=0$ is called C.-R. Singularity\cite{bi}   of the Real Submanifold $M\subset \mathbb{C}^{N+1}$, if $p=0$ is a jumping point for the   mapping $M\ni q\mapsto \dim_{\mathbb{C}}T^{c}_{q}M$ defined near $p=0$. In particular, 
let the coordinates $$\left(w,z\right)=\left(w,z_{1},z_{2},\dots,z_{N}\right)\in \mathbb{C}^{N+1}.$$ 
 
 The  first   result  is the following: 
\bt\label{t1} Let   $M\subset\mathbb{C}^{N+1}$ be the following real-analytic submanifold  
\begin{equation}w=z_{1}\overline{z}_{1}+z_{2}\overline{z}_{2}+\dots+z_{N}\overline{z}_{N}+\lambda_{1}\left(z_{1}^{2}+\overline{z}_{1}^{2}\right)+\lambda_{2}\left(z_{2}^{2}+\overline{z}_{2}^{2}\right)+\dots+\lambda_{N}\left(z_{N}^{2}+\overline{z}_{N}^{2}\right)+\mbox{O}(3), \label{var1}
\end{equation}
where the not all vanishing numbers $\lambda_{1},\lambda_{2},\dots,\lambda_{N}$ are chosen such that  
\begin{equation}  \lambda_{1},\lambda_{2},\dots, \lambda_{N}\in \left[0,\frac{1}{2}\right),\label{lambida}
\end{equation}
such that $M$ is formally  equivalent to the following  model
\begin{equation}M_{\infty}:\hspace{0.1 cm}w=z_{1}\overline{z}_{1}+z_{2}\overline{z}_{2}+\dots+z_{N}\overline{z}_{N}+\lambda_{1}\left(z_{1}^{2}+\overline{z}_{1}^{2}\right)+\lambda_{2}\left(z_{2}^{2}+\overline{z}_{2}^{2}\right)+\dots+\lambda_{N}\left(z_{N}^{2}+\overline{z}_{N}^{2}\right).\label{mmodel1}\end{equation}

Then   $M$ is  holomorphicaly equivalent to $M_{\infty}$.
\et 

This first result has been proven in $\mathbb{C}^{2}$ by Moser\cite{mo} for $\lambda_{1}=0$, being known  as the Theorem of Moser\cite{mo}. More generally,  Huang-Yin\cite{huyi1} obtained this result for $\lambda_{1}=\dots=\lambda_{N}=0$ in $\mathbb{C}^{N+1}$, and also the author in   more general situations by adapting their approaches. It  follows also in $\mathbb{C}^{2}$ from    Moser-Webster\cite{mowe}.    The author proved\cite{bu4} a similar result  when the model (\ref{mmodel1}) is perturbated by  terms of degree $3$, being motivated by Gong-Lebl\cite{gole}. His approach\cite{bu4} was based on trying to understand the C.-R. structure existent near a C.-R. Singularity, but that approach\cite{bu4} does not apply in this situation defined by (\ref{var1}).

This first result  is proven according to the following  approach. We  construct  a partial normal form by imposing convenient normalizations derived from Fischer Decompositions\cite{sh}. More precisely, it is required to consider  certain Fischer Normalization Spaces recalling the strategy from \cite{bu2}. This procedure  leaves underdetermined an infinite number of parameters making  the formal (holomorphic) equivalence  possibly divergent. These undetermined parameters are eliminated using compositions with formal automorphisms of the   model (\ref{mmodel1}).   Then the convergence  of the normalized formal (holomorphic) equivalence  is concluded  by adapting\cite{bu4} the  iterating approach of Moser\cite{mo} after there are made  suitable estimations using the     natural presence in the local computations of  the following matrix
\begin{equation}\Lambda=\begin{pmatrix}\lambda_{1} & 0     &\dots & 0  \\ 0 & \lambda_{2}   &\dots & 0  \\  \vdots & \vdots &\ddots & \vdots \\ 0  & 0  &   \dots & \lambda_{N} 
\end{pmatrix}. \label{Lambda}
\end{equation}

 Then  Moser's method\cite{mo} applies and the result follows recalling procedures from \cite{bu4}.   Moreover, the considered approach allows  to develop a formal normal form for a certain class of real-formal submanifolds.  In particular, we obtain a formal normal form in $\mathbb{C}^{2}$ different than Moser-Webster's Normal Form\cite{mowe} by  applying indirectly  the construction procedure of Huang-Yin's Normal Form\cite{huyi2}. However, our considered situation is   different because to the very complicated interactions of terms in the local defining equations.  Again, the Fischer Decompositions\cite{sh} simplify the local defining equations by imposing convenient normalization conditions. 
 
 Actually, our   constructions depends on these so-called Spaces of Fischer Normalizations, which are the hidden core subject of this paper. These spaces are fueled  by the approach from \cite{bu3} according to the author's initiation\cite{bu1}. Their definitions  use generalized Fischer Decompositions\cite{sh}   as in \cite{bu4}, being based on an approach depending on    complicated systems of equations depending on very non-trivial interactions of matrices and many tedious computations. 
 
  These methods are adapted in order to study an  analogue  of The Theorem of Moser\cite{mo} using Formal Holomorphic Segre Equivalences. We recall from \cite{V77} that the Formal Holomorphic Segre Equivalences are just    Formal Holomorphic Segre Preserving Mappings\cite{An} with non-vanishing Jacobian.

  Let $M,\hspace{0.1 cm}N\subset \mathbb{C}^{N+1}$ two real-formal submanifolds   defined near $0\in\mathbb{C}^{N+1}$ as follows
\begin{equation}M:\rho_{1}\left(w,z,\overline{z},\overline{w}\right)=0,\quad N:\rho_{2}\left(w,z,\overline{z},\overline{w}\right)=0,
\end{equation}
with respect to the following considerations
\begin{equation}\begin{split}&\rho_{1}\left(w,z,\overline{z},\overline{w}\right)=-w+z_{1}\overline{z}_{1}+z_{2}\overline{z}_{2}+\dots+z_{N}\overline{z}_{N}+\lambda_{1}\left(z_{1}^{2}+\overline{z}_{1}^{2}\right)+\lambda_{2}\left(z_{2}^{2}+\overline{z}_{2}^{2}\right)+\dots+\lambda_{N}\left(z_{N}^{2}+\overline{z}_{N}^{2}\right)+\mbox{O}(3),\\&\rho_{2}\left(w,z,\overline{z},\overline{w}\right)=-w+z_{1}\overline{z}_{1}+z_{2}\overline{z}_{2}+\dots+z_{N}\overline{z}_{N}+\lambda_{1}\left(z_{1}^{2}+\overline{z}_{1}^{2}\right)+\lambda_{2}\left(z_{2}^{2}+\overline{z}_{2}^{2}\right)+\dots+\lambda_{N}\left(z_{N}^{2}+\overline{z}_{N}^{2}\right)+\mbox{O}(3),\end{split}\label{cete2}
\end{equation}
where $\lambda_{1},\lambda_{2},\dots,\lambda_{N}\geq 0$.

Then, any (Formal) Holomorphic Segre Preserving Mapping between $M$ and $N$ is defined by \cite{An} as follows
\begin{equation}\mathcal{H}:\mathbb{C}^{2N+2}\rightarrow \mathbb{C}^{2N+2}\quad\mbox{such that $\mathcal{H}\left(w,z,\zeta,\nu\right)=\left(H\left(z,w\right),\tilde{H}\left(\zeta,\nu\right)\right)$,} \label{cete}
\end{equation} 
 where $H,\tilde{H}:\mathbb{C}^{N+1}\rightarrow \mathbb{C}^{N+1}$ are formal holomorphic mappings such that  
\begin{equation}\rho_{2}\left(H\left(z,w\right),\tilde{H}\left(\zeta,\nu\right) \right)=0,\quad\mbox{ for all $\left(w,z,\zeta,\nu\right)\in \mathbb{C}^{2N+2}$ such that $\rho_{1}\left(w,z,\zeta,\nu\right)=0$,}
\end{equation}
where the variables $$\zeta=\left(\zeta_{1},\zeta_{2},\dots,\zeta_{N},\right)$$ 
replace the variables $$\overline{z}=\left(\overline{z}_{1},\overline{z}_{2},\dots,\overline{z}_{N},\right)$$ in the formal expansions from (\ref{cete2}), and respectively where $\nu$ replaces $\overline{w}$ in the formal expansions from (\ref{cete2}). 

Such mappings have the obvious property that preserve the Segre varieties. The  second main result  is the following:  
  
\bt\label{t2} Let   $M\subset\mathbb{C}^{N+1}$ be the following real-analytic submanifold  
\begin{equation}w=z_{1}\overline{z}_{1}+z_{2}\overline{z}_{2}+\dots+z_{N}\overline{z}_{N}+\lambda_{1}\left(z_{1}^{2}+\overline{z}_{1}^{2}\right)+\lambda_{2}\left(z_{2}^{2}+\overline{z}_{2}^{2}\right)+\dots+\lambda_{N}\left(z_{N}^{2}+\overline{z}_{N}^{2}\right)+\mbox{O}(3), \label{var11}
\end{equation}
where the numbers $\lambda_{1},\lambda_{2},\dots,\lambda_{N}$ are chosen such that  
\begin{equation}\lambda_{1}>0,\quad  \lambda_{2},\dots, \lambda_{N}\in \left[0,\infty\right),\label{lambidalala}
\end{equation}
such  that $M$ is      Formally Holomorphic Segre equivalent  to the following  model
\begin{equation}M_{\infty}:\hspace{0.1 cm}w=z_{1}\overline{z}_{1}+z_{2}\overline{z}_{2}+\dots+z_{N}\overline{z}_{N}+\lambda_{1}\left(z_{1}^{2}+\overline{z}_{1}^{2}\right)+\lambda_{2}\left(z_{2}^{2}+\overline{z}_{2}^{2}\right)+\dots+\lambda_{N}\left(z_{N}^{2}+\overline{z}_{N}^{2}\right).\label{mmodel1se}\end{equation}

Then   $M$ is  Segre-Holomorphically  equivalent to $M_{\infty}$.
\et 
 
This second result  is proven by applying a similar approach. We consider  a construction of partial normal form type by imposing convenient normalizations derived from Fischer Decompositions\cite{sh} as previously. On the other hand, it is required to consider Double-Spaces  of Fischer Normalizations  according to the previous   strategy. This leads to the solving of very non-trivial systems of equations with respect to (\ref{Lambda}). Now, the convergence  of the   formal  Segre-holomorphic  preserving equivalence  is concluded    by adapting The Theorem of Artin\cite{A} according to the approach from \cite{V77}. In particular, the methods of Mir\cite{mi1},\cite{mi2} provide the convergence of this formal mapping.     
 
Regardless of the  complicated computational aspects, these approaches are   very useful for obtaining other results related to the local equivalence problem \cite{Salla} in Complex Analysis. More precisely, it is studied when two CR Singular Real-Formal Submanifolds are equivalent  without considering classical  assumptions as their formal  equivalence: 
 
  Let   $M_{\infty}\subset\mathbb{C}^{N+1}$ be the following model 
\begin{equation}w=z_{1}\overline{z}_{1}+z_{2}\overline{z}_{2}+\dots+z_{N}\overline{z}_{N}+\lambda_{1}\left(z_{1}^{2}+\overline{z}_{1}^{2}\right)+\lambda_{2}\left(z_{2}^{2}+\overline{z}_{2}^{2}\right)+\dots+\lambda_{N}\left(z_{N}^{2}+\overline{z}_{N}^{2}\right), \label{var1A}
\end{equation}
where the positive real numbers $\lambda_{1},\lambda_{2},\dots,\lambda_{N}$ are chosen such that
\begin{equation}   \lambda_{1},\lambda_{2},\dots, \lambda_{N}\in \left(0,\frac{1}{2}\right).\label{lambida1}
\end{equation}

  Let $N<N'$  in respect to the coordinates $$\left(w',z'_{1},z'_{2},\dots,z'_{N'}\right)\in\mathbb{C}^{N'+1}.$$
 
  Let   $M'_{\infty}\subset\mathbb{C}^{N'+1}$ be the following model 
\begin{equation}w'={z'}_{1}\overline{z'}_{1}+{z'}_{2}\overline{z'}_{2}+\dots+{z'}_{N}\overline{z'}_{N'}+\lambda'_{1}\left({z'}_{1}^{2}+\overline{z'}_{1}^{2}\right)+\lambda'_{2}\left({z'}_{2}^{2}+\overline{z'}_{2}^{2}\right)+\dots+{\lambda'}_{N'}\left({z'}_{N}^{2}+\overline{z'}_{N'}^{2}\right), \label{var1B}
\end{equation} 
where the positive real numbers ${\lambda'}_{1},{\lambda'}_{2},\dots,{\lambda'}_{N'}$ are chosen   such that  
\begin{equation}\lambda_{1}={\lambda'}_{1},\lambda_{2}={\lambda'}_{2},\dots,\lambda_{N}={\lambda'}_{N}  .\label{maro}\end{equation}

 An application is the following: 
 \bt\label{tA1} Let $\left(G,F\right): M_{\infty}\rightarrow {M'}_{\infty}$ be a non-constant formal   mapping  respecting (\ref{var1A}), (\ref{lambida1}), (\ref{var1B}), (\ref{maro}). Then, there  exist $\psi \in \mbox{Aut}_{0}\left(M_{\infty}\right)$ and $\varphi \in  \mbox{Aut}_{0}\left(M'_{\infty}\right)$ such that 
\begin{equation}
\varphi \circ \left(G,F\right) \circ \psi\left(w,z_{1},z_{2},\dots,z_{N}\right)=\left(w,z_{1},z_{2},\dots,z_{N},0,\dots,0\right).\label{bebe1}
\end{equation}
\et  

This   result is interesting in the light of the classifications  for  proper holomorphic mappings initiated by Webster\cite{we}. It is well-known the relationship between proper holomorphic mappings and Cauchy-Riemann Mappings\cite{DA}.  The classification  of the proper  holomorphic mappings, between unit balls situated in complex spaces of different dimensions, is reduced to the study  of  the CR mappings between  hyperquadrics\cite{DA}, making interesting the studies \cite{BH},\cite{BEH},\cite{CJX},\cite{Far},\cite{H},\cite{huang1},\cite{huang2},\cite{HJ},\cite{JX},\cite{L}  concerning   mappings  between hyperquadrics situated in Complex Spaces of different dimensions.  Similar classifications have been concluded also by the author\cite{bu5} in other situations. Thus, it is interesting to study formal (holomorphic) mappings defined between Models situated in Complex Spaces of different dimensions.

This result defines  an analogue   to the classifications  of Huang\cite{huang1} and  Kim-Zaitsev\cite{kz}. The approach of Moser\cite{mo} naturally applies  obtaining surprisingly just a class of equivalence, which is given by the standard linear embedding (\ref{bebe1}). Therefore, there are not encountered gapping phenomenons\cite{HJ} as it normally happens, because of the existence of non-vanishing Bishop invariants. The Spaces of Fischer Normalizations  provide simple conclusions in such particular case, regardless of the    non-trivial interactions of terms, including that 
\begin{equation}
 \frac{\partial G}{\partial w}(0)\neq 0,\label{transver}
 \end{equation}
which is an expected transversality conclusion in the light of Baouendi-Huang's Conjecture\cite{BHC}. 
 
Again surprisingly, the analogue of the this result holds also by considering Segre-Holomorphic Preserving Mappings:
\bt\label{tA1sec} Let $\left(G,\tilde{G},F,\tilde{F}\right)$ be a non-constant formal Segre-holomorphic preserving  mapping from $M_{\infty}$ into ${M'}_{\infty}$  respecting (\ref{var1A}), (\ref{lambida1}), (\ref{var1B}), (\ref{maro}). Then, there  exist $\tilde{\psi} \in \tilde{\mbox{Aut}_{0}}\left(M_{\infty}\right)$ and $\tilde{\varphi} \in \tilde{\mbox{Aut}_{0}}\left(M'_{\infty}\right)$ such that 
\begin{equation}
\tilde{\varphi} \circ \left(G,\tilde{G},F,\tilde{F}\right) \circ \tilde{\psi}\left(w,z_{1},z_{2},\dots,z_{N},\nu,\xi_{1},\xi_{2},\dots,\xi_{N}\right)=\left(w,z_{1},z_{2},\dots,z_{N},0,\dots,0,\nu,\xi_{1},\xi_{2},\dots,\xi_{N},0,\dots,0\right),\label{bebe1se}
\end{equation}
where $\tilde{\mbox{Aut}_{0}}\left(M_{\infty}\right)$ represents the group of formal Segre-holomorphic preserving automorphism of $M_{\infty}$, and respectively $\tilde{\mbox{Aut}_{0}}\left(M'_{\infty}\right)$  represents the group of formal Segre-holomorphic preserving automorphism of $ {M'}_{\infty}$.
\et

This result defines  an analogue    to the classification   of Zhang\cite{Zha}. The previous computational approach is adapted   obtaining  just a class of equivalence, again surprisingly, which is given by the standard linear Segre-embedding (\ref{bebe1se}), because Zhang's Classification\cite{Zha} is different. Similarly as previously, the Double-Spaces  of Fischer Normalizations provide simple and unexpcted conclusions, regardless of the    non-trivial interactions of terms, including that 
\begin{equation}
 \det\frac{\partial \left(G,\tilde{G}\right)}{\partial \left(w,\nu\right)}(0)\neq 0,\label{transverse}
 \end{equation}
because   Zhang\cite{Zha} proved that the transversality of the Segre-Holomorphic preserving mappings may not occur generally.

The main purpose, of these applications, is to understand   local equivalence problems for Real-Formal Submanifolds embedded in model manifolds in Complex Spaces.  The above results may be interesting from this point of view, according to the studies initiated by Webster\cite{we1} and Ebenfelt-Huang-Zaitsev\cite{ebhuza}. More precisely, we consider the following the Real-Formal Submanifold 
\begin{equation}M:\quad w=z_{1}\overline{z}_{1}+z_{2}\overline{z}_{2}+\dots+z_{N}\overline{z}_{N}+\lambda_{1}\left(z_{1}^{2}+\overline{z}_{1}^{2}\right)+\lambda_{2}\left(z_{2}^{2}+\overline{z}_{2}^{2}\right)+\dots+\lambda_{N}\left(z_{N}^{2}+\overline{z}_{N}^{2}\right)+\mbox{O}(3)\subset\mathbb{C}^{N+1}, \label{var1Ase}
\end{equation}
and respectively, we consider   the Real-Formal Submanifold
\begin{equation}M':\quad w'={z'}_{1}\overline{z'}_{1}+{z'}_{2}\overline{z'}_{2}+\dots+{z'}_{N}\overline{z'}_{N'}+\lambda'_{1}\left({z'}_{1}^{2}+\overline{z'}_{1}^{2}\right)+\lambda'_{2}\left({z'}_{2}^{2}+\overline{z'}_{2}^{2}\right)+\dots+{\lambda'}_{N'}\left({z'}_{N}^{2}+\overline{z'}_{N'}^{2}\right)+\mbox{O}(3)\subset\mathbb{C}^{N'+1}. \label{var1Bse}
\end{equation}

  Let $N<\tilde{N}$ in respect to the coordinates
 $$\left(\tilde{w},\tilde{z}_{1},\tilde{z}_{2},\dots,\tilde{z}_{N}\right)\in\mathbb{C}^{\tilde{N}+1}.$$

 Let   $\tilde{M}_{\infty}\subset\mathbb{C}^{\tilde{N}+1}$ be the following model 
\begin{equation}\tilde{w}=\tilde{z}_{1}\overline{\tilde{z}}_{1}+\tilde{z}_{2}\overline{\tilde{z}}_{2}+\dots+\tilde{z}_{\tilde{N}}\overline{\tilde{z}}_{\tilde{N}}+\tilde{\lambda}_{1}\left(\tilde{z}_{1}^{2}+\overline{\tilde{z}}_{1}^{2}\right)+\tilde{\lambda}_{2}\left(\tilde{z}_{2}^{2}+\overline{\tilde{z}}_{2}^{2}\right)+\dots+\tilde{\lambda}_{\tilde{N}}\left(\tilde{z}_{\tilde{N}}^{2}+\overline{\tilde{z}}_{\tilde{N}}^{2}\right), \label{var1AX}
\end{equation}
where the positive real numbers $\tilde{\lambda}_{1},\tilde{\lambda}_{2},\dots,\tilde{\lambda}_{\tilde{N}}$ are chosen such that
\begin{equation}\tilde{\lambda}_{1}={\lambda'}_{1},\tilde{\lambda}_{2}={\lambda'}_{2},\dots, \tilde{\lambda}_{N'}={\lambda'}_{N'}.\label{lambida1X}
\end{equation}

 Let $N'<\tilde{N'}$ in respect to the coordinates
 $$\left(\tilde{w'},\tilde{z'}_{1},\tilde{z'}_{2},\dots,\tilde{z'}_{\tilde{N'}}\right)\in\mathbb{C}^{\tilde{N'}+1}.$$

Let   $\tilde{M'}_{\infty}\subset\mathbb{C}^{\tilde{N'}+1}$ be the following model 
\begin{equation}\tilde{w'}=\tilde{z'}_{1}\overline{\tilde{z'}}_{1}+\tilde{z'}_{2}\overline{\tilde{z'}}_{2}+\dots+\tilde{z'}_{\tilde{N'}}\overline{\tilde{z'}}_{\tilde{N'}}+\tilde{\lambda'}_{1}\left(\tilde{z'}_{1}^{2}+\overline{\tilde{z'}}_{1}^{2}\right)+\tilde{\lambda'}_{2}\left(\tilde{z'}_{2}^{2}+\overline{\tilde{z'}}_{2}^{2}\right)+\dots+\tilde{\lambda'}_{N}\left(\tilde{z'}_{N}^{2}+\overline{\tilde{z'}}_{N}^{2}\right), \label{var1BX}
\end{equation}
where the positive real numbers $\tilde{\lambda'}_{1},\tilde{\lambda'}_{2},\dots,\tilde{\lambda'}_{\tilde{N'}}$ are chosen such that
\begin{equation}\tilde{\lambda'}_{1}={\lambda'}_{1},\tilde{\lambda'}_{2}={\lambda'}_{2},\dots, \tilde{\lambda'}_{N'}={\lambda'}_{N'}.\label{lambida1X1}
\end{equation}

Moreover, we assume the following
\begin{equation}   {\lambda'}_{1},{\lambda'}_{2},\dots, {\lambda'}_{N'}\in \left(0,\frac{1}{2}\right).\label{lambida1111}
\end{equation}

We say that $M$ is (formally holomorphically) embendable in $\tilde{M}_{\infty}$ when there exist   formal (holomorphic) embeddings of $M$ into $\tilde{M}_{\infty}$. Respectively, we say that $M$ is  formally Segre-holomorphically  embendable in $\tilde{M}_{\infty}$ when there exist   formal  Segre-holomorphic preserving embeddings of $M$ into $\tilde{M}_{\infty}$. This   terminology is analogously considered  in the case of  $M'$ and $\tilde{M'}_{\infty}$ in the standard case as well as in the Segre-case. Moreover, by using abuses of notations, we write as follows
\begin{equation*}
  M\mapsto\tilde{M}_{\infty},\quad\quad\quad   M'\mapsto \tilde{M'}_{\infty}. 
\end{equation*}
  
Clearly, these definitions and notations are locally considered.  Of course, it remains to study the uniqueness of such embenddings.  We will show  that two such embeddings are  always unique up to a composition with an automorphism of the corresponding model in the standard case as well as in the Segre-case. Then, we can establish similar results  towards to Webster\cite{we1} and Ebenfelt-Huang-Zaitsev\cite{ebhuza}. We have 
\bt\label{tA2}   Consider the following diagram
 \begin{equation}\begin{array}[c]{ccc}
\mathbb{C}^{\tilde{N}+1}\supset\tilde{M}_{\infty}& \ni\left(w,z_{1},z_{2},\dots,z_{\tilde{N}} \right) \longmapsto\left(z_{1},z_{2},\dots,z_{\tilde{N}},0,\dots,0\right)\in &\tilde{M'}_{\infty}\subset \mathbb{C}^{\tilde{N'}+1}\\
\uparrow\scriptstyle{ }&&\uparrow\scriptstyle{}\\\mathbb{C}^{N+1}
\supset M & & M'\subset \mathbb{C}^{N'+1}
\end{array}. \label{diag1}\end{equation} 
 Assume    that (\ref{var1Ase}), (\ref{var1Bse}), (\ref{var1AX}),   (\ref{var1BX}) hold and    that $ \tilde{N}\leq\tilde{N'}$ with   $N',N>1$.  Then   $M$ is embedded   in $M'$ according to the standard linear embedding in the sense of (\ref{bebe1}). Moreover, any formal     mapping,  between $M$ and $M'$ respecting (\ref{transver}), is equivalent the standard linear embedding in the sense of (\ref{bebe1}). \et 
 
Respectively, we have 
\bt\label{tA2sec}   Consider the following diagram
 \begin{equation}\begin{array}[c]{ccc}
\mathbb{C}^{\tilde{N}+1}\supset\tilde{M}_{\infty}& \ni\left(w,z_{1},z_{2},\dots,z_{\tilde{N}} \right) \longmapsto\left(w,z_{1},z_{2},\dots,z_{\tilde{N}},0,\dots,0\right)\in &\tilde{M'}_{\infty}\subset \mathbb{C}^{\tilde{N'}+1}\\
\uparrow\scriptstyle{ }&&\uparrow\scriptstyle{}\\\mathbb{C}^{N+1}
\supset M & & M'\subset \mathbb{C}^{N'+1}
\end{array}. \label{diag1}\end{equation} 
 Assume    that (\ref{var1Ase}), (\ref{var1Bse}), (\ref{var1AX}),   (\ref{var1BX}) hold and    that $ \tilde{N}\leq\tilde{N'}$ with   $N',N>1$.  Then   $M$ is    embedded  in $M'$ according to     the standard linear Segre-holomorphic preserving embeddeding in the sense of (\ref{bebe1se}). Moreover, any non-constant formal  Segre-holomorphic preserving  mapping,  between $M$ and $M'$ respecting (\ref{transverse}), is equivalent   the standard linear Segre-holomorphic preserving embeddeding in the sense of (\ref{bebe1se}). \et  
 
 These results may be considered    non-equidimensional versions of The Theorem of Moser\cite{mo}. 
 The proofs depend  on the previously mentioned Spaces of Fischer Normalization  and on the  Double-Spaces of Fischer Normalizations, which are used in order to  make computations and then to  normalize   the (formal holomorphic) mappings using the corresponding groups of automorphisms of   models. Clearly, other models can be considered as in \cite{bu11}, but we use   usual models corresponding to the Special Elliptic CR Singularities \cite{do1},\cite{do2},\cite{dotoza1},\cite{dotoza2}. 
 
\subsection{Organization} The first part of this paper is organized as follows:  Section  $2$ introduces the preliminaries  in order to define the so-called Spaces of Normalizations in Sections $3$ and $4$, on which are based the normalization procedures developed in Section $5$. In Sections  $9$ and $10$ are adapted the computations    developed in Sections $3$ and $4$ in order to construct the so-called Double-Spaces of Normalizations, on which are based the normalization procedures developed in Section $11$ according to the ingredients from Section $8$. Then, the proof of Theorem \cite{t1} is completed in Section respecting the normalizations developed in Section $5$ according to a modification of Moser's Iterative Procedure\cite{mo}. The proof of Theorem \ref{t2} is completed in Section respecting the normalizations developed in Section $11$ according to Artin's Approximation Theorem\cite{A}. The proofs of Theorems \ref{tA1} and \ref{tA2} are completed in Section $7$. The proofs of Theorems \ref{tA1sec} and \ref{tA2sec} are completed in Section $12$.

\section{Acknowledgements} 

Special Thanks to Capes for generous funding while I was working in The Federal University of Santa Catarina, Brazil, including some travel funding and mentioning  that I did not receive the scholarship (of value $4100$ BRL) for the first month  for while I had been working in The Federal University of Minas Gerais, in the conditions that  I had lost $6000$ BRL before arriving in the light of several obstacles. Special Thanks to   Science Foundation of Ireland for  enormous investments in my  doctoral formation while I was working in Trinity College Dublin under the supervision of Prof. Dmitri Zaitsev.   Special Thanks also to Dr. J.Lebl for stimulating conversations. 

 The Financial Thanks are organized like the relevance of each Funding.  I also make clear the main part of my doctoral thesis \cite{bu1} was fully supported by Science Foundation Ireland Grant 06/RFP/MAT 018.

 \newpage
\section{Ingredients and Computations}
  We make computations according to indications from my professor\cite{za}. In particular, we   study how a formal (holomorphic)  equivalence occurs in the local defining equations:
\subsection{Changes of Coordinates}Let $M\subset\mathbb{C}^{N+1}$ be a Real-Formal
Submanifold defined near  $p=0$ as follows
\begin{equation}
w=Q\left(z,\overline{z}\right)+\displaystyle\sum
_{k\geq 3}\varphi_{k}(z,\overline{z}), \label{ecuatie}
\end{equation}
where $Q\left(z,\overline{z}\right)$ is a real-quadratic form, and respectively $\varphi _{k}(z,\overline{z})$ is a  polynomial of    degree $k$ in $(z,\overline{z})$,  for all  $k\geq 3$. 
 
 We consider   another Real-Formal
Submanifoldd $M'\subset\mathbb{C}^{N+1}$ defined near $p=0$ as follows
\begin{equation}w'= Q\left(z',\overline{z'}\right)+\displaystyle\sum _{k\geq
3}\varphi'_{k}\left(z',\overline{z'}\right),\label{ecuatie1}
\end{equation}
where $\varphi' _{k}\left(z',\overline{z'}\right)$ is a
 polynomial of   degree $k$ in $\left(z',\overline{z'}\right)$,  for all  $k\geq 3$. 

Let  a formal equivalence between  $M$ and $M'$, which fixes the point
$0\in\mathbb{C}^{N+1}$, denoted as  
$$\left(z',w'\right)=\left(F(z,w),G(z,w)\right).$$ 

Then, for  $w$   defined by (\ref{ecuatie}), we have  
\begin{equation} G(z,w)=Q\left(
F(z,w),F(z,w)\right)+\displaystyle\sum _{k\geq
3}\varphi'_{k}\left(F(z,w),\overline{F(z,w)}\right). \label{ec}
\end{equation}

In order to understand better the interactions of terms in (\ref{ec}), we write the (formal) equivalence as follows
\begin{equation}\left(F(z,w),G(z,w)\right)
=\left(\displaystyle\sum_{m,n\geq 0}F_{m,n}(z)w^{n},\displaystyle\sum_{m,n\geq 0}G_{m,n}(z)w^{n}\right),
\label{map}\end{equation}
where we deal with  homogeneous polynomials of degree
$m$ in $z$ for all  $ n\in\mathbb{N}$, denoted as
\begin{equation}
G_{m,n}(z), \label{polin1}
\end{equation}
and respectively, denoted as 
\begin{equation}
F_{m,n}(z),\label{polin2}
\end{equation}
concluding by (\ref{ec}) and (\ref{map}) that
\begin{equation}\left.\begin{split}\displaystyle &  \quad\quad\quad\quad\quad\quad\quad\quad\quad\quad\quad \sum
_{m,n\geq0}G_{m,n}(z)\left(Q\left(z,\overline{z}\right) +\displaystyle\sum
_{k\geq 3}\varphi_{k}(z,\overline{z})\right)^{n}\\& \quad\quad\quad\quad\quad   \quad\quad\quad\quad\quad\quad\quad \quad\quad\quad\quad\quad   \quad\quad\quad\quad\quad   \begin{tabular}{l} \rotatebox[origin=c]{270}{$=$}\end{tabular} \\&  Q\left(\displaystyle\sum _{m,n \geq 0}
F_{m,n}(z)\left(Q\left(z,\overline{z}\right)+\displaystyle\sum
_{k\geq 3}\varphi_{k}(z,\overline{z})\right)^{n},\overline{\displaystyle\sum _{m,n \geq 0}
F_{m,n}(z)\left(Q\left(z,\overline{z}\right)+\displaystyle\sum
_{k\geq 3}\varphi_{k}(z,\overline{z})\right)^{n}}\right)  \\& \quad\quad\quad\quad\quad   \quad\quad\quad\quad\quad\quad\quad\quad\quad\quad\quad\quad   \quad\quad\quad\quad\quad     \begin{tabular}{l} \rotatebox[origin=c]{270}{$+$}\end{tabular} \\& \displaystyle\sum
_{k\geq 3}\varphi_{k}'
\left(\displaystyle\sum _{m,n \geq
0}F_{m,n}(z)\left(Q\left(z,\overline{z}\right)+\displaystyle\sum
_{k\geq 3}\varphi_{k}(z,\overline{z})\right)^{n},\overline{\displaystyle\sum _{m,n \geq
0}F_{m,n}(z)\left(Q\left(z,\overline{z}\right)+\displaystyle\sum
_{k\geq 3}\varphi_{k}(z,\overline{z})\right)^{n}}\right).
\end{split}\right.
\label{ecuatiegenerala}\end{equation}
 
  We obtain $G_{0,0}(z)=0$ and
 $F_{0,0}(z)=0$, because the formal equivalence (\ref{map})   fixes  the
point $0\in\mathbb{C}^{N+1}$. Moreover, collecting the terms of bidegree $(1,0)$ and $(1,1)$ in
$(z,\overline{z})$ from (\ref{ecuatiegenerala}), we obtain   
$G_{1,0}(z)=0$  and   $\Im G_{0,1}=0$, because
  \begin{equation}G_{0,1}Q(z,\overline{z})=Q\left(
F_{1,0}(z),F_{1,0}(z)\right).\label{ge}\end{equation} 

We move forward with these computations assuming  
\begin{equation}
Q(z,\overline{z})=z_{1}\overline{z}_{1}+z_{2}\overline{z}_{2}+\dots+z_{N}\overline{z}_{N}+\lambda_{1}\left(z_{1}^{2}+\overline{z}_{1}^{2}\right)+\lambda_{2}\left(z_{2}^{2}+\overline{z}_{2}^{2}\right)+\dots+\lambda_{N}\left(z_{N}^{2}+\overline{z}_{N}^{2}\right),\label{model}
\end{equation}
where $\lambda_{1},\lambda_{2},\dots,\lambda_{N}$ satisfy (\ref{lambida}).

We can thus assume $G_{0,1}(z)=1$ and $F_{1,0}(z)=z$  by   composing (\ref{map})
with a   linear holomorphic automorphism of the following  model 
\begin{equation}
\Re w=z_{1}\overline{z}_{1}+z_{2}\overline{z}_{2}+\dots+z_{N}\overline{z}_{N}+\lambda_{1}\left(z_{1}^{2}+\overline{z}_{1}^{2}\right)+\lambda_{2}\left(z_{2}^{2}+\overline{z}_{2}^{2}\right)+\dots+\lambda_{N}\left(z_{N}^{2}+\overline{z}_{N}^{2}\right) ,\label{model1}
\end{equation}
where $\lambda_{1},\lambda_{2},\dots,\lambda_{N}$ satisfy (\ref{lambida}), because such automorphism induces trivially a linear holomorphic automorphism of the following  model 
$$  w=z_{1}\overline{z}_{1}+z_{2}\overline{z}_{2}+\dots+z_{N}\overline{z}_{N}+\lambda_{1}\left(z_{1}^{2}+\overline{z}_{1}^{2}\right)+\lambda_{2}\left(z_{2}^{2}+\overline{z}_{2}^{2}\right)+\dots+\lambda_{N}\left(z_{N}^{2}+\overline{z}_{N}^{2}\right), 
 $$
preserving thus the quadratic form in (\ref{ecuatie}). 

The interactions of the homogeneous terms of the formal    equivalence (\ref{map})  is more complicated in this situation from (\ref{ecuatiegenerala}) than in \cite{bu1}. However, these  computational obstacles are eliminated using iterative Fischer Decompositions\cite{sh} as follows:

\subsection{Fischer Decompositions\cite{sh}}Recalling the strategy from \cite{bu2}, we define by (\ref{model}) the following differential operator
\begin{equation} \tr =\frac{\partial^{2}}{\partial z_{1}\partial\overline{z}_{1}}+\frac{\partial^{2}}{\partial z_{2}\partial\overline{z}_{2}}+\dots+\frac{\partial^{2}}{\partial z_{N}\partial\overline{z}_{N}}+\lambda_{1}\left( \frac{\partial^{2} }{\partial z_{1}^{2}}+\frac{\partial^{2} }{ \partial\overline{z}_{1}^{2}} \right)+\lambda_{2}\left( \frac{\partial^{2} }{\partial z_{2}^{2}}+\frac{\partial^{2} }{ \partial\overline{z}_{2}^{2}} \right)+\dots+\lambda_{N}\left(\frac{\partial^{2} }{\partial z_{N}^{2}}+\frac{\partial^{2} }{\partial\overline{z}_{N}^{2}}\right),\label{tracce}
\end{equation}
where $\lambda_{1},\lambda_{2},\dots,\lambda_{N}$ satisfy (\ref{lambida}).

This differential operator (\ref{tracce}) is just the Fischer differential operator associated to the polynomial (\ref{model}). It may be seen as an analogue of the trace operator considered by the author\cite{bu1} under the supervision of his   professor\cite{za}  in order to construct (formal) normal forms\cite{bu1}.

 Recalling the Fischer Decomposition from Shapiro\cite{sh}, we   write by (\ref{model}) uniquely as follows
\begin{equation}P(z)=A(z,\overline{z})Q(z,\overline{z})+C(z,\overline{z}),\quad\mbox{where $\tilde{\tr}\left(C(z,\overline{z})\right)=0$, }\label{V1} 
\end{equation}     
for any given homogeneous polynomial  $P(z)$, concluding  the following Fischer Decompositions
\begin{equation} \begin{split}& \quad\quad\quad\quad\quad\quad \hspace{0.05 cm} z^{I} =A(z,\overline{z})Q(z,\overline{z})+C(z,\overline{z}),\quad\quad\mbox{where $ \tr \left(C(z,\overline{z})\right)=0$, }\\&\left(  \overline{z}_{k}+2\lambda_{k}z_{k}\right)   z^{J} =C_{k}(z,\overline{z})Q(z,\overline{z})+D_{k}(z,\overline{z}),\quad\mbox{where $ \tr \left(D_{k}(z,\overline{z})\right)=0$,  }\end{split} 
\label{V2}
\end{equation}  
for all $k=1,\dots,N$, where we have used  the following notations
\begin{equation} \begin{split}&
z^{I}=z^{i_{1}}\cdot\dots\cdot z^{i_{N}},\quad\hspace{0.13 cm}\mbox{where $I=\left(i_{1},\dots,i_{N}\right)\in\mathbb{N}^{N}$ such that $\left|I\right|=i_{1}+\dots+i_{N}\geq 3$,}\\& z^{J}=z^{j_{1}}\cdot\dots\cdot z^{j_{N}},\quad\mbox{where $J=\left(j_{1},\dots,j_{N}\right)\in\mathbb{N}^{N}$ such that $\left|J\right|=j_{1}+\dots+j_{N}\geq 2$.} \end{split} \label{yy}
\end{equation} 

Clearly, it may occur the following
\begin{equation} \begin{split}&\quad\quad\quad\quad\quad \quad \tr  \left(z^{I}\right)=0,\quad\mbox{for $I=\left(i_{1},i_{2},\dots,i_{N}\right) \in\mathbb{N}^{N} $ with $i_{1}+i_{2}+\dots+i_{N}\geq 3$,}\\& \tr  \left(\left(  \overline{z}_{k}+2\lambda_{k}z_{k}\right) z^{J}\right)=0,\quad\mbox{for $J=\left(j_{1},j_{2},\dots,j_{N}\right) \in\mathbb{N}^{N} $ with $j_{1}+j_{2}+\dots+j_{N}\geq 2$, for all $k= 1,\dots,N$.}\end{split} \label{lap}
\end{equation}

In particular, we can choose $N=3$   and  $I=\left(1,1,1\right)$ or $k=1$ and $J=\left(0,1,1\right)$  in order to have (\ref{lap}) satisfied. Then, it is naturally required to define  the following  sets of multi-indexes
\begin{equation} \begin{split}&\hspace{0.1 cm}\mathcal{S}=\left\{ \mbox{$I=\left(i_{1},i_{2},\dots,i_{N}\right) \in\mathbb{N}^{N} $ such that $ \tr  \left(z^{I}\right)=0$ and $i_{1}+i_{2}+\dots+i_{N}\geq 3$ }\right\},\\&\mathcal{T}_{k}=\left\{ \mbox{$J=\left(j_{1},j_{2},\dots,j_{N}\right) \in\mathbb{N}^{N} $ such that $ \tr \left(\left(  \overline{z}_{k}+2\lambda_{k}z_{k}\right) z^{J}\right)=0$ and $j_{1}+j_{2}+\dots+j_{N}\geq 2$ }\right\},\end{split} \label{788}
\end{equation}
for all $k=1,\dots,N$.

These sets  (\ref{788}) play important roles in order to   define  Fischer Normalization Spaces\cite{bu2} that shape the normalizations in the local defining equation. The generalized version of the Fischer Decomposition\cite{sh} is widely applied by separating the real parts and the imaginary parts of the local defining equation at each degree level. We proceed as follows:  

\section{Fischer $G$-Decompositions\cite{bu2},\cite{bu3}}For any multi-index $\tilde{I}\not\in \mathcal{S}$ of length $p$, we consider by (\ref{model}) the following Fischer Decomposition 
\begin{equation}z^{\tilde{I}}=A(z,\overline{z}) Q(z,\overline{z})+C(z,\overline{z}),\hspace{0.1 cm}  \tr  \left(C(z,\overline{z})\right)=0,\quad\mbox{where  $\tilde{I}=\left(\tilde{i}_{1},\tilde{i}_{2},\dots,\tilde{i}_{N}\right)\in\mathbb{N}^{N}$.}\label{op}
\end{equation}

These homogeneous polynomials $A(z,\overline{z})$ and $C(z,\overline{z}$) are uniquely determined according to Shapiro\cite{sh}. We compute them straightforwardly  in (\ref{op}),  writing   as follows
\begin{equation} A(z,\overline{z})=\displaystyle\sum_{I,J\in\mathbb{N}^{N}\atop\left|I\right|+\left|J\right|=p-2} a_{I;J}z^{I}\overline{z}^{J}, \quad C(z,\overline{z})=\displaystyle\sum_{I,J\in\mathbb{N}^{N}\atop\left|I\right|+\left|J\right|=p} c_{I;J}z^{I}\overline{z}^{J}.\label{oppv}
\end{equation}

We have to understand the  contributions of the terms from (\ref{oppv}) in (\ref{op}). We  observe by (\ref{yy}) that
\begin{equation} \begin{split}& \quad\hspace{0.15 cm} \frac{\partial^{2}}{\partial z_{k}^{2}}\left(z^{I}\overline{z}^{J}\right)=i_{k}\left(i_{k}-1\right) z_{1}^{i_{1}}\dots z_{k}^{i_{k}-2}\dots z_{N}^{i_{N}}\overline{z}_{1}^{j_{1}}\dots  \overline{z}_{N}^{j_{N}},  \\&\quad\hspace{0.15 cm} \frac{\partial^{2}}{\partial \overline{z}_{k}^{2}}\left(z^{I}\overline{z}^{J}\right)=j_{k}\left(j_{k}-1\right) z_{1}^{i_{1}}\dots  z_{N}^{i_{N}} \overline{z}_{1}^{j_{1}}\dots \overline{z}_{k}^{j_{k}-2}\dots \overline{z}_{N}^{j_{N}}, 
\\&   \frac{\partial^{2}}{\partial z_{k}\partial \overline{z}_{k}}\left(z^{I}\overline{z}^{J}\right)= i_{k}j_{k} z_{1}^{i_{1}}\dots z_{k}^{i_{k}-1} \dots z_{N}^{i_{N}} \overline{z}_{1}^{j_{1}}\dots \overline{z}_{k}^{j_{k}-1}\dots \overline{z}_{N}^{j_{N}}, \end{split}  \label{909se}
\end{equation}
for all $k=1,\dots,N$, where $I=\left(i_{1},\dots,i_{N}\right)$, $J=\left(j_{1},\dots,j_{N}\right)\in\mathbb{N}^{N}$ are chosen such that $i_{1}+j_{1}+\dots+i_{N}+j_{N}=p-2$.   
 
Let's compute. We have 
 \begin{equation} \begin{split}& \displaystyle\sum_{k=1}^{N} z_{k}\frac{\partial}{\partial \overline{z}_{k}}\left(A\left(z,\overline{z}\right)\right)=\displaystyle\sum_{k=1}^{N}\displaystyle\sum_{I,J\in\mathbb{N}^{N}\atop\left|I\right|+\left|J\right|=p-2} a_{I;J} j_{k}  z_{1}^{i_{1}} \dots z_{k}^{i_{k}+1}\dots  z_{N}^{i_{N}} \overline{z}_{1}^{j_{1}}\dots \overline{z}_{k}^{j_{k}-1}\dots  \overline{z}_{N}^{j_{N}},\\& \displaystyle\sum_{k=1}^{N}z_{k} \frac{\partial}{\partial z_{k}}\left(A\left(z,\overline{z}\right)\right) =\displaystyle\sum_{k=1}^{N}\displaystyle\sum_{I,J\in\mathbb{N}^{N}\atop\left|I\right|+\left|J\right|=p-2} a_{I;J} i_{k}  z_{1}^{i_{1}} \dots  z_{N}^{i_{N}} \overline{z}_{1}^{j_{1}}\dots   \overline{z}_{N}^{j_{N}},\\&\displaystyle\sum_{k=1}^{N}\overline{z}_{k}\frac{\partial}{\partial \overline{z}_{k}}\left( A\left(z,\overline{z}\right)\right)=\displaystyle\sum_{k=1}^{N}\displaystyle\sum_{I,J\in\mathbb{N}^{N}\atop\left|I\right|+\left|J\right|=p-2} a_{I;J} j_{k}  z_{1}^{i_{1}} \dots  z_{N}^{i_{N}} \overline{z}_{1}^{j_{1}}\dots   \overline{z}_{N}^{j_{N}},\\&\displaystyle\sum_{k=1}^{N} \overline{z}_{k}\frac{\partial}{\partial z_{k}}\left( A\left(z,\overline{z}\right)\right)=\displaystyle\sum_{k=1}^{N}\displaystyle\sum_{I,J\in\mathbb{N}^{N}\atop\left|I\right|+\left|J\right|=p-2} a_{I;J} i_{k}  z_{1}^{i_{1}} \dots z_{1}^{i_{k}-1} \dots z_{N}^{i_{N}} \overline{z}_{1}^{j_{1}}\dots\overline{z}_{k}^{j_{k}+1}   \dots\overline{z}_{N}^{j_{N}}.\end{split} \label{B11}\end{equation} 
 
Let's make more computations using (\ref{model}), (\ref{oppv}) and (\ref{tracce}). We have  
\begin{equation} \begin{split}&\quad\hspace{0.18 cm}  \frac{\partial^{2}}{\partial z_{k}^{2}}\left(A(z,\overline{z}) Q(z,\overline{z})\right)=\frac{\partial^{2} }{\partial z_{k}\partial z_{k}}\left(A (z,\overline{z})\right) Q(z,\overline{z})+\frac{\partial  }{\partial z_{k} }\left(A (z,\overline{z})\right)\left(2\overline{z}_{k}+4\lambda_{k} z_{k}\right)+2\lambda_{k}A(z,\overline{z}),  \\&    \frac{\partial^{2}}{\partial z_{k}\partial \overline{z}_{k}}\left(A(z,\overline{z}) Q(z,\overline{z})\right)=  \frac{\partial^{2} }{\partial z_{k}\partial \overline{z}_{k}}\left(A (z,\overline{z})\right)  Q(z,\overline{z})+\frac{\partial  }{\partial \overline{z}_{k} }\left(A (z,\overline{z})\right)\left( \overline{z}_{k}+2\lambda_{k}z_{k}\right)+\frac{\partial  }{\partial z_{k} }\left(A (z,\overline{z})\right)\left( z_{k}+2\lambda_{k}\overline{z}_{k}\right)+A(z,\overline{z}),   \\& \quad\hspace{0.18 cm} \frac{\partial^{2}}{\partial \overline{z}_{k}^{2}}\left(A(z,\overline{z}) Q(z,\overline{z})\right)= \frac{\partial^{2} }{\partial \overline{z}_{k}\partial \overline{z}_{k}}\left(A (z,\overline{z})\right) Q(z,\overline{z})+\frac{\partial  }{\partial \overline{z}_{k} }\left(A (z,\overline{z})\right)\left(2z_{k}+4\lambda_{k} \overline{z}_{k}\right)+2\lambda_{k}A(z,\overline{z}), \end{split} \label{909se11} 
\end{equation}
for all $k=1,\dots,N$, concluding by (\ref{tracce}) and (\ref{op})  the following important aspect
\begin{equation}\begin{split}& \tilde{\tr}\left(A(z,\overline{z}) \right)Q(z,\overline{z})=\lambda_{1}\tilde{i}_{1}\left(\tilde{i}_{1}-1\right)z_{1}^{\tilde{i}_{1}-2}\dots z_{N}^{\tilde{i}_{N}}+\dots+\lambda_{k}\tilde{i}_{k}\left(\tilde{i}_{k}-1\right)z_{1}^{\tilde{i}_{1}}\dots z_{k}^{\tilde{i}_{k}-2}\dots z_{N}^{\tilde{i}_{N}}+\dots+\lambda_{N}\tilde{i}_{N}\left(\tilde{i}_{N}-1\right)z_{1}^{\tilde{i}_{1}}\dots z_{N}^{\tilde{i}_{N}-2}   \\& \quad\quad   \begin{tabular}{l} \rotatebox[origin=c]{270}{$+$}\end{tabular} \\&  \displaystyle\sum_{k=1}^{N}\frac{\partial  }{\partial z_{k} }\left(A (z,\overline{z})\right)\left(z_{k}\left(1+4\lambda_{k}^{2}\right)+4\lambda_{k}\overline{z}_{k} \right)+\displaystyle\sum_{k=1}^{N}\frac{\partial  }{\partial \overline{z}_{k} }\left(A (z,\overline{z})\right)\left(\overline{z}_{k}\left(1+4\lambda_{k}^{2}\right)+4\lambda_{k}z_{k}\right)+ A(z,\overline{z})\left(N+4\displaystyle\sum_{k=1}^{N}\lambda_{k}^{2}\right).\end{split}\label{VVV1}
\end{equation} 

On the other hand, we combine (\ref{yy}), (\ref{oppv}) and (\ref{909se}). We obtain 
 \begin{equation}\begin{split}&    \quad\quad\quad \quad\quad\quad\quad\quad\quad\quad\quad\quad\quad\quad\quad\quad\quad\quad\quad \quad\quad\quad\quad\tr \left(A(z,\overline{z})\right)Q(z,\overline{z})\\& \quad\quad\quad\quad \quad\quad\quad\quad\quad\quad\quad\quad  \quad\quad\quad\quad\quad\quad\quad \quad\quad\quad\quad\quad   \quad\quad\hspace{0.1 cm}  \begin{tabular}{l} \rotatebox[origin=c]{270}{$=$}\end{tabular} \\&    \quad\left(z_{1}\overline{z}_{1}+ \dots+z_{N}\overline{z}_{N}+\lambda_{1}\left(z_{1}^{2}+\overline{z}_{1}^{2}\right)+ \dots+\lambda_{N}\left(z_{N}^{2}+\overline{z}_{N}^{2}\right)\right)\cdot   \left(\displaystyle\sum_{k=1}^{N}\displaystyle\sum_{I,J\in\mathbb{N}^{N}\atop\left|I\right|+\left|J\right|=p-2} a_{I;J}i_{k}j_{k} z_{1}^{i_{1}}\dots z_{k}^{i_{k}-1} \dots z_{N}^{i_{N}} \overline{z}_{1}^{j_{1}}\dots \overline{z}_{k}^{j_{k}-1}\dots \overline{z}_{N}^{j_{N}} \right.\\&\quad\quad\quad\quad \quad\quad \quad\quad  \quad\quad  \quad\quad\quad\quad \quad\quad \quad\quad  \quad\quad \quad\quad  \quad\quad \quad\quad \hspace{0.1 cm}  + \\&     \left.\displaystyle\sum_{k=1}^{N}\displaystyle\sum_{I,J\in\mathbb{N}^{N}\atop\left|I\right|+\left|J\right|=p-2} a_{I;J}\lambda_{k} i_{k}\left(i_{k}-1\right) z_{1}^{i_{1}}\dots z_{k}^{i_{k}-2}\dots z_{N}^{i_{N}}\overline{z}_{1}^{j_{1}}\dots  \overline{z}_{N}^{j_{N}}+\displaystyle\sum_{k=1}^{N}\displaystyle\sum_{I,J\in\mathbb{N}^{N}\atop\left|I\right|+\left|J\right|=p-2} a_{I;J}\lambda_{k} j_{k}\left(j_{k}-1\right) z_{1}^{i_{1}}\dots  z_{N}^{i_{N}} \overline{z}_{1}^{j_{1}}\dots \overline{z}_{k}^{j_{k}-2}\dots \overline{z}_{N}^{j_{N}}\right) .\end{split}
\label{B1}\end{equation}

The interactions of homogeneous terms  are very complicated in (\ref{B1}). We have thus to better organize all these terms  depending on their contributions at each  degree level in (\ref{B1}), observing the following
 \begin{equation}  \begin{split}&    \left(z_{1}\overline{z}_{1}+\dots+z_{N}\overline{z}_{N}+\lambda_{1}\left(z_{1}^{2}+\overline{z}_{1}^{2}\right) +\dots+\lambda_{N}\left(z_{N}^{2}+\overline{z}_{N}^{2}\right)\right) \cdot   \left( \displaystyle\sum_{I,J\in\mathbb{N}^{N}\atop\left|I\right|+\left|J\right|=p-2} b_{I;J} z_{1}^{i_{1}}\dots z_{k}^{i_{k}-1} \dots z_{N}^{i_{N}} \overline{z}_{1}^{j_{1}}\dots \overline{z}_{k}^{j_{k}-1}\dots \overline{z}_{N}^{j_{N}}\right.\\&\quad\quad\quad\quad \quad\quad \quad\quad  \quad\quad\quad\quad  \quad\quad\quad\quad\quad\quad\quad\quad\hspace{0.2 cm}   + \\&  \left.\displaystyle\sum_{I,J\in\mathbb{N}^{N}\atop\left|I\right|+\left|J\right|=p-2} {b'}_{I;J} z_{1}^{i_{1}}\dots z_{k}^{i_{k}-2}\dots z_{N}^{i_{N}}\overline{z}_{1}^{j_{1}}\dots  \overline{z}_{N}^{j_{N}}+  \displaystyle\sum_{I,J\in\mathbb{N}^{N}\atop\left|I\right|+\left|J\right|=p-2} {b''}_{I;J} z_{1}^{i_{1}}\dots  z_{N}^{i_{N}} \overline{z}_{1}^{j_{1}}\dots \overline{z}_{k}^{j_{k}-2}\dots \overline{z}_{N}^{j_{N}}\right)=\tr \left(A(z,\overline{z})\right)Q(z,\overline{z}) ,\end{split}
\label{B111}\end{equation}
where we have used by (\ref{yy}) the following notations
\begin{equation} \begin{split}& {b}_{I;J}=\left(i_{1}+1\right)\left(j_{1}+1\right)a_{\left(i_{1}+1,\dots,i_{k}-1,\dots,i_{N};j_{1}+1,\dots,j_{k}-1,\dots,j_{N}\right)}+\dots+i_{k}j_{k}a_{\left(i_{1},\dots,i_{k},\dots,i_{N};j_{1},\dots,j_{k},\dots,j_{N}\right)}\\&\quad\quad\quad\quad\quad\quad\quad\quad\quad\quad\quad\quad\quad\quad\quad\quad\quad\quad\quad\quad\quad\quad\quad\quad\quad\quad\quad\quad\quad\quad\quad\quad\quad\quad\quad\quad\quad\quad\hspace{0.1 cm} +\\& \quad\quad\quad\quad\quad\quad\quad\quad\quad\quad\quad\quad\quad\quad\quad\quad\quad\quad\quad\quad\quad\quad\quad\quad\quad\quad\quad\quad\quad\quad\quad\quad\quad\quad\quad\quad\quad\quad\hspace{0.1 cm}\begin{tabular}{l} \rotatebox[origin=c]{270}{$\dots$}\end{tabular}\\&\quad\quad\quad\quad\quad\quad\quad\quad\quad\quad\quad\quad\quad\quad\quad\quad\quad\quad\quad\quad\quad\quad\quad\quad\quad\quad\quad\quad\quad\quad\quad\quad\quad\quad\quad\quad\quad\quad\hspace{0.1 cm} +\\&\quad\quad\quad\quad\quad\quad\quad\quad\quad\quad\quad\quad\quad\quad\quad\quad\quad\quad\quad\quad\quad\quad\quad\left(i_{N}+1\right)\left(j_{N}+1\right)a_{\left(i_{1},\dots,i_{k}-1,\dots,i_{N}+1;j_{1},\dots,j_{k}-1,\dots,j_{N}+1\right)},  \\&{b'}_{I;J}=\lambda_{1}  \left(i_{1}+1\right)\left(i_{1}+2\right) a_{\left(i_{1}+2,\dots,i_{k}-2,\dots,i_{N};j_{1},\dots,j_{N}\right)}+ \dots+\lambda_{k} i_{k}\left(i_{k}-1\right)a_{\left(i_{1},\dots,i_{k},\dots,i_{N};j_{1},\dots,j_{N}\right)}\\&\quad\quad\quad\quad\quad\quad\quad\quad\quad\quad\quad\quad\quad\quad\quad\quad\quad\quad\quad\quad\quad\quad\quad\quad\quad\quad\quad\quad\quad\quad\quad\quad\quad\quad\quad\quad\quad\quad\hspace{0.1 cm} +\\& \quad\quad\quad\quad\quad\quad\quad\quad\quad\quad\quad\quad\quad\quad\quad\quad\quad\quad\quad\quad\quad\quad\quad\quad\quad\quad\quad\quad\quad\quad\quad\quad\quad\quad\quad\quad\quad\quad\hspace{0.1 cm}\begin{tabular}{l} \rotatebox[origin=c]{270}{$\dots$}\end{tabular}\\&\quad\quad\quad\quad\quad\quad\quad\quad\quad\quad\quad\quad\quad\quad\quad\quad\quad\quad\quad\quad\quad\quad\quad\quad\quad\quad\quad\quad\quad\quad\quad\quad\quad\quad\quad\quad\quad\quad\hspace{0.1 cm} +\\&\quad\quad\quad\quad\quad\quad\quad\quad\quad\quad\quad\quad\quad\quad\quad\quad\quad\quad\quad\quad\quad\quad\quad\lambda_{N} \left(i_{N}+1\right)\left(i_{N}+2\right)a_{\left(i_{1},\dots,i_{k}-2,\dots,i_{N}+2;j_{1},\dots,j_{N}\right)},   \\&{b''}_{I;J}=\lambda_{1}\left(j_{1}+1\right)\left(j_{1}+2\right) a_{\left(i_{1},\dots,i_{N};j_{1}+2,\dots,j_{k}-2,\dots,j_{N}\right)}+ \dots+\lambda_{k}j_{k}\left(j_{k}-1\right) a_{\left(i_{1},\dots,i_{N};j_{1},\dots,j_{k}\dots,j_{N}\right)}\\&\quad\quad\quad\quad\quad\quad\quad\quad\quad\quad\quad\quad\quad\quad\quad\quad\quad\quad\quad\quad\quad\quad\quad\quad\quad\quad\quad\quad\quad\quad\quad\quad\quad\quad\quad\quad\quad\quad\hspace{0.1 cm} +\\& \quad\quad\quad\quad\quad\quad\quad\quad\quad\quad\quad\quad\quad\quad\quad\quad\quad\quad\quad\quad\quad\quad\quad\quad\quad\quad\quad\quad\quad\quad\quad\quad\quad\quad\quad\quad\quad\quad\hspace{0.1 cm}\begin{tabular}{l} \rotatebox[origin=c]{270}{$\dots$}\end{tabular}\\&\quad\quad\quad\quad\quad\quad\quad\quad\quad\quad\quad\quad\quad\quad\quad\quad\quad\quad\quad\quad\quad\quad\quad\quad\quad\quad\quad\quad\quad\quad\quad\quad\quad\quad\quad\quad\quad\quad\hspace{0.1 cm} +\\&\quad\quad\quad\quad\quad\quad\quad\quad\quad\quad\quad\quad\quad\quad\quad\quad\quad\quad\quad\quad\quad\quad\quad\lambda_{N}\left(j_{N}+1\right)\left(j_{N}+2\right)a_{\left(i_{1},\dots,i_{N};j_{1},\dots,j_{k}-2\dots,j_{N}+2\right)}.   \end{split} \label{99}
\end{equation}

In order to simplify the computations regarding (\ref{VVV1}), we introduce by (\ref{yy}) and (\ref{99})  the following vectors
\begin{equation} \begin{split}&
X_{1}=X_{1}\left[I;J\right]=\begin{pmatrix}a_{\left(i_{1}+1,\dots,i_{k}-1,\dots,i_{N};j_{1}+1,\dots,j_{k}-1,\dots,j_{N}\right)} \\ \vdots \\ a_{\left(i_{1},\dots,i_{k},\dots,i_{N};j_{1},\dots,j_{k},\dots,j_{N}\right)} \\ \vdots \\  a_{\left(i_{1},\dots,i_{k}-1,\dots,i_{N}+1;j_{1},\dots,j_{k}-1,\dots,j_{N}+1\right)},\end{pmatrix},\\&    X_{2}=X_{2}\left[I;J\right]=\begin{pmatrix}a_{\left(i_{1}+2,\dots,i_{k}-2,\dots,i_{N};j_{1},\dots,j_{N}\right)}\\ \vdots \\ a_{\left(i_{1},\dots,i_{k},\dots,i_{N};j_{1},\dots,j_{N}\right)} \\ \vdots \\  a_{\left(i_{1},\dots,i_{k}-2,\dots,i_{N}+2;j_{1},\dots,j_{N}\right)}\end{pmatrix},\\&  X_{3}=X_{3}\left[I;J\right]=\begin{pmatrix}a_{\left(i_{1},\dots,i_{N};j_{1}+2,\dots,j_{k}-2,\dots,j_{N}\right)} \\ \vdots \\ a_{\left(i_{1},\dots,i_{N};j_{1},\dots,j_{k}\dots,j_{N}\right)} \\ \vdots \\  a_{\left(i_{1},\dots,i_{N};j_{1},\dots,j_{k}-2\dots,j_{N}+2\right)}\label{vectori}\end{pmatrix}.
\end{split} \end{equation}

The first very consistent sum of homogeneous terms  is obviously multiplied in  (\ref{B111}) by
$$\mbox{$\lambda_{1}z_{1}^{2},\lambda_{2}z_{2}^{2},\dots,\lambda_{N}z_{N}^{2}$, $z_{1}\overline{z}_{1},z_{2}\overline{z}_{2},\dots,z_{N}\overline{z}_{N}$ and $\lambda_{1}\overline{z}_{1}^{2},\lambda_{2}\overline{z}_{2}^{2},\dots,\lambda_{N}\overline{z}_{N}^{2}$.}$$ 

This sum    generates obviously  by  (\ref{Lambda}), (\ref{99}) and (\ref{vectori}) in (\ref{VVV1})  the following terms
\begin{equation}\mathcal{A}X_{1},\quad  \left(\Lambda\mathcal{A}\right)X_{1},\quad \left(\Lambda\mathcal{A} \Lambda\right)X_{1}, 
  \label{idiot1}
\end{equation}
where we have used by (\ref{yy})  the following matrix
\begin{equation}  \mathcal{A}=\mathcal{A}\left[I;J\right]= 
\begin{pmatrix} \left(i_{1}+1\right)\left(j_{1}+1\right) & \dots & \left(i_{k-1}+1\right)\left(j_{k-1}+1\right)  & i_{k}j_{k} & \left(i_{k+1}+1\right)\left(j_{k+1}+1\right)  &\dots & \left(i_{N}+1\right)\left(j_{N}+1\right)  \\    \vdots & \ddots & \vdots & \vdots & \vdots & \ddots & \vdots    \\ \left(i_{1}+1\right)\left(j_{1}+1\right) & \dots & \left(i_{k-1}+1\right)\left(j_{k-1}+1\right)  & i_{k}j_{k} & \left(i_{k+1}+1\right)\left(j_{k+1}+1\right)  &\dots & \left(i_{N}+1\right)\left(j_{N}+1\right) \\ \vdots & \ddots & \vdots & \vdots & \vdots & \ddots & \vdots  \\ \left(i_{1}+1\right)\left(j_{1}+1\right) & \dots & \left(i_{k-1}+1\right)\left(j_{k-1}+1\right)  & i_{k}j_{k} & \left(i_{k+1}+1\right)\left(j_{k+1}+1\right)  &\dots & \left(i_{N}+1\right)\left(j_{N}+1\right)
\end{pmatrix}. \label{idiot1se}
\end{equation}

The second very consistent sum of homogeneous terms  is obviously multiplied (\ref{B111}) by $$\mbox{$\lambda_{1}z_{1}^{2},\lambda_{2}z_{2}^{2},\dots,\lambda_{N}z_{N}^{2}$, $z_{1}\overline{z}_{1},z_{2}\overline{z}_{2},\dots,z_{N}\overline{z}_{N}$ and $\lambda_{1}\overline{z}_{1}^{2},\lambda_{2}\overline{z}_{2}^{2},\dots,\lambda_{N}\overline{z}_{N}^{2}$.}$$  

This sum     generates obviously  by  (\ref{Lambda}), (\ref{99}) and (\ref{vectori})  in (\ref{VVV1}) the following terms
\begin{equation}
\left( \mathcal{A'}\Lambda\right)X_{2} ,\quad \left(\Lambda\mathcal{A'}\Lambda\right)X_{2} ,\quad \left(\Lambda\mathcal{A'}\Lambda\right)X_{2}, \label{idiot2}
\end{equation} 
where we have used by (\ref{yy})  the following matrix
\begin{equation}\mathcal{A'}=\mathcal{A'}\left[I;J\right]= 
\begin{pmatrix} \left(i_{1}+1\right)\left(i_{1}+2\right) & \dots   & \left(i_{k}-1\right)i_{k} & \dots & \left(i_{N}+1\right)\left(i_{N}+2\right) \\    \vdots & \ddots &   \vdots &   \ddots & \vdots    \\ \left(i_{1}+1\right)\left(i_{1}+2\right) & \dots & \left(i_{k}-1\right)i_{k}  & \dots & \left(i_{N}+1\right)\left(i_{N}+2\right) \\ \vdots & \ddots &   \vdots &  \ddots & \vdots  \\ \left(i_{1}+1\right)\left(i_{1}+2\right) & \dots &   \left(i_{k}-1\right)i_{k} &\dots & \left(i_{N}+1\right)\left(i_{N}+2\right)
\end{pmatrix}.\label{idiot2se} 
\end{equation}

The third very consistent sum of homogeneous terms  is obviously multiplied in  (\ref{B111}) by $$\mbox{$\lambda_{1}z_{1}^{2},\lambda_{2}z_{2}^{2},\dots,\lambda_{N}z_{N}^{2}$, $z_{1}\overline{z}_{1},z_{2}\overline{z}_{2},\dots,z_{N}\overline{z}_{N}$ and $\lambda_{1}\overline{z}_{1}^{2},\lambda_{2}\overline{z}_{2}^{2},\dots,\lambda_{N}\overline{z}_{N}^{2}$.}$$ 

 This sum     generates obviously   by  (\ref{Lambda}), (\ref{99}) and (\ref{vectori}) in (\ref{VVV1}) the following terms
\begin{equation}
\left(\ \mathcal{A''}\Lambda\right)X_{3} ,\quad \left(\Lambda\mathcal{A''}\Lambda\right)X_{3} ,\quad \left(\Lambda\mathcal{A''}\Lambda\right)X_{3}, \label{idiot3}
\end{equation}
where we have used by (\ref{yy})  the following matrix
\begin{equation} \mathcal{A''}=\mathcal{A''}\left[I;J\right]= 
\begin{pmatrix} \left(j_{1}+1\right)\left(j_{1}+2\right) & \dots   & \left(j_{k}-1\right)j_{k}  & \dots & \left(j_{N}+1\right)\left(j_{N}+2\right) \\    \vdots & \ddots &   \vdots &   \ddots & \vdots    \\ \left(j_{1}+1\right)\left(j_{1}+2\right) & \dots &  \left(j_{k}-1\right)j_{k}   &\dots & \left(j_{N}+1\right)\left(j_{N}+2\right) \\ \vdots & \ddots &   \vdots &  \ddots & \vdots  \\ \left(j_{1}+1\right)\left(j_{1}+2\right) & \dots &   \left(j_{k}-1\right)j_{k}  &  \dots & \left(j_{N}+1\right)\left(j_{N}+2\right)
\end{pmatrix}. \label{idiot3se}
\end{equation}

 The role of the  matrix $\Lambda$  is crucial   in making convenient estimates of the radius of convergence from the local defining equation (\ref{ecuatiegenerala}). It is obvious that   $\Lambda$ does not generally commute with any of the matrices $\mathcal{A}$, $\mathcal{A'}$,  $\mathcal{A''}$ from (\ref{idiot1}), (\ref{idiot2}) and (\ref{idiot3}),  but (\ref{idiot1}), (\ref{idiot2}) and (\ref{idiot3}) provide a better  understanding of the very complicated interactions of of homogeneous terms in (\ref{VVV1}) as follows.

\subsection{Systems of Equations} It is imposed  the standard lexicografic order corresponding to \begin{equation}\left(z_{1},z_{2},\dots,z_{N}; \overline{z}_{1},\overline{z}_{2},\dots,\overline{z}_{N}\right)\in\mathbb{C}^{N}\times\mathbb{C}^{N},
\label{ODIN}
\end{equation}   considering     the following vectors   
\begin{equation}\begin{split}& Y_{1}^{t}=\left\{\left(a_{I;0}\right)\right\}_{I\in\mathbb{N}^{N}\atop \left|I\right|=p},\\& Y_{2}^{t}=\left\{\left(a_{I;J}\right)\right\}_{I,J\in\mathbb{N}^{N}\atop \left|I\right|=p-1, \left|J\right|=1},\\&\quad\quad\quad\quad \vdots   \\& Y_{k}^{t}=\left\{\left(a_{I;J}\right)\right\}_{I,J\in\mathbb{N}^{N}\atop \left|I\right|=p-k+1, \left|J\right|=k-1}\\& \quad\quad\quad\quad\vdots\\&  Y_{p+1}^{t}=\left\{\left(a_{0;J}\right)\right\}_{J\in\mathbb{N}^{N}\atop   \left|J\right|=p},\end{split}
\label{91A}
\end{equation}
where  $p=\tilde{i}_{1}+\tilde{i}_{2}+\dots+\tilde{i}_{N}$, in order to  
 construct   systems of equations,   extracting homogeneous terms in (\ref{VVV1})  recalling (\ref{op}),   using
 \begin{equation}  V_{1} =\begin{pmatrix}   0\\  \vdots \\0 \\   \lambda_{1}\tilde{i}_{1}\left(\tilde{i}_{1}-1\right)  \\0\\ \dots \\0 \\ \lambda_{k}\tilde{i}_{k}\left(\tilde{i}_{k}-1\right) \\0\\ \vdots \\ \lambda_{N}\tilde{i}_{N}\left(\tilde{i}_{N}-1\right) \\0\\ \vdots  \\ 0 
\end{pmatrix}
,\quad V_{0} =\begin{pmatrix}  0\\ \vdots \\0\\0\\0\\  \vdots\\  0\\0\\0\\  \vdots \\0\\0\\0\\  \vdots  
\\ 0\end{pmatrix}. \label{special1}
\end{equation}

We obtain    the following non-trivial system of equations
 \begin{equation}\begin{pmatrix} \mathcal{M}_{1,1} & \mathcal{M}_{1,2} & \mathcal{M}_{1,3} & \mbox{O}_{N^{p}} & \mbox{O}_{N^{p}} &\dots & \mbox{O}_{N^{p}} & \mbox{O}_{N^{p}} & \mbox{O}_{N^{p}} & \mbox{O}_{N^{p}} \\ \mathcal{M}_{2,1} & \mathcal{M}_{2,2} & \mathcal{M}_{2,3} & \mathcal{M}_{2,4} & \mbox{O}_{N^{p}} &\dots & \mbox{O}_{N^{p}} & \mbox{O}_{N^{p}} & \mbox{O}_{N^{p}} & \mbox{O}_{N^{p}} \\ \mathcal{M}_{3,1}  & \mathcal{M}_{3,2}  & \mathcal{M}_{3,3}  & \mathcal{M}_{3,4}  & \mathcal{M}_{3,5}  &\dots & \mbox{O}_{N^{p}} & \mbox{O}_{N^{p}} & \mbox{O}_{N^{p}} & \mbox{O}_{N^{p}} \\ \mbox{O}_{N^{p}} & \mathcal{M}_{4,2}  & \mathcal{M}_{4,3}  & \mathcal{M}_{4,4}  & \mathcal{M}_{4,5}  &\dots & \mbox{O}_{N^{p}} & \mbox{O}_{N^{p}} & \mbox{O}_{N^{p}}& \mbox{O}_{N^{p}} \\ \mbox{O}_{N^{p}} & \mbox{O}_{N^{p}} & \mathcal{M}_{5,3} & \mathcal{M}_{5,4}& \mathcal{M}_{5,5} &\dots & \mbox{O}_{N^{p}} &\mbox{O}_{N^{p}} & \mbox{O}_{N^{p}} & \mbox{O}_{N^{p}} \cr \mbox{O}_{N^{p}} & \mbox{O}_{N^{p}} & \mbox{O}_{N^{p}} & \mathcal{M}_{6,4} & \mathcal{M}_{6,5} &\dots & \mbox{O}_{N^{p}} & \mbox{O}_{N^{p}} & \mbox{O}_{N^{p}} & \mbox{O}_{N^{p}}  \\ \vdots & \vdots & \vdots & \vdots & \vdots & \ddots & \vdots & \vdots & \vdots & \vdots \\ \mbox{O}_{N^{p}} & \mbox{O}_{N^{p}} & \mbox{O}_{N^{p}} & \mbox{O}_{N^{p}} & \mbox{O}_{N^{p}} &\dots & \mathcal{M}_{p-2,p-2}  & \mathcal{M}_{p-2,p-1}  & \mathcal{M}_{p-2,p} & \mbox{O}_{N^{p}} \\ \mbox{O}_{N^{p}} & \mbox{O}_{N^{p}} & \mbox{O}_{N^{p}} & \mbox{O}_{N^{p}} & \mbox{O}_{N^{p}} &\dots & \mathcal{M}_{p-1,p-2}  & \mathcal{M}_{p-1,p-1}  &\mathcal{M}_{p-1,p}  & \mathcal{M}_{p-1,p+1}  \\ \mbox{O}_{N^{p}} & \mbox{O}_{N^{p}} & \mbox{O}_{N^{p}} & \mbox{O}_{N^{p}} & \mbox{O}_{N^{p}} &\dots & \mathcal{M}_{p,p-2}  & \mathcal{M}_{p,p-1}  & \mathcal{M}_{p,p}  & \mathcal{M}_{p,p+1}  \\ \mbox{O}_{N^{p}} & \mbox{O}_{N^{p}} & \mbox{O}_{N^{p}} & \mbox{O}_{N^{p}} &\mbox{O}_{N^{p}} &\dots & \mbox{O}_{N^{p}} & \mathcal{M}_{p+1,p-1}  & \mathcal{M}_{p+1,p}  & \mathcal{M}_{p+1,p+1}           
\end{pmatrix}\begin{pmatrix} Y_{1}\\Y_{2}\\Y_{3}\\Y_{4}\\Y_{5}\\ \vdots\\ Y_{p-3}\\ Y_{p-2} \\ Y_{p-1}\\ Y_{p}\\ Y_{p+1}
\end{pmatrix}=\begin{pmatrix} V_{1} \\ V_{0} \\ V_{0} \\ V_{0}\\ V_{0}  \\   \vdots \\ V_{0} \\ V_{0}\\ V_{0}\\ V_{0}\\ V_{0}
\end{pmatrix}, \label{beb1}
\end{equation}
where its   elements are matrices, which  are defined as follows
\begin{equation} \begin{split}& \mathcal{M}_{k-2,k}=\mathcal{W}_{k}^{'},\quad\quad\quad\quad  \quad\quad\quad\quad\quad\quad\quad\quad\quad\quad\quad\quad\quad\quad\quad\quad\quad\quad\hspace{0.1 cm}  \mbox{for all $k=3,\dots,p+1$,} \\&\mathcal{M}_{k-1,k}=\mathcal{V}_{k}^{(0)}+\mathcal{W}_{k}+\mathcal{O}_{k}^{'} ,\quad\quad\quad  \quad\quad\quad\quad\quad\quad\quad\quad\quad\quad\quad\quad\quad\hspace{0.2 cm}\mbox{for all $k=2,\dots,p+1 $,}\\&\quad \mathcal{M}_{k,k}= \mathcal{O}_{k}^{(0)}+\mathcal{O}_{k}+\mathcal{V}_{k}^{'}+\mathcal{W}_{k}^{''}+\left(N+4\lambda_{1}^{2}+\dots+4\lambda_{N}^{2}\right)I_{N^{p}} ,\quad\mbox{for all $k=1,\dots,p+1$,}  \\& \mathcal{M}_{k+1,k}=\mathcal{W}_{k}^{(0)}+\mathcal{V}_{k}+\mathcal{O}_{k}^{''} ,\quad\quad\quad  \quad\quad\quad\quad\quad\quad\quad\quad\quad\quad\quad\quad\hspace{0.23 cm}\quad\mbox{for all $k=1,\dots,p$,} \\&  \mathcal{M}_{k+2,k}=\mathcal{V}_{k}^{''},\quad\quad\quad\quad  \quad\quad\quad\quad\quad\quad\quad\quad\quad\quad\quad\quad\quad\quad\quad\quad\quad\quad\hspace{0.26 cm}  \mbox{for all $k=1,\dots,p-1$.} \end{split}  \label{lil}
\end{equation} 

These occurring  matrices   are defined as follows:

$\bf{The\hspace{0.1 cm}matrices:}$
$$\bf{\left\{\mathcal{O}_{k}^{(0)}\right\}_{k=1,\dots,p+1}\hspace{0.1 cm}, \left\{\mathcal{V}_{k}^{(0)}\right\}_{k=2,\dots,p+1},\hspace{0.1 cm} \left\{\mathcal{W}_{k}^{(0)}\right\}_{k=1,\dots,p}.}$$

 Occurring from (\ref{VVV1}), we have used in (\ref{lil}) the following matrices 
\begin{equation}\mathcal{O}_{k}^{(0)}\left[I;J\right]=\begin{pmatrix} 1& \dots & 0 & \dots & 0 & \dots & 0\\ \vdots & \ddots& \vdots & \ddots & \vdots& \ddots & \vdots \\ 0&\dots & \left(i_{1}+ j_{1}+1\right)\left(1+4\lambda_{1}^{2}\right)    &\dots & 0   &   \dots & 0 &  \\ \vdots & \ddots& \vdots & \ddots & \vdots& \ddots & \vdots  \\ 0&\dots & 0& \dots &  \left(i_{N}+ j_{N}+1\right)\left(1+4\lambda_{N}^{2}\right)    & \dots  & 0  \\   \vdots & \ddots& \vdots & \ddots & \vdots& \ddots & \vdots\\ 0&\dots &0 &  \dots & 0 & \dots &  1  \end{pmatrix}\in \mathcal{M}_{N^{p}\times N^{p}}\left(\mathbb{C}\right), \label{lili1}
\end{equation} 
such that   
\begin{equation}j_{1}+j_{2}+\dots+j_{N}+i_{1}+i_{2}+\dots+i_{N}=p ,\quad j_{1}+j_{2}+\dots+j_{N}=k-1.\label{lili}
\end{equation}

Then, (\ref{lili}) defines in (\ref{lil}) also the following   matrices
\begin{equation} \mathcal{V}_{k}^{(0)}\left[I;J\right]=\begin{pmatrix} 1& \dots & 0 & \dots & 0 & \dots & 0\\ \vdots & \ddots& \vdots & \ddots & \vdots& \ddots & \vdots \\ 0&\dots & 4\lambda_{1} i_{1}    &\dots & 0   &   \dots & 0 &  \\ \vdots & \ddots& \vdots & \ddots & \vdots& \ddots & \vdots  \\ 0&\dots & 0& \dots & 4\lambda_{N} i_{N}     & \dots  & 0  \\   \vdots & \ddots& \vdots & \ddots & \vdots& \ddots & \vdots\\ 0&\dots &0 &  \dots & 0 & \dots &  1  \end{pmatrix}\in \mathcal{M}_{N^{p}\times N^{p}}\left(\mathbb{C}\right),\quad\mbox{for all $k=2,\dots,p+1 $,}\label{lili2} 
\end{equation} 
and respectively the following matrices
\begin{equation} \mathcal{W}_{k}^{(0)}\left[I;J\right]=\begin{pmatrix} 1& \dots & 0 & \dots & 0 & \dots & 0\\ \vdots & \ddots& \vdots & \ddots & \vdots& \ddots & \vdots \\ 0&\dots & 4\lambda_{1} j_{1}    &\dots & 0   &   \dots & 0 &  \\ \vdots & \ddots& \vdots & \ddots & \vdots& \ddots & \vdots  \\ 0&\dots & 0& \dots & 4\lambda_{N} j_{N}     & \dots  & 0  \\   \vdots & \ddots& \vdots & \ddots & \vdots& \ddots & \vdots\\ 0&\dots &0 &  \dots & 0 & \dots &  1  \end{pmatrix}\in \mathcal{M}_{N^{p}\times N^{p}}\left(\mathbb{C}\right),\quad\mbox{for all $k=1,\dots,p $,}\label{lili3} 
\end{equation}

Thus  in (\ref{lil}), the following matrices are by (\ref{lili1}), (\ref{lili2}) and (\ref{lili3})  defined as follows
\begin{equation} \begin{split}&\hspace{0.05 cm}
\mathcal{O}_{k}^{(0)}=\mathcal{O}_{k}^{(0)}\left[p-2,\dots,0;0,\dots,0\right]\dots\mathcal{O}_{k}^{(0)}\left[p-3,\dots,0;1,\dots,0\right]\dots\mathcal{O}_{k}^{(0)}\left[0,\dots,0;0,\dots,p-2\right],\quad\hspace{0.22 cm}\mbox{for all $k=1,\dots,p+1$,}  \\&\hspace{0.1 cm}\mathcal{V}_{k}^{(0)}=\mathcal{V}_{k}^{(0)}\left[p-2,\dots,0;0,\dots,0\right]\dots\mathcal{V}_{k}^{(0)}\left[p-3,\dots,0;1,\dots,0\right]\dots \mathcal{V}_{k}^{(0)}\left[0,\dots,0;0,\dots,p-2\right],\quad\quad \mbox{for all $k=2,\dots,p+1$,}  \\& \mathcal{W}_{k}^{(0)}=\mathcal{W}_{k}^{(0)}\left[p-2,\dots,0;0,\dots,0\right]\dots\mathcal{W}_{k}^{(0)}\left[p-3,\dots,0;1,\dots,0\right]\dots \mathcal{W}_{k}^{(0)}\left[0,\dots,0;0,\dots,p-2\right],\quad\mbox{for all $k=1,\dots,p$.} \end{split} \label{90000se1extra}
\end{equation}

$\bf{The\hspace{0.1 cm}matrices:}$
$$\bf{\left\{\mathcal{O}_{k}\right\}_{k=1,\dots,p+1},\hspace{0.1 cm}  \left\{\mathcal{V}_{k}\right\}_{k=1,\dots,p},\hspace{0.1 cm} \left\{\mathcal{W}_{k}\right\}_{k=2,\dots,p+1}.}$$

 The first   matrix from (\ref{idiot1}) induces by (\ref{yy}), (\ref{B1}),   (\ref{lili})  the following matrices  
\begin{equation}   \mathcal{O}_{k}\left[I;J\right]=
\begin{pmatrix}  1 & \dots   & 0 & \dots & 0 & \dots & 0& \dots & 0  \\ \vdots & \ddots & \vdots & \ddots & \vdots & \ddots & \vdots& \ddots & \vdots   \\ 0 & \dots & \left(i_{1}+1\right)\left(j_{1}+1\right) & \dots & i_{k}j_{k} & \dots & \left(i_{N}+1\right)\left(j_{N}+1\right) & \dots & 0 
\\ \vdots & \ddots & \vdots  & \ddots & \vdots  & \ddots & \vdots  & \ddots & \vdots  \\ 0 & \dots & \left(i_{1}+1\right)\left(j_{1}+1\right) & \dots & i_{k}j_{k} & \dots & \left(i_{N}+1\right)\left(j_{N}+1\right) & \dots & 0 \\ \vdots  & \ddots & \vdots  & \ddots & \vdots  & \ddots & \vdots  & \ddots & \vdots  \\ 0 & \dots & \left(i_{1}+1\right)\left(j_{1}+1\right) & \dots & i_{k}j_{k} & \dots & \left(i_{N}+1\right)\left(j_{N}+1\right) & \dots & 0 \\ \vdots  & \ddots & \vdots  & \ddots & \vdots  & \ddots & \vdots  & \ddots & \vdots  \\ 0  & \dots & 0  & \dots &0  & \dots & 0  & \dots & 1  
\end{pmatrix}\in \mathcal{M}_{N^{p}\times N^{p}}\left(\mathbb{C}\right),  \label{suc1}
\end{equation}
for all $k=1,\dots,p+1$.  

This matrix (\ref{suc1}) has the characteristic that $i_{k}j_{k}$    is  the diagonal entry on the following row $$\left(i_{1},\dots,i_{N};j_{1},\dots,j_{N}\right)\in\mathbb{N}^{N}\times\mathbb{N}^{N},$$  according to the    lexicografic order  related to  (\ref{ODIN}), otherwise having   $1$ as diagonal entry.

Similarly, we consider by (\ref{yy})  and (\ref{idiot1}) the following matrices
 \begin{equation}   \mathcal{V}_{k}\left[I;J\right]=
\begin{pmatrix}  1 & \dots   & 0 & \dots & 0 & \dots & 0& \dots & 0  \\ \vdots & \ddots & \vdots & \ddots & \vdots & \ddots & \vdots& \ddots & \vdots   \\ 0 & \dots & \left(i_{1}+1\right)\left(j_{1}+1\right)
 \lambda_{1} & \dots & i_{k}j_{k}
 \lambda_{k} & \dots & \left(i_{N}+1\right)\left(j_{N}+1\right)
 \lambda_{N} & \dots & 0 
\\ \vdots & \ddots & \vdots  & \ddots & \vdots  & \ddots & \vdots  & \ddots & \vdots  \\ 0 & \dots & \left(i_{1}+1\right)\left(j_{1}+1\right)
 \lambda_{1} & \dots & i_{k}j_{k}
 \lambda_{k} & \dots & \left(i_{N}+1\right)\left(j_{N}+1\right)
 \lambda_{N} & \dots & 0 \\ \vdots  & \ddots & \vdots  & \ddots & \vdots  & \ddots & \vdots  & \ddots & \vdots  \\ 0 & \dots & \left(i_{1}+1\right)\left(j_{1}+1\right)
 \lambda_{1} & \dots & i_{k}j_{k}
 \lambda_{k} & \dots & \left(i_{N}+1\right)\left(j_{N}+1\right)
 \lambda_{N} & \dots & 0 \\ \vdots  & \ddots & \vdots  & \ddots & \vdots  & \ddots & \vdots  & \ddots & \vdots  \\ 0  & \dots & 0  & \dots &0  & \dots & 0  & \dots & 1  
\end{pmatrix}\in \mathcal{M}_{N^{p}\times N^{p}}\left(\mathbb{C}\right),  \label{suc2}
\end{equation}
for all $k=1,\dots,p$.  

This matrix  (\ref{suc2}) has the characteristic that $i_{k}j_{k}\lambda_{k}$    is  the diagonal entry  on the  following row $$ \left(i_{1},\dots,i_{k}+1,\dots,i_{N};j_{1},\dots,j_{k}-1,\dots,j_{N}\right)\in\mathbb{N}^{N}\times\mathbb{N}^{N},$$  according to the corresponding  lexicografic order  related to  (\ref{ODIN}), otherwise having   $1$ as diagonal entry. 

Then (\ref{suc2}) induces by (\ref{yy})  and (\ref{idiot1}) analogously another matrices denoted  as follows
\begin{equation}\mathcal{W}_{k}\left[I;J\right]\in \mathcal{M}_{N^{p}\times N^{p}}\left(\mathbb{C}\right),\quad\mbox{for all $k=2,\dots,p+1$,}\label{suc33}\end{equation} having the characteristic that $i_{k}j_{k}\lambda_{k}$    is  the diagonal entry  on the  following row   $$ \left(i_{1},\dots,i_{k}-1,\dots,i_{N};j_{1},\dots,j_{k}+1,\dots,j_{N}\right)\in\mathbb{N}^{N}\times\mathbb{N}^{N},$$ according to the    lexicografic order  related to  (\ref{ODIN}), otherwise having   $1$ as diagonal entry.  

 Thus in (\ref{lil}),  the following matrices are by (\ref{suc1}), (\ref{suc2})  and (\ref{suc33}) naturally defined  as follows
\begin{equation} \begin{split}&\hspace{0.05 cm}
\mathcal{O}_{k}=\mathcal{O}_{k}\left[p-2,\dots,0;0,\dots,0\right]\dots\mathcal{O}_{k}\left[p-3,\dots,0;1,\dots,0\right]\dots \mathcal{O}_{k}\left[0,\dots,0;0,\dots,p-2\right], \hspace{0.16 cm}\quad\mbox{for all $k=1,\dots,p+1$,} \\& \hspace{0.1 cm}\mathcal{V}_{k}=\mathcal{V}_{k}\left[p-2,\dots,0;0,\dots,0\right]\dots\mathcal{V}_{k}\left[p-3,\dots,0;1,\dots,0\right]\dots \mathcal{V}_{k}\left[0,\dots,0;0,\dots,p-2\right], \quad\quad\mbox{for all $k=1,\dots,p$,}  \\& \mathcal{W}_{k}=\mathcal{W}_{k}\left[p-2,\dots,0;0,\dots,0\right]\dots\mathcal{W}_{k}\left[p-3,\dots,0;1,\dots,0\right]\dots \mathcal{W}_{k}\left[0,\dots,0;0,\dots,p-2\right],\quad\mbox{for all $k=2,\dots,p+1$.}\end{split} \label{90000}
\end{equation}

It is important to better explain these products from (\ref{90000}). These matrices are commuting, because there are no commune non-trivial rows. Thus, (\ref{90000}) is clear because   (\ref{suc1}), (\ref{suc2}) and (\ref{suc33}) define classes of commuting matrices.

$\bf{The\hspace{0.1 cm}matrices:}$
$$\bf{\left\{\mathcal{O}_{k}^{'}\right\}_{k=2,\dots,p+1},\hspace{0.1 cm}  \left\{\mathcal{V}_{k}^{'}\right\}_{k=1,\dots,p+1},\hspace{0.1 cm} \left\{\mathcal{W}_{k}^{'}\right\}_{k=3,\dots,p+1}.}$$ 

The first  matrix from (\ref{idiot2}) induces   by   (\ref{yy}),   (\ref{B1}) and (\ref{lili})    the following matrices 
 \begin{equation}      \mathcal{O'}_{k}\left[I;J\right]= \begin{pmatrix}  1 & \dots   & 0 & \dots & 0 & \dots & 0& \dots & 0  \\ \vdots & \ddots & \vdots & \ddots & \vdots & \ddots & \vdots& \ddots & \vdots   \\ 0 & \dots & \left(i_{1}+1\right)\left(i_{1}+2\right)\lambda_{1}  & \dots & \left(i_{k}-1\right)i_{k}\lambda_{k} & \dots & \left(i_{N}+1\right)\left(i_{N}+2\right)\lambda_{N} & \dots & 0 
\\ \vdots & \ddots & \vdots  & \ddots & \vdots  & \ddots & \vdots  & \ddots & \vdots  \\ 0 & \dots & \left(i_{1}+1\right)\left(i_{1}+2\right)\lambda_{1} & \dots & \left(i_{k}-1\right)i_{k}\lambda_{k} & \dots & \left(i_{N}+1\right)\left(i_{N}+2\right)\lambda_{N} & \dots & 0 \\ \vdots  & \ddots & \vdots  & \ddots & \vdots  & \ddots & \vdots  & \ddots & \vdots  \\ 0 & \dots & \left(i_{1}+1\right)\left(i_{1}+2\right)\lambda_{1}  & \dots & \left(i_{k}-1\right)i_{k}\lambda_{k} & \dots & \left(i_{N}+1\right)\left(i_{N}+2\right)\lambda_{N} & \dots & 0 \\ \vdots  & \ddots & \vdots  & \ddots & \vdots  & \ddots & \vdots  & \ddots & \vdots  \\ 0  & \dots & 0  & \dots &0  & \dots & 0  & \dots & 1  
\end{pmatrix}\in \mathcal{M}_{N^{p}\times N^{p}}\left(\mathbb{C}\right),\label{sucL1} 
\end{equation} 
for all  $k=2,\dots,p+1$.  

This matrix   (\ref{sucL1}) has the characteristic that  $\left(i_{k}-1\right)i_{k}\lambda_{k}$ is  the diagonal entry  of the following  row  $$\left(i_{1},\dots,i_{k}-1,\dots,i_{N};j_{1},\dots,j_{k}+1,\dots,j_{N}\right)\in\mathbb{N}^{N}\times\mathbb{N}^{N},$$  according to the  lexicografic order  related to  (\ref{ODIN}), otherwise having   $1$ as diagonal entry. 

Then, the second matrix from (\ref{idiot2}) induces   by (\ref{yy})   the following matrices
  \begin{equation}     \mathcal{V'}_{k}\left[I;J\right]= \begin{pmatrix}  1 & \dots   & 0 & \dots & 0 & \dots & 0& \dots & 0  \\ \vdots & \ddots & \vdots & \ddots & \vdots & \ddots & \vdots& \ddots & \vdots   \\ 0 & \dots & \left(i_{1}+1\right)\left(i_{1}+2\right)\lambda_{1}^{2} & \dots & \left(i_{k}-1\right)i_{k}\lambda_{1} \lambda_{k}   & \dots & \left(i_{N}+1\right)\left(i_{N}+2\right)\lambda_{1}\lambda_{N}   & \dots & 0 
\\ \vdots & \ddots & \vdots  & \ddots & \vdots  & \ddots & \vdots  & \ddots & \vdots  \\ 0 & \dots & \left(i_{1}+1\right)\left(i_{1}+2\right)\lambda_{1}\lambda_{k}   & \dots & \left(i_{k}-1\right)i_{k}\lambda_{k}^{2} & \dots & \left(i_{N}+1\right)\left(i_{N}+2\right)\lambda_{N}\lambda_{k}  & \dots & 0 \\ \vdots  & \ddots & \vdots  & \ddots & \vdots  & \ddots & \vdots  & \ddots & \vdots  \\ 0 & \dots & \left(i_{1}+1\right)\left(i_{1}+2\right)\lambda_{1}\lambda_{N}   & \dots & \left(i_{k}-1\right)i_{k}\lambda_{k}\lambda_{N}  & \dots & \left(i_{N}+1\right)\left(i_{N}+2\right)\lambda_{N}^{2} & \dots & 0 \\ \vdots  & \ddots & \vdots  & \ddots & \vdots  & \ddots & \vdots  & \ddots & \vdots  \\ 0  & \dots & 0  & \dots &0  & \dots & 0  & \dots & 1  
\end{pmatrix}\in \mathcal{M}_{N^{p}\times N^{p}}\left(\mathbb{C}\right),\label{sucL11} 
\end{equation} 
for all  $k=1,\dots,p+1$.  

This matrix (\ref{sucL11}) has the characteristic that  $\left(i_{k}-1\right)i_{k}\lambda_{k}^{2}$ is  the diagonal entry  on the  following row   $$\left(i_{1}, \dots,i_{N};j_{1},\dots,j_{N}\right)\in\mathbb{N}^{N}\times\mathbb{N}^{N},$$  according to the   lexicografic order  related to  (\ref{ODIN}), otherwise having   $1$ as diagonal entry.  

Also (\ref{sucL11}) induces by (\ref{yy})  and (\ref{idiot2}) another matrices denoted   as follows 
\begin{equation}\mathcal{W'}_{k}\left[I;J\right]\in \mathcal{M}_{N^{p}\times N^{p}}\left(\mathbb{C}\right),\quad\mbox{for all $k=3,\dots,p+1$.}\label{sucL111}\end{equation}

 This matrix (\ref{sucL111}) has the characteristic that $\left(i_{k}-1\right) i_{k}  \lambda_{k}^{2}$ is  the diagonal entry  on the  following row  $$\left(i_{1},\dots,i_{k}-2,\dots,i_{N};j_{1},\dots,j_{k}+2, \dots,j_{N}\right)\in\mathbb{N}^{N}\times\mathbb{N}^{N},$$  according to the    lexicografic order  related to  (\ref{ODIN}), otherwise having   $1$ as diagonal entry. 
 
Thus in (\ref{lil}), the following matrices are by (\ref{sucL1}), (\ref{sucL11}) and (\ref{sucL111}) naturally defined as follows
\begin{equation} \begin{split}&\hspace{0.05 cm}
\mathcal{O'}_{k}=\mathcal{O'}_{k}\left[p-2,\dots,0;0,\dots,0\right]\dots\mathcal{O'}_{k}\left[p-3,\dots,0;1,\dots,0\right]\dots\mathcal{O'}_{k}\left[0,\dots,0;0,\dots,p-2\right], \quad\hspace{0.19 cm}\mbox{for all $k=2,\dots,p+1$,} \\& \hspace{0.1 cm}\mathcal{V'}_{k}=\mathcal{V'}_{k}\left[p-2,\dots,0;0,\dots,0\right]\dots\mathcal{V'}_{k}\left[p-3,\dots,0;1,\dots,0\right]\dots\mathcal{V'}_{k}\left[0,\dots,0;0,\dots,p-2\right],\quad\quad\mbox{for all $k=1,\dots,p+1$,}  \\& \mathcal{W'}_{k}=\mathcal{W'}_{k}\left[p-2,\dots,0;0,\dots,0\right]\dots\mathcal{W'}_{k}\left[p-3,\dots,0;1,\dots,0\right]\dots\mathcal{W'}_{k}\left[0,\dots,0;0,\dots,p-2\right],\quad\mbox{for all $k=3,\dots,p+1$.}\end{split} \label{90000se}
\end{equation}

$\bf{The\hspace{0.1 cm}matrices:}$
$$ \bf{\left\{\mathcal{O}_{k}^{''}\right\}_{k=1,\dots,p},\hspace{0.1 cm}  \left\{\mathcal{V}_{k}^{''}\right\}_{k=1,\dots,p-1},\hspace{0.1 cm} \left\{\mathcal{W}_{k}^{''}\right\}_{k=1,\dots,p+1}.} $$

The first  matrix from (\ref{idiot3}) induces by    (\ref{yy}), (\ref{B1}) and  (\ref{lili}) the following matrices 
 \begin{equation}  \mathcal{O''}_{k}\left[I;J\right]= \begin{pmatrix}  1 & \dots   & 0 & \dots & 0 & \dots & 0& \dots & 0  \\ \vdots & \ddots & \vdots & \ddots & \vdots & \ddots & \vdots& \ddots & \vdots   \\ 0 & \dots & \left(j_{1}+1\right)\left(j_{1}+2\right)\lambda_{1}  & \dots &  \left(j_{k}-1\right)j_{k}\lambda_{k} & \dots & \left(j_{N}+1\right)\left(j_{N}+2\right)\lambda_{N} & \dots & 0 
\\ \vdots & \ddots & \vdots  & \ddots & \vdots  & \ddots & \vdots  & \ddots & \vdots  \\ 0 & \dots & \left(j_{1}+1\right)\left(j_{1}+2\right)\lambda_{1} & \dots &  \left(j_{k}-1\right)j_{k}\lambda_{k} & \dots & \left(j_{N}+1\right)\left(j_{N}+2\right)\lambda_{N} & \dots & 0 \\ \vdots  & \ddots & \vdots  & \ddots & \vdots  & \ddots & \vdots  & \ddots & \vdots  \\ 0 & \dots & \left(j_{1}+1\right)\left(j_{1}+2\right)\lambda_{1}  & \dots &  \left(j_{k}-1\right)j_{k}\lambda_{k} & \dots & \left(j_{N}+1\right)\left(j_{N}+2\right)\lambda_{N} & \dots & 0 \\ \vdots  & \ddots & \vdots  & \ddots & \vdots  & \ddots & \vdots  & \ddots & \vdots  \\ 0  & \dots & 0  & \dots &0  & \dots & 0  & \dots & 1  
\end{pmatrix}\in \mathcal{M}_{N^{p}\times N^{p}}\left(\mathbb{C}\right), \label{sucL1sese}
\end{equation}
for all $k=1,\dots,p$.  
 
This matrix   (\ref{sucL1sese}) has the characteristic that $\left(j_{k}-1\right)j_{k}\lambda_{k}$    is  the diagonal entry  on the  following row $$\left(i_{1},\dots,i_{k}+1,\dots,i_{N};j_{1},\dots,j_{k}-1, \dots,j_{N}\right)\in\mathbb{N}^{N}\times\mathbb{N}^{N},$$    according to the  lexicografic order  related to  (\ref{ODIN}), otherwise having   $1$ as diagonal entry. 

Similarly as previously, we consider by (\ref{yy})   the following matrices
 \begin{equation} \mathcal{V''}_{k}\left[I;J\right]= \begin{pmatrix}  1 & \dots   & 0 & \dots & 0 & \dots & 0& \dots & 0  \\ \vdots & \ddots & \vdots & \ddots & \vdots & \ddots & \vdots& \ddots & \vdots   \\ 0 & \dots & \left(j_{1}+1\right)\left(j_{1}+2\right)\lambda_{1}^{2} & \dots & \left(j_{k}-1\right)j_{k}\lambda_{1}\lambda_{k}  & \dots & \left(j_{N}+1\right)\left(j_{N}+2\right)\lambda_{1}\lambda_{N}  & \dots & 0 
\\ \vdots & \ddots & \vdots  & \ddots & \vdots  & \ddots & \vdots  & \ddots & \vdots  \\ 0 & \dots & \left(j_{1}+1\right)\left(j_{1}+2\right)\lambda_{1}\lambda_{k}   & \dots & \left(j_{k}-1\right)j_{k}\lambda_{k}^{2} & \dots & \left(j_{N}+1\right)\left(j_{N}+2\right)\lambda_{k} \lambda_{N} & \dots & 0 \\ \vdots  & \ddots & \vdots  & \ddots & \vdots  & \ddots & \vdots  & \ddots & \vdots  \\ 0 & \dots & \left(j_{1}+1\right)\left(j_{1}+2\right)\lambda_{1}\lambda_{N}   & \dots & \left(j_{k}-1\right)j_{k}\lambda_{k}\lambda_{N}   & \dots & \left(j_{N}+1\right)\left(j_{N}+2\right)\lambda_{N}^{2} & \dots & 0 \\ \vdots  & \ddots & \vdots  & \ddots & \vdots  & \ddots & \vdots  & \ddots & \vdots  \\ 0  & \dots & 0  & \dots &0  & \dots & 0  & \dots & 1  
\end{pmatrix}\in \mathcal{M}_{N^{p}\times N^{p}}\left(\mathbb{C}\right) , \label{sucL11sese}
\end{equation}
for all $k=1,\dots,p-1$.  

This matrix   (\ref{sucL11sese}) has the characteristic that $\left(j_{k}-1\right)j_{k}\lambda_{k}^{2}$    is  the diagonal entry  on the  following row   
$$\left(i_{1},\dots,i_{k}+2,\dots,i_{N};j_{1},\dots,j_{k}-2, \dots,j_{N}\right)\in\mathbb{N}^{N}\times\mathbb{N}^{N},$$  
 according to the   lexicografic order related to (\ref{ODIN}), otherwise having   $1$ as diagonal entry.
 
Then (\ref{sucL11sese}) induces by (\ref{yy})   another matrices denoted   as follows 
\begin{equation}\mathcal{W''}_{k}\left[I;J\right]\in\mathcal{M}_{N^{p}\times N^{p}}\left(\mathbb{C}\right) ,\quad\mbox{ for all $k=1,\dots,p+1$}.\label{sucL111se}\end{equation}  

This matrix (\ref{sucL111se}) has the characteristic that $\left(j_{k}-1\right)j_{k} \lambda_{k}^{2}$   is  the diagonal entry  on the  following row 
$$\left(i_{1}, i_{2},\dots,i_{N};j_{1}, j_{2}, \dots,j_{N}\right)\in\mathbb{N}^{N}\times\mathbb{N}^{N},$$ 
according to the   lexicografic order  related to  (\ref{ODIN}), otherwise having  $1$ as diagonal entry.

Thus, in (\ref{lil}) the following matrices are by (\ref{suc1}), (\ref{suc2}) and (\ref{sucL111se})  defined as follows
\begin{equation} \begin{split}&\hspace{0.05 cm}
\mathcal{O''}_{k}=\mathcal{O''}_{k}\left[p-2,\dots,0;0,\dots,0\right]\dots\mathcal{O''}_{k}\left[p-3,\dots,0;1,\dots,0\right]\dots\mathcal{O''}_{k}\left[0,\dots,0;0,\dots,p-2\right],\quad\hspace{0.1 cm} \hspace{0.1 cm}\mbox{for all $k=1,\dots,p$,}  \\&\hspace{0.1 cm} \mathcal{V''}_{k}=\mathcal{V''}_{k}\left[p-2,\dots,0;0,\dots,0\right]\dots\mathcal{V''}_{k}\left[p-3,\dots,0;1,\dots,0\right]\dots \mathcal{V''}_{k}\left[0,\dots,0;0,\dots,p-2\right],\quad\quad \mbox{for all $k=1,\dots,p-1$,}  \\& \mathcal{W''}_{k}=\mathcal{W''}_{k}\left[p-2,\dots,0;0,\dots,0\right]\dots\mathcal{W''}_{k}\left[p-3,\dots,0;1,\dots,0\right]\dots \mathcal{W''}_{k}\left[0,\dots,0;0,\dots,p-2\right],\quad\mbox{for all $k=1,\dots,p+1$.}\end{split} \label{90000se1}
\end{equation}
 
We were using the following obvious observation 
\begin{equation}   \#\left\{\left(i_{1},i_{2},\dots,i_{N}:j_{1},j_{2},\dots,j_{N}\right)\in \mathbb{N}^{N}\times \mathbb{N}^{N};\quad  \mbox{   $i_{1}+i_{2}+\dots+i_{N}=p-k$, $j_{1}+j_{2}+\dots+j_{N}=k$}\right\}=N^{p},\hspace{0.1 cm} \mbox{for all $k=0,\dots,p$.}\label{set} \end{equation}
\subsection{Remarks} It is crucial    to prove  the    invertibility  of each of the following matrices
\begin{equation}\mathcal{M}_{k,k}=\mathcal{O}_{k}^{(0)}+\mathcal{O}_{k}+\mathcal{V}_{k}^{'}+\mathcal{W}_{k}^{''}+\left(N+4\lambda_{1}^{2}+\dots+4\lambda_{N}^{2}\right)I_{N^{p}},\quad\mbox{for all $k= 1,\dots,p+1$}.\label{hihi}\end{equation}

We recall (\ref{lili1}), (\ref{suc1}),  (\ref{sucL11}),   (\ref{sucL111se}) and (\ref{90000}), (\ref{90000se}), (\ref{90000se1}). We obtain
 \begin{equation} \mathcal{M}_{k,k}\left[I;J\right]=\begin{pmatrix}  1 & \dots   & 0 & \dots & 0 & \dots & 0& \dots & 0  \\ \vdots & \ddots & \vdots & \ddots & \vdots & \ddots & \vdots& \ddots & \vdots   \\ 0 & \dots & \gamma_{1,1}\left[I;J\right]  & \dots & \gamma_{1,k}\left[I;J\right]  & \dots & \gamma_{1,N}\left[I;J\right]& \dots& 0 
\\ \vdots & \ddots & \vdots  & \ddots & \vdots  & \ddots & \vdots  & \ddots & \vdots  \\ 0 & \dots &  \gamma_{k,1}\left[I;J\right]  & \dots &  \gamma_{k,k}\left[I;J\right]    & \dots &  \gamma_{k,N}\left[I;J\right]& \dots & 0 \\ \vdots  & \ddots & \vdots  & \ddots & \vdots  & \ddots & \vdots  & \ddots & \vdots  \\ 0 & \dots &  \gamma_{N,1}\left[I,J\right]  & \dots &  \gamma_{N,k}\left[I;J\right]   & \dots &  \gamma_{N,N}\left[I;J\right]  & \dots & 0 \\ \vdots  & \ddots & \vdots  & \ddots & \vdots  & \ddots & \vdots  & \ddots & \vdots  \\ 0  & \dots & 0  & \dots &0  & \dots & 0  & \dots & 1  
\end{pmatrix}\in \mathcal{M}_{N^{p}\times N^{p}}\left(\mathbb{C}\right), \label{shobi}
\end{equation}
where  (\ref{lili}) is satisfied, for all $k=1,\dots,p+1$, using by (\ref{yy})  new notations starting with the first column, which is defined as follows
\begin{equation*} \begin{split}&\gamma_{1,1}\left[I;J\right]=\left(j_{1}+1\right)\left(i_{1}+1\right)+ \left(j_{1}+1\right)\left(j_{1}+2\right)\lambda_{1}^{2} +  \left(i_{1}+1\right)\left(i_{1}+2\right)\lambda_{1}^{2}\\&\quad\quad\quad\quad\quad\quad\quad\quad\quad\quad\quad\quad\hspace{0.1 cm} +  \left(i_{1}+j_{1}+1\right)\left(1+4\lambda_{1}^{2}\right)+\left(N+4\lambda_{1}^{2}+\dots+4\lambda_{N}^{2}\right),\\&\quad\quad\quad\quad\quad\vdots \\& \gamma_{1,k}\left[I;J\right]=  {\left(j_{1}+1\right)\left(i_{1}+1\right)+ \left(j_{1}+1\right)\left(j_{1}+2\right)\lambda_{1}\lambda_{k} +\left(i_{1}+1\right)\left(i_{1}+2\right)\lambda_{1}\lambda_{k}}, \\& \quad\quad\quad\quad\quad\vdots \\& \gamma_{1,N}\left[I;J\right]=\left(j_{1}+1\right)\left(i_{1}+1\right)+ \left(j_{1}+1\right)\left(j_{1}+2\right)\lambda_{1} \lambda_{N}+\left(i_{1}+1\right)\left(i_{1}+2\right)\lambda_{1}\lambda_{N},\end{split} 
\end{equation*}
up to the last column, which is defined as follows
\begin{equation*} \begin{split}&\gamma_{N,1}\left[I;J\right]=\left(j_{N}+1\right)\left(i_{N}+1\right)+ \left(j_{N}+1\right)\left(j_{N}+2\right)\lambda_{1}\lambda_{N} +  \left(i_{N}+1\right)\left(i_{N}+2\right)\lambda_{1}\lambda_{N},\\&\quad\quad\quad\quad\quad\vdots \\& \gamma_{N,k}\left[I;J\right]=\left(j_{N}+1\right)\left(i_{N}+1\right)+\left(j_{1}+1\right)\left(j_{1}+2\right)\lambda_{1}\lambda_{k} + \left(i_{1}+1\right)\left(i_{1}+2\right)\lambda_{1}\lambda_{k},   \\& \quad\quad\quad\quad\quad\vdots \\& \gamma_{N,N}\left[I;J\right]=\left(j_{N}+1\right)\left(i_{N}+1\right)+ \left(j_{N}+1\right)\left(j_{N}+2\right)\lambda_{N}^{2} +  \left(i_{N}+1\right)\left(i_{N}+2\right)\lambda_{N}^{2}\\&\quad\quad\quad\quad\quad\quad\quad\quad\quad\quad\quad\quad\quad\hspace{0.18 cm}+  \left(i_{N}+j_{N}+1\right)\left(1+4\lambda_{N}^{2}\right)+\left(N+4\lambda_{1}^{2}+\dots+4\lambda_{N}^{2}\right),\end{split} 
\end{equation*} 
but excepting the $k$-column, which is defined as follows
\begin{equation*} \begin{split}&\gamma_{k,1}\left[I;J\right]=j_{k}i_{k}+  \left(j_{k}-1\right)j_{k}\lambda_{1}\lambda_{k} +  \left(i_{k}-1\right)i_{k}\lambda_{1}\lambda_{k},\\&\quad\quad\quad\quad\quad\vdots \\& \gamma_{k,k}\left[I;J\right]=  j_{k}i_{k}  +\left(j_{k}-1\right)j_{k} \lambda_{k}^{2} +   \left(i_{k}-1\right)i_{k}  \lambda_{k}^{2}+  \left(i_{k}+j_{k}+1\right)\left(1+4\lambda_{k}^{2}\right)+\left(N+4\lambda_{1}^{2}+\dots+4\lambda_{N}^{2}\right),  \\& \quad\quad\quad\quad\quad\vdots \\& \gamma_{k,N}\left[I;J\right]=j_{k}i_{k}+ \left(j_{k}-1\right)j_{k}\lambda_{k}\lambda_{N} +  \left(i_{k}-1\right)i_{k}\lambda_{k}\lambda_{N}.\end{split} 
\end{equation*}
 
It suffices to show the  invertibility of the matrices from (\ref{shobi}), because we can make simple permutations between its rows or its columns.  Then, we obtain that the invertibility of  (\ref{shobi}) is equivalent to the invertibility of the following matrix
\begin{equation}Q\left[I,J\right]=\begin{pmatrix} \gamma_{1,1}\left[I;J\right] & \dots &\gamma_{1,k}\left[I;J\right] & \dots &\gamma_{1,N}\left[I;J\right] \\ \vdots & \ddots &\vdots& \ddots &\vdots \\ \gamma_{k,1}\left[I;J\right] & \dots &\gamma_{k,k}\left[I;J\right] & \dots &\gamma_{k,N}\left[I;J\right]\\ \vdots & \ddots &\vdots& \ddots &\vdots \\ \gamma_{N,1}\left[I;J\right] & \dots &\gamma_{N,k}\left[I;J\right] & \dots &\gamma_{N,N}\left[I;J\right]
\end{pmatrix},\label{shobi1}
\end{equation}
which is just assumed to be  invertible, because its invertibility   will be proven  in Section $2.6$.  

Respecting to the previous computations,  we introduce  new notations in order to make look   simpler  the further computations. Thus, let's study  the following matrices
\begin{equation}\mathcal{M}_{k+2,k}=\mathcal{V}_{k}^{''},\quad\mbox{for all $k= 1,\dots,p-1$}.\label{hihi1} \end{equation}

We recall      (\ref{sucL11sese}) and (\ref{90000se1}). We obtain
 \begin{equation}\mathcal{M}_{k+2,k}\left[I;J\right]=\begin{pmatrix}  1 & \dots   & 0 & \dots & 0 & \dots & 0& \dots & 0  \\ \vdots & \ddots & \vdots & \ddots & \vdots & \ddots & \vdots& \ddots & \vdots   \\ 0 & \dots & \gamma_{1,1}^{''}\left[I;J\right]  & \dots & \gamma_{1,k}^{''}\left[I;J\right]  & \dots & \gamma_{1,N}^{''}\left[I;J\right]& \dots& 0 
\\ \vdots & \ddots & \vdots  & \ddots & \vdots  & \ddots & \vdots  & \ddots & \vdots  \\ 0 & \dots &  \gamma_{k,1}^{''}\left[I;J\right]  & \dots &  \gamma_{k,k}^{''}\left[I;J\right]    & \dots &  \gamma_{k,N}^{''}\left[I;J\right]& \dots & 0 \\ \vdots  & \ddots & \vdots  & \ddots & \vdots  & \ddots & \vdots  & \ddots & \vdots  \\ 0 & \dots &  \gamma_{N,1}^{''}\left[I,J\right]  & \dots &  \gamma_{N,k}^{''}\left[I;J\right]   & \dots &  \gamma_{N,N}^{''}\left[I;J\right]  & \dots & 0 \\ \vdots  & \ddots & \vdots  & \ddots & \vdots  & \ddots & \vdots  & \ddots & \vdots  \\ 0  & \dots & 0  & \dots &0  & \dots & 0  & \dots & 1  
\end{pmatrix}\in \mathcal{M}_{N^{p}\times N^{p}}\left(\mathbb{C}\right), \label{shobix1}
\end{equation}
 for all $k=1,\dots,p-1$,  using by (\ref{yy})  new notations starting with the first column, which is defined as follows
\begin{equation*} \begin{split}&\gamma_{1,1}^{''}\left[I;J\right]=\left(j_{1}+1\right)\left(j_{1}+2\right)\lambda_{1}^{2} ,\\&\quad\quad\quad\quad\quad\vdots \\& \gamma_{1,k}^{''}\left[I;J\right]= \left(j_{k}-1\right) j_{k} \lambda_{1}\lambda_{k}   , \\& \quad\quad\quad\quad\quad\vdots \\& \gamma_{1,N}^{''}\left[I;J\right]=\left(j_{N}+1\right)\left(j_{N}+2\right)\lambda_{1}\lambda_{N} ,\end{split} \end{equation*}  
going to the $k$-column, which is defined as follows
\begin{equation*} \begin{split}&\gamma_{k,1}^{''}\left[I;J\right]=\left(j_{1}+1\right)\left(j_{1}+2\right)\lambda_{1}\lambda_{k} ,\\&\quad\quad\quad\quad\quad\vdots \\& \gamma_{k,k}^{''}\left[I;J\right]=   \left(j_{k}-1\right) j_{k}  \lambda_{k}^{2} , \\& \quad\quad\quad\quad\quad\vdots \\& \gamma_{k,N}^{''}\left[I;J\right]=\left(j_{N}+1\right) \left(j_{N}+2\right) \lambda_{N}\lambda_{k} ,\end{split} \end{equation*}
and ending with the last column, which is defined as follows
\begin{equation*} \begin{split}&\gamma_{N,1}^{''}\left[I;J\right]=\left(j_{1}+1\right)\left(j_{1}+2\right)\lambda_{1}\lambda_{N} ,\\&\quad\quad\quad\quad\quad\vdots \\& \gamma_{N,k}^{''}\left[I;J\right]= \left(j_{k}-1\right) j_{k} \lambda_{N}\lambda_{k}   , \\& \quad\quad\quad\quad\quad\vdots \\& \gamma_{N,N}^{''}\left[I;J\right]=\left(j_{N}+1\right)\left(j_{N}+2\right) \lambda_{N}^{2}.\end{split} \end{equation*} 
 
Next, we continue to reshape the previous matrices having a look on the following matrices
\begin{equation}\mathcal{M}_{k+1,k}=\mathcal{W}_{k}^{(0)}+\mathcal{V}_{k}+\mathcal{O}_{k}^{''},\quad\mbox{for all $k= 1,\dots,p$}. \label{hihi2}\end{equation}

We recall   (\ref{lili3}), (\ref{suc2}), (\ref{sucL1sese}) and (\ref{90000}), (\ref{90000se1}). We obtain
 \begin{equation} \mathcal{M}_{k+1,k}\left[I;J\right]=\begin{pmatrix}  1 & \dots   & 0 & \dots & 0 & \dots & 0& \dots & 0  \\ \vdots & \ddots & \vdots & \ddots & \vdots & \ddots & \vdots& \ddots & \vdots   \\ 0 & \dots & \gamma_{1,1}^{'}\left[I;J\right]  & \dots & \gamma_{1,k}^{'}\left[I;J\right]  & \dots & \gamma_{1,N}^{'}\left[I;J\right]& \dots& 0 
\\ \vdots & \ddots & \vdots  & \ddots & \vdots  & \ddots & \vdots  & \ddots & \vdots  \\ 0 & \dots &  \gamma_{k,1}^{'}\left[I;J\right]  & \dots &  \gamma_{k,k}^{'}\left[I;J\right]    & \dots &  \gamma_{k,N}^{'}\left[I;J\right]& \dots & 0 \\ \vdots  & \ddots & \vdots  & \ddots & \vdots  & \ddots & \vdots  & \ddots & \vdots  \\ 0 & \dots &  \gamma_{N,1}^{'}\left[I,J\right]  & \dots &  \gamma_{N,k}^{'}\left[I;J\right]   & \dots &  \gamma_{N,N}^{'}\left[I;J\right]  & \dots & 0 \\ \vdots  & \ddots & \vdots  & \ddots & \vdots  & \ddots & \vdots  & \ddots & \vdots  \\ 0  & \dots & 0  & \dots &0  & \dots & 0  & \dots & 1  
\end{pmatrix}\in \mathcal{M}_{N^{p}\times N^{p}}\left(\mathbb{C}\right), \label{shobix2}
\end{equation}
 for all $k=1,\dots,p$,   using by (\ref{yy}) new notations starting with the first column, which is defined as follows
 \begin{equation*}\left.\begin{split}&\gamma_{1,1}^{'}\left[I;J\right]= \left(j_{1}+1\right)\left(j_{1}+2\right)\lambda_{1}+\left(i_{1}+1\right)\left(j_{1}+1\right)
 \lambda_{1}+4\lambda_{1}j_{1}  ,\\&\quad\quad\quad\quad\quad\vdots \\& \gamma_{1,k}^{'}\left[I;J\right]=  \left(j_{k}-1\right)j_{k}\lambda_{k}+i_{k}j_{k}
 \lambda_{k} , & \\& \quad\quad\quad\quad\quad\vdots \\& \gamma_{1,N}^{'}\left[I;J\right]= \left(j_{N}+1\right)\left(j_{N}+2\right)\lambda_{N}
 +\left(i_{N}+1\right)\left(j_{N}+1\right)\lambda_{N} ,\end{split}\right.\end{equation*}
ending with the last column, which is defined as follows
  \begin{equation*}\left.\begin{split}&\gamma_{N,1}^{'}\left[I;J\right]=\left(j_{1}+1\right)\left(j_{1}+2\right)\lambda_{1}+\left(i_{1}+1\right)\left(j_{1}+1\right)
 \lambda_{1}  ,\\&\quad\quad\quad\quad\quad\vdots \\& \gamma_{N,k}^{'}\left[I;J\right]= \left(j_{k}-1\right)j_{k}\lambda_{k}+i_{k}j_{k}
 \lambda_{k},&    \\& \quad\quad\quad\quad\quad\vdots \\& \gamma_{N,N}^{'}\left[I;J\right]=\left(j_{N}+1\right)\left(j_{N}+2\right)\lambda_{N}+\left(i_{N}+1\right)\left(j_{N}+1\right)
 \lambda_{N}+4\lambda_{N}j_{N}.\end{split}\right.\end{equation*}
  
Next, we continue to reshape the previous matrices having a look on the following matrices
\begin{equation}\mathcal{M}_{k-1,k}=\mathcal{V}_{k}^{(0)}+\mathcal{W}_{k}+\mathcal{O}_{k}^{'},\quad\mbox{for all $k= 2,\dots,p+1$}. \label{hihi1x}\end{equation}

We recall   (\ref{lili2}), (\ref{suc33}), (\ref{sucL1}) and (\ref{90000}), (\ref{90000se}). We obtain
 \begin{equation} \mathcal{M}_{k-1,k}\left[I;J\right]=\begin{pmatrix}  1 & \dots   & 0 & \dots & 0 & \dots & 0& \dots & 0  \\ \vdots & \ddots & \vdots & \ddots & \vdots & \ddots & \vdots& \ddots & \vdots   \\ 0 & \dots & \tilde{\gamma}_{1,1}^{'}\left[I;J\right]  & \dots & \tilde{\gamma}_{1,k}^{'}\left[I;J\right]  & \dots & \tilde{\gamma}_{1,N}^{'}\left[I;J\right]& \dots& 0 
\\ \vdots & \ddots & \vdots  & \ddots & \vdots  & \ddots & \vdots  & \ddots & \vdots  \\ 0 & \dots &  \tilde{\gamma}_{k,1}^{'}\left[I;J\right]  & \dots &  \tilde{\gamma}_{k,k}^{'}\left[I;J\right]    & \dots &  \tilde{\gamma}_{k,N}^{'}\left[I;J\right]& \dots & 0 \\ \vdots  & \ddots & \vdots  & \ddots & \vdots  & \ddots & \vdots  & \ddots & \vdots  \\ 0 & \dots &  \tilde{\gamma}_{N,1}^{'}\left[I,J\right]  & \dots &  \gamma_{N,k}^{'}\left[I;J\right]   & \dots &  \tilde{\gamma}_{N,N}^{'}\left[I;J\right]  & \dots & 0 \\ \vdots  & \ddots & \vdots  & \ddots & \vdots  & \ddots & \vdots  & \ddots & \vdots  \\ 0  & \dots & 0  & \dots &0  & \dots & 0  & \dots & 1  
\end{pmatrix}\in \mathcal{M}_{N^{p}\times N^{p}}\left(\mathbb{C}\right), \label{shobiy1}
\end{equation}
 for all $k=2,\dots,p+1$,   using by (\ref{yy}) new notations  starting with the first column,  which is defined as follows
\begin{equation*}  \left. \begin{split}& \tilde{\gamma}_{1,1}^{'}\left[I;J\right]= \left(i_{1}+1\right)\left(i_{1}+2\right)\lambda_{1}+\left(i_{1}+1\right)\left(j_{1}+1\right)
  \lambda_{1}+4\lambda_{1} i_{1},\\& \quad\quad\quad\quad\quad\vdots \\& \tilde{\gamma}_{1,k}^{'}\left[I;J\right]=  \left(i_{k}-1\right)i_{k}\lambda_{k}+i_{k}j_{k}
 \lambda_{k} ,  \\& \quad\quad\quad\quad\quad\vdots \\& \tilde{\gamma}_{1,N}^{'}\left[I;J\right]= \left(i_{N}+1\right)\left(i_{N}+2\right)\lambda_{N}+\left(i_{N}+1\right)\left(j_{N}+1\right) \lambda_{N}\end{split}\right. \end{equation*} 
ending with the last column,  which is defined as follows
 \begin{equation*}\left.\begin{split}& \tilde{\gamma}_{N,1}^{'}\left[I;J\right]= \left(i_{1}+1\right)\left(i_{1}+2\right)\lambda_{1}+\left(i_{1}+1\right)\left(j_{1}+1\right)\lambda_{1},\\&\quad\quad\quad\quad\quad\vdots \\&
 \tilde{\gamma}_{N,k}^{'}\left[I;J\right]= \left(i_{k}-1\right)i_{k}\lambda_{k}+i_{k}j_{k}\lambda_{k},    \\& \quad\quad\quad\quad\quad\vdots \\& \tilde{\gamma}_{N,N}^{'}\left[I;J\right]=\left(i_{N}+1\right)\left(i_{N}+2\right)\lambda_{N}+\left(i_{N}+1\right)\left(j_{N}+1\right)\lambda_{N}+4\lambda_{N} i_{N}.\end{split}\right.   
\end{equation*} 
 
We study now the following matrices
\begin{equation}\mathcal{M}_{k-2,k}=\mathcal{W}_{k}^{'},\quad\mbox{for all $k= 3,\dots,p+1$}.\label{hihi2x}\end{equation}

We recall  (\ref{sucL111}) and (\ref{90000se}). We obtain
 \begin{equation}  \mathcal{M}_{k-2,k}\left[I;J\right]=\begin{pmatrix}  1 & \dots   & 0 & \dots & 0 & \dots & 0& \dots & 0  \\ \vdots & \ddots & \vdots & \ddots & \vdots & \ddots & \vdots& \ddots & \vdots   \\ 0 & \dots & \tilde{\gamma}_{1,1}^{''}\left[I;J\right]  & \dots & \tilde{\gamma}_{1,k}^{''}\left[I;J\right]  & \dots & \tilde{\gamma}_{1,N}^{''}\left[I;J\right]& \dots& 0 
\\ \vdots & \ddots & \vdots  & \ddots & \vdots  & \ddots & \vdots  & \ddots & \vdots  \\ 0 & \dots &  \tilde{\gamma}_{k,1}^{''}\left[I;J\right]  & \dots &  \tilde{\gamma}_{k,k}^{''}\left[I;J\right]    & \dots &  \tilde{\gamma}_{k,N}^{''}\left[I;J\right]& \dots & 0 \\ \vdots  & \ddots & \vdots  & \ddots & \vdots  & \ddots & \vdots  & \ddots & \vdots  \\ 0 & \dots &  \tilde{\gamma}_{N,1}^{''}\left[I,J\right]  & \dots &  \tilde{\gamma}_{N,k}^{''}\left[I;J\right]   & \dots & \tilde{\gamma}_{N,N}^{''}\left[I;J\right]  & \dots & 0 \\ \vdots  & \ddots & \vdots  & \ddots & \vdots  & \ddots & \vdots  & \ddots & \vdots  \\ 0  & \dots & 0  & \dots &0  & \dots & 0  & \dots & 1  
\end{pmatrix}\in \mathcal{M}_{N^{p}\times N^{p}}\left(\mathbb{C}\right), \label{shobiy2}
\end{equation}
for all $k=3,\dots,p+1$, using by (\ref{yy})  new notations  starting with the first column, which is defined as follows
\begin{equation*}\left.\begin{split}&\tilde{\gamma}_{1,1}^{''}\left[I;J\right]=\left(i_{1}+1\right)\left(i_{1}+2\right)\lambda_{1}^{2} ,\\&\quad\quad\quad\quad\quad\vdots \\& \tilde{\gamma}_{1,k}^{''}\left[I;J\right]= \left(i_{k}-1\right) i_{k} \lambda_{1}\lambda_{k}   , \\& \quad\quad\quad\quad\quad\vdots \\& \tilde{\gamma}_{1,N}^{''}\left[I;J\right]=\left(i_{N}+1\right)\left(i_{N}+2\right)\lambda_{1}\lambda_{N} ,\end{split}\right.\end{equation*}
going to the $k$-column, which is defined as follows
\begin{equation*}\left.\begin{split}&\tilde{\gamma}_{k,1}^{''}\left[I;J\right]=\left(i_{1}+1\right)\left(i_{1}+2\right)\lambda_{1}\lambda_{k} ,\\&\quad\quad\quad\quad\quad\vdots \\& \tilde{\gamma}_{k,k}^{''}\left[I;J\right]=   \left(i_{k}-1\right) i_{k}  \lambda_{k}^{2} , \\& \quad\quad\quad\quad\quad\vdots \\& \tilde{\gamma}_{k,N}^{''}\left[I;J\right]=\left(i_{N}+1\right) \left(i_{N}+2\right) \lambda_{N}\lambda_{k} ,\end{split}\right.\end{equation*}
and ending with the last column, which is defined as follows
\begin{equation*}\left.\begin{split}&\tilde{\gamma}_{N,1}^{''}\left[I;J\right]=\left(i_{1}+1\right)\left(i_{1}+2\right)\lambda_{1}\lambda_{N} ,\\&\quad\quad\quad\quad\quad\vdots \\& \tilde{\gamma}_{N,k}^{''}\left[I;J\right]= \left(i_{k}-1\right) i_{k} \lambda_{N}\lambda_{k}   , \\& \quad\quad\quad\quad\quad\vdots \\& \tilde{\gamma}_{N,N}^{''}\left[I;J\right]=\left(i_{N}+1\right)\left(i_{N}+2\right) \lambda_{N}^{2}.\end{split}\right.\end{equation*}
\subsection{Simplified Systems of Equations} We  simplify the system of equations (\ref{beb1}), which  is  equivalent to
 \begin{equation}\begin{pmatrix} I_{N^{p}} & \frac{\mathcal{M}_{1,2}}{\mathcal{M}_{1,1}} & \frac{\mathcal{M}_{1,3}}{\mathcal{M}_{1,1}} & \mbox{O}_{N^{p}} & \mbox{O}_{N^{p}} &\dots & \mbox{O}_{N^{p}} & \mbox{O}_{N^{p}} & \mbox{O}_{N^{p}} & \mbox{O}_{N^{p}} \\ \frac{\mathcal{M}_{2,1}}{\mathcal{M}_{2,2}} & I_{N^{p}} & \frac{\mathcal{M}_{2,3}}{\mathcal{M}_{2,2}} & \frac{\mathcal{M}_{2,4}}{\mathcal{M}_{2,2}} & \mbox{O}_{N^{p}} &\dots & \mbox{O}_{N^{p}} & \mbox{O}_{N^{p}} & \mbox{O}_{N^{p}} & \mbox{O}_{N^{p}} \\ \frac{\mathcal{M}_{3,1}}{\mathcal{M}_{3,3}}  & \frac{\mathcal{M}_{3,2}}{\mathcal{M}_{3,3}}  & I_{N^{p}}  & \frac{\mathcal{M}_{3,4}}{\mathcal{M}_{3,3}}  & \frac{\mathcal{M}_{3,5}}{\mathcal{M}_{3,3}}  &\dots & \mbox{O}_{N^{p}} & \mbox{O}_{N^{p}} & \mbox{O}_{N^{p}} & \mbox{O}_{N^{p}}  \\ \mbox{O}_{N^{p}} & \frac{\mathcal{M}_{4,2}}{\mathcal{M}_{4,4}}  & \frac{\mathcal{M}_{4,3}}{\mathcal{M}_{4,4}}  & I_{N^{p}} & \frac{\mathcal{M}_{4,5}}{\mathcal{M}_{4,4}}  &\dots & \mbox{O}_{N^{p}}& \mbox{O}_{N^{p}} &\mbox{O}_{N^{p}}& \mbox{O}_{N^{p}} \\ \mbox{O}_{N^{p}}& \mbox{O}_{N^{p}} & \frac{\mathcal{M}_{5,3}}{\mathcal{M}_{5,5}} & \frac{\mathcal{M}_{5,4}}{\mathcal{M}_{5,5}}& I_{N^{p}} &\dots & \mbox{O}_{N^{p}} & \mbox{O}_{N^{p}} & \mbox{O}_{N^{p}} &\mbox{O}_{N^{p}}   \\ \vdots & \vdots & \vdots & \vdots & \vdots & \ddots & \vdots & \vdots & \vdots & \vdots \\ \mbox{O}_{N^{p}} & \mbox{O}_{N^{p}} & \mbox{O}_{N^{p}} & \mbox{O}_{N^{p}} & \mbox{O}_{N^{p}} &\dots & I_{N^{p}} & \frac{\mathcal{M}_{p-2,p-1}}{\mathcal{M}_{p-2,p-2}}  & \frac{\mathcal{M}_{p-2,p}}{\mathcal{M}_{p-2,p-2}} & \mbox{O}_{N^{p}} \\ \mbox{O}_{N^{p}} & \mbox{O}_{N^{p}} & \mbox{O}_{N^{p}} & \mbox{O}_{N^{p}} & \mbox{O}_{N^{p}} &\dots & \frac{\mathcal{M}_{p-1,p-2}}{\mathcal{M}_{p-1,p-1}}  & I_{N^{p}} & \frac{\mathcal{M}_{p-1,p}}{\mathcal{M}_{p-1,p-1}}  & \frac{\mathcal{M}_{p-1,p+1}}{\mathcal{M}_{p-1,p-1}}   \\ \mbox{O}_{N^{p}} & \mbox{O}_{N^{p}} & \mbox{O}_{N^{p}} & \mbox{O}_{N^{p}} & \mbox{O}_{N^{p}} &\dots & \frac{\mathcal{M}_{p,p-2}}{\mathcal{M}_{p,p}}  & \frac{\mathcal{M}_{p,p-1}}{\mathcal{M}_{p,p}}  & I_{N^{p}}  & \frac{\mathcal{M}_{p,p+1}}{\mathcal{M}_{p,p}}  \\ \mbox{O}_{N^{p}} & \mbox{O}_{N^{p}} & \mbox{O}_{N^{p}} & \mbox{O}_{N^{p}} & \mbox{O}_{N^{p}} &\dots & \mbox{O}_{N^{p}} & \frac{\mathcal{M}_{p+1,p-1}}{\mathcal{M}_{p+1,p+1}}  & \frac{\mathcal{M}_{p+1,p}}{\mathcal{M}_{p+1,p+1}}  & I_{N^{p}}        
\end{pmatrix}\begin{pmatrix} Y_{1}\\Y_{2}\\Y_{3}\\Y_{4}\\Y_{5}\\ \vdots\\ Y_{p-3}\\ Y_{p-2} \\ Y_{p-1}\\ Y_{p}\\ Y_{p+1}
\end{pmatrix}=\begin{pmatrix} \frac{1}{\mathcal{M}_{11}}V_{1}\\ V_{0} \\ V_{0} \\ V_{0}\\ V_{0}  \\   \vdots \\ V_{0} \\ V_{0}\\ V_{0}\\ V_{0}\\ V_{0}
\end{pmatrix}, \label{beb120}
 \end{equation}
where its  elements are matrices   defined by (\ref{lil}). 

We  rewrite its first equation  as follows
\begin{equation}\mbox{I}_{N^{p}} Y_{1}+ \mathcal{N}_{1,2}Y_{2}+ \mathcal{N}_{1,3}Y_{3}+\displaystyle\sum_{k=2}^{p+1}\mbox{O}_{N^{p}}Y_{k}=\tilde{V}_{1},\label{gringo1}
\end{equation}
where we have used the following matrices
\begin{equation}\mathcal{N}_{1,2}=\frac{\mathcal{M}_{1,2}}{\mathcal{M}_{1,1}},\quad \mathcal{N}_{1,3}=\frac{\mathcal{M}_{1,3}}{\mathcal{M}_{1,1}},\quad \tilde{V}_{1}=\frac{1}{\mathcal{M}_{1,1}}V_{1}.\label{nene1} \end{equation}

Making   substraction from the second equation in (\ref{beb120}) using  (\ref{gringo1}), we obtain 
\begin{equation} \left(\mbox{I}_{N^{p}} -\frac{\mathcal{M}_{2,1}}{\mathcal{M}_{2,2}}\frac{\mathcal{M}_{1,2}}{\mathcal{M}_{1,1}} \right)Y_{2}+\left(\frac{\mathcal{M}_{2,3}}{\mathcal{M}_{2,2}}-\frac{\mathcal{M}_{2,1}}{\mathcal{M}_{2,2}}\frac{\mathcal{M}_{1,3}}{\mathcal{M}_{1,1}}\right)Y_{3}+\frac{\mathcal{M}_{2,4}}{\mathcal{M}_{2,2}}Y_{4}+\displaystyle\sum_{k=5}^{p+1}\mbox{O}_{N^{p}}Y_{k}= -\frac{\mathcal{M}_{2,1}}{\mathcal{M}_{2,2}}\frac{1}{\mathcal{M}_{1,1}}V_{1},\label{ggrin1} \end{equation}
or  equivalently by (\ref{nene1}) to
\begin{equation}\mbox{I}_{N^{p}} Y_{2}+\mathcal{N}_{2,3}Y_{3}+\mathcal{N}_{2,4}Y_{4}+\mbox{O}_{N^{p}}Y_{5}+\dots+\mbox{O}_{N^{p}}Y_{p+1}=\tilde{V}_{2},\label{gringo2}
\end{equation}
where we have used the following matrices
\begin{equation}\mathcal{N}_{2,3}=\frac{ \frac{\mathcal{M}_{2,3}}{\mathcal{M}_{2,2}}-\frac{\mathcal{M}_{2,1}}{\mathcal{M}_{2,2}} \mathcal{N}_{1,3} }{\mbox{I}_{N^{p}} -\frac{\mathcal{M}_{2,1}}{\mathcal{M}_{2,2}}\mathcal{N}_{1,2}},\quad \mathcal{N}_{2,4}=\frac{\frac{\mathcal{M}_{2,4}}{\mathcal{M}_{2,2}}}{\mbox{I}_{N^{p}} -\frac{\mathcal{M}_{2,1}}{\mathcal{M}_{2,2}}\mathcal{N}_{1,2}},\quad \tilde{V}_{2}= \frac{-\frac{\mathcal{M}_{2,1}}{\mathcal{M}_{2,2}}\frac{1}{\mathcal{M}_{1,1}}}{\mbox{I}_{N^{p}} -\frac{\mathcal{M}_{2,1}}{\mathcal{M}_{2,2}}\mathcal{N}_{1,2}}V_{1}. \label{nene2}\end{equation}

Making   substraction from the third equation in (\ref{beb120}) using  (\ref{gringo1}), we obtain 
\begin{equation} \left(\frac{\mathcal{M}_{3,2}}{\mathcal{M}_{3,3}}-\frac{\mathcal{M}_{3,1}}{\mathcal{M}_{3,3}}\frac{\mathcal{M}_{1,2}}{\mathcal{M}_{1,1}}\right)Y_{2}+\left(\mbox{I}_{N^{p}} -\frac{\mathcal{M}_{3,1}}{\mathcal{M}_{3,3}}\frac{\mathcal{M}_{1,3}}{\mathcal{M}_{1,1}} \right)Y_{3}+ \frac{\mathcal{M}_{3,4}}{\mathcal{M}_{3,3}}Y_{4}+ \displaystyle\sum_{k=5}^{p+1}\mbox{O}_{N^{p}}Y_{k}= -\frac{\mathcal{M}_{3,1}}{\mathcal{M}_{3,3}}\frac{1}{\mathcal{M}_{1,1}}V_{1}, \label{gringo3}\end{equation}

Making   substraction from    (\ref{gringo3}) using  (\ref{gringo2}), we obtain 
\begin{equation}\begin{split}&\left(\mbox{I}_{N^{p}} -\frac{\mathcal{M}_{3,1}}{\mathcal{M}_{3,3}}\frac{\mathcal{M}_{1,3}}{\mathcal{M}_{1,1}} -\left(\frac{\mathcal{M}_{3,2}}{\mathcal{M}_{3,3}}-\frac{\mathcal{M}_{3,1}}{\mathcal{M}_{3,3}}\frac{\mathcal{M}_{1,2}}{\mathcal{M}_{1,1}}\right)\mathcal{N}_{2,3}\right)Y_{3}\\& \quad\quad\quad\quad\quad \quad   \quad\quad\quad\quad\quad\quad      \begin{tabular}{l} \rotatebox[origin=c]{270}{$+$}\end{tabular} \\&\left( \frac{\mathcal{M}_{3,4}}{\mathcal{M}_{3,3}}-\left(\frac{\mathcal{M}_{3,2}}{\mathcal{M}_{3,3}}-\frac{\mathcal{M}_{3,1}}{\mathcal{M}_{3,3}}\frac{\mathcal{M}_{1,2}}{\mathcal{M}_{1,1}}\right)\mathcal{N}_{2,4}\right)Y_{4}+\frac{\mathcal{M}_{3,5}}{\mathcal{M}_{3,3}}Y_{5}  +\displaystyle\sum_{k=6}^{p+1}\mbox{O}_{N^{p}}Y_{k}=  -\frac{\mathcal{M}_{3,1}}{\mathcal{M}_{3,3}}\frac{1}{\mathcal{M}_{1,1}} V_{1}-\left(\frac{\mathcal{M}_{3,2}}{\mathcal{M}_{3,3}}-\frac{\mathcal{M}_{3,1}}{\mathcal{M}_{3,3}}\frac{\mathcal{M}_{1,2}}{\mathcal{M}_{1,1}}\right)\tilde{V}_{2},\end{split}\label{ggrin2}\end{equation}
or equivalently, by (\ref{nene1}) and (\ref{nene2}), to
\begin{equation}  \mbox{I}_{N^{p}} Y_{3}+\mathcal{N}_{3,4}Y_{4}+ \mathcal{N}_{3,5}Y_{5}+\displaystyle\sum_{k=6}^{p+1}\mbox{O}_{N^{p}}Y_{k}=\tilde{V}_{3}, \label{gringo4} 
\end{equation} 
where we have used the following  matrices
\begin{equation} \begin{split}&  \mathcal{N}_{3,4}=\frac{ \frac{\mathcal{M}_{3,4}}{\mathcal{M}_{3,3}}-\left(\frac{\mathcal{M}_{3,2}}{\mathcal{M}_{3,3}}-\frac{\mathcal{M}_{3,1}}{\mathcal{M}_{3,3}}\mathcal{N}_{1,2}\right)\mathcal{N}_{2,4}}{\mbox{I}_{N^{p}} -\frac{\mathcal{M}_{3,1}}{\mathcal{M}_{3,3}}\mathcal{N}_{1,3}   -\left(\frac{\mathcal{M}_{3,2}}{\mathcal{M}_{3,3}}-\frac{\mathcal{M}_{3,1}}{\mathcal{M}_{3,3}} \mathcal{N}_{1,2} \right)\mathcal{N}_{2,3}}, \\& \mathcal{N}_{3,5}=\frac{\frac{\mathcal{M}_{3,5}}{\mathcal{M}_{3,3}}}{\mbox{I}_{N^{p}} -\frac{\mathcal{M}_{3,1}}{\mathcal{M}_{3,3}}\mathcal{N}_{1,3}   -\left(\frac{\mathcal{M}_{3,2}}{\mathcal{M}_{3,3}}-\frac{\mathcal{M}_{3,1}}{\mathcal{M}_{3,3}} \mathcal{N}_{1,2} \right)\mathcal{N}_{2,3}},\\&  \quad \tilde{V}_{3}=\frac{-\frac{\mathcal{M}_{3,1}}{\mathcal{M}_{3,3}}\tilde{V}_{1}-\left(\frac{\mathcal{M}_{3,2}}{\mathcal{M}_{3,3}}-\frac{\mathcal{M}_{3,1}}{\mathcal{M}_{3,3}}\mathcal{N}_{1,2}\right)\tilde{V}_{2}}{\mbox{I}_{N^{p}} -\frac{\mathcal{M}_{3,1}}{\mathcal{M}_{3,3}}\mathcal{N}_{1,3}   -\left(\frac{\mathcal{M}_{3,2}}{\mathcal{M}_{3,3}}-\frac{\mathcal{M}_{3,1}}{\mathcal{M}_{3,3}} \mathcal{N}_{1,2} \right)\mathcal{N}_{2,3}}.\label{nene3}\end{split} \end{equation}

Making   substraction from the fourth equation in (\ref{beb120}) using  (\ref{gringo2}), we obtain 
\begin{equation}\left(\frac{\mathcal{M}_{4,3}}{\mathcal{M}_{4,4}}-\frac{\mathcal{M}_{4,2}}{\mathcal{M}_{4,4}}\mathcal{N}_{2,3}\right)Y_{3}+\left(\mbox{I}_{N^{p}} -\frac{\mathcal{M}_{4,2}}{\mathcal{M}_{4,4}}\mathcal{N}_{2,4}\right)Y_{4}+  \frac{\mathcal{M}_{4,5}}{\mathcal{M}_{4,4}}Y_{5}+ \frac{\mathcal{M}_{4,6}}{\mathcal{M}_{4,4}}Y_{6}+\displaystyle\sum_{k=7}^{p+1}\mbox{O}_{N^{p}}Y_{k}= -\frac{\mathcal{M}_{4,2}}{\mathcal{M}_{4,4}}\tilde{V}_{2}.\label{gringo5}\end{equation}

Making   substraction  from  (\ref{gringo5}) using  (\ref{gringo4}), we obtain 
\begin{equation}  \mbox{I}_{N^{p}} Y_{4}+\mathcal{N}_{4,5}Y_{5}+ \mathcal{N}_{4,6}Y_{6}+\displaystyle\sum_{k=7}^{p+1}\mbox{O}_{N^{p}}Y_{k}=\tilde{V}_{4}, \label{gringo4se} 
\end{equation}
where we have used the following  matrices
\begin{equation} \begin{split}&  \mathcal{N}_{4,5}=\frac{ \frac{\mathcal{M}_{4,5}}{\mathcal{M}_{4,4}}-\left(\frac{\mathcal{M}_{4,3}}{\mathcal{M}_{4,4}}-\frac{\mathcal{M}_{4,2}}{\mathcal{M}_{4,4}}\mathcal{N}_{2,3}\right)\mathcal{N}_{3,5}}{\mbox{I}_{N^{p}} -\frac{\mathcal{M}_{4,2}}{\mathcal{M}_{4,4}}\mathcal{N}_{2,4}   -\left(\frac{\mathcal{M}_{4,3}}{\mathcal{M}_{4,4}}-\frac{\mathcal{M}_{4,2}}{\mathcal{M}_{4,4}}\mathcal{N}_{2,3}\right)\mathcal{N}_{3,4}}, \\& \mathcal{N}_{4,6}=\frac{\frac{\mathcal{M}_{4,6}}{\mathcal{M}_{4,4}}}{\mbox{I}_{N^{p}} -\frac{\mathcal{M}_{4,2}}{\mathcal{M}_{4,4}}\mathcal{N}_{2,4}   -\left(\frac{\mathcal{M}_{4,3}}{\mathcal{M}_{4,4}}-\frac{\mathcal{M}_{4,2}}{\mathcal{M}_{4,4}}\mathcal{N}_{2,3}\right)\mathcal{N}_{3,4}},\\&  \quad\tilde{V}_{4}=\frac{-\frac{\mathcal{M}_{4,2}}{\mathcal{M}_{4,4}}\tilde{V}_{2}-\left(\frac{\mathcal{M}_{4,3}}{\mathcal{M}_{4,4}}-\frac{\mathcal{M}_{4,2}}{\mathcal{M}_{4,4}}\mathcal{N}_{2,3}\right)\tilde{V}_{3}}{\mbox{I}_{N^{p}} -\frac{\mathcal{M}_{4,2}}{\mathcal{M}_{4,4}}\mathcal{N}_{2,4}   -\left(\frac{\mathcal{M}_{4,3}}{\mathcal{M}_{4,4}}-\frac{\mathcal{M}_{4,2}}{\mathcal{M}_{4,4}}\mathcal{N}_{2,3}\right)\mathcal{N}_{3,4}}.\label{nene4}\end{split} \end{equation}

We continue the computations considering an induction process depending on  $ 3\leq l  \leq p$. Assume 
\begin{equation}  \mbox{I}_{N^{p}} Y_{l-2}+\mathcal{N}_{l-2,l-1}Y_{l-1}+ \mathcal{N}_{l-2,l}Y_{l}+\displaystyle\sum_{k=l+1}^{p+1}\mbox{O}_{N^{p}}Y_{k}=\tilde{V}_{l-2}, \label{gringo44se} 
\end{equation}  
and respectively 
\begin{equation}\mbox{I}_{N^{p}} Y_{l-1}+\mathcal{N}_{l-1,l}Y_{l}+ \mathcal{N}_{l-1,l+1}Y_{l+1}+\displaystyle\sum_{k=l+2}^{p+1}\mbox{O}_{N^{p}}Y_{k}=\tilde{V}_{l-1}, \label{gringo44se1} 
\end{equation}  

We write again the equation $l$ from (\ref{beb120}): 
\begin{equation} \frac{\mathcal{M}_{l,l-2}}{\mathcal{M}_{l,l}}Y_{l-2}+\frac{\mathcal{M}_{l,l-1}}{\mathcal{M}_{l,l}}Y_{l-1}+ \mbox{I}_{N^{p}}  Y_{l}+ \frac{\mathcal{M}_{l,l+1}}{\mathcal{M}_{l,l}} Y_{l+1}+ \frac{\mathcal{M}_{l,l+2}}{\mathcal{M}_{l,l}} Y_{l+2}+\displaystyle\sum_{k=l+3}^{p+1}\mbox{O}_{N^{p}}Y_{k}=V_{0}.\label{gringoo1}
\end{equation}  

Making substraction from (\ref{gringoo1}) using (\ref{gringo44se}), we obtain  
\begin{equation}    \left(\frac{\mathcal{M}_{l,l-1}}{\mathcal{M}_{l,l}}-\frac{\mathcal{M}_{l,l-2}}{\mathcal{M}_{l,l}}\mathcal{N}_{l-2,l-1}\right)Y_{l-1}+\left(\mbox{I}_{N^{p}} -\frac{\mathcal{M}_{l,l-2}}{\mathcal{M}_{l,l}}\mathcal{N}_{l-2,l}\right)Y_{l}+  \frac{\mathcal{M}_{l,l+1}}{\mathcal{M}_{l,l}}Y_{l+1}+ \frac{\mathcal{M}_{l,l+2}}{\mathcal{M}_{l,l}}Y_{l+2}+\displaystyle\sum_{k=l+3}^{p+1}\mbox{O}_{N^{p}}Y_{k}=  -\frac{\mathcal{M}_{l,l-2}}{\mathcal{M}_{l,l}}\tilde{V}_{l-2}. \label{gringo55}\end{equation}  

Making substraction from  (\ref{gringo55}) using (\ref{gringo44se1}), we obtain  
 \begin{equation}\mbox{I}_{N^{p}} Y_{l}+ \mathcal{N}_{l,l+1}Y_{l+1}+ \mathcal{N}_{l,l+2}Y_{l+2}+\displaystyle\sum_{k=l+3}^{p+1}\mbox{O}_{N^{p}}Y_{k}=\tilde{V}_{l},\label{gringooo}
\end{equation} 
where we have used the following  matrices
\begin{equation} \begin{split}&  \mathcal{N}_{l,l+1}=\frac{\frac{\mathcal{M}_{l,l+1}}{\mathcal{M}_{l,l}}-\left(\frac{\mathcal{M}_{l,l-1}}{\mathcal{M}_{l,l}}-\frac{\mathcal{M}_{l,l-2}}{\mathcal{M}_{l,l}}\mathcal{N}_{l-2,l-1}\right)\mathcal{N}_{l-1,l+1} }{\mbox{I}_{N^{p}} -\frac{\mathcal{M}_{l,l-2}}{\mathcal{M}_{l,l}}\mathcal{N}_{l-2,l}-\left(\frac{\mathcal{M}_{l,l-1}}{\mathcal{M}_{l,l}}-\frac{\mathcal{M}_{l,l-2}}{\mathcal{M}_{l,l}}\mathcal{N}_{l-2,l-1}\right)\mathcal{N}_{l-1,l} } ,\\&   \mathcal{N}_{l,l+2}=\frac{\frac{\mathcal{M}_{l,l+2}}{\mathcal{M}_{l,l}} }{\mbox{I}_{N^{p}} -\frac{\mathcal{M}_{l,l-2}}{\mathcal{M}_{l,l}}\mathcal{N}_{l-2,l}-\left(\frac{\mathcal{M}_{l,l-1}}{\mathcal{M}_{l,l}}-\frac{\mathcal{M}_{l,l-2}}{\mathcal{M}_{l,l}}\mathcal{N}_{l-2,l-1}\right)\mathcal{N}_{l-1,l} },\\& \quad\hspace{0.25 cm}   \tilde{V}_{l}=\frac{-\frac{\mathcal{M}_{l,l-2}}{\mathcal{M}_{l,l}}\tilde{V}_{l-2}-\left(\frac{\mathcal{M}_{l,l-1}}{\mathcal{M}_{l,l}}-\frac{\mathcal{M}_{l,l-2}}{\mathcal{M}_{l,l}}\mathcal{N}_{l-2,l-1}\right)\tilde{V}_{l-1} }{\mbox{I}_{N^{p}} -\frac{\mathcal{M}_{l,l-2}}{\mathcal{M}_{l,l}}\mathcal{N}_{l-2,l}-\left(\frac{\mathcal{M}_{l,l-1}}{\mathcal{M}_{l,l}}-\frac{\mathcal{M}_{l,l-2}}{\mathcal{M}_{l,l}}\mathcal{N}_{l-2,l-1}\right)\mathcal{N}_{l-1,l} }.\label{nene5}\end{split} \end{equation} 
 
Choosing $l=p-1$, we have
\begin{equation}\mbox{I}_{N^{p}} Y_{p-1}+ \mathcal{N}_{p-1,p}Y_{p}+ \mathcal{N}_{p-1,p+1}Y_{p+1}=\tilde{V}_{p-1},\label{gringooo1}
\end{equation}  
where we have  used the following matrices
\begin{equation} \begin{split}& \quad \mathcal{N}_{p-1,p}=\frac{\frac{\mathcal{M}_{p-1,p}}{\mathcal{M}_{p-1,p-1}}-\left(\frac{\mathcal{M}_{p-1,p-2}}{\mathcal{M}_{p-1,p-1}}-\frac{\mathcal{M}_{p-1,p-3}}{\mathcal{M}_{p-1,p-1}}\mathcal{N}_{p-3,p-2}\right)\mathcal{N}_{p-2,p} }{\mbox{I}_{N^{p}} -\frac{\mathcal{M}_{p-1,p-3}}{\mathcal{M}_{p-1,p-1}}\mathcal{N}_{p-3,p-1}-\left(\frac{\mathcal{M}_{p-1,p-2}}{\mathcal{M}_{p-1,p-1}}-\frac{\mathcal{M}_{p-1,p-3}}{\mathcal{M}_{p-1,p-1}}\mathcal{N}_{p-3,p-2}\right)\mathcal{N}_{p-2,p-1} } ,\\&   \mathcal{N}_{p-1,p+1}=\frac{\frac{\mathcal{M}_{p-1,p+1}}{\mathcal{M}_{p-1,p-1}} }{\mbox{I}_{N^{p}} -\frac{\mathcal{M}_{p-1,p-3}}{\mathcal{M}_{p-1,p-1}}\mathcal{N}_{p-3,p-1}-\left(\frac{\mathcal{M}_{p-1,p-2}}{\mathcal{M}_{p-1,p-1}}-\frac{\mathcal{M}_{p-1,p-3}}{\mathcal{M}_{p-1,p-1}}\mathcal{N}_{p-3,p-2}\right)\mathcal{N}_{p-2,p-1} },\\&  \quad\quad  \tilde{V}_{p-1}=\frac{-\frac{\mathcal{M}_{p-1,p-3}}{\mathcal{M}_{p-1,p-1}}\tilde{V}_{p-3}-\left(\frac{\mathcal{M}_{p-1,p-2}}{\mathcal{M}_{p-1,p-1}}-\frac{\mathcal{M}_{p-1,p-3}}{\mathcal{M}_{p-1,p-1}}\mathcal{N}_{p-3,p-2}\right)\tilde{V}_{p-2} }{\mbox{I}_{N^{p}} -\frac{\mathcal{M}_{p-1,p-3}}{\mathcal{M}_{p-1,p-1}}\mathcal{N}_{p-3,p-1}-\left(\frac{\mathcal{M}_{p-1,p-2}}{\mathcal{M}_{p-1,p-1}}-\frac{\mathcal{M}_{p-1,p-3}}{\mathcal{M}_{p-1,p-1}}\mathcal{N}_{p-3,p-2}\right)\mathcal{N}_{p-2,p-1} }.\label{nene4441}\end{split} \end{equation} 
 
Choosing  $l=p$, we obtain
 \begin{equation}\mbox{I}_{N^{p}} Y_{p}+ \mathcal{N}_{p,p+1}Y_{p+1}=\tilde{V}_{p},\label{gringooo2}
\end{equation}   
where we have used the following matrices
\begin{equation} \begin{split}&  \mathcal{N}_{p,p+1}=\frac{\frac{\mathcal{M}_{p,p+1}}{\mathcal{M}_{p,p}}-\left(\frac{\mathcal{M}_{p,p-1}}{\mathcal{M}_{p,p}}-\frac{\mathcal{M}_{p,p-2}}{\mathcal{M}_{p,p}}\mathcal{N}_{p-2,p-1}\right)\mathcal{N}_{p-1,p+1} }{\mbox{I}_{N^{p}} -\frac{\mathcal{M}_{p,p-2}}{\mathcal{M}_{p,p}}\mathcal{N}_{p-2,p}-\left(\frac{\mathcal{M}_{p,p-1}}{\mathcal{M}_{p,p}}-\frac{\mathcal{M}_{p,p-2}}{\mathcal{M}_{p,p}}\mathcal{N}_{p-2,p-1}\right)\mathcal{N}_{p-1,p} } ,\\&  \quad\quad     \tilde{V}_{p}=\frac{-\frac{\mathcal{M}_{p,p-2}}{\mathcal{M}_{p,p}}\tilde{V}_{p-2}-\left(\frac{\mathcal{M}_{p,p-1}}{\mathcal{M}_{p,p}}-\frac{\mathcal{M}_{p,p-2}}{\mathcal{M}_{p,p}}\mathcal{N}_{p-2,p-1}\right)\tilde{V}_{p-1} }{\mbox{I}_{N^{p}} -\frac{\mathcal{M}_{p,p-2}}{\mathcal{M}_{p,p}}\mathcal{N}_{p-2,p}-\left(\frac{\mathcal{M}_{p,p-1}}{\mathcal{M}_{p,p}}-\frac{\mathcal{M}_{p,p-2}}{\mathcal{M}_{p,p}}\mathcal{N}_{p-2,p-1}\right)\mathcal{N}_{p-1,p} }.\label{nene4442}\end{split} \end{equation} 

Choosing  $l=p+1$, we obtain
 \begin{equation}\mbox{I}_{N^{p}} Y_{p+1}=\tilde{V}_{p+1},\label{gringooo3}
\end{equation} 
where we have used the following matrix
 \begin{equation}\tilde{V}_{p+1}=\frac{-\frac{\mathcal{M}_{p+1,p-1}}{\mathcal{M}_{p+1,p+1}}\tilde{V}_{p-1}-\left(\frac{\mathcal{M}_{p+1,p}}{\mathcal{M}_{p+1,p+1}}-\frac{\mathcal{M}_{p+1,p-1}}{\mathcal{M}_{p+1,p+1}}\mathcal{N}_{p-1,p}\right)\tilde{V}_{p} }{\mbox{I}_{N^{p}} -\frac{\mathcal{M}_{p+1,p-1}}{\mathcal{M}_{p+1,p+1}}\mathcal{N}_{p-1,p+1}-\left(\frac{\mathcal{M}_{p+1,p}}{\mathcal{M}_{p+1,p+1}}-\frac{\mathcal{M}_{p+1,p-1}}{\mathcal{M}_{p+1,p+1}}\mathcal{N}_{p-1,p}\right)\mathcal{N}_{p,p+1} }.\label{gringooo31}
\end{equation}
 
 We summarize  (\ref{gringo1}), (\ref{nene1}), (\ref{gringo2}),  (\ref{nene2}), (\ref{gringo4se}), (\ref{nene4}),  (\ref{gringooo}), (\ref{gringooo1}), (\ref{gringooo3}) and  (\ref{gringooo31}) together with the recurrences  (\ref{nene5}), (\ref{nene4441}) and (\ref{nene4442}). These iterative  computations  simplifies  (\ref{beb120}), which is equivalent to 
 \begin{equation}\begin{pmatrix} \mbox{I}_{N^{p}} & \mathcal{N}_{1,2} & \mathcal{N}_{1,3} & \mbox{O}_{N^{p}} & \mbox{O}_{N^{p}} &\dots & \mbox{O}_{N^{p}} & \mbox{O}_{N^{p}} & \mbox{O}_{N^{p}} & \mbox{O}_{N^{p}} \\ \mbox{O}_{N^{p}} & \mbox{I}_{N^{p}} & \mathcal{N}_{2,3} &\mathcal{N}_{2,4} & \mbox{O}_{N^{p}} &\dots & \mbox{O}_{N^{p}} & \mbox{O}_{N^{p}} & \mbox{O}_{N^{p}} & \mbox{O}_{N^{p}} \\ \mbox{O}_{N^{p}} & \mbox{O}_{N^{p}} & \mbox{I}_{N^{p}} & \mathcal{N}_{3,4} & \mathcal{N}_{3,5} &\dots & \mbox{O}_{N^{p}} & \mbox{O}_{N^{p}} & \mbox{O}_{N^{p}} & \mbox{O}_{N^{p}} \\ \mbox{O}_{N^{p}} & \mbox{O}_{N^{p}} & \mbox{O}_{N^{p}} & \mbox{I}_{N^{p}} & \mathcal{N}_{4,5} &\dots & \mbox{O}_{N^{p}} & \mbox{O}_{N^{p}} & \mbox{O}_{N^{p}} & \mbox{O}_{N^{p}}   \\ \vdots & \vdots & \vdots & \vdots & \vdots & \ddots & \vdots & \vdots & \vdots & \vdots \\ \mbox{O}_{N^{p}} & \mbox{O}_{N^{p}} & \mbox{O}_{N^{p}} & \mbox{O}_{N^{p}} & \mbox{O}_{N^{p}} &\dots & \mbox{I}_{N^{p}} & \mathcal{N}_{p-2,p-1} & \mathcal{N}_{p-2,p} & \mbox{O}_{N^{p}} \\ \mbox{O}_{N^{p}} & \mbox{O}_{N^{p}} & \mbox{O}_{N^{p}} & \mbox{O}_{N^{p}} & \mbox{O}_{N^{p}} &\dots & \mbox{O}_{N^{p}} & \mbox{I}_{N^{p}} & \mathcal{N}_{p-1,p} & \mathcal{N}_{p-1,p+1} \\ \mbox{O}_{N^{p}} & \mbox{O}_{N^{p}} & \mbox{O}_{N^{p}} & \mbox{O}_{N^{p}} & \mbox{O}_{N^{p}} &\dots & \mbox{O}_{N^{p}} & \mbox{O}_{N^{p}} & \mbox{I}_{N^{p}} & \mathcal{N}_{p,p+1} \\ \mbox{O}_{N^{p}} & \mbox{O}_{N^{p}} & \mbox{O}_{N^{p}} & \mbox{O}_{N^{p}} & \mbox{O}_{N^{p}} &\dots & \mbox{O}_{N^{p}} & \mbox{O}_{N^{p}} & \mbox{O}_{N^{p}} & \mbox{I}_{N^{p}}
\end{pmatrix}\begin{pmatrix} Y_{1}\\ Y_{2}\\ Y_{3}\\ Y_{4}\\  \vdots \\ Y_{p-2} \\ Y_{p-1}\\ Y_{p}\\ Y_{p+1}
\end{pmatrix}=\begin{pmatrix} \tilde{V}_{1} \\ \tilde{V}_{2} \\ \tilde{V}_{3} \\ \tilde{V}_{4}\\    \vdots    \\ \tilde{V}_{p-2}\\ \tilde{V}_{p-1}\\ \tilde{V}_{p}\\ \tilde{V}_{p+1}
\end{pmatrix}.\label{beb1se1}
\end{equation} 

Now, we  show that  (\ref{beb1se1}) has unique solution. Moreover, we  compute this solution inductively    in (\ref{beb1se1}) as follows.    The last equation in (\ref{beb1se1}) is equivalent to (\ref{gringooo3}). Making substraction from the next equation in (\ref{beb1se1}), we obtain
\begin{equation}
Y_{p}=\tilde{V}_{p}-\mathcal{N}_{p,p+1}\tilde{V}_{p+1}.\label{yx1}
\end{equation}

Substracting (\ref{gringooo3}) from the next equation in (\ref{beb1se1}), we obtain
\begin{equation}
Y_{p-1}=\tilde{V}_{p-1}-\mathcal{N}_{p-1,p}\left( \tilde{V}_{p}-\mathcal{N}_{p,p+1}\tilde{V}_{p+1}\right)-   \mathcal{N}_{p-1,p+1}\tilde{V}_{p+1}.\label{yx2}
\end{equation}

Walking backwards among the equations  of (\ref{beb1se1}), we obtain
\begin{equation}
Y_{1}=\tilde{V}_{1}-\mathcal{N}_{1,2}Y_{2}-\mathcal{N}_{1,3}Y_{3}.\label{yx3}
\end{equation}

Now, the system of equations (\ref{beb1}) is solved.
\subsection{Remarks} It remains    to better explain   the  invertibility of those matrices    occurring   in (\ref{nene1}), (\ref{nene2}), (\ref{nene3}), (\ref{nene4}), (\ref{nene5}), (\ref{nene4441}), (\ref{nene4442}),  (\ref{gringooo31}), because this fact was skipped throughout the previous computations. This is achieved by making convenient  estimations for the solution of (\ref{beb1})  with respect to  the matrix $\Lambda$,  from where we will obtain the intertibilities of those matrices. 

Before beginning, it is   introduced the following norm
\begin{equation}\left\|V\right\|, \label{norma}
\end{equation}
which is just the standard definition   of the   norm of $V$  as operator.  

  We introduce also by (\ref{yy}) the following matrix
\begin{equation} \tilde{Q}\left[I,J\right]=\begin{pmatrix} \frac{\gamma_{1,1}\left[I;J\right]}{\gamma_{k,k}\left[I;J\right]}-1 & \dots & \frac{\gamma_{1,k}\left[I;J\right]}{\gamma_{k,k}\left[I;J\right]} & \dots & \frac{\gamma_{1,N}\left[I;J\right]}{\gamma_{k,k}\left[I;J\right]} \\ \vdots & \ddots &\vdots& \ddots &\vdots \\ \frac{\gamma_{k,1}\left[I;J\right]}{\gamma_{k,k}\left[I;J\right]} & \dots & 0 & \dots &\frac{\gamma_{k,N}\left[I;J\right]}{\gamma_{k,k}\left[I;J\right]}\\ \vdots & \ddots &\vdots& \ddots &\vdots \\ \frac{\gamma_{N,1}\left[I;J\right]}{\gamma_{k,k}\left[I;J\right]} & \dots & \frac{\gamma_{N,k}\left[I;J\right]}{\gamma_{k,k}\left[I;J\right]} & \dots & \frac{\gamma_{N,N}\left[I;J\right]}{\gamma_{k,k}\left[I;J\right]}-1
\end{pmatrix}.\label{shobi11}
\end{equation}

According to (\ref{shobi11}), we clearly have
$$Q\left[I,J\right]=\gamma_{k,k}\left[I;J\right]I_{N}\left(I_{N}+\tilde{Q}\left[I,J\right] \right), 
$$
when by (\ref{yy}) the following holds
$$\gamma_{k,k}\left[I;J\right]=\displaystyle\max_{1\leq j,t\leq N}\gamma_{l,t}\left[I;J\right],\quad\mbox{where $k\in 1,\dots,N$.}$$

Because the entries, of this matrix $\tilde{Q}\left[I,J\right]$, are strictly positive numbers  less than $1$, excepting the middle entry that is $0$, we obtain
\begin{equation}\left\|\tilde{Q}\left[I,J\right]\right\|<1,\label{chiu}
\end{equation}
concluding the  invertibility of the matrix from (\ref{shobi1}), and  also the invertibility of the matrix $\mathcal{M}_{1,1}$. 

 It follows  that (\ref{nene1}) has sense, since the other remaining situations may be similarly approached eventually switching its columns, but it remains to show the following
\bl\label{lima}Let $A=\left(a_{ij}\right)_{1\leq i,j\leq N}$ and $B=\left(b_{ij}\right)_{1\leq i,j\leq N}$ two matrices with positive non-vanishing entries such that
\begin{equation*}\det B\neq 0,\quad \frac{a_{ij}}{b_{ij}}<x,\quad\mbox{for some $x>0$ and for all  $i,j=1,\dots, N$.}\label{ix1}
\end{equation*}

Then, we have
\begin{equation}
\left\|\frac{A}{B}\right\|<x.\label{ix}
\end{equation}
 \el
\begin{proof} Clearly,   (\ref{ix}) is concluded, because (\ref{ix1}) implies
 \begin{equation*}\frac{A}{B}   \leq xI_{n},\quad\mbox{because $A \leq xI_{n} B$,}
\end{equation*}
which is obvious and thus (\ref{ix}) becomes clear. 
\end{proof}

We have now to show that also (\ref{nene2}) has sense. It remains just to  conclude    invertibility for the following matrix
\begin{equation} I_{N^{p}}-\frac{\mathcal{M}_{2,1}}{\mathcal{M}_{2,2}}\mathcal{N}_{1,2},\label{vulep1}
\end{equation}
or equivalently  existence for the following matrix
\begin{equation}\frac{1}{I_{N^{p}}-\frac{\mathcal{M}_{2,1}}{\mathcal{M}_{2,2}}\mathcal{N}_{1,2}}.\label{tigru}
\end{equation}
 
Obviously, any of our  matrices may be written as product of more ,,simple'' matrices similarly as we were arguing  around (\ref{90000se1extra}),(\ref{90000}),(\ref{90000se}) and (\ref{90000se1}). This fact is useful in order to make convenient estimations by (\ref{norma})  of the norm     for  each of the following matrices
\begin{equation}\frac{\mathcal{M}_{2,1}}{\mathcal{M}_{2,2}}\mathcal{N}_{1,2},\quad\frac{\mathcal{M}_{2,1}}{\mathcal{M}_{2,2}},\quad\mathcal{N}_{1,2},\quad\mathcal{M}_{2,1},\quad \mathcal{M}_{2,2}.\label{girafa}
\end{equation}

 Recalling now (\ref{lil}) and (\ref{nene1}), it suffices to consider the following matrices
\begin{equation*} I_{N^{p}}-\left[\frac{\mathcal{W}^{(0)}_{1}+\mathcal{V}_{1}+\mathcal{O}_{1}''}{\mathcal{O}^{(0)}_{2}+\mathcal{O}_{2}+\mathcal{V}_{2}'+\mathcal{W}_{2}''+\left(N+4\lambda_{1}^{2}+\dots\right)I_{N^{p}} }\frac{\mathcal{V}^{(0)}_{2}+\mathcal{W}_{2}+\mathcal{O}_{2}'}{\mathcal{O}^{(0)}_{1}+\mathcal{O}_{1}+\mathcal{V}_{1}'+\mathcal{W}_{1}''+\left(N+4\lambda_{1}^{2}+\dots\right)I_{N^{p}} }\right]\left[i_{1},\dots,i_{N};j_{1},\dots,j_{N}\right], 
\end{equation*}
where $i_{1}+\dots+i_{N}+j_{1}+\dots+j_{N}=p\geq 3$.

The  product of  all these matrices defines (\ref{vulep1}) similarly as   there are defined  (\ref{90000se1extra}), (\ref{90000}), (\ref{90000se}), (\ref{90000se1})  recalling (\ref{lil}), (\ref{lili}), (\ref{lili2}), (\ref{lili3}), (\ref{suc1}), (\ref{suc2}), (\ref{sucL1}), (\ref{sucL11}), (\ref{sucL11sese}), concluding the following
\begin{equation}
\left\|\frac{\mathcal{M}_{2,1}}{\mathcal{M}_{2,2}}\mathcal{N}_{1,2}\right\|<\frac{1}{2}. \label{2211}
\end{equation}
which gives sense to (\ref{tigru}) in the light of the following expansion
\begin{equation} \frac{1}{I_{N^{p}}-\frac{\mathcal{M}_{2,1}}{\mathcal{M}_{2,2}}\mathcal{N}_{1,2}}=I_{N^{p}}+\frac{\mathcal{M}_{2,1}}{\mathcal{M}_{2,2}}\mathcal{N}_{1,2}+\left(\frac{\mathcal{M}_{2,1}}{\mathcal{M}_{2,2}}\mathcal{N}_{1,2}\right)^{2}+\dots.\label{urs1}
\end{equation}

We avoid to explain (\ref{2211})  at this moment, because this will be explained throughout further computations, showing that (\ref{nene5}) has sense by concluding invertibility for the following matrices
\begin{equation} I_{N^{p}}-\frac{\mathcal{M}_{l,l-2}}{\mathcal{M}_{l,l}}\mathcal{N}_{l-2,l}-\left(\frac{\mathcal{M}_{l,l-1}}{\mathcal{M}_{l,l}}-\frac{\mathcal{M}_{l,l-2}}{\mathcal{M}_{l,l}}\mathcal{N}_{l-2,l-1}\right)\mathcal{N}_{l-1,l},\quad\mbox{for all $l=3,\dots,p-1$.}   \label{vulep2}
\end{equation}

We want to make convenient estimations of the norm (\ref{norma})     for each of the following matrices
\begin{equation}\frac{\mathcal{M}_{l,l-1}}{\mathcal{M}_{l,l}},\quad \frac{\mathcal{M}_{l,l-2}}{\mathcal{M}_{l,l}},\quad\mathcal{M}_{l,l-2},\quad\mathcal{M}_{l-1,l},\quad  \mathcal{M}_{l,l},\quad  \mathcal{N}_{l-1,l},\quad\mathcal{N}_{l-2,l-1},\quad\mathcal{N}_{l-2,l},\quad\mbox{for all $l=3,\dots,p-1$.} \label{asia}
\end{equation}

It is required again to recall the factorizations using the simple matrices from (\ref{90000se1extra}), (\ref{90000}), (\ref{90000se}), (\ref{90000se1}), concluding that each of  the above matrices can be written as   products of simple matrices as in (\ref{90000se1extra}), (\ref{90000}), (\ref{90000se}), (\ref{90000se1}). 

Recalling now (\ref{lil}) and (\ref{nene1}), it suffices to consider the following matrices
\begin{equation}  I_{N^{p}}-\left[\frac{ \mathcal{V}_{l-2}''\mathcal{N}_{l-2,l}-\left(\mathcal{W}^{(0)}_{l-1}+\mathcal{V}_{l-1}+\mathcal{O}_{l-1}''\right)\mathcal{N}_{l-1,l}+\mathcal{V}_{l-2}''\mathcal{N}_{l-2,l-1}\mathcal{N}_{l-1,l}}{\mathcal{O}^{(0)}_{l}+\mathcal{O}_{l}+\mathcal{V}_{l}'+\mathcal{W}_{l}''+\left(N+4\lambda_{1}^{2}+\dots+4\lambda_{N}^{2}\right)I_{N^{p}} }  \right]\left[i_{1},\dots,i_{N};j_{1},\dots,j_{N}\right], \label{bib}
\end{equation}
where $i_{1}+\dots+i_{N}+j_{1}+\dots+j_{N}=p\geq 3$. 

The  product of all these matrices (\ref{bib})  defines (\ref{vulep2}) similarly as there are defined (\ref{90000se1extra}), (\ref{90000}), (\ref{90000se}), (\ref{90000se1}). We recall also (\ref{lil}), (\ref{lili}), (\ref{lili2}), (\ref{lili3}), (\ref{suc1}), (\ref{suc2}), (\ref{sucL1}), (\ref{sucL11}), (\ref{sucL11sese}) and (\ref{lili}), dealing with certain sets of  positive integers and notations
\begin{equation}\left\{\begin{split}&\left\{i_{l}\right\}_{l=1,\dots,N},\quad\hspace{0.12 cm}\left\{j_{l}\right\}_{l=1,\dots,N}
,\\& \left\{\tilde{i}_{l}\right\}_{l=1,\dots,N},\quad\hspace{0.07 cm}\left\{\tilde{j}_{l}\right\}_{l=1,\dots,N},\\& \left\{\tilde{\tilde{i}}_{l}\right\}_{l=1,\dots,N},\quad\left\{\tilde{\tilde{j}}_{l}\right\}_{l=1,\dots,N},\end{split}\right.\quad\quad\mbox{such that:}\quad \left\{\begin{split}&\displaystyle\sum_{l=1}^{N} i_{l}=p+1-l,\hspace{0.25 cm}\displaystyle\sum_{l=1}^{N} j_{l}=l-1 ,\\& \displaystyle\sum_{l=1}^{N} \tilde{i}_{l}=p-l ,\quad\hspace{0.1 cm}\hspace{0.1 cm}\quad\displaystyle\sum_{l=1}^{N} \tilde{j}_{l}=l,\\&  \displaystyle\sum_{l=1}^{N} \tilde{\tilde{i}}_{l}=p-l-1 ,\quad\displaystyle\sum_{l=1}^{N} \tilde{\tilde{j}}_{l}=l+1.
\end{split}\right.  \label{liliwww1}
\end{equation}

More precisely, we consider by (\ref{yy})  the following  multi-indexes \begin{equation}
\begin{split}& I=\left(i_{1},i_{2},\dots,i_{N}\right)\in\mathbb{N}^{N},\quad\hspace{0.12 cm} \tilde{I}=\left(\tilde{i}_{1},\tilde{i}_{2},\dots,\tilde{i}_{N}\right)\in\mathbb{N}^{N},\quad\hspace{0.09 cm} \tilde{\tilde{I}}=\left(\tilde{\tilde{i}}_{1},\tilde{\tilde{i}}_{2},\dots,\tilde{\tilde{i}}_{N}\right)\in\mathbb{N}^{N}, \\& J=\left(j_{1},j_{2},\dots,j_{N}\right)\in\mathbb{N}^{N},\quad  \tilde{J}=\left(\tilde{j}_{1},\tilde{j}_{2},\dots,\tilde{j}_{N}\right)\in\mathbb{N}^{N}, \quad \tilde{\tilde{J}}=\left(\tilde{\tilde{j}}_{1},\tilde{\tilde{j}}_{2},\dots,\tilde{\tilde{j}}_{N}\right)\in\mathbb{N}^{N},   \end{split}\label{indexes}
\end{equation}
respecting appropriately (\ref{liliwww1}).

We recall (\ref{shobi}), (\ref{shobix1})  and (\ref{shobiy2}), wishing to better understand the following   matrices
 $$\frac{\mathcal{M}_{1,3}}{\mathcal{M}_{1,1}},\quad\frac{\mathcal{M}_{3,1}}{\mathcal{M}_{3,3}},\quad\mbox{or}\quad \frac{\mathcal{M}_{1,3}}{\mathcal{M}_{3,3}},\quad\frac{\mathcal{M}_{3,1}}{\mathcal{M}_{1,1}}.$$
 
 We   reconsider the strategy applied in the light of (\ref{90000se1extra}), (\ref{90000}),  (\ref{90000se}),  and (\ref{90000se1}). Thus, according to  (\ref{liliwww1})  and  (\ref{indexes}), we have to focus by (\ref{yy})  on the following matrices
 \begin{equation} \mathcal{M}_{N^{p}\times N^{p}}\left(\mathbb{C}\right)\ni\frac{ \begin{pmatrix}  1 & \dots   & 0 & \dots & 0 & \dots & 0& \dots & 0  \\ \vdots & \ddots & \vdots & \ddots & \vdots & \ddots & \vdots& \ddots & \vdots   \\ 0 & \dots &  \tilde{\gamma}_{1,1}^{''}\left[I;J\right]   & \dots &  \tilde{\gamma}_{1,k}^{''}\left[I;J\right]   & \dots &  \tilde{\gamma}_{1,N}^{''}\left[I;J\right] & \dots& 0 
\\ \vdots & \ddots & \vdots  & \ddots & \vdots  & \ddots & \vdots  & \ddots & \vdots  \\ 0 & \dots &   \tilde{\gamma}_{k,1}^{''}\left[I;J\right]   & \dots &  \tilde{\gamma}_{k,k}^{''}\left[I;J\right]    & \dots &   \tilde{\gamma}_{k,N}^{''}\left[I;J\right] & \dots & 0 \\ \vdots  & \ddots & \vdots  & \ddots & \vdots  & \ddots & \vdots  & \ddots & \vdots  \\ 0 & \dots &   \tilde{\gamma}_{N,1}^{''}\left[I,J\right]   & \dots &   \tilde{\gamma}_{N,k}^{''}\left[I;J\right]   & \dots &   \tilde{\gamma}_{N,N}^{''}\left[I;J\right]  & \dots & 0 \\ \vdots  & \ddots & \vdots  & \ddots & \vdots  & \ddots & \vdots  & \ddots & \vdots  \\ 0  & \dots & 0  & \dots &0  & \dots & 0  & \dots & 1  
\end{pmatrix}}{ \begin{pmatrix}  1 & \dots   & 0 & \dots & 0 & \dots & 0& \dots & 0  \\ \vdots & \ddots & \vdots & \ddots & \vdots & \ddots & \vdots& \ddots & \vdots   \\ 0 & \dots &  \gamma_{1,1}\left[I;J\right]   & \dots &  \gamma_{1,k}\left[I;J\right]   & \dots &  \gamma_{1,N}\left[I;J\right] & \dots& 0 
\\ \vdots & \ddots & \vdots  & \ddots & \vdots  & \ddots & \vdots  & \ddots & \vdots  \\ 0 & \dots &   \gamma_{k,1}\left[I;J\right]   & \dots &   \gamma_{k,k}\left[I;J\right]      & \dots &   \gamma_{k,N}\left[I;J\right] & \dots & 0 \\ \vdots  & \ddots & \vdots  & \ddots & \vdots  & \ddots & \vdots  & \ddots & \vdots  \\ 0 & \dots &   \gamma_{N,1}\left[I,J\right]   & \dots &   \gamma_{N,k}\left[I;J\right]    & \dots &   \gamma_{N,N}\left[I;J\right]    & \dots & 0 \\ \vdots  & \ddots & \vdots  & \ddots & \vdots  & \ddots & \vdots  & \ddots & \vdots  \\ 0  & \dots & 0  & \dots &0  & \dots & 0  & \dots & 1  
\end{pmatrix}}, \label{shobiV1}
\end{equation}
where  (\ref{lili}) is satisfied, for all $k=3,\dots,p+1$.

Respectively, we have to focus by (\ref{yy})  on the following matrices
 \begin{equation} \mathcal{M}_{N^{p}\times N^{p}}\left(\mathbb{C}\right)\ni \frac{\begin{pmatrix}  1 & \dots   & 0 & \dots & 0 & \dots & 0& \dots & 0  \\ \vdots & \ddots & \vdots & \ddots & \vdots & \ddots & \vdots& \ddots & \vdots   \\ 0 & \dots &  \gamma_{1,1}^{''}\left[\tilde{\tilde{I}};\tilde{\tilde{J}}\right]  & \dots & \gamma_{1,k}^{''}\left[\tilde{\tilde{I}};\tilde{\tilde{J}}\right]   & \dots &  \gamma_{1,N}^{''}\left[\tilde{\tilde{I}};\tilde{\tilde{J}}\right] & \dots& 0 
\\ \vdots & \ddots & \vdots  & \ddots & \vdots  & \ddots & \vdots  & \ddots & \vdots  \\ 0 & \dots &   \gamma_{k,1}^{''}\left[\tilde{\tilde{I}};\tilde{\tilde{J}}\right]   & \dots &   \gamma_{k,k}^{''}\left[\tilde{\tilde{I}};\tilde{\tilde{J}}\right]   & \dots &   \gamma_{k,N}^{''}\left[\tilde{\tilde{I}};\tilde{\tilde{J}}\right] & \dots & 0 \\ \vdots  & \ddots & \vdots  & \ddots & \vdots  & \ddots & \vdots  & \ddots & \vdots  \\ 0 & \dots &  \gamma_{N,1}^{''}\left[\tilde{\tilde{I}};\tilde{\tilde{J}}\right]   & \dots &   \gamma_{N,k}^{''}\left[\tilde{\tilde{I}};\tilde{\tilde{J}}\right]    & \dots &   \gamma_{N,N}^{''}\left[\tilde{\tilde{I}};\tilde{\tilde{J}}\right]  & \dots & 0 \\ \vdots  & \ddots & \vdots  & \ddots & \vdots  & \ddots & \vdots  & \ddots & \vdots  \\ 0  & \dots & 0  & \dots &0  & \dots & 0  & \dots & 1  
\end{pmatrix}}{\begin{pmatrix}  1 & \dots   & 0 & \dots & 0 & \dots & 0& \dots & 0  \\ \vdots & \ddots & \vdots & \ddots & \vdots & \ddots & \vdots& \ddots & \vdots   \\ 0 & \dots &  \gamma_{1,1}\left[\tilde{\tilde{I}};\tilde{\tilde{J}}\right]   & \dots &  \gamma_{1,k}\left[\tilde{\tilde{I}};\tilde{\tilde{J}}\right]   & \dots &  \gamma_{1,N}\left[\tilde{\tilde{I}};\tilde{\tilde{J}}\right] & \dots& 0 
\\ \vdots & \ddots & \vdots  & \ddots & \vdots  & \ddots & \vdots  & \ddots & \vdots  \\ 0 & \dots &   \gamma_{k,1}\left[\tilde{\tilde{I}};\tilde{\tilde{J}}\right]   & \dots &   \gamma_{k,k}\left[\tilde{\tilde{I}};\tilde{\tilde{J}}\right]      & \dots &   \gamma_{k,N}\left[\tilde{\tilde{I}};\tilde{\tilde{J}}\right] & \dots & 0 \\ \vdots  & \ddots & \vdots  & \ddots & \vdots  & \ddots & \vdots  & \ddots & \vdots  \\ 0 & \dots &   \gamma_{N,1}\left[\tilde{\tilde{I}};\tilde{\tilde{J}}\right]   & \dots &   \gamma_{N,k}\left[\tilde{\tilde{I}};\tilde{\tilde{J}}\right]    & \dots &   \gamma_{N,N}\left[\tilde{\tilde{I}};\tilde{\tilde{J}}\right]    & \dots & 0 \\ \vdots  & \ddots & \vdots  & \ddots & \vdots  & \ddots & \vdots  & \ddots & \vdots  \\ 0  & \dots & 0  & \dots &0  & \dots & 0  & \dots & 1  
\end{pmatrix}}, \label{shobiV2}
\end{equation}
where  (\ref{lili}) is satisfied, for all $k=1,\dots,p-1$.

We recall again (\ref{shobi}), (\ref{shobix1})  and (\ref{shobiy2}), wishing to better understand the following two matrices
 $$ \frac{\mathcal{M}_{2,3}}{\mathcal{M}_{2,2}},\quad \frac{\mathcal{M}_{3,2}}{\mathcal{M}_{3,3}},\quad\mbox{or}\quad \frac{\mathcal{M}_{2,3}}{\mathcal{M}_{3,3}},\quad \frac{\mathcal{M}_{3,2}}{\mathcal{M}_{2,2}}.$$
 
 We   reconsider the strategy applied in the light of (\ref{90000se1extra}), (\ref{90000}),  (\ref{90000se}),  and (\ref{90000se1}). Thus, according to  (\ref{liliwww1})  and  (\ref{indexes}), we have to focus by (\ref{yy})  on the following matrices
 \begin{equation} \mathcal{M}_{N^{p}\times N^{p}}\left(\mathbb{C}\right)\ni \frac{\begin{pmatrix}  1 & \dots   & 0 & \dots & 0 & \dots & 0& \dots & 0  \\ \vdots & \ddots & \vdots & \ddots & \vdots & \ddots & \vdots& \ddots & \vdots   \\ 0 & \dots &  \tilde{\gamma}_{1,1}^{'}\left[\tilde{I};\tilde{J}\right]   & \dots &  \tilde{\gamma}_{1,k}^{'}\left[\tilde{I};\tilde{J}\right]   & \dots &  \tilde{\gamma}_{1,N}^{'}\left[\tilde{I};\tilde{J}\right] & \dots& 0 
\\ \vdots & \ddots & \vdots  & \ddots & \vdots  & \ddots & \vdots  & \ddots & \vdots  \\ 0 & \dots &   \tilde{\gamma}_{k,1}^{'}\left[\tilde{I};\tilde{J}\right]   & \dots &   \tilde{\gamma}_{k,k}^{'}\left[\tilde{I};\tilde{J}\right]     & \dots &   \tilde{\gamma}_{k,N}^{'}\left[\tilde{I};\tilde{J}\right] & \dots & 0 \\ \vdots  & \ddots & \vdots  & \ddots & \vdots  & \ddots & \vdots  & \ddots & \vdots  \\ 0 & \dots &   \tilde{\gamma}_{N,1}^{'}\left[\tilde{I};\tilde{J}\right]   & \dots &  \tilde{\gamma}_{N,k}^{'}\left[\tilde{I};\tilde{J}\right]   & \dots &   \tilde{\gamma}_{N,N}^{'}\left[\tilde{I};\tilde{J}\right]   & \dots & 0 \\ \vdots  & \ddots & \vdots  & \ddots & \vdots  & \ddots & \vdots  & \ddots & \vdots  \\ 0  & \dots & 0  & \dots &0  & \dots & 0  & \dots & 1  
\end{pmatrix}}{\begin{pmatrix}  1 & \dots   & 0 & \dots & 0 & \dots & 0& \dots & 0  \\ \vdots & \ddots & \vdots & \ddots & \vdots & \ddots & \vdots& \ddots & \vdots   \\ 0 & \dots &  \gamma_{1,1}\left[\tilde{I};\tilde{J}\right]   & \dots &  \gamma_{1,k}\left[\tilde{I};\tilde{J}\right]   & \dots &  \gamma_{1,N}\left[\tilde{I};\tilde{J}\right] & \dots& 0 
\\ \vdots & \ddots & \vdots  & \ddots & \vdots  & \ddots & \vdots  & \ddots & \vdots  \\ 0 & \dots &  \gamma_{k,1}\left[\tilde{I};\tilde{J}\right]   & \dots &   \gamma_{k,k}\left[\tilde{I};\tilde{J}\right]      & \dots &   \gamma_{k,N}\left[\tilde{I};\tilde{J}\right] & \dots & 0 \\ \vdots  & \ddots & \vdots  & \ddots & \vdots  & \ddots & \vdots  & \ddots & \vdots  \\ 0 & \dots &   \gamma_{N,1}\left[\tilde{I};\tilde{J}\right]   & \dots &   \gamma_{N,k}\left[\tilde{I};\tilde{J}\right]    & \dots &   \gamma_{N,N}\left[\tilde{I};\tilde{J}\right]   & \dots & 0 \\ \vdots  & \ddots & \vdots  & \ddots & \vdots  & \ddots & \vdots  & \ddots & \vdots  \\ 0  & \dots & 0  & \dots &0  & \dots & 0  & \dots & 1  
\end{pmatrix}}, \label{shobiV11}
\end{equation}
where  (\ref{lili}) is satisfied, for all $k=2,\dots,p+1$.

Respectively, we have to focus by (\ref{yy})  on the following matrices
 \begin{equation} \mathcal{M}_{N^{p}\times N^{p}}\left(\mathbb{C}\right)\ni \frac{\begin{pmatrix}  1 & \dots   & 0 & \dots & 0 & \dots & 0& \dots & 0  \\ \vdots & \ddots & \vdots & \ddots & \vdots & \ddots & \vdots& \ddots & \vdots   \\ 0 & \dots &  \gamma_{1,1}^{'}\left[\tilde{\tilde{I}};\tilde{\tilde{J}}\right]   & \dots &  \gamma_{1,k}^{'}\left[\tilde{\tilde{I}};\tilde{\tilde{J}}\right]   & \dots &  \gamma_{1,N}^{'}\left[\tilde{\tilde{I}};\tilde{\tilde{J}}\right] & \dots& 0 
\\ \vdots & \ddots & \vdots  & \ddots & \vdots  & \ddots & \vdots  & \ddots & \vdots  \\ 0 & \dots &   \gamma_{k,1}^{'}\left[\tilde{\tilde{I}};\tilde{\tilde{J}}\right]   & \dots &   \gamma_{k,k}^{'}\left[\tilde{\tilde{I}};\tilde{\tilde{J}}\right]     & \dots &   \gamma_{k,N}^{'}\left[\tilde{\tilde{I}};\tilde{\tilde{J}}\right] & \dots & 0 \\ \vdots  & \ddots & \vdots  & \ddots & \vdots  & \ddots & \vdots  & \ddots & \vdots  \\ 0 & \dots &   \gamma_{N,1}^{'}\left[\tilde{\tilde{I}};\tilde{\tilde{J}}\right]   & \dots &   \gamma_{N,k}^{'}\left[\tilde{\tilde{I}};\tilde{\tilde{J}}\right]    & \dots &   \gamma_{N,N}^{'}\left[\tilde{\tilde{I}};\tilde{\tilde{J}}\right]  & \dots & 0 \\ \vdots  & \ddots & \vdots  & \ddots & \vdots  & \ddots & \vdots  & \ddots & \vdots  \\ 0  & \dots & 0  & \dots &0  & \dots & 0  & \dots & 1  
\end{pmatrix}}{\begin{pmatrix}  1 & \dots   & 0 & \dots & 0 & \dots & 0& \dots & 0  \\ \vdots & \ddots & \vdots & \ddots & \vdots & \ddots & \vdots& \ddots & \vdots   \\ 0 & \dots &  \gamma_{1,1}\left[\tilde{\tilde{I}};\tilde{\tilde{J}}\right]   & \dots &  \gamma_{1,k}\left[\tilde{\tilde{I}};\tilde{\tilde{J}}\right]   & \dots &  \gamma_{1,N}\left[\tilde{\tilde{I}};\tilde{\tilde{J}}\right] & \dots& 0 
\\ \vdots & \ddots & \vdots  & \ddots & \vdots  & \ddots & \vdots  & \ddots & \vdots  \\ 0 & \dots &   \gamma_{k,1}\left[\tilde{\tilde{I}};\tilde{\tilde{J}}\right]   & \dots &   \gamma_{k,k}\left[\tilde{\tilde{I}};\tilde{\tilde{J}}\right]      & \dots &   \gamma_{k,N}\left[\tilde{\tilde{I}};\tilde{\tilde{J}}\right] & \dots & 0 \\ \vdots  & \ddots & \vdots  & \ddots & \vdots  & \ddots & \vdots  & \ddots & \vdots  \\ 0 & \dots &   \gamma_{N,1}\left[\tilde{\tilde{I}};\tilde{\tilde{J}}\right]   & \dots &   \gamma_{N,k}\left[\tilde{\tilde{I}};\tilde{\tilde{J}}\right]    & \dots &  \gamma_{N,N}\left[\tilde{\tilde{I}};\tilde{\tilde{J}}\right]    & \dots & 0 \\ \vdots  & \ddots & \vdots  & \ddots & \vdots  & \ddots & \vdots  & \ddots & \vdots  \\ 0  & \dots & 0  & \dots &0  & \dots & 0  & \dots & 1  
\end{pmatrix}}, \label{shobiV21}
\end{equation}
where  (\ref{lili}) is satisfied, for all $k=1,\dots,p$.

Now, we are ready to analyse (\ref{shobiV1}), (\ref{shobiV2}), (\ref{shobiV11}) and  (\ref{shobiV21}) in the light of Lemma \ref{lima}. Thus, we have to study by (\ref{yy})  the interactions of the following numbers
\begin{equation} \displaystyle\max_{1\leq i,j\leq N}\frac{\tilde{\gamma}_{i,j}^{''}\left[I;J\right]}{ \gamma_{i,j}\left[I;J\right]},\quad  \displaystyle\max_{1\leq i,j\leq N}\frac{ \gamma_{i,j}^{''}\left[\tilde{\tilde{I}};\tilde{\tilde{J}}\right]}{ \gamma_{i,j}\left[\tilde{\tilde{I}};\tilde{\tilde{J}}\right]},\quad  \displaystyle\max_{1\leq i,j\leq N}\frac{ \tilde{\gamma}_{i,j}^{'}\left[\tilde{I};\tilde{J}\right]}{ \gamma_{i,j}\left[\tilde{I};\tilde{J}\right]},\quad  \displaystyle\max_{1\leq i,j\leq N}\frac{ \tilde{\gamma}_{i,j}^{'}\left[\tilde{\tilde{I}};\tilde{\tilde{J}}\right]}{ \gamma_{i,j}\left[\tilde{\tilde{I}};\tilde{\tilde{J}}\right]}.\label{detalii}
\end{equation}

For instance, we can encounter the following case
\begin{equation}\displaystyle\max_{1\leq i,j\leq N}\frac{\tilde{\gamma}_{i,j}^{''}\left[I;J\right]}{ \gamma_{i,j}\left[I;J\right]}=\frac{\left(j_{1}+1\right)\left(j_{1}+2\right)\lambda_{1}^{2}}{\left(j_{1}+1\right)\left(i_{1}+1\right)+ \left(j_{1}+1\right)\left(j_{1}+2\right)\lambda_{1}^{2} +  \left(i_{1}+1\right)\left(i_{1}+2\right)\lambda_{1}^{2}}<1,\label{ixrr1} 
\end{equation}
or the following case  
\begin{equation*}\begin{split} \displaystyle\max_{1\leq i,j\leq N}\frac{\tilde{\gamma}_{i,j}^{''}\left[I;J\right]}{ \gamma_{i,j}\left[I;J\right]}& =\frac{\left(j_{k}-1\right)j_{k} \lambda_{k}^{2}}{j_{k}i_{k}  +\left(j_{k}-1\right)j_{k} \lambda_{k}^{2} +   \left(i_{k}-1\right)i_{k}  \lambda_{k}^{2}+ \left(i_{k}+j_{k}+1\right)\left(1+4\lambda_{k}^{2}\right)+\left(N+4\lambda_{1}^{2}+\dots+4\lambda_{N}^{2}\right) }\\&< \frac{\left(j_{k}-1\right)j_{k} \lambda_{k}^{2}}{j_{k}i_{k}  +\left(j_{k}-1\right)j_{k} \lambda_{k}^{2} +   \left(i_{k}-1\right)i_{k}  \lambda_{k}^{2}  } <1 .\end{split}   
\end{equation*}

We can also encounter the following case   
\begin{equation}\displaystyle\max_{1\leq i,j\leq N}\frac{ \gamma_{i,j}^{''}\left[\tilde{\tilde{I}};\tilde{\tilde{J}}\right]}{ \gamma_{i,j}\left[\tilde{\tilde{I}};\tilde{\tilde{J}}\right]}=\frac{\left(\tilde{\tilde{i}}_{1}+1\right)\left(\tilde{\tilde{i}}_{1}+2\right)\lambda_{1}^{2}}{\left(\tilde{\tilde{j}}_{1}+1\right)\left(\tilde{\tilde{i}}_{1}+1\right)+ \left(\tilde{\tilde{j}}_{1}+1\right)\left(\tilde{\tilde{j}}_{1}+2\right)\lambda_{1}^{2} +  \left(\tilde{\tilde{i}}_{1}+1\right)\left(\tilde{\tilde{i}}_{1}+2\right)\lambda_{1}^{2}}<1, \label{ixrr2} 
\end{equation}
or the following case
 \begin{equation*}\begin{split} \displaystyle\max_{1\leq i,j\leq N}\frac{ \gamma_{i,j}^{''}\left[\tilde{\tilde{I}};\tilde{\tilde{J}}\right]}{ \gamma_{i,j}\left[\tilde{\tilde{I}};\tilde{\tilde{J}}\right]} &=\frac{\left(\tilde{\tilde{i}}_{k}-1\right)\tilde{\tilde{i}}_{k} \lambda_{k}^{2}}{\tilde{\tilde{j}}_{k}\tilde{\tilde{i}}_{k}  +\left(\tilde{\tilde{j}}_{k}-1\right)\tilde{\tilde{j}}_{k} \lambda_{k}^{2} +   \left(\tilde{\tilde{i}}_{k}-1\right)\tilde{\tilde{i}}_{k}  \lambda_{k}^{2}+   \left(\tilde{\tilde{i}}_{k}+\tilde{\tilde{j}}_{k}+1\right)\left(1+4\lambda_{k}^{2}\right)+\left(N+4\lambda_{1}^{2}+\dots+4\lambda_{N}^{2}\right)}\\& < \frac{\left(\tilde{\tilde{i}}_{k}-1\right)\tilde{\tilde{i}}_{k} \lambda_{k}^{2}}{\tilde{\tilde{j}}_{k}\tilde{\tilde{i}}_{k}  +\left(\tilde{\tilde{j}}_{k}-1\right)\tilde{\tilde{j}}_{k} \lambda_{k}^{2} +   \left(\tilde{\tilde{i}}_{k}-1\right)\tilde{\tilde{i}}_{k}  \lambda_{k}^{2}} <1 .\end{split}   
\end{equation*}

We  prove that
\begin{equation}  \displaystyle\max_{1\leq i,j\leq N}\frac{\tilde{\gamma}_{i,j}^{''}\left[I;J\right]}{ \gamma_{i,j}\left[I;J\right]}\displaystyle\max_{1\leq i,j\leq N}\frac{ \gamma_{i,j}^{''}\left[\tilde{\tilde{I}};\tilde{\tilde{J}}\right]}{ \gamma_{i,j}\left[\tilde{\tilde{I}};\tilde{\tilde{J}}\right]}  < \frac{1}{7}. \label{viena1} \end{equation}
 
It suffices to prove our claim (\ref{viena1}) assuming that (\ref{ixrr1}) and (\ref{ixrr2}) hold respecting (\ref{liliwww1}), because otherwise (\ref{viena1}) follows by similar explanations  dealing with   the  simple inequalities. In this regard, we introduce the following expression

\begin{equation*}\begin{split}  E\left[I,J,\tilde{\tilde{I}},\tilde{\tilde{J}}\right]=& \left(\left(j_{1}+1\right)\left(i_{1}+1\right)+ \left(j_{1}+1\right)\left(j_{1}+2\right)\lambda_{1}^{2} +  \left(i_{1}+1\right)\left(i_{1}+2\right)\lambda_{1}^{2}\right)\left(\left(\tilde{\tilde{j}}_{1}+1\right)\left(\tilde{\tilde{i}}_{1}+1\right)\right.\\&\left.\quad\quad\quad\quad\quad\quad\quad\quad\quad\quad\quad\quad\quad\quad + \left(\tilde{\tilde{j}}_{1}+1\right)\left(\tilde{\tilde{j}}_{1}+2\right)\lambda_{1}^{2} +  \left(\tilde{\tilde{i}}_{1}+1\right)\left(\tilde{\tilde{i}}_{1}+2\right)\lambda_{1}^{2}\right),\end{split}
\end{equation*}
which may be written and developed as follows
\begin{equation}\begin{split}  E\left[I,J,\tilde{\tilde{I}},\tilde{\tilde{J}}\right]=& \left(j_{1}+1\right)\left(i_{1}+1\right)\left(\tilde{\tilde{j}}_{1}+1\right)\left(\tilde{\tilde{i}}_{1}+1\right)+\left(j_{1}+1\right)\left(i_{1}+1\right)\left(\tilde{\tilde{j}}_{1}+1\right)\left(\tilde{\tilde{j}}_{1}+2\right)\lambda_{1}^{2} +\left(j_{1}+1\right)\left(i_{1}+1\right)\left(\tilde{\tilde{i}}_{1}+1\right)\left(\tilde{\tilde{i}}_{1}+2\right)\lambda_{1}^{2}\\&   +\left(j_{1}+1\right)\left(j_{1}+2\right)\left(\tilde{\tilde{j}}_{1}+1\right)\left(\tilde{\tilde{i}}_{1}+1\right)\lambda_{1}^{2} +\left(j_{1}+1\right)\left(j_{1}+2\right) \left(\tilde{\tilde{j}}_{1}+1\right)\left(\tilde{\tilde{j}}_{1}+2\right)\lambda_{1}^{4}\\&+\left(j_{1}+1\right)\left(j_{1}+2\right) \left(\tilde{\tilde{i}}_{1}+1\right)\left(\tilde{\tilde{i}}_{1}+2\right)\lambda_{1}^{4} +\left(i_{1}+1\right)\left(i_{1}+2\right)\left(\tilde{\tilde{j}}_{1}+1\right)\left(\tilde{\tilde{i}}_{1}+1\right)\lambda_{1}^{2}\\&+\left(i_{1}+1\right)\left(i_{1}+2\right)\left(\tilde{\tilde{j}}_{1}+1\right)\left(\tilde{\tilde{j}}_{1}+2\right)\lambda_{1}^{4}  +\left(i_{1}+1\right)\left(i_{1}+2\right)\lambda_{1}^{2}\left(\tilde{\tilde{i}}_{1}+1\right)\left(\tilde{\tilde{i}}_{1}+2\right)\lambda_{1}^{2}.\end{split}\label{eee}
\end{equation}

In order to prove (\ref{viena1}), we have to recall (\ref{liliwww1})   and (\ref{indexes}). Then,  we deal with the following cases:
\begin{equation}\begin{split}&\quad\hspace{0.2 cm}\left(i_{1},j_{1}\right)=\left(\tilde{\tilde{i}}_{1},\tilde{\tilde{j}}_{1}\right),\quad\quad\hspace{0.2 cm} \left(i_{1}-1,j_{1}\right)=\left(\tilde{\tilde{i}}_{1},\tilde{\tilde{j}}_{1}\right),\quad\quad\hspace{0.23 cm} \left(i_{1}-2,j_{1}\right)=\left(\tilde{\tilde{i}}_{1},\tilde{\tilde{j}}_{1}\right), \\& \left(i_{1},j_{1}+1\right)=\left(\tilde{\tilde{i}}_{1},\tilde{\tilde{j}}_{1}\right),\quad \left(i_{1}-1,j_{1}+1\right)=\left(\tilde{\tilde{i}}_{1},\tilde{\tilde{j}}_{1}\right),\quad \left(i_{1}-2,j_{1}+1\right)=\left(\tilde{\tilde{i}}_{1},\tilde{\tilde{j}}_{1}\right), \\& \left(i_{1},j_{1}+2\right)=\left(\tilde{\tilde{i}}_{1},\tilde{\tilde{j}}_{1}\right),\quad\left(i_{1}-1,j_{1}+2\right)=\left(\tilde{\tilde{i}}_{1},\tilde{\tilde{j}}_{1}\right),\quad  \left(i_{1}-2,j_{1}+1\right)=\left(\tilde{\tilde{i}}_{1},\tilde{\tilde{j}}_{1}\right). \end{split}\label{cazuri}
\end{equation}

It suffices to show the following inequalities. We have 
\begin{equation*} \left(j_{1}+1\right)\left(j_{1}+2\right) \left(\tilde{\tilde{i}}_{1}+1\right)\left(\tilde{\tilde{i}}_{1}+2\right)\lambda_{1}^{4}\leq \frac{4\left(j_{1}+1\right)\left(i_{1}+1\right)\left(\tilde{\tilde{j}}_{1}+1\right)\left(\tilde{\tilde{i}}_{1}+1\right)}{16},
\end{equation*}
or equivalently
\begin{equation}  \left(j_{1}+2\right) \left(\tilde{\tilde{i}}_{1}+2\right)\lambda_{1}^{4}\leq \frac{4\left(i_{1}+1\right)\left(\tilde{\tilde{j}}_{1}+1\right)}{16},\label{gege1}
\end{equation}
and we  have
\begin{equation}\begin{split}\left(j_{1}+1\right)\left(j_{1}+2\right) \left(\tilde{\tilde{i}}_{1}+1\right)\left(\tilde{\tilde{i}}_{1}+2\right) \leq &  \left(j_{1}+1\right)\left(j_{1}+2\right)\left(\tilde{\tilde{j}}_{1}+1\right)\left(\tilde{\tilde{i}}_{1}+1\right) +\left(i_{1}+1\right)\left(i_{1}+2\right)\left(\tilde{\tilde{j}}_{1}+1\right)\left(\tilde{\tilde{i}}_{1}+1\right) \\&+ \left(j_{1}+1\right)\left(i_{1}+1\right)\left(\tilde{\tilde{j}}_{1}+1\right)\left(\tilde{\tilde{j}}_{1}+2\right)      +\left(j_{1}+1\right)\left(i_{1}+1\right)\left(\tilde{\tilde{i}}_{1}+1\right)\left(\tilde{\tilde{i}}_{1}+2\right)  ,\end{split}\label{gege2}\end{equation}
and we clearly have
$$\left(j_{1}+1\right)\left(j_{1}+2\right) \left(\tilde{\tilde{i}}_{1}+1\right)\left(\tilde{\tilde{i}}_{1}+2\right)\lambda_{1}^{4}\leq  \left(j_{1}+1\right)\left(j_{1}+2\right) \left(\tilde{\tilde{i}}_{1}+1\right)\left(\tilde{\tilde{i}}_{1}+2\right)\lambda_{1}^{4}  ,$$
and we have
\begin{equation}\begin{split}
 \left(j_{1}+1\right)\left(j_{1}+2\right) \left(\tilde{\tilde{i}}_{1}+1\right)\left(\tilde{\tilde{i}}_{1}+2\right)\lambda_{1}^{4}\leq&  \left(i_{1}+1\right)\left(i_{1}+2\right)\left(\tilde{\tilde{j}}_{1}+1\right)\left(\tilde{\tilde{j}}_{1}+2\right)\lambda_{1}^{4}  +\left(i_{1}+1\right)\left(i_{1}+2\right) \left(\tilde{\tilde{i}}_{1}+1\right)\left(\tilde{\tilde{i}}_{1}+2\right)\lambda_{1}^{4} \\&\quad\quad\quad\quad\quad\quad\quad\quad\quad\quad\quad\quad\quad\quad\quad\quad + \left(j_{1}+1\right)\left(j_{1}+2\right) \left(\tilde{\tilde{j}}_{1}+1\right)\left(\tilde{\tilde{j}}_{1}+2\right)\lambda_{1}^{4},\end{split}  \label{gege3}
\end{equation}
which shows by summing that (\ref{viena1}) holds. 

The last inequalities (\ref{gege2}) and  (\ref{gege3}) may be easily checked using (\ref{cazuri}) providing not very sharp, but very useful estimates which will be applied later. 
The other estimates hold as well providing (\ref{viena1}). Indeed, when $\left(i_{1},j_{1}\right)=\left(\tilde{\tilde{i}}_{1},\tilde{\tilde{j}}_{1}\right)$, we have that (\ref{gege1}) is equivalent to
$$\left(j_{1}+2\right) \left(i_{1}+2\right)\lambda_{1}^{4}\leq \frac{4\left(i_{1}+1\right)\left(j_{1}+1\right)}{16}, $$
which is obviously true, confirming (\ref{gege1}) is this case, because $4\left(i_{1}+1\right)\left(j_{1}+1\right)\geq \left(j_{1}+2\right) \left(i_{1}+2\right)$, which holds because   $$3\left(i_{1}+1\right)\left(j_{1}+1\right)\geq \left(j_{1}+1\right)\left(i_{1}+1\right)+\left(j_{1}+1\right)+\left(i_{1}+1\right)+1.$$

Respectively, when $\left(i_{1}-1,j_{1}\right)=\left(\tilde{\tilde{i}}_{1},\tilde{\tilde{j}}_{1}\right)$, we have that (\ref{gege1}) is equivalent to
$$\left(j_{1}+2\right) \left(i_{1}+1\right)\lambda_{1}^{4}\leq \frac{4\left(i_{1}+1\right)\left(j_{1}+1\right)}{16}, $$
which is obviously true, confirming (\ref{gege1})  because $4\left(i_{1}+1\right)\left( j_{1}+1\right)\geq \left(j_{1}+2\right) \left(i_{1}+1\right)$, which holds because    $$ 2\left(i_{1}+1\right)\geq \left(i_{1}+1\right) \quad\mbox{and}\quad 2\left( j_{1}+1\right) \geq\left(j_{1}+2\right). $$ 

Respectively, when $\left(i_{1}-2,j_{1}\right)=\left(\tilde{\tilde{i}}_{1},\tilde{\tilde{j}}_{1}\right)$, we have that (\ref{gege1}) is equivalent to
$$\left(j_{1}+2\right)  i_{1} \lambda_{1}^{4}\leq \frac{4\left(i_{1}+1\right)\left(j_{1}+1\right)}{16}, $$
which is obviously true, confirming (\ref{gege1}) because $4\left(i_{1}+1\right)\left(j_{1}+1\right)\geq \left(j_{1}+2\right) i_{1} $, which holds because
$ 4\left(j_{1}+1\right) \geq \left(j_{1}+2\right)$.
 
Respectively, when $ \left(i_{1},j_{1}+1\right)=\left(\tilde{\tilde{i}}_{1},\tilde{\tilde{j}}_{1}\right)$,  we have that (\ref{gege1}) is equivalent to
$$\left(j_{1}+3\right) \left(i_{1}+2\right)\lambda_{1}^{4}\leq \frac{4\left(i_{1}+1\right)\left(j_{1}+2\right)}{16},  
$$
which is obviously true, confirming (\ref{gege1}) because $4\left(i_{1}+1\right)\left( j_{1}+2\right)\geq \left(j_{1}+3\right) \left(i_{1}+2\right)$, which holds because
$$ 2\left( j_{1}+2\right)\geq j_{1}+3\quad\mbox{and}\quad 2\left(i_{1}+1\right)\geq  i_{1}+2 .$$ 
  
Respectively, when $ \left(i_{1}-1,j_{1}+1\right)=\left(\tilde{\tilde{i}}_{1},\tilde{\tilde{j}}_{1}\right),$  we have that (\ref{gege1}) is equivalent to
$$\left(j_{1}+3\right) \left(i_{1}+1\right)\lambda_{1}^{4}\leq \frac{4 i_{1}  \left(j_{1}+2\right)}{16},  
$$
which is obviously true, confirming (\ref{gege1}) because $4 i_{1} \left( j_{1}+2\right)\geq \left(j_{1}+3\right) \left(i_{1}+1\right)$, which holds because
$ 2\ i_{1} \geq  i_{1}+1$ and $2\left( j_{1}+2\right)\geq j_{1}+3$. Clearly, we deal with $i_{1}\in\mathbb{N}^{\star}$, because $\tilde{\tilde{i}}_{1}=i_{1}-1\in\mathbb{N}$, explaining the previous arguments which show why (\ref{gege1}) holds in this case. 

Respectively, when $\left(i_{1}-2,j_{1}+1\right)=\left(\tilde{\tilde{i}}_{1},\tilde{\tilde{j}}_{1}\right)$, we have that (\ref{gege1}) is equivalent to
$$\left(j_{1}+2\right)  i_{1} \lambda_{1}^{4}\leq \frac{4\left(i_{1}-1\right)\left(j_{1}+1\right)}{16},  
$$
which is obviously true, confirming (\ref{gege1}) because $4\left(i_{1}-1\right)\left( j_{1}+2\right)\geq \left(j_{1}+3\right)  i_{1} $, which holds because
$ 2\left( i_{1}-1\right) \geq  i_{1}$ and $2\left( j_{1}+2\right)\geq j_{1}+3$. Clearly, we deal with $i_{1}-1\in\mathbb{N}^{\star}$, because $\tilde{\tilde{i}}_{1}=i_{1}-2\in\mathbb{N}$, explaining the previous arguments which show why (\ref{gege1}) holds in this case.

Respectively, when $\left(i_{1},j_{1}+2\right)=\left(\tilde{\tilde{i}}_{1},\tilde{\tilde{j}}_{1} \right)$, we have that (\ref{gege1}) is equivalent to
$$\left(j_{1}+4\right) \left(i_{1}+2\right)\lambda_{1}^{4}\leq \frac{4\left(i_{1}+1\right)\left(j_{1}+3\right)}{16},  
$$
which is obviously true, confirming (\ref{gege1}) because $4\left(i_{1}+1\right)\left( j_{1}+3\right)\geq \left(j_{1}+4\right) \left(i_{1}+2\right)$, which holds because
$$ 2\left( j_{1}+3\right)\geq j_{1}+4\quad\mbox{and}\quad 2\left(i_{1}+1\right)\geq  i_{1}+2 .$$  

Respectively, when $\left(i_{1}-1,j_{1}+2\right)=\left(\tilde{\tilde{i}}_{1},\tilde{\tilde{j}}_{1} \right)$, we have that (\ref{gege1}) is equivalent to
$$\left(j_{1}+4\right) \left(i_{1}+1\right)\lambda_{1}^{4}\leq \frac{4 i_{1} \left(j_{1}+3\right)}{16},  
$$
which is obviously true, confirming (\ref{gege1}) because $4 i_{1} \left( j_{1}+3\right)\geq \left(j_{1}+4\right) \left(i_{1}+1\right)$, which holds because
$ 2\left( j_{1}+3\right)\geq j_{1}+4$ and  $2 i_{1} \geq  i_{1}+1$. Clearly, we deal with $i_{1}\in\mathbb{N}^{\star}$, because  $\tilde{\tilde{i}}_{1}=i_{1}-1\in\mathbb{N}$. 

Respectively, when $\left(i_{1}-2,j_{1}+2\right)=\left(\tilde{\tilde{i}}_{1},\tilde{\tilde{j}}_{1} \right)$, we have that (\ref{gege1}) is equivalent to
$$\left(j_{1}+4\right)  i_{1} \lambda_{1}^{4}\leq \frac{4\left(i_{1}-1\right)\left(j_{1}+3\right)}{16},  
$$
which is obviously true, confirming (\ref{gege1}) because $4\left(i_{1}-1\right)\left( j_{1}+3\right)\geq \left(j_{1}+4\right)  i_{1} $, which holds because
$ 2\left( j_{1}+3\right)\geq j_{1}+4$ and  $2 \left(i_{1}-1\right) \geq  i_{1}$. Clearly, we deal with $i_{1}-1\in\mathbb{N}^{\star}$, because  $\tilde{\tilde{i}}_{1}=i_{1}-2\in\mathbb{N}$.  We move forward:

We have  
\begin{equation}\displaystyle\max_{1\leq i,j\leq N}\frac{ \tilde{\gamma}_{i,j}^{'}\left[\tilde{I};\tilde{J}\right]}{ \gamma_{i,j}\left[\tilde{I};\tilde{J}\right]}=\frac{\left(\tilde{j}_{1}+1\right)\left(\tilde{j}_{1}+2\right)\lambda_{1}+\left(\tilde{i}_{1}+1\right)\left(\tilde{j}_{1}+1\right)
 \lambda_{1} }{\left(\tilde{j}_{1}+1\right)\left(\tilde{i}_{1}+1\right)+ \left(\tilde{j}_{1}+1\right)\left(\tilde{j}_{1}+2\right)\lambda_{1}^{2} +  \left(\tilde{i}_{1}+1\right)\left(\tilde{i}_{1}+2\right)\lambda_{1}^{2} },\label{ixrr3}
\end{equation}
or the following
\begin{equation} \displaystyle\max_{1\leq i,j\leq N}\frac{ \tilde{\gamma}_{i,j}^{'}\left[\tilde{I};\tilde{J}\right]}{ \gamma_{i,j}\left[\tilde{I};\tilde{J}\right]}= \frac{  4\lambda_{k} \tilde{j}_{k} + \left(\tilde{j}_{k}-1\right)\tilde{j}_{k}\lambda_{k}+\tilde{i}_{k}\tilde{j}_{k}
 \lambda_{k}}{\tilde{j}_{k}\tilde{i}_{k}  +\left(\tilde{j}_{k}-1\right)\tilde{j}_{k} \lambda_{k}^{2} +   \left(\tilde{i}_{k}-1\right)\tilde{i}_{k}  \lambda_{k}^{2}+   \left(\tilde{i}_{k}+\tilde{j}_{k}+1\right)\left(1+4\lambda_{k}^{2}\right)+\left(N+4\lambda_{1}^{2}+\dots+4\lambda_{N}^{2}\right) }. \label{ixrr3se}
\end{equation} 

Also, we  have 
\begin{equation}\displaystyle\max_{1\leq i,j\leq N}\frac{ \tilde{\gamma}_{i,j}^{'}\left[\tilde{\tilde{I}};\tilde{\tilde{J}}\right]}{ \gamma_{i,j}\left[\tilde{\tilde{I}};\tilde{\tilde{J}}\right]}=\frac{\left(\tilde{\tilde{i}}_{1}+1\right)\left(\tilde{\tilde{i}}_{1}+2\right)\lambda_{1}+\left(\tilde{\tilde{i}}_{1}+1\right)\left(\tilde{\tilde{j}}_{1}+1\right)
 \lambda_{1} }{\left(\tilde{\tilde{j}}_{1}+1\right)\left(\tilde{\tilde{i}}_{1}+1\right)+ \left(\tilde{\tilde{j}}_{1}+1\right)\left(\tilde{\tilde{j}}_{1}+2\right)\lambda_{1}^{2} +  \left(\tilde{\tilde{i}}_{1}+1\right)\left(\tilde{\tilde{i}}_{1}+2\right)\lambda_{1}^{2} },\label{ixrr4}
\end{equation}
or the following  
\begin{equation} \displaystyle\max_{1\leq i,j\leq N}\frac{ \tilde{\gamma}_{i,j}^{'}\left[\tilde{\tilde{I}};\tilde{\tilde{J}}\right]}{ \gamma_{i,j}\left[\tilde{\tilde{I}};\tilde{\tilde{J}}\right]} = \frac{  4\lambda_{k} \tilde{\tilde{i}}_{k} + \left(\tilde{\tilde{i}}_{k}-1\right)\tilde{\tilde{i}}_{k}\lambda_{k}+\tilde{\tilde{i}}_{k}\tilde{\tilde{j}}_{k}
 \lambda_{k}}{\tilde{\tilde{j}}_{k}\tilde{\tilde{i}}_{k}  +\left(\tilde{\tilde{j}}_{k}-1\right)\tilde{\tilde{j}}_{k} \lambda_{k}^{2} +   \left(\tilde{\tilde{i}}_{k}-1\right)\tilde{\tilde{i}}_{k}  \lambda_{k}^{2}+ \displaystyle  \left(\tilde{\tilde{i}}_{k}+\tilde{\tilde{j}}_{k}+1\right)\left(1+4\lambda_{k}^{2}\right)+\left(N+4\lambda_{1}^{2}+\dots+4\lambda_{N}^{2}\right) }. \label{ixrr4se}
\end{equation}  

We prove that
\begin{equation}  \displaystyle\max_{1\leq i,j\leq N}\frac{ \tilde{\gamma}_{i,j}^{'}\left[\tilde{I};\tilde{J}\right]}{ \gamma_{i,j}\left[\tilde{I};\tilde{J}\right]}\displaystyle\max_{1\leq i,j\leq N}\frac{ \tilde{\gamma}_{i,j}^{'}\left[\tilde{\tilde{I}};\tilde{\tilde{J}}\right]}{ \gamma_{i,j}\left[\tilde{\tilde{I}};\tilde{\tilde{J}}\right]}   < \frac{1}{2} . \label{viena2} \end{equation} 

Assume that (\ref{ixrr3}) and hold (\ref{ixrr4}).  Then, we have
\begin{equation}\begin{split}& \left(\left(\tilde{j}_{1}+1\right)\left(\tilde{j}_{1}+2\right)\lambda_{1}+\left(\tilde{i}_{1}+1\right)\left(\tilde{j}_{1}+1\right)
 \lambda_{1}\right)\left(\left(\tilde{\tilde{i}}_{1}+1\right)\left(\tilde{\tilde{i}}_{1}+2\right)\lambda_{1}+\left(\tilde{\tilde{i}}_{1}+1\right)\left(\tilde{\tilde{j}}_{1}+1\right)
 \lambda_{1}\right)=\left(\tilde{j}_{1}+1\right)\left(\tilde{j}_{1}+2\right)\left(\tilde{\tilde{i}}_{1}+1\right)\left(\tilde{\tilde{i}}_{1}+2\right)\lambda_{1}^{2}\\&\quad\quad\quad\quad\quad  +\left(\tilde{j}_{1}+1\right)\left(\tilde{j}_{1}+2\right)\left(\tilde{\tilde{i}}_{1}+1\right)\left(\tilde{\tilde{j}}_{1}+1\right)\lambda_{1}^{2}+    \left(\tilde{i}_{1}+1\right)\left(\tilde{j}_{1}+1\right)\left(\tilde{\tilde{i}}_{1}+1\right)\left(\tilde{\tilde{i}}_{1}+2\right)
 \lambda_{1}^{2}+\left(\tilde{i}_{1}+1\right)\left(\tilde{j}_{1}+1\right)\left(\tilde{\tilde{i}}_{1}+1\right)\left(\tilde{\tilde{j}}_{1}+1\right) 
 \lambda_{1}^{2}.\end{split}\label{bbb}
\end{equation}

Also, rewriting (\ref{eee}) in this case, we obtain
\begin{equation}\begin{split}  E\left[\tilde{I},\tilde{J},\tilde{\tilde{I}},\tilde{\tilde{J}}\right]=&\left(\tilde{j}_{1}+1\right)\left(\tilde{i}_{1}+1\right)\left(\tilde{\tilde{j}}_{1}+1\right)\left(\tilde{\tilde{i}}_{1}+1\right)+\left(\tilde{j}_{1}+1\right)\left(\tilde{i}_{1}+1\right)\left(\tilde{\tilde{j}}_{1}+1\right)\left(\tilde{\tilde{j}}_{1}+2\right)\lambda_{1}^{2}   +\left(\tilde{j}_{1}+1\right)\left(\tilde{i}_{1}+1\right)\left(\tilde{\tilde{i}}_{1}+1\right)\left(\tilde{\tilde{i}}_{1}+2\right)\lambda_{1}^{2}\\& +\left(\tilde{j}_{1}+1\right)\left(\tilde{j}_{1}+2\right)\left(\tilde{\tilde{j}}_{1}+1\right)\left(\tilde{\tilde{i}}_{1}+1\right)\lambda_{1}^{2}
 +\left(\tilde{j}_{1}+1\right)\left(\tilde{j}_{1}+2\right) \left(\tilde{\tilde{j}}_{1}+1\right)\left(\tilde{\tilde{j}}_{1}+2\right)\lambda_{1}^{4}\\&+\left(\tilde{j}_{1}+1\right)\left(\tilde{j}_{1}+2\right) \left(\tilde{\tilde{i}}_{1}+1\right)\left(\tilde{\tilde{i}}_{1}+2\right)\lambda_{1}^{4} +\left(\tilde{i}_{1}+1\right)\left(\tilde{i}_{1}+2\right)\left(\tilde{\tilde{j}}_{1}+1\right)\left(\tilde{\tilde{i}}_{1}+1\right)\lambda_{1}^{2}\\&+\left(\tilde{i}_{1}+1\right)\left(\tilde{i}_{1}+2\right)\left(\tilde{\tilde{j}}_{1}+1\right)\left(\tilde{\tilde{j}}_{1}+2\right)\lambda_{1}^{4}  +\left(\tilde{i}_{1}+1\right)\left(\tilde{i}_{1}+2\right)\lambda_{1}^{2}\left(\tilde{\tilde{i}}_{1}+1\right)\left(\tilde{\tilde{i}}_{1}+2\right)\lambda_{1}^{2}.\end{split}\label{eee1}
\end{equation}

Then, according to the assumption (\ref{liliwww1}), we can encounter only the following cases
\begin{equation} \left(\tilde{i}_{1},\tilde{j}_{1}\right)=\left(\tilde{\tilde{i}}_{1},\tilde{\tilde{j}}_{1}\right),\quad \left(\tilde{i}_{1}-1,\tilde{j}_{1}\right)=\left(\tilde{\tilde{i}}_{1},\tilde{\tilde{j}}_{1}\right), \quad \left(\tilde{i}_{1},\tilde{j}_{1}+1\right)=\left(\tilde{\tilde{i}}_{1},\tilde{\tilde{j}}_{1}\right),\quad  \left(\tilde{i}_{1}-1,\tilde{j}_{1}+1\right)=\left(\tilde{\tilde{i}}_{1},\tilde{\tilde{j}}_{1}\right).  \label{cazuri1}
\end{equation}

Then, we  observe by (\ref{lambida}) the following
\begin{equation}2\left(\tilde{j}_{1}+1\right)\left(\tilde{j}_{1}+2\right)\left(\tilde{\tilde{i}}_{1}+1\right)\left(\tilde{\tilde{j}}_{1}+1\right)\lambda_{1}^{2}   \leq \left(\tilde{j}_{1}+1\right)\left(\tilde{i}_{1}+1\right)\left(\tilde{\tilde{j}}_{1}+1\right)\left(\tilde{\tilde{j}}_{1}+2\right)\lambda_{1}^{2}+\left(\tilde{j}_{1}+1\right)\left(\tilde{j}_{1}+2\right)\left(\tilde{\tilde{j}}_{1}+1\right)\left(\tilde{\tilde{i}}_{1}+1\right)\lambda_{1}^{2}.\label{ega1}
\end{equation}

Indeed, (\ref{ega1}) happens because of the following simple inequalities. When $\left(\tilde{i}_{1},\tilde{j}_{1}\right)=\left(\tilde{\tilde{i}}_{1},\tilde{\tilde{j}}_{1}\right)$, we have that (\ref{ega1}) is equivalent to
\begin{equation*}2\left(\tilde{j}_{1}+1\right)\left(\tilde{j}_{1}+2\right)\left(\tilde{i}_{1}+1\right)\left(\tilde{j}_{1}+1\right)\lambda_{1}^{2}   \leq \left(\tilde{j}_{1}+1\right)\left(\tilde{i}_{1}+1\right)\left(\tilde{j}_{1}+1\right)\left(\tilde{j}_{1}+2\right)\lambda_{1}^{2}+\left(\tilde{j}_{1}+1\right)\left(\tilde{j}_{1}+2\right)\left(\tilde{j}_{1}+1\right)\left(\tilde{i}_{1}+1\right)\lambda_{1}^{2}, 
\end{equation*}
which is obviously true, concluding that (\ref{ega1})  holds under this assumption.  

Respectively, when $\left(\tilde{i}_{1}-1,\tilde{j}_{1}\right)=\left(\tilde{\tilde{i}}_{1},\tilde{\tilde{j}}_{1}\right)$, we have that (\ref{ega1}) is equivalent to
\begin{equation*}2\left(\tilde{j}_{1}+1\right)\left(\tilde{j}_{1}+2\right)\tilde{i}_{1}\left(\tilde{j}_{1}+1\right)\lambda_{1}^{2}   \leq \left(\tilde{j}_{1}+1\right)\left(\tilde{i}_{1}+1\right)\left(\tilde{j}_{1}+1\right)\left(\tilde{j}_{1}+2\right)\lambda_{1}^{2}+\left(\tilde{j}_{1}+1\right)\left(\tilde{j}_{1}+2\right)\left(\tilde{j}_{1}+1\right)\tilde{i}_{1}\lambda_{1}^{2}, 
\end{equation*}
which is obviously true, concluding that (\ref{ega1})  holds under this assumption. 

Respectively, when $\left(\tilde{i}_{1},\tilde{j}_{1}+1\right)=\left(\tilde{\tilde{i}}_{1},\tilde{\tilde{j}}_{1}\right)$, we have that (\ref{ega1}) is equivalent to
\begin{equation*}2\left(\tilde{j}_{1}+1\right)\left(\tilde{j}_{1}+2\right)\left(\tilde{i}_{1}+1\right)\left(\tilde{j}_{1}+2\right)\lambda_{1}^{2}   \leq \left(\tilde{j}_{1}+1\right)\left(\tilde{i}_{1}+1\right)\left(\tilde{j}_{1}+2\right)\left(\tilde{j}_{1}+3\right)\lambda_{1}^{2}+\left(\tilde{j}_{1}+1\right)\left(\tilde{j}_{1}+2\right)\left(\tilde{j}_{1}+2\right)\left(\tilde{i}_{1}+1\right)\lambda_{1}^{2}, 
\end{equation*}
which is obviously true, concluding that (\ref{ega1})  holds under this assumption. 

Respectively, when $\left(\tilde{i}_{1}-1,\tilde{j}_{1}+1\right)=\left(\tilde{\tilde{i}}_{1},\tilde{\tilde{j}}_{1}\right)$, we have that (\ref{ega1}) is equivalent to
\begin{equation*}2\left(\tilde{j}_{1}+1\right)\left(\tilde{j}_{1}+2\right)\tilde{i}_{1}\left(\tilde{j}_{1}+2\right)\lambda_{1}^{2}   \leq \left(\tilde{j}_{1}+1\right)\left(\tilde{i}_{1}+1\right)\left(\tilde{j}_{1}+2\right)\left(\tilde{j}_{1}+3\right)\lambda_{1}^{2}+\left(\tilde{j}_{1}+1\right)\left(\tilde{j}_{1}+2\right)\left(\tilde{j}_{1}+2\right)\tilde{i}_{1}\lambda_{1}^{2}, 
\end{equation*}
which is obviously true, concluding that (\ref{ega1})  holds under this assumption, because we have
$$\left(\tilde{j}_{1}+1\right)\left(\tilde{j}_{1}+2\right)\tilde{i}_{1}\left(\tilde{j}_{1}+2\right)\lambda_{1}^{2}   \leq \left(\tilde{j}_{1}+1\right)\left(\tilde{i}_{1}+1\right)\left(\tilde{j}_{1}+2\right)\left(\tilde{j}_{1}+3\right)\lambda_{1}^{2}, $$
in the light of the following obvious inequality
$$ \left(\tilde{j}_{1}+2\right)\tilde{i}_{1}   \leq  \left(\tilde{i}_{1}+1\right) \left(\tilde{j}_{1}+3\right) , $$
and respectively, because we obviously have
$$\left(\tilde{j}_{1}+1\right)\left(\tilde{j}_{1}+2\right)\tilde{i}_{1}\left(\tilde{j}_{1}+2\right)\lambda_{1}^{2}   \leq \left(\tilde{j}_{1}+1\right)\left(\tilde{j}_{1}+2\right)\left(\tilde{j}_{1}+2\right)\tilde{i}_{1}\lambda_{1}^{2}, $$
concluding that (\ref{ega1})  holds under this assumption.

Then, we  observe by (\ref{lambida}) the following
\begin{equation}2\left(\tilde{i}_{1}+1\right)\left(\tilde{j}_{1}+1\right)\left(\tilde{\tilde{i}}_{1}+1\right)\left(\tilde{\tilde{i}}_{1}+2\right)\lambda_{1}^{2}\leq \left(\tilde{j}_{1}+1\right)\left(\tilde{i}_{1}+1\right)\left(\tilde{\tilde{i}}_{1}+1\right)\left(\tilde{\tilde{i}}_{1}+2\right)\lambda_{1}^{2}    +\left(\tilde{i}_{1}+1\right)\left(\tilde{i}_{1}+2\right)\left(\tilde{\tilde{j}}_{1}+1\right)\left(\tilde{\tilde{i}}_{1}+1\right)\lambda_{1}^{2}.\label{ega2}
\end{equation}

Indeed, (\ref{ega2}) happens because of the following simple inequalities. When $\left(\tilde{i}_{1},\tilde{j}_{1}\right)=\left(\tilde{\tilde{i}}_{1},\tilde{\tilde{j}}_{1}\right)$, we have that (\ref{ega2}) is equivalent to
\begin{equation*}2\left(\tilde{i}_{1}+1\right)\left(\tilde{j}_{1}+1\right)\left(\tilde{i}_{1}+1\right)\left(\tilde{i}_{1}+2\right)\lambda_{1}^{2}\leq \left(\tilde{j}_{1}+1\right)\left(\tilde{i}_{1}+1\right)\left(\tilde{i}_{1}+1\right)\left(\tilde{i}_{1}+2\right)\lambda_{1}^{2}    +\left(\tilde{i}_{1}+1\right)\left(\tilde{i}_{1}+2\right)\left(\tilde{j}_{1}+1\right)\left(\tilde{i}_{1}+1\right)\lambda_{1}^{2} , 
\end{equation*}
which is obviously true, concluding that (\ref{ega2})  holds under this assumption. 
  
Respectively, when $\left(\tilde{i}_{1}-1,\tilde{j}_{1}\right)=\left(\tilde{\tilde{i}}_{1},\tilde{\tilde{j}}_{1}\right)$, we have that (\ref{ega2}) is equivalent to
\begin{equation*}2\left(\tilde{i}_{1}+1\right)\left(\tilde{j}_{1}+1\right)\tilde{i}_{1}\left(\tilde{i}_{1}+1\right)\lambda_{1}^{2}\leq \left(\tilde{j}_{1}+1\right)\left(\tilde{i}_{1}+1\right)\tilde{i}_{1}\left(\tilde{i}_{1}+1\right)\lambda_{1}^{2}    +\left(\tilde{i}_{1}+1\right)\left(\tilde{i}_{1}+2\right)\left(\tilde{j}_{1}+1\right)\tilde{i}_{1}\lambda_{1}^{2} , 
\end{equation*}
which is obviously true, concluding that (\ref{ega2})  holds under this assumption.   
  
Respectively, when $\left(\tilde{i}_{1},\tilde{j}_{1}+1\right)=\left(\tilde{\tilde{i}}_{1},\tilde{\tilde{j}}_{1}\right)$, we have that  (\ref{ega2}) is equivalent to
\begin{equation*}2\left(\tilde{i}_{1}+1\right)\left(\tilde{j}_{1}+1\right)\left(\tilde{i}_{1}+1\right)\left(\tilde{i}_{1}+2\right)\lambda_{1}^{2}\leq \left(\tilde{j}_{1}+1\right)\left(\tilde{i}_{1}+1\right)\left(\tilde{i}_{1}+1\right)\left(\tilde{i}_{1}+2\right)\lambda_{1}^{2}    +\left(\tilde{i}_{1}+1\right)\left(\tilde{i}_{1}+2\right)\left(\tilde{j}_{1}+2\right)\left(\tilde{i}_{1}+1\right)\lambda_{1}^{2} , 
\end{equation*}
which is obviously true, concluding that (\ref{ega2})  holds under this assumption.  
  
Respectively, when $\left(\tilde{i}_{1}-1,\tilde{j}_{1}+1\right)=\left(\tilde{\tilde{i}}_{1},\tilde{\tilde{j}}_{1}\right)$, we have that  (\ref{ega2}) is equivalent to
\begin{equation*}2\left(\tilde{i}_{1}+1\right)\left(\tilde{j}_{1}+1\right)\tilde{i}_{1}\left(\tilde{i}_{1}+1\right)\lambda_{1}^{2}\leq \left(\tilde{j}_{1}+1\right)\left(\tilde{i}_{1}+1\right)\tilde{i}_{1}\left(\tilde{i}_{1}+1\right)\lambda_{1}^{2}    +\left(\tilde{i}_{1}+1\right)\left(\tilde{i}_{1}+2\right)\left(\tilde{j}_{1}+2\right)\tilde{i}_{1}\lambda_{1}^{2} , 
\end{equation*}
which is obviously true, because we have
$$\left(\tilde{i}_{1}+1\right)\left(\tilde{j}_{1}+1\right)\tilde{i}_{1}\left(\tilde{i}_{1}+1\right)\lambda_{1}^{2}\leq \left(\tilde{i}_{1}+1\right)\left(\tilde{i}_{1}+2\right)\left(\tilde{j}_{1}+2\right)\tilde{i}_{1}\lambda_{1}^{2},$$
which holds in the light of the following obvious inequality
$$ \left(\tilde{j}_{1}+1\right) \left(\tilde{i}_{1}+1\right) \leq  \left(\tilde{i}_{1}+2\right)\left(\tilde{j}_{1}+2\right)  ,$$
and respectively, because we obviously have
$$\left(\tilde{i}_{1}+1\right)\left(\tilde{j}_{1}+1\right)\tilde{i}_{1}\left(\tilde{i}_{1}+1\right)\lambda_{1}^{2}\leq \left(\tilde{j}_{1}+1\right)\left(\tilde{i}_{1}+1\right)\tilde{i}_{1}\left(\tilde{i}_{1}+1\right)\lambda_{1}^{2}  ,$$
concluding that (\ref{ega2})  holds under this assumption.   
  
Then, we  observe by (\ref{lambida}) the following
\begin{equation}\begin{split}&2\left(\tilde{j}_{1}+1\right)\left(\tilde{j}_{1}+2\right)\left(\tilde{\tilde{i}}_{1}+1\right)\left(\tilde{\tilde{i}}_{1}+2\right)\lambda_{1}^{2}+2\left(\tilde{j}_{1}+1\right)\left(\tilde{\tilde{j}}_{1}+1\right)\left( \tilde{i} _{1}+1\right)\left(\tilde{\tilde{i}}_{1}+1\right) \lambda_{1}^{2}\\& \leq \left(\tilde{j}_{1}+1\right)\left(\tilde{i}_{1}+1\right)\left(\tilde{\tilde{j}}_{1}+1\right)\left(\tilde{\tilde{i}}_{1}+1\right)+\left(\left(\tilde{\tilde{j}}_{1}+1\right)\left(\tilde{\tilde{j}}_{1}+2\right)  +  \left(\tilde{\tilde{i}}_{1}+1\right)\left(\tilde{\tilde{i}}_{1}+2\right) \right)\left(\left(\tilde{j}_{1}+1\right)\left(\tilde{j}_{1}+2\right)  +  \left(\tilde{i}_{1}+1\right)\left(\tilde{i}_{1}+2\right) \right) \lambda_{1}^{4},\end{split}\label{ega3}
\end{equation}
which, in the light of the following inequality
\begin{equation*}  2\left(\tilde{j}_{1}+1\right)\left(\tilde{\tilde{j}}_{1}+1\right)\left( \tilde{i} _{1}+1\right)\left(\tilde{\tilde{i}}_{1}+1\right) \lambda_{1}^{2}\leq \frac{  \left(\tilde{j}_{1}+1\right)\left(\tilde{i}_{1}+1\right)\left(\tilde{\tilde{j}}_{1}+1\right)\left(\tilde{\tilde{i}}_{1}+1\right)}{2},
\end{equation*}
is equivalent to
\begin{equation}\begin{split} \left(\tilde{j}_{1}+1\right)\left(\tilde{j}_{1}+2\right)\left(\tilde{\tilde{i}}_{1}+1\right)\left(\tilde{\tilde{i}}_{1}+2\right)4\lambda_{1}^{2}\leq&   \left(\tilde{j}_{1}+1\right)\left(\tilde{i}_{1}+1\right)\left(\tilde{\tilde{j}}_{1}+1\right)\left(\tilde{\tilde{i}}_{1}+1\right)  \\&   +\left(\left(\tilde{\tilde{j}}_{1}+1\right)\left(\tilde{\tilde{j}}_{1}+2\right)  +  \left(\tilde{\tilde{i}}_{1}+1\right)\left(\tilde{\tilde{i}}_{1}+2\right) \right)\left(\left(\tilde{j}_{1}+1\right)\left(\tilde{j}_{1}+2\right)  +  \left(\tilde{i}_{1}+1\right)\left(\tilde{i}_{1}+2\right) \right) 2\lambda_{1}^{4}.\end{split}\label{ega3se}
\end{equation}  
    
Then (\ref{ega3se}) happens in the light of the following inequalities. When $\left(\tilde{i}_{1},\tilde{j}_{1}\right)=\left(\tilde{\tilde{i}}_{1},\tilde{\tilde{j}}_{1}\right)$, we have
\begin{equation*}\begin{split} \left(\tilde{j}_{1}+1\right)\left(\tilde{j}_{1}+2\right)\left(\tilde{i}_{1}+1\right)\left(\tilde{i}_{1}+2\right)4\lambda_{1}^{2}\leq&   \left(\tilde{j}_{1}+1\right)\left(\tilde{i}_{1}+1\right)\left(\tilde{j}_{1}+1\right)\left(\tilde{i}_{1}+1\right)  \\&   +\left(\left(\tilde{j}_{1}+1\right)\left(\tilde{j}_{1}+2\right)  +  \left(\tilde{i}_{1}+1\right)\left(\tilde{i}_{1}+2\right) \right)\left(\left(\tilde{j}_{1}+1\right)\left(\tilde{j}_{1}+2\right)  +  \left(\tilde{i}_{1}+1\right)\left(\tilde{i}_{1}+2\right) \right) 2\lambda_{1}^{4}.\end{split} 
\end{equation*}

We have
\begin{equation}\left(\tilde{j}_{1}+1\right)\left(\tilde{j}_{1}+2\right)\left(\tilde{i}_{1}+1\right)\left(\tilde{i}_{1}+2\right)\left(2\lambda_{1}^{2}-1\right)\leq  0 \leq \frac{\left(\tilde{j}_{1}+1\right)\left(\tilde{i}_{1}+1\right)\left(\tilde{j}_{1}+1\right)\left(\tilde{i}_{1}+1\right)}{2},   \label{vv}
\end{equation}
and, we also have
\begin{equation}\left(\tilde{j}_{1}+1\right)\left(\tilde{j}_{1}+2\right)\left(\tilde{i}_{1}+1\right)\left(\tilde{i}_{1}+2\right) \leq  \left(\left(\tilde{j}_{1}+1\right)\left(\tilde{j}_{1}+2\right)  +  \left(\tilde{i}_{1}+1\right)\left(\tilde{i}_{1}+2\right) \right)^{2} \lambda_{1}^{4},  \label{vv1}
\end{equation} 
because it is equivalent to
\begin{equation} 0 \leq   \left(\frac{\left(\tilde{j}_{1}+1\right)\left(\tilde{j}_{1}+2\right)}{\left(\tilde{i}_{1}+1\right)\left(\tilde{i}_{1}+2\right)} \right)^{2}  +  \frac{\left(\tilde{j}_{1}+1\right)\left(\tilde{j}_{1}+2\right)}{\left(\tilde{i}_{1}+1\right)\left(\tilde{i}_{1}+2\right)}\frac{2\lambda_{1}^{4}-1}{\lambda_{1}^{4}}  +1,\label{i1}
\end{equation}
which holds, because we have 
\begin{equation}\frac{2\lambda_{1}^{4}-1}{\lambda_{1}^{4}} \leq 1\Longleftrightarrow \lambda_{1}^{4} \leq 1.\label{i2}
\end{equation}

Now, the summing of (\ref{vv}) with (\ref{vv1})  confirm  that (\ref{ega3se}) holds. Thus   (\ref{ega3}) holds.

Respectively, when $\left(\tilde{i}_{1}-1,\tilde{j}_{1}\right)=\left(\tilde{\tilde{i}}_{1},\tilde{\tilde{j}}_{1}\right)$, we have that (\ref{ega3se}) is equivalent to
\begin{equation*}\begin{split} \left(\tilde{j}_{1}+1\right)\left(\tilde{j}_{1}+2\right)\tilde{i}_{1}\left(\tilde{i}_{1}+1\right)4\lambda_{1}^{2}\leq&   \left(\tilde{j}_{1}+1\right)\left(\tilde{i}_{1}+1\right)\left(\tilde{j}_{1}+1\right)\tilde{i}_{1}  \\&   +\left(\left(\tilde{j}_{1}+1\right)\left(\tilde{j}_{1}+2\right)  +  \tilde{i}_{1}\left(\tilde{i}_{1}+1\right) \right)\left(\left(\tilde{j}_{1}+1\right)\left(\tilde{j}_{1}+2\right)  +  \left(\tilde{i}_{1}+1\right)\left(\tilde{i}_{1}+2\right) \right) 2\lambda_{1}^{4}.\end{split}
\end{equation*} 

We have
\begin{equation}\left(\tilde{j}_{1}+1\right)\left(\tilde{j}_{1}+2\right)\tilde{i}_{1}\left(\tilde{i}_{1}+1\right)\left(2\lambda_{1}^{2}-1\right)\leq  0 \leq \frac{\left(\tilde{j}_{1}+1\right)\left(\tilde{i}_{1}+1\right)\left(\tilde{j}_{1}+1\right)\tilde{i}_{1}}{2}   \label{vv11}
\end{equation}
and, we also have
\begin{equation}\left(\tilde{j}_{1}+1\right)\left(\tilde{j}_{1}+2\right)\tilde{i}_{1}\left(\tilde{i}_{1}+1\right) \leq  \left(\left(\tilde{j}_{1}+1\right)\left(\tilde{j}_{1}+2\right)  +  \tilde{i}_{1}\left(\tilde{i}_{1}+1\right) \right)\left(\left(\tilde{j}_{1}+1\right)\left(\tilde{j}_{1}+2\right)  +  \left(\tilde{i}_{1}+1\right)\left(\tilde{i}_{1}+2\right) \right)  \lambda_{1}^{4},  \label{vv12}
\end{equation} 
because we have
\begin{equation} \left(\tilde{j}_{1}+1\right)\left(\tilde{j}_{1}+2\right)\tilde{i}_{1}\left(\tilde{i}_{1}+1\right) \leq  \left(\left(\tilde{j}_{1}+1\right)\left(\tilde{j}_{1}+2\right)  +  \tilde{i}_{1}\left(\tilde{i}_{1}+1\right) \right)\left(\left(\tilde{j}_{1}+1\right)\left(\tilde{j}_{1}+2\right)  +  \left(\tilde{i}_{1}+1\right) \tilde{i}_{1}  \right)  \lambda_{1}^{4} ,\label{3445}
\end{equation}
which holds in the light of (\ref{i1}) and (\ref{i2}).

Now, the summing of (\ref{vv11}) with (\ref{vv12})  confirm  that (\ref{ega3se})  holds. Thus (\ref{ega3}) holds.

Respectively, when $\left(\tilde{i}_{1},\tilde{j}_{1}+1\right)=\left(\tilde{\tilde{i}}_{1},\tilde{\tilde{j}}_{1}\right)$, we have that  (\ref{ega3se}) is equivalent to
\begin{equation*}\begin{split} \left(\tilde{j}_{1}+1\right)\left(\tilde{j}_{1}+2\right)\left(\tilde{i}_{1}+1\right)\left(\tilde{i}_{1}+2\right)4\lambda_{1}^{2}\leq&   \left(\tilde{j}_{1}+1\right)\left(\tilde{i}_{1}+1\right)\left(\tilde{j}_{1}+2\right)\left(\tilde{i}_{1}+1\right)  \\&   +\left(\left(\tilde{j}_{1}+2\right)\left(\tilde{j}_{1}+3\right)  +  \left(\tilde{i}_{1}+1\right)\left(\tilde{i}_{1}+2\right) \right)\left(\left(\tilde{j}_{1}+1\right)\left(\tilde{j}_{1}+2\right)  +  \left(\tilde{i}_{1}+1\right)\left(\tilde{i}_{1}+2\right) \right) 2\lambda_{1}^{4}.\end{split}
\end{equation*} 

We have
\begin{equation}\left(\tilde{j}_{1}+1\right)\left(\tilde{j}_{1}+2\right)\left(\tilde{i}_{1}+1\right)\left(\tilde{i}_{1}+2\right)\left(2\lambda_{1}^{2}-1\right)\leq  0 \leq \frac{\left(\tilde{j}_{1}+1\right)\left(\tilde{i}_{1}+1\right)\left(\tilde{j}_{1}+2\right)\left(\tilde{i}_{1}+1\right)}{2} ,  \label{vv3}
\end{equation}
and, we also have
\begin{equation}\left(\tilde{j}_{1}+1\right)\left(\tilde{j}_{1}+2\right)\left(\tilde{i}_{1}+1\right)\left(\tilde{i}_{1}+2\right) \leq  \left(\left(\tilde{j}_{1}+2\right)\left(\tilde{j}_{1}+3\right)  +  \left(\tilde{i}_{1}+1\right)\left(\tilde{i}_{1}+2\right) \right)\left(\left(\tilde{j}_{1}+1\right)\left(\tilde{j}_{1}+2\right)  +  \left(\tilde{i}_{1}+1\right)\left(\tilde{i}_{1}+2\right) \right) \lambda_{1}^{4},  \label{vv14}
\end{equation} 
because we have
\begin{equation*} \left(\tilde{j}_{1}+1\right)\left(\tilde{j}_{1}+2\right)\left(\tilde{i}_{1}+1\right)\left(\tilde{i}_{1}+2\right) \leq  \left(\left(\tilde{j}_{1}+1\right)\left(\tilde{j}_{1}+2\right)  +  \left(\tilde{i}_{1}+1\right)\left(\tilde{i}_{1}+2\right) \right)\left(\left(\tilde{j}_{1}+1\right)\left(\tilde{j}_{1}+2\right)  +  \left(\tilde{i}_{1}+1\right)\left(\tilde{i}_{1}+2\right) \right) \lambda_{1}^{4} ,
\end{equation*}
which holds in the light of (\ref{i1}) and (\ref{i2}).

Now, the summing of (\ref{vv3}) with (\ref{vv14})  confirm  that (\ref{ega3se})  holds. Thus  (\ref{ega3}) holds.

Respectively, when $\left(\tilde{i}_{1}-1,\tilde{j}_{1}+1\right)=\left(\tilde{\tilde{i}}_{1},\tilde{\tilde{j}}_{1}\right)$, we have that  (\ref{ega3se}) is equivalent to
 \begin{equation*}\begin{split} \left(\tilde{j}_{1}+1\right)\left(\tilde{j}_{1}+2\right)\left(\tilde{i}_{1}+1\right)\tilde{i}_{1}4\lambda_{1}^{2}\leq&   \left(\tilde{j}_{1}+1\right)\left(\tilde{i}_{1}+1\right)\left(\tilde{j}_{1}+2\right) \tilde{i}_{1}   \\&   +\left(\left(\tilde{j}_{1}+2\right)\left(\tilde{j}_{1}+3\right)  +  \left(\tilde{i}_{1}+1\right)\tilde{i}_{1}\right)\left(\left(\tilde{j}_{1}+1\right)\left(\tilde{j}_{1}+2\right)  +  \left(\tilde{i}_{1}+1\right)\left(\tilde{i}_{1}+2\right) \right) 2\lambda_{1}^{4}.\end{split}
\end{equation*} 

We have
\begin{equation}\left(\tilde{j}_{1}+1\right)\left(\tilde{j}_{1}+2\right)\left(\tilde{i}_{1}+1\right)\tilde{i}_{1}\left(2\lambda_{1}^{2}-1\right)\leq  0 \leq \frac{\left(\tilde{j}_{1}+1\right)\left(\tilde{i}_{1}+1\right)\left(\tilde{j}_{1}+2\right) \tilde{i}_{1}}{2},   \label{vv8}
\end{equation}
and, we also have
\begin{equation}\left(\tilde{j}_{1}+1\right)\left(\tilde{j}_{1}+2\right)\left(\tilde{i}_{1}+1\right)\tilde{i}_{1} \leq  \left(\left(\tilde{j}_{1}+2\right)\left(\tilde{j}_{1}+3\right)  +  \left(\tilde{i}_{1}+1\right)\tilde{i}_{1}\right)\left(\left(\tilde{j}_{1}+1\right)\left(\tilde{j}_{1}+2\right)  +  \left(\tilde{i}_{1}+1\right)\left(\tilde{i}_{1}+2\right) \right) \lambda_{1}^{4},  \label{vv18}
\end{equation} 
because (\ref{3445}) holds. 
 
Now, the summing of (\ref{vv8}) with (\ref{vv18})  confirm  that (\ref{ega3se}) holds. Thus  (\ref{ega3}) holds. Moreover,  (\ref{viena2}) becomes clear in the light of (\ref{ega1}), (\ref{ega2}),  (\ref{ega3}) and (\ref{bbb}), (\ref{eee1}).
 
Assume that (\ref{ixrr3se}) and  (\ref{ixrr4se}) hold. Then, we have
\begin{equation}\begin{split}& \frac{  4\lambda_{k} \tilde{j}_{k} + \left(\tilde{j}_{k}-1\right)\tilde{j}_{k}\lambda_{k}+\tilde{i}_{k}\tilde{j}_{k}
 \lambda_{k}}{\tilde{j}_{k}\tilde{i}_{k}  +\left(\tilde{j}_{k}-1\right)\tilde{j}_{k} \lambda_{k}^{2} +   \left(\tilde{i}_{k}-1\right)\tilde{i}_{k}  \lambda_{k}^{2}+   \left(\tilde{i}_{k}+\tilde{j}_{k}+1\right)\left(1+4\lambda_{k}^{2}\right)+\left(N+4\lambda_{1}^{2}+\dots+4\lambda_{N}^{2}\right) }\\&\quad\quad \frac{  4\lambda_{k} \tilde{\tilde{i}}_{k} + \left(\tilde{\tilde{i}}_{k}-1\right)\tilde{\tilde{i}}_{k}\lambda_{k}+\tilde{\tilde{i}}_{k}\tilde{\tilde{j}}_{k}
 \lambda_{k}}{\tilde{\tilde{j}}_{k}\tilde{\tilde{i}}_{k}  +\left(\tilde{\tilde{j}}_{k}-1\right)\tilde{\tilde{j}}_{k} \lambda_{k}^{2} +   \left(\tilde{\tilde{i}}_{k}-1\right)\tilde{\tilde{i}}_{k}  \lambda_{k}^{2}+ \displaystyle  \left(\tilde{\tilde{i}}_{k}+\tilde{\tilde{j}}_{k}+1\right)\left(1+4\lambda_{k}^{2}\right)+\left(N+4\lambda_{1}^{2}+\dots+4\lambda_{N}^{2}\right) }\\&\leq \frac{  4\lambda_{k} \tilde{j}_{k} + \left(\tilde{j}_{k}+1\right)\tilde{j}_{k}\lambda_{k}+\tilde{i}_{k}\tilde{j}_{k}
 \lambda_{k}}{\tilde{j}_{k}\tilde{i}_{k}  +\left(\tilde{j}_{k}+1\right)\tilde{j}_{k} \lambda_{k}^{2} +   \left(\tilde{i}_{k}+1\right)\tilde{i}_{k}  \lambda_{k}^{2}+   \left(\tilde{i}_{k}+\tilde{j}_{k}+1\right)\left(1+2\lambda_{k}^{2}\right)+\left(N+4\lambda_{1}^{2}+\dots+4\lambda_{N}^{2}\right) }\\&\quad\quad \frac{  4\lambda_{k} \tilde{\tilde{i}}_{k} + \left(\tilde{\tilde{i}}_{k}+1\right)\tilde{\tilde{i}}_{k}\lambda_{k}+\tilde{\tilde{i}}_{k}\tilde{\tilde{j}}_{k}
 \lambda_{k}}{\tilde{\tilde{j}}_{k}\tilde{\tilde{i}}_{k}  +\left(\tilde{\tilde{j}}_{k}+1\right)\tilde{\tilde{j}}_{k} \lambda_{k}^{2} +   \left(\tilde{\tilde{i}}_{k}+1\right)\tilde{\tilde{i}}_{k}  \lambda_{k}^{2}+ \displaystyle  \left(\tilde{\tilde{i}}_{k}+\tilde{\tilde{j}}_{k}+1\right)\left(1+2\lambda_{k}^{2}\right)+\left(N+4\lambda_{1}^{2}+\dots+4\lambda_{N}^{2}\right) }\\&\leq \frac{  4\lambda_{k} \tilde{j}_{k} + \left(\tilde{j}_{k}+1\right)\tilde{j}_{k}\lambda_{k}+\tilde{i}_{k}\tilde{j}_{k}
 \lambda_{k}}{\tilde{j}_{k}\tilde{i}_{k}  +\left(\tilde{j}_{k}+1\right)\tilde{j}_{k} \lambda_{k}^{2} +   \left(\tilde{i}_{k}+1\right)\tilde{i}_{k}  \lambda_{k}^{2} } \frac{  \left(\tilde{\tilde{i}}_{k}+1\right)\tilde{\tilde{i}}_{k}\lambda_{k}+\tilde{\tilde{i}}_{k}\tilde{\tilde{j}}_{k}
 \lambda_{k}}{\tilde{\tilde{j}}_{k}\tilde{\tilde{i}}_{k}  +\left(\tilde{\tilde{j}}_{k}+1\right)\tilde{\tilde{j}}_{k} \lambda_{k}^{2} +   \left(\tilde{\tilde{i}}_{k}+1\right)\tilde{\tilde{i}}_{k}  \lambda_{k}^{2}}<\frac{1}{2},\end{split}\label{911}
\end{equation}
according to the previous computations,  finalizing thus (\ref{viena2}).
  
Similarly, we have  \begin{equation}\begin{split}&  \displaystyle\max_{1\leq i,j\leq N}\frac{ \tilde{\gamma}_{i,j}^{'}\left[\tilde{I};\tilde{J}\right]}{ \gamma_{i,j}\left[\tilde{I};\tilde{J}\right]}\displaystyle\max_{1\leq i,j\leq N}\frac{\tilde{\gamma}_{i,j}^{''}\left[I;J\right]}{ \gamma_{i,j}\left[I;J\right]}    < \frac{1}{4},\\& \displaystyle\max_{1\leq i,j\leq N}\frac{ \gamma_{i,j}^{''}\left[\tilde{\tilde{I}};\tilde{\tilde{J}}\right]}{ \gamma_{i,j}\left[\tilde{\tilde{I}};\tilde{\tilde{J}}\right]}\displaystyle\max_{1\leq i,j\leq N}\frac{ \tilde{\gamma}_{i,j}^{'}\left[\tilde{\tilde{I}};\tilde{\tilde{J}}\right]}{ \gamma_{i,j}\left[\tilde{\tilde{I}};\tilde{\tilde{J}}\right]}  < \frac{1}{4}.\end{split}\label{viena3} \end{equation}     
These inequalities are clear in the light of   the  computations  previously worked.  
 
We are ready now to move forward with these heavy computations. We firstly prove that (\ref{asia}) holds for $l=3$ starting to perform induction on $l=3,\dots,p-1$. Recalling (\ref{lil}), (\ref{nene1}) and (\ref{nene2}), we show the invertibility of following matrix
\begin{equation}\mbox{I}_{N^{p}} -\frac{\mathcal{M}_{3,1}}{\mathcal{M}_{3,3}}\mathcal{N}_{1,3}   -\left(\frac{\mathcal{M}_{3,2}}{\mathcal{M}_{3,3}}-\frac{\mathcal{M}_{3,1}}{\mathcal{M}_{3,3}} \mathcal{N}_{1,2} \right)\mathcal{N}_{2,3}.\label{congo}
\end{equation}
 
Recalling also (\ref{lil}), (\ref{nene1}) and (\ref{shobi}), we have  
\begin{equation}0<\epsilon_{1}< \left\|\frac{\mathcal{M}_{1,2}}{\mathcal{M}_{1,1}}\right\|, \hspace{0.1 cm}\left\|\frac{\mathcal{M}_{2,1}}{\mathcal{M}_{2,2}}\right\| < \epsilon_{2}<1,\quad\mbox{where $\epsilon_{1}$, $\epsilon_{2}\in (0,1)$ are chosen such that $\epsilon_{2}^{2}<\left\|\frac{\mathcal{M}_{2,1}}{\mathcal{M}_{2,2}}\right\| < \epsilon_{2}$,}\label{bobo}
\end{equation} 
resulting by (\ref{urs1}), (\ref{viena1}), (\ref{viena2}), (\ref{viena3}) and (\ref{nene2})  the following evaluation
\begin{equation}\begin{split} \left\| \frac{\mathcal{M}_{3,1}}{\mathcal{M}_{3,3}}\mathcal{N}_{1,3}   +\left(\frac{\mathcal{M}_{3,2}}{\mathcal{M}_{3,3}}-\frac{\mathcal{M}_{3,1}}{\mathcal{M}_{3,3}} \mathcal{N}_{1,2} \right)\mathcal{N}_{2,3}  \right\|&= \left\|\frac{\mathcal{M}_{3,1}}{\mathcal{M}_{3,3}} \frac{\mathcal{M}_{1,3}}{\mathcal{M}_{1,1}}+\left(\frac{ \mathcal{M}_{3,2}}{\mathcal{M}_{3,3}} -   \frac{\mathcal{M}_{3,1}}{\mathcal{M}_{3,3}}  \frac{\mathcal{M}_{1,2}}{\mathcal{M}_{1,1}}  \right)\frac{ \frac{\mathcal{M}_{2,3}}{\mathcal{M}_{2,2}}-\frac{\mathcal{M}_{2,1}}{\mathcal{M}_{2,2}} \frac{\mathcal{M}_{1,3}}{\mathcal{M}_{1,1}} }{\mbox{I}_{N^{p}} -\frac{\mathcal{M}_{2,1}}{\mathcal{M}_{2,2}}\frac{\mathcal{M}_{1,2}}{\mathcal{M}_{1,1}}}\right\|  \\&  \leq  \frac{1}{3}+\frac{1}{2} \left(1- \epsilon_{1}-\epsilon_{2}^{2}+\epsilon_{2}^{2}\right)  < 1, \end{split}\label{yor}
\end{equation}
according to (\ref{bobo}), concluding similarly as in (\ref{tigru}) and (\ref{urs1}), the existence of the following matrix
\begin{equation} \frac{1}{\mbox{I}_{N^{p}} -\frac{\mathcal{M}_{3,1}}{\mathcal{M}_{3,3}}\mathcal{N}_{1,3}   -\left(\frac{\mathcal{M}_{3,2}}{\mathcal{M}_{3,3}}-\frac{\mathcal{M}_{3,1}}{\mathcal{M}_{3,3}} \mathcal{N}_{1,2} \right)\mathcal{N}_{2,3}}.\label{urs11}
\end{equation}

Let's  prove now that (\ref{asia}) holds for $l=4$ continuing to perform induction on $l=3,4,5,\dots,p-1$. Recalling (\ref{lil}), (\ref{nene1}) and (\ref{nene2}), we apply the above procedures   on the following matrix 
\begin{equation} \begin{split}&\quad\quad\quad \mbox{I}_{N^{p}} -\frac{\mathcal{M}_{4,2}}{\mathcal{M}_{4,4}}\mathcal{N}_{2,4}   -\left(\frac{\mathcal{M}_{4,3}}{\mathcal{M}_{4,4}}-\frac{\mathcal{M}_{4,2}}{\mathcal{M}_{4,4}}\mathcal{N}_{2,3}\right)\mathcal{N}_{3,4}=\mbox{I}_{N^{p}} -\frac{\mathcal{M}_{4,2}}{\mathcal{M}_{4,4}}\frac{\frac{\mathcal{M}_{2,4}}{\mathcal{M}_{2,2}}}{\mbox{I}_{N^{p}} -\frac{\mathcal{M}_{3,1}}{\mathcal{M}_{3,3}}\mathcal{N}_{1,3}   -\left(\frac{\mathcal{M}_{3,2}}{\mathcal{M}_{3,3}}-\frac{\mathcal{M}_{3,1}}{\mathcal{M}_{3,3}} \mathcal{N}_{1,2} \right)\mathcal{N}_{2,3}}\\&   -\left(\frac{\mathcal{M}_{4,3}}{\mathcal{M}_{4,4}}-\frac{\mathcal{M}_{4,2}}{\mathcal{M}_{4,4}}\frac{ \frac{\mathcal{M}_{2,3}}{\mathcal{M}_{2,2}}-\frac{\mathcal{M}_{2,1}}{\mathcal{M}_{2,2}} \frac{\mathcal{M}_{1,3}}{\mathcal{M}_{1,1}} }{\mbox{I}_{N^{p}} -\frac{\mathcal{M}_{2,1}}{\mathcal{M}_{2,2}}\frac{\mathcal{M}_{1,2}}{\mathcal{M}_{1,1}}}\right) \frac{ \frac{\mathcal{M}_{3,4}}{\mathcal{M}_{3,3}}-\left(\frac{\mathcal{M}_{3,2}}{\mathcal{M}_{3,3}}-\frac{\mathcal{M}_{3,1}}{\mathcal{M}_{3,3}}\mathcal{N}_{1,2}\right)\frac{\mathcal{M}_{4,2}}{\mathcal{M}_{4,4}}\frac{\frac{\mathcal{M}_{2,4}}{\mathcal{M}_{2,2}}}{\mbox{I}_{N^{p}} -\frac{\mathcal{M}_{3,1}}{\mathcal{M}_{3,3}}\mathcal{N}_{1,3}   -\left(\frac{\mathcal{M}_{3,2}}{\mathcal{M}_{3,3}}-\frac{\mathcal{M}_{3,1}}{\mathcal{M}_{3,3}} \mathcal{N}_{1,2} \right)\mathcal{N}_{2,3}}}{\mbox{I}_{N^{p}} -\frac{\mathcal{M}_{3,1}}{\mathcal{M}_{3,3}}\mathcal{N}_{1,3}   -\left(\frac{\mathcal{M}_{3,2}}{\mathcal{M}_{3,3}}-\frac{\mathcal{M}_{3,1}}{\mathcal{M}_{3,3}} \mathcal{N}_{1,2} \right)\mathcal{N}_{2,3}}.\end{split}\label{york}  
\end{equation}

Similarly to the arguing on which  (\ref{viena1}), (\ref{viena2}) and (\ref{viena3})    were concluded, we obtain
\begin{equation}\begin{split}&  \left\| \frac{\mathcal{M}_{4,2}}{\mathcal{M}_{4,4}}\frac{\mathcal{M}_{2,4}}{\mathcal{M}_{2,2}}\dots\right\|<\frac{1}{6}\left(1+\frac{5}{18}+\dots\right)=\frac{3}{13}<\frac{1}{4},\\& \left\| \frac{\mathcal{M}_{4,3}}{\mathcal{M}_{4,4}}\frac{\mathcal{M}_{3,4}}{\mathcal{M}_{3,3}}\dots\right\| <\frac{1}{2}\left(1+\frac{1}{3}+\frac{1}{3^{2}} +\dots\right)=\frac{3}{4} ,\label{ioc}\end{split}
\end{equation}
where  ,,$\dots$'' are defined by terms coming from following formal expansions:
\begin{equation}\begin{split}&\quad\quad\quad\quad\quad\quad\quad\quad\quad\quad\quad\quad\hspace{0.1 cm}      \frac{1}{ \mbox{I}_{N^{p}}-\frac{\mathcal{M}_{2,1}}{\mathcal{M}_{2,2}}\frac{\mathcal{M}_{1,2}}{\mathcal{M}_{1,1}} }=\mbox{I}_{N^{p}}+\frac{\mathcal{M}_{2,1}}{\mathcal{M}_{2,2}}\frac{\mathcal{M}_{1,2}}{\mathcal{M}_{1,1}}+\frac{\mathcal{M}_{2,1}}{\mathcal{M}_{2,2}}\frac{\mathcal{M}_{1,2}}{\mathcal{M}_{1,1}}\frac{\mathcal{M}_{2,1}}{\mathcal{M}_{2,2}}\frac{\mathcal{M}_{1,2}}{\mathcal{M}_{1,1}}+\dots,\\& \frac{1}{\mbox{I}_{N^{p}} -\frac{\mathcal{M}_{3,1}}{\mathcal{M}_{3,3}}\mathcal{N}_{1,3}   -\left(\frac{\mathcal{M}_{3,2}}{\mathcal{M}_{3,3}}-\frac{\mathcal{M}_{3,1}}{\mathcal{M}_{3,3}} \mathcal{N}_{1,2} \right)\mathcal{N}_{2,3}}=\mbox{I}_{N^{p}}+\frac{\mathcal{M}_{3,1}}{\mathcal{M}_{3,3}}\mathcal{N}_{1,3}   -\left(\frac{\mathcal{M}_{3,2}}{\mathcal{M}_{3,3}}+\frac{\mathcal{M}_{3,1}}{\mathcal{M}_{3,3}} \mathcal{N}_{1,2} \right)\mathcal{N}_{2,3}\\&\quad\quad\quad\quad\quad\quad\quad\quad\quad\quad\quad\quad\quad\quad\quad\quad\quad\quad\quad\quad\quad\hspace{0.1 cm}      +\left(\frac{\mathcal{M}_{3,1}}{\mathcal{M}_{3,3}}\mathcal{N}_{1,3}   -\left(\frac{\mathcal{M}_{3,2}}{\mathcal{M}_{3,3}}+\frac{\mathcal{M}_{3,1}}{\mathcal{M}_{3,3}} \mathcal{N}_{1,2} \right)\mathcal{N}_{2,3}\right)^{2}+\dots,\end{split}\label{expan}
\end{equation} 
in the light of the following evaluations
\begin{equation*}\begin{split}&\left\| \frac{\mathcal{M}_{4,3}}{\mathcal{M}_{4,4}}\frac{\mathcal{M}_{3,4}}{\mathcal{M}_{3,3}} \right\|<\frac{1}{2}\quad \left\| \frac{\mathcal{M}_{4,3}}{\mathcal{M}_{4,4}}\frac{\mathcal{M}_{3,4}}{\mathcal{M}_{3,3}} \frac{\mathcal{M}_{3,1}}{\mathcal{M}_{3,3}}\mathcal{N}_{1,3} \right\|<\frac{1}{24}, \quad  \left\| \frac{\mathcal{M}_{4,3}}{\mathcal{M}_{4,4}}\frac{\mathcal{M}_{3,4}}{\mathcal{M}_{3,3}}      \left(\frac{\mathcal{M}_{3,2}}{\mathcal{M}_{3,3}}-\frac{\mathcal{M}_{3,1}}{\mathcal{M}_{3,3}} \mathcal{N}_{1,2} \right)\mathcal{N}_{2,3}  \right\|<\frac{1}{8},\quad \dots, \\& \left\| \frac{\mathcal{M}_{4,2}}{\mathcal{M}_{4,4}}\frac{\mathcal{M}_{2,4}}{\mathcal{M}_{2,2}} \right\|<\frac{1}{6},\quad \left\| \frac{\mathcal{M}_{4,2}}{\mathcal{M}_{4,4}}\frac{\mathcal{M}_{2,4}}{\mathcal{M}_{2,2}} \frac{\mathcal{M}_{3,1}}{\mathcal{M}_{3,3}}\mathcal{N}_{1,3}\right\|<\frac{1}{72},\quad \left\| \frac{\mathcal{M}_{4,2}}{\mathcal{M}_{4,4}}\frac{\mathcal{M}_{2,4}}{\mathcal{M}_{2,2}}  \left(\frac{\mathcal{M}_{3,2}}{\mathcal{M}_{3,3}}-\frac{\mathcal{M}_{3,1}}{\mathcal{M}_{3,3}} \mathcal{N}_{1,2} \right)\mathcal{N}_{2,3} \right\|<\frac{1}{24},\quad \dots,   \end{split}
\end{equation*}
in the light of the estimations related to (\ref{viena1}), (\ref{viena2}), (\ref{viena3}), observing that the resulted denominador may be multiplied by $2$ according to similar computations, improving thus the initial estimations.

Respectively, we  analogously obtain
\begin{equation}\begin{split}& \quad\quad\quad\quad\hspace{0.1 cm}\quad\quad  \quad  \left\| \frac{\mathcal{M}_{4,2}}{\mathcal{M}_{4,4}}\left( \frac{\mathcal{M}_{2,3}}{\mathcal{M}_{2,2}}-\frac{\mathcal{M}_{2,1}}{\mathcal{M}_{2,2}} \frac{\mathcal{M}_{1,3}}{\mathcal{M}_{1,1}} \right) \frac{\mathcal{M}_{3,4}}{\mathcal{M}_{3,3}}\dots\right\|<\frac{3}{4}\left(1-\epsilon_{1}\right) ,\\& \left\|\frac{\mathcal{M}_{4,2}}{\mathcal{M}_{4,4}}\left( \frac{\mathcal{M}_{2,3}}{\mathcal{M}_{2,2}}-\frac{\mathcal{M}_{2,1}}{\mathcal{M}_{2,2}} \frac{\mathcal{M}_{1,3}}{\mathcal{M}_{1,1}} \right) \left(\frac{\mathcal{M}_{3,2}}{\mathcal{M}_{3,3}}-\frac{\mathcal{M}_{3,1}}{\mathcal{M}_{3,3}}\mathcal{N}_{1,2}\right) \dots\right\|<\frac{3}{4}\left(1-\epsilon_{1}\right)^{2} ,\\& \quad\quad\quad \quad \quad \quad   \quad \quad \quad \quad \quad \left\| \frac{\mathcal{M}_{4,3}}{\mathcal{M}_{4,4}}\left(\frac{\mathcal{M}_{3,2}}{\mathcal{M}_{3,3}}-\frac{\mathcal{M}_{3,1}}{\mathcal{M}_{3,3}}\mathcal{N}_{1,2}\right)\dots \right\| <\frac{3}{4}\left(1-\epsilon_{1}\right) ,\end{split}\label{iocc}
\end{equation}
recalling (\ref{bobo}), concluding   by (\ref{urs1}), (\ref{urs11}) ,(\ref{york}) and (\ref{ioc})   the following evaluation
  \begin{equation}\left\|  \frac{\mathcal{M}_{4,2}}{\mathcal{M}_{4,4}}\mathcal{N}_{2,4}+\left(\frac{\mathcal{M}_{4,3}}{\mathcal{M}_{4,4}}-\frac{\mathcal{M}_{4,2}}{\mathcal{M}_{4,4}}\mathcal{N}_{2,3}\right)\mathcal{N}_{3,4}  \right\|\leq \frac{1}{2}+\frac{1}{3}\left(1-2\left(1-\epsilon_{1}\right)+\left(1-\epsilon_{1}\right)^{2}\right)=\frac{3}{13}+\frac{3}{4}\left(1-\epsilon_{1}-1\right)^{4}<\frac{51}{52}<1.\label{yor2}
\end{equation}

Exactly as previously, we obtain the existence of the following matrix
\begin{equation} \frac{1}{\mbox{I}_{N^{p}} -\frac{\mathcal{M}_{4,2}}{\mathcal{M}_{4,4}}\mathcal{N}_{2,4}-\left(\frac{\mathcal{M}_{4,3}}{\mathcal{M}_{4,4}}-\frac{\mathcal{M}_{4,2}}{\mathcal{M}_{4,4}}\mathcal{N}_{2,3}\right)\mathcal{N}_{3,4}} .\label{urs11se}
\end{equation} 

 Now, we are prepared to continue the induction process.  Let $l\in 3,\dots,p-2$ and assume existence for each of the following matrices 
\begin{equation}  \frac{1 }{I_{N^{p}}-\frac{\mathcal{M}_{l,l-2}}{\mathcal{M}_{l,l}}\mathcal{N}_{l-2,l}-\left(\frac{\mathcal{M}_{l,l-1}}{\mathcal{M}_{l,l}}-\frac{\mathcal{M}_{l,l-2}}{\mathcal{M}_{l,l}}\mathcal{N}_{l-2,l-1}\right)\mathcal{N}_{l-1,l} } ,\hspace{0.1 cm} \dots   ,\hspace{0.1 cm}  \frac{1}{I_{N^{p}}-\frac{\mathcal{M}_{2,1}}{\mathcal{M}_{2,2}}\mathcal{N}_{1,2}} . \label{urs11sese}
\end{equation}

Moreover, we assume that we can write as follows
\begin{equation}\begin{split}&\frac{1 }{I_{N^{p}}-\frac{\mathcal{M}_{l,l-2}}{\mathcal{M}_{l,l}}\mathcal{N}_{l-2,l}-\left(\frac{\mathcal{M}_{l,l-1}}{\mathcal{M}_{l,l}}-\frac{\mathcal{M}_{l,l-2}}{\mathcal{M}_{l,l}}\mathcal{N}_{l-2,l-1}\right)\mathcal{N}_{l-1,l} }=I_{N^{p}}+\frac{\mathcal{M}_{l,l-2}}{\mathcal{M}_{l,l}}\mathcal{N}_{l-2,l}\\&\quad\quad\quad\quad\quad\quad\quad\quad\quad\quad\quad\quad\quad\quad\quad\quad\quad\quad\quad\quad\quad\quad\quad\quad\quad\quad\quad        +\left(\frac{\mathcal{M}_{l,l-1}}{\mathcal{M}_{l,l}}-\frac{\mathcal{M}_{l,l-2}}{\mathcal{M}_{l,l}}\mathcal{N}_{l-2,l-1}\right)\mathcal{N}_{l-1,l}+\dots ,\\&\quad\quad\quad\quad\quad\quad\quad\quad \quad\quad\quad\quad\quad\quad\quad\quad\quad\quad\quad\quad\quad\quad\quad\quad\quad\quad\quad  \vdots \\&\frac{1}{I_{N^{p}}-\frac{\mathcal{M}_{2,1}}{\mathcal{M}_{2,2}}\mathcal{N}_{1,2}}=I_{N^{p}}+\frac{\mathcal{M}_{2,1}}{\mathcal{M}_{2,2}}\mathcal{N}_{1,2}+\left(\frac{\mathcal{M}_{2,1}}{\mathcal{M}_{2,2}}\mathcal{N}_{1,2}\right)^{2}+\dots.\end{split}\label{vulpe}
\end{equation}

Similarly as previously, there are recalled (\ref{gringooo}), (\ref{nene5}) and (\ref{gringooo1})  in order to
study the following matrix
 \begin{equation}\begin{split} &\quad\quad\quad\quad\hspace{0.2 cm} \quad\quad\quad\quad\quad\quad\quad \quad I_{N^{p}}-\frac{\mathcal{M}_{l+1,l-1}}{\mathcal{M}_{l+1,l+1}}\mathcal{N}_{l-1,l+1}-\left(\frac{\mathcal{M}_{l+1,l}}{\mathcal{M}_{l+1,l+1}}-\frac{\mathcal{M}_{l+1,l-1}}{\mathcal{M}_{l+1,l+1}}\mathcal{N}_{l-1,l}\right)\mathcal{N}_{l,l+1}=\\& \hspace{0.1 cm}\quad\quad\quad\quad\quad I_{N^{p}}-\frac{\mathcal{M}_{l+1,l-1}}{\mathcal{M}_{l+1,l+1}}\frac{\frac{\mathcal{M}_{l-1,l+1}}{\mathcal{M}_{l-1,l-1}} }{\mbox{I}_{N^{p}} -\frac{\mathcal{M}_{l-1,l-3}}{\mathcal{M}_{l-1,l-1}}\mathcal{N}_{l-3,l-1}-\left(\frac{\mathcal{M}_{l-1,l-2}}{\mathcal{M}_{l-1,l-1}}-\frac{\mathcal{M}_{l-1,l-3}}{\mathcal{M}_{l-1,l-1}}\mathcal{N}_{l-3,l-2}\right)\mathcal{N}_{l-2,l-1} }\\&-\left(\frac{\mathcal{M}_{l+1,l}}{\mathcal{M}_{l+1,l+1}} -\frac{\mathcal{M}_{l+1,l-1}}{\mathcal{M}_{l+1,l+1}}\frac{\frac{\mathcal{M}_{l-1,l }}{\mathcal{M}_{l-1,l-1}}-\left(\frac{\mathcal{M}_{l-1,l-2}}{\mathcal{M}_{l-1,l-1}}-\frac{\mathcal{M}_{l-1,l-3}}{\mathcal{M}_{l-1,l-1}}\mathcal{N}_{l-3,l-2}\right)\mathcal{N}_{l-2,l } }{\mbox{I}_{N^{p}} -\frac{\mathcal{M}_{l-1,l-1}}{\mathcal{M}_{l-1,l-1}}\mathcal{N}_{l-2,l-1}-\left(\frac{\mathcal{M}_{l-1,l-2}}{\mathcal{M}_{l-1,l-1}}-\frac{\mathcal{M}_{l-1,l-3}}{\mathcal{M}_{l-1,l-1}}\mathcal{N}_{l-3,l-3}\right)\mathcal{N}_{l-2,l-1} }\right)\\&\quad\quad\quad\quad\quad\quad\quad\quad\quad\quad\quad\quad\quad\quad\quad\quad\hspace{0.1 cm}\quad\quad\frac{\frac{\mathcal{M}_{l,l+1}}{\mathcal{M}_{l,l}}-\left(\frac{\mathcal{M}_{l,l-1}}{\mathcal{M}_{l,l}}-\frac{\mathcal{M}_{l,l-2}}{\mathcal{M}_{l,l}}\mathcal{N}_{l-2,l-1}\right)\mathcal{N}_{l-1,l+1} }{\mbox{I}_{N^{p}} -\frac{\mathcal{M}_{l,l-2}}{\mathcal{M}_{l,l}}\mathcal{N}_{l-2,l}-\left(\frac{\mathcal{M}_{l,l-1}}{\mathcal{M}_{l,l}}-\frac{\mathcal{M}_{l,l-2}}{\mathcal{M}_{l,l}}\mathcal{N}_{l-2,l-1}\right)\mathcal{N}_{l-1,l} },\end{split}\label{9909}
\end{equation} 
 
Similarly to the arguing on which (\ref{ioc}) and (\ref{iocc})  are concluded, we obtain
 similar estimations also in this case defined by (\ref{9909}) having in mind (\ref{vulpe}). It remains to make suitable estimations for the remaining  matrices in (\ref{9909}) using the recurrences (\ref{nene5}), (\ref{nene4}), (\ref{nene3}), (\ref{nene2}) and (\ref{nene1}). There are occurring similar estimates as in (\ref{ioc}), but for more general products of matrices.  We make expansions using these recurrences and we arrive at similar estimations  to (\ref{yor2}). Similarly as in   (\ref{yor}) and (\ref{yor2}), we obtain
$$\left\|\frac{\mathcal{M}_{l+1,l-1}}{\mathcal{M}_{l+1,l+1}}\mathcal{N}_{l-1,l+1}+\left(\frac{\mathcal{M}_{l+1,l}}{\mathcal{M}_{l+1,l+1}}-\frac{\mathcal{M}_{l+1,l-1}}{\mathcal{M}_{l+1,l+1}}\mathcal{N}_{l-1,l}\right)\mathcal{N}_{l,l+1}\right\|  <\frac{51}{52}<1,
$$
  concluding similarly as in (\ref{tigru}) and (\ref{urs1})  the invertibility of the following matrix 
\begin{equation} \frac{1 }{I_{N^{p}}-\frac{\mathcal{M}_{l+1,l-1}}{\mathcal{M}_{l+1,l+1}}\mathcal{N}_{l-1,l+1}-\left(\frac{\mathcal{M}_{l+1,l}}{\mathcal{M}_{l+1,l+1}}-\frac{\mathcal{M}_{l+1,l-1}}{\mathcal{M}_{l+1,l+1}}\mathcal{N}_{l-1,l}\right)\mathcal{N}_{l,l+1}}.\label{urs11sesese}
\end{equation}

It is clear how the condition (\ref{lambida}), used within  the previous computations, appears  in doing relevant computations.   The estimations (\ref{viena1}), (\ref{viena2}) and (\ref{viena3}) are crucial in order to show  the existence of the matrices from (\ref{urs11}), (\ref{urs11se}) and (\ref{urs11sese}). The condition (\ref{lambida})  is the  analogue to the natural condition imposed by Gong\cite{go1},\cite{go2}, in the light of which we can  conjecture that our result may not hold when (\ref{lambida}) does not hold. Thus, it is difficult to hope to find  a much better and simpler assumption than (\ref{lambida}), because the computations are infernal.

The existences of these matrices (\ref{urs11}), (\ref{urs11se}) and (\ref{urs11sese}), which are concluded using evaluations (\ref{ioc}) and (\ref{iocc}), are fundamental   in order to reach a unique solution for the system of equations from (\ref{beb1}). These evaluations (\ref{ioc}) and (\ref{iocc})  are not very sharp, but are enough for the establishing of convenient bounds. Thus, this condition (\ref{lambida}) may be relaxed, but it is preferred to be kept in that form (\ref{lambida}), because otherwise the computations may get out of any form of control.  We move forward:
\subsection{Fischer Normalization $G$-Spaces\cite{bu4},\cite{bu5}}The Fischer Decomposition from (\ref{op}) gives
\begin{equation}z^{I}=A_{I}(z,\overline{z}) Q(z,\overline{z})+C_{I}(z,\overline{z}),\quad\mbox{where}\hspace{0.1 cm} \tilde{\tr}\left(C_{I}(z,\overline{z})\right)=0,\label{opt}
\end{equation}
where $I\not\in \mathcal{S}$ having   length $p$.
 
This set of homogeneous polynomials, derived from (\ref{opt}),   gives  normalizations based on  certain   Normalization Spaces, which is defined using the generalized version of  the Fischer Decomposition\cite{sh}: for a given  real homogeneous polynomial of degree $p\geq 1$ in $(z,\overline{z})$ denoted by $P(z,\overline{z})$, we  have
\begin{equation} \begin{split}& P(z,\overline{z})=P_{1}(z,\overline{z})Q(z,\overline{z})+R_{1}(z,\overline{z}),\quad\mbox{where $\tr\left(R_{1}(z,\overline{z})\right)=0$ and:}\\& \quad\quad\quad\quad\hspace{0.1 cm}  R_{1}(z,\overline{z})=\displaystyle\sum_{I\in\mathbb{N}^{N}\atop {\left|I\right|=p}} \left(a_{I}C_{I}(z,\overline{z})+b_{I}\overline{C_{I}(z,\overline{z})}\right)+R_{1,0}(z,\overline{z}),\quad\mbox{such that:} \\&\quad\quad\quad\quad\quad\quad \quad\quad\quad R_{1,0}(z,\overline{z})\in \left(\displaystyle\bigcap_{I\in\mathbb{N}^{N},\hspace{0.1 cm} I\not\in \mathcal{S} \atop {\left|I\right|=p}}\left( \ker  C^{\star}_{I}  \cap  \ker  \overline{C}^{\star}_{I} \right)\right) \bigcap\left( \displaystyle\bigcap_{I\in\mathbb{N}^{N},\hspace{0.1 cm} I \in \mathcal{S} \atop {\left|I\right|=p}}\left(\ker  \left(z^{I}\right)^{\star}   \cap  \ker  \left(\overline{z}^{I}\right)^{\star}\right)\right) ,
\\& P_{1}(z,\overline{z})=P_{2}(z,\overline{z})Q(z,\overline{z})+R_{2}(z,\overline{z}),\quad\mbox{where $\tr\left(R_{2}(z,\overline{z})\right)=0$ and:}\\& \quad\quad\quad\quad\hspace{0.1 cm}  R_{2}(z,\overline{z})=\displaystyle\sum_{I\in\mathbb{N}^{N}\atop {\left|I\right|=p-2}} \left(a_{I}C_{I}(z,\overline{z})+b_{I}\overline{C_{I}(z,\overline{z})}\right)+R_{2,0}(z,\overline{z}),\quad\mbox{such that:}\\&\quad\quad\quad\quad \quad\quad \quad\quad\quad R_{2,0}(z,\overline{z})\in \left(\displaystyle\bigcap_{I\in\mathbb{N}^{N},\hspace{0.1 cm} I\not\in \mathcal{S} \atop {\left|I\right|=p-2}} \ker  C^{\star}_{I}  \cap  \ker  \overline{C}^{\star}_{I} \right) \bigcap \left(\displaystyle\bigcap_{I\in\mathbb{N}^{N},\hspace{0.1 cm} I \in \mathcal{S} \atop {\left|I\right|=p-2}}\ker  \left(z^{I}\right)^{\star}   \cap  \ker  \left(\overline{z}^{I}\right)^{\star}\right),   \\&\quad\quad\quad\quad\vdots\quad\quad\quad\quad\quad\quad\quad\quad\vdots\quad\quad\quad\quad\quad\quad\quad\quad\vdots\quad\quad\quad\quad\quad\quad\quad\quad\vdots \quad\quad\quad\quad\quad\quad\quad\quad\quad\quad\quad\quad\quad\quad\vdots
\\& P_{k}(z,\overline{z})=P_{k+1}(z,\overline{z})Q(z,\overline{z})+R_{k+1}(z,\overline{z}),\quad\mbox{where $\tr\left(R_{k+1}(z,\overline{z})\right)=0$ and:}\\& \quad\quad\quad \quad \hspace{0.1 cm} R_{k+1}(z,\overline{z})=\displaystyle\sum_{I\in\mathbb{N}^{N}\atop {\left|I\right|=p-2k}} \left(a_{I}C_{I}(z,\overline{z})+b_{I}\overline{C_{I}(z,\overline{z})}\right)+R_{k+1,0}(z,\overline{z}),\quad\mbox{such that:}\\&\quad\quad\quad\quad \quad\quad \quad\quad\quad R_{k+1,0}(z,\overline{z})\in \left(\displaystyle\bigcap_{I\in\mathbb{N}^{N}, \hspace{0.1 cm}I\not\in \mathcal{S} \atop {\left|I\right|=p-2k}} \left(\ker  C^{\star}_{I}  \cap  \ker  \overline{C}^{\star}_{I}\right)\right)\bigcap\left(   \displaystyle\bigcap_{I\in\mathbb{N}^{N},\hspace{0.1 cm} I \in \mathcal{S} \atop {\left|I\right|=p-2k}}\left(\ker  \left(z^{I}\right)^{\star}   \cap  \ker \left(\overline{z}^{I}\right)^{\star}\right) \right),\\&\quad\quad\quad\quad\vdots\quad\quad\quad\quad\quad\quad\quad\quad\vdots\quad\quad\quad\quad\quad\quad\quad\quad\vdots\quad\quad\quad\quad\quad\quad\quad\quad\vdots \quad\quad\quad\quad\quad\quad\quad\quad\quad\quad\quad\quad\quad\quad\vdots\end{split} \label{new1}
\end{equation}
where   these occurring polynomials 
\begin{equation}\left\{P_{k}(z,\overline{z})\right\}_{k=1,\dots,\left[\frac{p}{2}\right]},\quad \left\{R_{k}(z,\overline{z})\right\}_{k=1,\dots,\left[\frac{p}{2}\right]},
\label{poll1}
\end{equation}
are iteratively obtained  using the generalized version  of the Fischer Decomposition\cite{sh}.  

Recalling   strategies from \cite{bu4} and \cite{bu5}, we define 
\begin{equation}\mathcal{G}_{p},\quad p\in\mathbb{N}^{\star},\label{spartiuG}
\end{equation} 
which consist in real-valued polynomials $P(z,\overline{z})$  of degree $p\geq 1$ in $(z,\overline{z})$ satisfying the normalizations:
$$ P_{k}^{\left(p\right)}(z,\overline{z})=P_{k+1}^{\left(p\right)}(z,\overline{z})Q(z,\overline{z})+R_{k+1}^{\left(p\right)}(z,\overline{z}),\hspace{0.1 cm}\mbox{for all  $k=0,\dots, \left[\frac{p}{2}\right]$ and given $P_{0}^{\left(p\right)}(z,\overline{z})=P(z,\overline{z})$,}  $$ 
such that 
$$ R_{k+1,0}^{\left(p\right)}(z,\overline{z})\in \left(\displaystyle\bigcap_{I\in\mathbb{N}^{N},\hspace{0.1 cm}I\not\in\mathcal{S}\atop {\left|I\right|=p-2k}}\left( \ker  C^{\star}_{I}  \cap  \ker  \overline{C}^{\star}_{I}\cap \ker  \tr\right)\right)\bigcap\left(\displaystyle\bigcap_{I\in\mathbb{N}^{N},\hspace{0.1 cm} I \in \mathcal{S} \atop {\left|I\right|=p-2k}}\left(\ker  \left(z^{I}\right)^{\star}   \bigcap   \ker \left(\overline{z}^{I}\right)^{\star}\right)\right),\hspace{0.1 cm}\mbox{for all  $k=0,\dots, \left[\frac{p}{2}\right]$.}$$

We consider the Fischer Decompositions (\ref{new1}) choosing 
\begin{equation}P(z,\overline{z})=\frac{\varphi_{p}(z,\overline{z})-\overline{\varphi_{p}(z,\overline{z})}}{2\sqrt{-1}},\quad\mbox{for given $p\in\mathbb{N}^{\star}$.}\label{kama1}\end{equation}

The   $G$-component of the formal transformation (\ref{map}) is  iteratively computed by  imposing the   normalizations   described by (\ref{spartiuG}). In particular, these  polynomials (\ref{poll1}) are uniquely determined by (\ref{spartiuG}), but it is not clear if these polynomials (\ref{poll1}) are real valued, because the reality of the polynomial $P(z,\overline{z})$ does not guarantee the reality of the homogeneous polynomials involved in the iterative Fischer Decompositions  described by (\ref{new1}).

In order to show that these  Spaces of Normalizations (\ref{spartiuG}) uniquely determine  the $G$-component of the transformation (\ref{map}), it is required also to show the linear independence, considering complex numbers, of the following set of polynomials 
\begin{equation}\left\{C_{I}(z,\overline{z}),\hspace{0.1 cm} \overline{C_{I}(z,\overline{z})}\right\}_{I\in\mathbb{N}^{N}\atop {{\left|I\right|=p}\atop I\not\in\mathcal{S}}},\quad\mbox{for all $p\in\mathbb{N}^{\star}$.}\label{330}
\end{equation}
concluding, in particular, the following
\begin{equation}a_{I}=\overline{b_{I}},\quad\mbox{for all $I\in\mathbb{N}^{N}$ having length $p\geq 3$.}\label{331}
\end{equation}

These computations are difficult  to conclude because of  the overlapping of the  homogeneous   polynomials from (\ref{330}). It is thus necessary to return to (\ref{beb1se1}) in order to study more carefully the computations of their solution recalling (\ref{yx1}), (\ref{yx2}) and (\ref{yx3}) together with   (\ref{gringooo3}), (\ref{gringooo2}), (\ref{gringooo1}), (\ref{gringooo}), (\ref{gringo4se}), (\ref{gringo4}), (\ref{gringo3}), (\ref{gringo2})  and (\ref{gringo1}).  

The non-triviality of the Fischer Decomposition\cite{sh} forces the use of (\ref{opt}) as follows: for any given multi-index 
\begin{equation}
  I=\left(i_{1},\dots,i_{k},\dots,i_{N}\right)\in\mathbb{N}^{N}\quad\mbox{such that $i_{1}+\dots+i_{k}+\dots+i_{N}=p\geq 3, $}\label{IJ1}
\end{equation}
similarly somehow to (\ref{vectori}),   we consider by (\ref{IJ1}) the following vector
\begin{equation*} Z\left[I \right]=\begin{pmatrix}a_{\left(i_{1}-2,\dots,i_{k},\dots,i_{N}\right) }\\ \vdots \\ a_{\left(i_{1},\dots,i_{k}-2,\dots,i_{N}\right) }\\ \vdots \\  a_{\left(i_{1},\dots,i_{k},\dots,i_{N}-2\right) }.\end{pmatrix}
\end{equation*}

Immediately from (\ref{opt}), we obtain
\begin{equation}z^{I}-A_{I}(z,\overline{z}) Q(z,\overline{z})=C_{I}(z,\overline{z}),\quad\mbox{for all $I\not\in \mathcal{S}$ having   length $p$.}\label{opt1}
\end{equation} 

Clearly (\ref{opt1}) naturally introduces the following matrix
\begin{equation} 
\mbox{Aux}=  \begin{pmatrix} 1& \dots &  1& \dots & 1 \\ \vdots & \ddots & \vdots & \ddots & \vdots    \\ 1& \dots &  1 & \dots & 1\\ \vdots & \ddots & \vdots & \ddots & \vdots   \\ 1 & \dots &  1& \dots &  1
\end{pmatrix}\label{aux1},\end{equation} 
because  each very consistent  sum of terms  multiplied by $z_{1}^{2},z_{2}^{2},\dots,z_{N}^{2}$ in  (\ref{opt1}), generates by  (\ref{Lambda})  and (\ref{aux1})  the following terms $$\left( \mbox{Aux}\Lambda  \right)Z\left[I \right].
$$

Similarly as previously to (\ref{suc1}), (\ref{suc2}), (\ref{sucL1}), (\ref{sucL11}), (\ref{sucL1sese}), (\ref{sucL11sese}), we consider  by (\ref{IJ1}) the following matrices
 \begin{equation} \mbox{Aux}\left[I\right]= \begin{pmatrix}  1 & \dots   & 0 & \dots & 0 & \dots & 0& \dots & 0  \\ \vdots & \ddots & \vdots & \ddots & \vdots & \ddots & \vdots& \ddots & \vdots   \\ 0 & \dots & \lambda_{1} & \dots & \lambda_{k} & \dots & \lambda_{N} & \dots & 0 
\\ \vdots & \ddots & \vdots  & \ddots & \vdots  & \ddots & \vdots  & \ddots & \vdots  \\ 0 & \dots & \lambda_{1} & \dots & \lambda_{k}& \dots & \lambda_{N}& \dots & 0 \\ \vdots  & \ddots & \vdots  & \ddots & \vdots  & \ddots & \vdots  & \ddots & \vdots  \\ 0 & \dots &\lambda_{1} & \dots & \lambda_{k} & \dots & \lambda_{N}& \dots & 0 \\ \vdots  & \ddots & \vdots  & \ddots & \vdots  & \ddots & \vdots  & \ddots & \vdots  \\ 0  & \dots & 0  & \dots &0  & \dots & 0  & \dots & 1  
\end{pmatrix} .\label{calcan1}
\end{equation}

This matrix (\ref{calcan1})  has the characteristic the middle entry containing $\lambda_{k}$  is  the diagonal entry  on the  following row   
$$\left(i_{1},\dots,i_{k},\dots,i_{N}\right),$$  
according to the  lexicografic order related to (\ref{ODIN}), otherwise having  $1$  as diagonal entry,   being important in order to consider  products of matrices as in (\ref{90000se1extra}), (\ref{90000}), (\ref{90000se}), (\ref{90000se1})  in order to make  approximations of the formal transformation (\ref{map}). 

Defining $$Z^{t}=\left\{a_{I}\right\}_{I\in\mathbb{N}^{N}\atop \left|I\right|=p},$$
which respects the lexicografic order related to (\ref{ODIN}), we obtain the following system    of  equations 
\begin{equation}\left(I_{N^{p}}-\mbox{Aux}_{p}A\right) Z+B\overline{Z}=V\left(z_{1},z_{2},\dots,z_{N}\right),\quad\mbox{where $A$, $B\in \mathcal{M}_{N^{p}\times N^{p}}\left(\mathbb{C}\right)$}, \label{sisi}
\end{equation} 
and $V\left(z_{1},z_{2},\dots,z_{N}\right)$ is a known vector homogeneous polynomial  of degree $p-2$, dealing with the following products of ,,simple'' matrices
\begin{equation} 
\mbox{Aux}_{p}=\mbox{Aux}\left[p-2,\dots,0\right]\mbox{Aux}\left[p-3,1,\dots,0\right]\dots\mbox{Aux}\left[0,\dots,1,p-3\right] \mbox{Aux}\left[0,\dots,p-2\right],\quad\mbox{for all $p \geq 3$.}\label{calcan2}
\end{equation}

In order to show  that (\ref{sisi}) has unique solution, it suffices to show  the existences of each of the following matrices
\begin{equation}\frac{1}{I_{N^{p}}-\mbox{Aux}_{p}A-B},\quad\frac{1}{I_{N^{p}}-\mbox{Aux}_{p}A+B}.\label{550}
\end{equation}

It remains now to explain the existence these two appearing matrices in (\ref{550})  recalling  (\ref{nene1}), (\ref{nene2}), (\ref{nene3}), (\ref{nene4}), (\ref{nene5}), (\ref{nene4441}), (\ref{nene4442}), (\ref{gringooo31}), which are matrices that clearly exist in the light of the invetibility of the matrices from (\ref{urs11}), (\ref{urs11se}), (\ref{urs11sese}) and (\ref{urs11sesese}). 

We have
\begin{equation}  \left\| \frac{\mathcal{M}_{2,1}}{\mathcal{M}_{2,2}}\frac{V_{1}}{\mathcal{M}_{1,1}}\dots\right\|<\frac{1}{2}\left(1+\frac{1}{4}+\dots\right)=\frac{2}{3} , \label{ioc1}
\end{equation}
in the light of (\ref{expan})  and in the light of the following evaluations
$$ \left\| \frac{\mathcal{M}_{4,3}}{\mathcal{M}_{4,4}}\frac{\mathcal{M}_{3,4}}{\mathcal{M}_{3,3}} \right\|<\frac{1}{2}\quad \left\| \frac{\mathcal{M}_{4,3}}{\mathcal{M}_{4,4}}\frac{\mathcal{M}_{3,4}}{\mathcal{M}_{3,3}} \frac{\mathcal{M}_{3,1}}{\mathcal{M}_{3,3}}\mathcal{N}_{1,3} \right\|<\frac{1}{24},  
$$
which are clear according to the estimations related to (\ref{viena1}), (\ref{viena2}), (\ref{viena3}), observing that the resulted denominador may be multiplied by $2$ by similar computations, improving thus the initial estimations in the light of similar expansions as (\ref{york}), (\ref{expan}),  (\ref{vulpe}) and (\ref{9909}). 

Moving forward in (\ref{beb120}) from left to right, we obtain   
\begin{equation}\left|B\right\|<\frac{1}{4},\label{iok1}
\end{equation}
because the matrix $B$ is defined via (\ref{gringooo31}). 

Similarly, returning from right to left in (\ref{beb1se1}) according to similar arguing, we obtain 
\begin{equation}\left\|\mbox{Aux}_{p}A\right\|<\frac{1}{4}.\label{iok2}\end{equation}
  
It becomes now clear  the existence of the  matrices from  (\ref{550}). Moreover, we  obtain
\begin{equation}\left\|\frac{1}{I_{N^{p}}-\mbox{Aux}_{p}A-B}\right\|<2,\quad\left\|\frac{1}{I_{N^{p}}-\mbox{Aux}_{p}A+B}\right\|<2.\label{550se}
\end{equation}
 
These estimations (\ref{550se}) are important in order to establish   bounds for the solutions of the systems of equations defined as in (\ref{beb1}), which give the computations of the $G$-component of the formal equivalence (\ref{map}) by considering products of matrices as in  (\ref{90000se1extra}), (\ref{90000}), (\ref{90000se}), (\ref{90000se1}) in (\ref{sisi}) using  the product of simple matrices (\ref{calcan2}). Then, the solution is computed, being unique in the light of the invertibility of the matrices (\ref{550}), after there are separated the real part and the imaginary in the local defining equation, which is defined by a real-analytic function, because the source real-manifold is assumed to be real-analytic.  Then the norms of these matrices (\ref{550}) naturally provide  suitable estimates   in order to conclude the convergence of the formal transformation (\ref{map})  using the rapid iterations arguments of Moser\cite{mo}. 

Moreover, the previous procedures   can be     generalized using the language of matrices, which gives  computations   described as follows:
\section{Fischer  $F$-Decompositions\cite{bu4},\cite{bu5}}For any multi-index $$\mbox{$\tilde{J}=\left(\tilde{j}_{1},\tilde{j}_{2},\dots,\tilde{j}_{N}\right)\not\in \mathcal{T}_{l}$, for all   $l\in 1,\dots, N$}$$   we consider by (\ref{788}) the following Fischer Decompositions 
\begin{equation} \left(\overline{z}_{l}+2\lambda_{l}z_{l}\right) z^{\tilde{J}}=A_{l}(z,\overline{z}) Q(z,\overline{z})+C_{l}(z,\overline{z}),\quad \tilde{\tr}\left(C_{l}(z,\overline{z})\right)=0,\quad\mbox{for all   $l\in 1,\dots, N$. } \label{opsec} 
\end{equation}

These homogeneous polynomials $$A_{1}(z,\overline{z}),A_{2}(z,\overline{z}),\dots,A_{N}(z,\overline{z}), C_{1}(z,\overline{z}),C_{2}(z,\overline{z}),\dots,C_{N}(z,\overline{z})$$ are uniquely determined according to Shapiro\cite{sh}, being directly  computed  from (\ref{opsec}) writing   as follows
\begin{equation} A_{l}(z,\overline{z})=\displaystyle\sum_{I,J\in\mathbb{N}^{N}\atop\left|I\right|+\left|J\right|=p-2} a_{I,J}^{\left(l\right)}z^{I}\overline{z}^{J}, \quad C_{l}(z,\overline{z})=\displaystyle\sum_{I,J\in\mathbb{N}^{N}\atop\left|I\right|+\left|J\right|=p} c_{I,J}^{\left(l\right)}z^{I}\overline{z}^{J},\quad\mbox{for all $l=1,\dots,N$.}\label{oppvsec}
\end{equation}

We apply the operator $\tilde{tr}$ in (\ref{opsec}) using (\ref{oppvsec}) and recalling (\ref{909se11}). We obtain 
\begin{equation}\begin{split}& \lambda_{1}^{2}\tilde{j}_{1}\left(\tilde{j}_{1}-1\right)z_{1}^{\tilde{j}_{1}-2}\dots z_{l}^{\tilde{j}_{l}+1} \dots z_{N}^{\tilde{j}_{N}}+\dots+\lambda_{l}\lambda_{1}\tilde{j}_{l}\left(\tilde{j}_{l}+1\right)z_{1}^{\tilde{j}_{1}}\dots z_{l}^{\tilde{j}_{l}-1} \dots z_{N}^{\tilde{j}_{N}}+\dots+\lambda_{N}\lambda_{1}\tilde{j}_{N}\left(\tilde{j}_{N}-1\right)z_{1}^{\tilde{j}_{1}}\dots z_{l}^{\tilde{j}_{l}+1} \dots z_{N}^{\tilde{j}_{N}-2} \\&      \overline{z}_{l}\left(\lambda_{1} \tilde{j}_{1}\left(\tilde{j}_{1}-1\right)z_{1}^{\tilde{j}_{1}-2}  \dots z_{N}^{\tilde{j}_{N}}+\dots+ \lambda_{N}\tilde{j}_{N}\left(\tilde{j}_{N}-1\right)z_{1}^{\tilde{j}_{1}}  \dots z_{N}^{\tilde{j}_{N}-2}\right)+j_{l} z_{1}^{\tilde{j}_{1}}\dots  z_{l}^{\tilde{j}_{l}-1} \dots z_{N}^{\tilde{j}_{N}}  \\& \quad\quad\quad\quad\quad   \quad\quad\quad\quad\quad\quad  \quad\quad\quad\quad\quad \quad\quad\quad  \quad\quad\quad\quad\quad\quad\quad\quad  \begin{tabular}{l} \rotatebox[origin=c]{270}{$=$}\end{tabular} \\& \tilde{\tr}\left(A_{l}(z,\overline{z}) \right)Q(z,\overline{z})+\displaystyle\sum_{k=1}^{N}\left(A_{l}\right)_{z_{k}}(z,\overline{z})\left(z_{k}\left(1+4\lambda_{k}^{2}\right)+4\lambda_{k}\overline{z}_{k} \right)+\displaystyle\sum_{k=1}^{N}\left(A_{l}\right)_{\overline{z}_{k}}(z,\overline{z})\left(\overline{z}_{k}\left(1+4\lambda_{k}^{2}\right)+4\lambda_{k}z_{k}\right) +A_{l}(z,\overline{z})\left(N+4\lambda_{1}^{2}+\dots+4\lambda_{N}^{2}\right),\end{split}\label{VVV1se}
\end{equation} 
for all $l=1,\dots,N$. 

On the other hand, we have    
 \begin{equation}\begin{split}& \quad\quad\quad\quad\quad\quad\quad\quad\quad\quad\quad\quad\quad\quad\quad\quad\quad\quad\quad\quad\quad\quad\quad\quad\quad \tilde{\tr}\left(A_{l}(z,\overline{z})\right)Q(z,\overline{z})\\& \quad\quad\quad\quad\quad   \quad\quad\quad\quad\quad\quad  \quad\quad\quad\quad\quad \quad\quad\quad  \quad\quad\quad\quad\quad\quad\quad\quad  \begin{tabular}{l} \rotatebox[origin=c]{270}{$=$}\end{tabular} \\&  \left(z_{1}\overline{z}_{1}+\dots+z_{N}\overline{z}_{N}+\lambda_{1}\left(z_{1}^{2}+\overline{z}_{1}^{2}\right)+\dots+\lambda_{N}\left(z_{N}^{2}+\overline{z}_{N}^{2}\right)\right)\left(\displaystyle\sum_{k=1}^{N}\displaystyle\sum_{I,J\in\mathbb{N}^{N}\atop\left|I\right|+\left|J\right|=p-2} a_{I,J}^{\left(l\right)}i_{k}j_{k} z_{1}^{i_{1}}\dots z_{k}^{i_{k}-1} \dots z_{N}^{i_{N}} \overline{z}_{1}^{j_{1}}\dots \overline{z}_{k}^{j_{k}-1}\dots \overline{z}_{N}^{j_{N}}+\right.\\& \quad\quad\quad\quad\quad   \quad\quad\quad\quad\quad\quad  \quad\quad\quad\quad\quad \quad\quad\quad  \quad\quad\quad\quad\quad\quad\quad\quad  \begin{tabular}{l} \rotatebox[origin=c]{270}{$+$}\end{tabular} \\& \left.\displaystyle\sum_{k=1}^{N}\displaystyle\sum_{I,J\in\mathbb{N}^{N}\atop\left|I\right|+\left|J\right|=p-2} a_{I,J}^{\left(l\right)}\lambda_{k} i_{k}\left(i_{k}-1\right) z_{1}^{i_{1}}\dots z_{k}^{i_{k}-2}\dots z_{N}^{i_{N}}\overline{z}_{1}^{j_{1}}\dots  \overline{z}_{N}^{j_{N}}+\displaystyle\sum_{k=1}^{N}\displaystyle\sum_{I,J\in\mathbb{N}^{N}\atop\left|I\right|+\left|J\right|=p-2} a_{I,J}^{\left(l\right)}\lambda_{k} j_{k}\left(j_{k}-1\right) z_{1}^{i_{1}}\dots  z_{N}^{i_{N}} \overline{z}_{1}^{j_{1}}\dots \overline{z}_{k}^{j_{k}-2}\dots \overline{z}_{N}^{j_{N}}\right) ,\end{split}
\label{B1se}\end{equation}
for all $l=1,\dots,N$. 

Clearly, the above interactions of terms are more complicated than the previous case. However, we impose the previous strategy according to the lexicografic order related to (\ref{ODIN}) on each component of the $F$-component of the formal equivalence (\ref{map}). More precisely, we  extract homogeneous terms in (\ref{VVV1se})  considering   the following vectors 
\begin{equation}\begin{split}&  \tilde{Y}_{1}^{t}=\left\{\left(a_{I,0}^{\left(1\right)}\right),\left(a_{I,0}^{(2)}\right),\dots,\left(a_{I,0}^{(N)}\right)\right\}_{I\in\mathbb{N}^{N}\atop \left|I\right|=p},\hspace{0.1 cm} \tilde{Y}_{2}^{t}=\left\{\left(a_{I,J}^{\left(1\right)}\right),\left(a_{I,J}^{(2)}\right),\dots,\left(a_{I,J}^{(N)}\right)\right\}_{I,J\in\mathbb{N}^{N}\atop \left|I\right|=p-1, \left|J\right|=1},\dots,\\& \tilde{Y}_{k}^{t}=\left\{\left(a_{I,J}^{\left(1\right)}\right),\left(a_{I,J}^{(2)}\right),\dots,\left(a_{I,J}^{(N)}\right)\right\}_{I,J\in\mathbb{N}^{N}\atop \left|I\right|=p-k+1, \left|J\right|=k-1},\dots, \tilde{Y}_{p+1}^{t}=\left\{\left(a_{0,J}^{\left(1\right)}\right),\left(a_{0,J}^{(2)}\right),\dots,\left(a_{0,J}^{(N)}\right)\right\}_{J\in\mathbb{N}^{N}\atop   \left|J\right|=p},\end{split}
\label{91ASee}
\end{equation}where $p-1=\tilde{j}_{1}+\dots+\tilde{j}_{N}$. 

This approach allows the construction of    systems of equations   extracting homogeneous
 terms in (\ref{VVV1se})  defining
 \begin{equation} W_{1}=\begin{pmatrix}0 \\ \vdots\\ 0\\ j_{l}\\0\\ \vdots\\ 0\\\lambda_{N}\lambda_{1}\tilde{j}_{N}\left(\tilde{j}_{N}-1\right)\\ 0\\\vdots\\ 0\\ \lambda_{l}\lambda_{1}\tilde{j}_{l}\left(\tilde{j}_{l}+1\right)\\0\\ \vdots \\0\\  \lambda_{1}^{2}\tilde{j}_{1}\left(\tilde{j}_{1}-1\right)\\ 0 \\\vdots\\0\\ \lambda_{N}\tilde{j}_{N}\left(\tilde{j}_{N}-1\right)  \\0\\ \vdots  \\0\\ \lambda_{1}\tilde{i}_{1}\left(\tilde{i}_{1}-1\right)\\ 0\\ \vdots\\0\end{pmatrix},\quad \quad W_{0}=\begin{pmatrix}
  0\\ \vdots\\  0\\0\\0\\ \vdots\\  0\\0\\0 \\\vdots \\ 0\\0\\0 \\ \vdots\\   0\\0\\0\\  \vdots \\ 0\\0\\0 \\ \vdots\\ 0\\0\\0\\  \vdots\\  0 
\end{pmatrix} .\label{special2} 
\end{equation}

We obtain   the following system of equations
 \begin{equation}\begin{pmatrix} \tilde{\mathcal{M}}_{1,1} & \tilde{\mathcal{M}}_{1,2} & \tilde{\mathcal{M}}_{1,3} & \mbox{O}_{p+1} & \mbox{O}_{p+1} &\dots & \mbox{O}_{p+1} & \mbox{O}_{p+1} & \mbox{O}_{p+1} & \mbox{O}_{p+1} \\ \tilde{\mathcal{M}}_{2,1} & \tilde{\mathcal{M}}_{2,2} & \tilde{\mathcal{M}}_{2,3} & \tilde{\mathcal{M}}_{2,4} & \mbox{O}_{p+1} &\dots & \mbox{O}_{p+1} & \mbox{O}_{p+1} & \mbox{O}_{p+1} & \mbox{O}_{p+1} \\ \tilde{\mathcal{M}}_{3,1}  & \tilde{\mathcal{M}}_{3,2}  & \tilde{\mathcal{M}}_{3,3}  & \tilde{\mathcal{M}}_{3,4}  & \tilde{\mathcal{M}}_{3,4}  &\dots & \mbox{O}_{p+1} & \mbox{O}_{p+1} & \mbox{O}_{p+1} & \mbox{O}_{p+1} \\ \mbox{O}_{p+1} & \tilde{\mathcal{M}}_{4,2}  & \tilde{\mathcal{M}}_{4,3}  & \tilde{\mathcal{M}}_{4,4}  & \tilde{\mathcal{M}}_{4,5}  &\dots & \mbox{O}_{p+1} & \mbox{O}_{p+1} & \mbox{O}_{p+1} & \mbox{O}_{p+1} \\ \mbox{O}_{p+1} & \mbox{O}_{p+1} & \tilde{\mathcal{M}}_{5,3} & \tilde{\mathcal{M}}_{5,4}& \tilde{\mathcal{M}}_{5,5} &\dots & \mbox{O}_{p+1} & \mbox{O}_{p+1} & \mbox{O}_{p+1} & \mbox{O}_{p+1} \cr \mbox{O}_{p+1} & \mbox{O}_{p+1} & \mbox{O}_{p+1} & \tilde{\mathcal{M}}_{6,4} & \tilde{\mathcal{M}}_{6,5} &\dots & \mbox{O}_{p+1} & \mbox{O}_{p+1} & \mbox{O}_{p+1} & \mbox{O}_{p+1}   \\ \vdots & \vdots & \vdots & \vdots & \vdots & \ddots & \vdots & \vdots & \vdots & \vdots \\ \mbox{O}_{p+1} & \mbox{O}_{p+1} & \mbox{O}_{p+1} & \mbox{O}_{p+1} & \mbox{O}_{p+1} &\dots & \tilde{\mathcal{M}}_{p-2,p-2}  & \tilde{\mathcal{M}}_{p-2,p-1}  & \tilde{\mathcal{M}}_{p-2,p} & \mbox{O}_{p+1} \\ \mbox{O}_{p+1} & \mbox{O}_{p+1} & \mbox{O}_{p+1} & \mbox{O}_{p+1} & \mbox{O}_{p+1} &\dots & \tilde{\mathcal{M}}_{p-1,p-2}  & \tilde{\mathcal{M}}_{p-1,p-1}  &\tilde{\mathcal{M}}_{p-1,p}  & \tilde{\mathcal{M}}_{p-1,p+1}  \\ \mbox{O}_{p+1} & \mbox{O}_{p+1} & \mbox{O}_{p+1} & \mbox{O}_{p+1} & \mbox{O}_{p+1} &\dots & \tilde{\mathcal{M}}_{p,p-2}  & \tilde{\mathcal{M}}_{p,p-1}  & \tilde{\mathcal{M}}_{p,p}  & \tilde{\mathcal{M}}_{p,p+1}  \\ \mbox{O}_{p+1} & \mbox{O}_{p+1} & \mbox{O}_{p+1} & \mbox{O}_{p+1} & \mbox{O}_{p+1} &\dots & \mbox{O}_{p+1} & \tilde{\mathcal{M}}_{p+1,p-1}  & \tilde{\mathcal{M}}_{p+1,p}  & \tilde{\mathcal{M}}_{p+1,p+1}           
\end{pmatrix}\begin{pmatrix} \tilde{Y}_{1}\\\tilde{Y}_{2}\\\tilde{Y}_{3}\\\tilde{Y}_{4}\\\tilde{Y}_{5}\\ \vdots\\ \tilde{Y}_{p-3}\\ \tilde{Y}_{p-2} \\ \tilde{Y}_{p-1}\\ \tilde{Y}_{p}\\ \tilde{Y}_{p+1}
\end{pmatrix}=\begin{pmatrix} W_{1} \\ W_{0} \\ W_{0} \\ W_{0}\\ W_{0}  \\   \vdots \\ W_{0} \\ W_{0}\\ W_{0}\\ W_{0}\\ W_{0}
\end{pmatrix}, \label{Africabeb1}
\end{equation}
where its  elements are matrices defined  similarly as in (\ref{lil}):
\begin{equation} \begin{split}& \tilde{\mathcal{M}}_{k-2,k}=\tilde{\mathcal{W}}_{k}^{'},\quad\quad\quad  \quad\quad\quad\quad\quad\quad\quad\quad\quad\quad\quad\quad\quad\quad\quad\quad\quad\quad\hspace{0.1 cm}  \mbox{for all $k=3,\dots,p$,} \\&\tilde{\mathcal{M}}_{k-1,k}=\tilde{\mathcal{V}}_{k}^{(0)}+\tilde{\mathcal{W}}_{k}+\tilde{\mathcal{O}}_{k}^{'},\quad\quad\quad  \quad\quad\quad\quad\quad\quad\quad\quad\quad\quad\quad\quad\hspace{0.18 cm}\mbox{for all $k=2,\dots,p $,}\\&\tilde{\mathcal{M}}_{k,k}= \tilde{\mathcal{O}}_{k}^{(0)}+\tilde{\mathcal{O}}_{k}+\tilde{\mathcal{V}}_{k}^{'}+\tilde{\mathcal{W}}_{k}^{''}+\left(N+4\lambda_{1}^{2}+\dots+4\lambda_{N}^{2}\right)\tilde{I}_{N^{p}} ,\quad\mbox{for all $k=1,\dots,p$,}  \\& \tilde{\mathcal{M}}_{k+1,k}=\tilde{\mathcal{W}}_{k}^{(0)}+\tilde{\mathcal{V}}_{k}+\tilde{\mathcal{O}}_{k}^{''} ,\quad\quad\quad  \quad\quad\quad\quad\quad\quad\quad\quad\quad\quad\quad\quad\hspace{0.18 cm}\mbox{for all $k=1,\dots,p-1$,} \\&  \tilde{\mathcal{M}}_{k+2,k}=\tilde{\mathcal{V}}_{k}^{''},\quad\quad\quad\quad  \quad\quad\quad\quad\quad\quad\quad\quad\quad\quad\quad\quad\quad\quad\quad\quad\quad \hspace{0.2 cm} \mbox{for all $k=1,\dots,p-2$.} \end{split}  \label{Africalil}
\end{equation} 
These occurring matrices  are defined as follows:

$\bf{The\hspace{0.1 cm}matrices:}$
$$\bf{\left\{\tilde{\mathcal{O}}_{k}^{(0)}\right\}_{k=1,\dots,p}\hspace{0.1 cm}, \left\{\tilde{\mathcal{V}}_{k}^{(0)}\right\}_{k=2,\dots,p},\hspace{0.1 cm} \left\{\tilde{\mathcal{W}}_{k}^{(0)}\right\}_{k=1,\dots,p-1}.}$$

Occurring from (\ref{VVV1se}), we have used in (\ref{Africalil}) the following matrices 
\begin{equation}\tilde{\mathcal{O}}_{k}^{(0)}=\begin{pmatrix}  \mathcal{O}_{k}^{(0)1} &  \mbox{O}_{N^{p-1}}    &   \dots & \mbox{O}_{N^{p-1}}  &  \\    \mbox{O}_{N^{p-1}} &   \mathcal{O}_{k}^{(0)2}& \dots  & \mbox{O}_{N^{p-1}}  \\       \vdots   &  \vdots &   \ddots & \vdots\\ \mbox{O}_{N^{p-1}}  &    \mbox{O}_{N^{p-1}}  & \dots & \mathcal{O}_{k}^{(0)N} \end{pmatrix}\in\mathcal{M}_{N^{p}\times N^{p} }\left(\mathbb{C}\right),\quad\mbox{for all $k=1,\dots,p$,} \label{Africalili1}
\end{equation} 
where (\ref{lili1}) is respected by each of the following matrices
\begin{equation}\mathcal{O}_{k}^{(0)1},\mathcal{O}_{k}^{(0)2},\dots,\mathcal{O}_{k}^{(0)N}\in\mathcal{M}_{N^{p-1}\times N^{p-1} }\left(\mathbb{C}\right).\label{Congo}
\end{equation}

Then (\ref{lili}) defines in (\ref{Africalil}) also the following   matrices
\begin{equation}\tilde{\mathcal{V}}_{k}^{(0)}=\begin{pmatrix}  \mathcal{V}_{k}^{(0)1} &  \mbox{O}_{N^{p-1}}    &   \dots & \mbox{O}_{N^{p-1}}  &  \\    \mbox{O}_{N^{p-1}} &   \mathcal{V}_{k}^{(0)2}& \dots  & \mbox{O}_{N^{p-1}}  \\       \vdots   &  \vdots &   \ddots & \vdots\\ \mbox{O}_{N^{p-1}}  &    \mbox{O}_{N^{p-1}}  & \dots & \mathcal{V}_{k}^{(0)N} \end{pmatrix}\in\mathcal{M}_{N^{p}\times N^{p} }\left(\mathbb{C}\right),\quad\mbox{for all $k=2,\dots,p$,} \label{Africalili2}
\end{equation}
where (\ref{lili2}) is respected by each of the following matrices
\begin{equation}\mathcal{V}_{k}^{(0)1},\mathcal{V}_{k}^{(0)2},\dots,\mathcal{V}_{k}^{(0)N}\in\mathcal{M}_{N^{p-1}\times N^{p-1} }\left(\mathbb{C}\right),\label{Congoo1}
\end{equation}
and respectively,  the following matrices
\begin{equation}\tilde{\mathcal{W}}_{k}^{(0)}=\begin{pmatrix}  \mathcal{W}_{k}^{(0)1} &  \mbox{O}_{N^{p-1}}    &   \dots & \mbox{O}_{N^{p-1}}  &  \\    \mbox{O}_{N^{p-1}} &   \mathcal{W}_{k}^{(0)2}& \dots  & \mbox{O}_{N^{p-1}}  \\       \vdots   &  \vdots &   \ddots & \vdots\\ \mbox{O}_{N^{p-1}}  &    \mbox{O}_{N^{p-1}}  & \dots & \mathcal{W}_{k}^{(0)N} \end{pmatrix}\in\mathcal{M}_{N^{p-1}\times N^{p-1} }\left(\mathbb{C}\right),\quad\mbox{for all $k=1,\dots,p-1$,} \label{Africalili3}
\end{equation}
where (\ref{lili3}) is respected by each of the following matrices
\begin{equation}\mathcal{W}_{k}^{(0)1},\mathcal{W}_{k}^{(0)2},\dots,\mathcal{W}_{k}^{(0)N}\in\mathcal{M}_{N^{p-1}\times N^{p-1} }\left(\mathbb{C}\right).\label{Congoo2}
\end{equation}

Clearly, we were using in (\ref{Africalil})   the following notation
 \begin{equation}   \tilde{I}_{N^{p}}=
\begin{pmatrix}  I_{N^{p-1}} &  \mbox{O}_{N^{p-1}}    &   \dots & \mbox{O}_{N^{p-1}}  &  \\    \mbox{O}_{N^{p-1}} &   I_{N^{p-1}}& \dots  & \mbox{O}_{N^{p-1}}  \\       \vdots   &  \vdots &   \ddots & \vdots\\ \mbox{O}_{N^{p-1}}  &    \mbox{O}_{N^{p-1}}  & \dots & I_{N^{p-1}} \end{pmatrix}\in\mathcal{M}_{N^{p}\times N^{p} }\left(\mathbb{C}\right).  \label{Africasuc222}
\end{equation}

Thus, in (\ref{Africalil}) the following matrices are by (\ref{Africalili1}), (\ref{Africalili2}) and (\ref{Africalili3})  defined as follows
\begin{equation} \begin{split}&\hspace{0.05 cm}
\tilde{\mathcal{O}}_{k}^{(0)}=\tilde{\mathcal{O}}_{k}^{(0)}\left[p-2,\dots,0;0,\dots,0\right]\dots\tilde{\mathcal{O}}_{k}^{(0)}\left[p-3,\dots,0;1,\dots,0\right]\dots\tilde{\mathcal{O}}_{k}^{(0)}\left[0,\dots,0;0,\dots,p-2\right],\quad\hspace{0.22 cm}\mbox{for all $k=1,\dots,p$,}  \\&\hspace{0.1 cm}\tilde{\mathcal{V}}_{k}^{(0)}=\tilde{\mathcal{V}}_{k}^{(0)}\left[p-2,\dots,0;0,\dots,0\right]\dots\tilde{\mathcal{V}}_{k}^{(0)}\left[p-3,\dots,0;1,\dots,0\right]\dots \tilde{\mathcal{V}}_{k}^{(0)}\left[0,\dots,0;0,\dots,p-2\right],\quad\quad\mbox{for all $k=2,\dots,p$,}  \\& \tilde{\mathcal{W}}_{k}^{(0)}=\tilde{\mathcal{W}}_{k}^{(0)}\left[p-2,\dots,0;0,\dots,0\right]\dots\tilde{\mathcal{W}}_{k}^{(0)}\left[p-3,\dots,0;1,\dots,0\right]\dots \tilde{\mathcal{W}}_{k}^{(0)}\left[0,\dots,0;0,\dots,p-2\right],\quad\mbox{for all $k=1,\dots,p-1$.} \end{split} \label{90000se1extra16}
\end{equation}

$\bf{The\hspace{0.1 cm}matrices:}$
$$\bf{\left\{\tilde{\mathcal{O}}_{k}\right\}_{k=1,\dots,p},\hspace{0.1 cm}  \left\{\tilde{\mathcal{V}}_{k}\right\}_{k=1,\dots,p-1},\hspace{0.1 cm} \left\{\tilde{\mathcal{W}}_{k}\right\}_{k=2,\dots,p}.}$$

Respecting the imposed lexicografic order, the first   matrix from (\ref{idiot1}) induces by (\ref{yy}), (\ref{B1}),   (\ref{lili}) and  (\ref{idiot1se}) the following matrices  
\begin{equation}   \tilde{\mathcal{O}}_{k}\left[I;J\right]=
\begin{pmatrix}  \mathcal{O}_{k}^{1}\left[I;J\right] &  \mbox{O}_{N^{p-1}}    &   \dots & \mbox{O}_{N^{p-1}}  &  \\    \mbox{O}_{N^{p-1}} &   \mathcal{O}_{k}^{2}\left[I;J\right]& \dots  & \mbox{O}_{N^{p-1}}  \\       \vdots   &  \vdots &   \ddots & \vdots\\ \mbox{O}_{N^{p-1}}  &    \mbox{O}_{N^{p-1}}  & \dots & \mathcal{O}_{k}^{N}\left[I;J\right] \end{pmatrix}\in\mathcal{M}_{N^{p}\times N^{p} }\left(\mathbb{C}\right),\quad\mbox{for all $k=1,\dots,p$,}  \label{Africasuc1}
\end{equation}
where  (\ref{suc1}) is respected by each of the following matrices
\begin{equation}\mathcal{O}_{k}^{1}\left[I;J\right] ,\mathcal{O}_{k}^{2}\left[I;J\right] ,\dots,\mathcal{O}_{k}^{N}\left[I;J\right] \in\mathcal{M}_{N^{p-1}\times N^{p-1} }\left(\mathbb{C}\right).\label{Congo1}
\end{equation}

Similarly, we consider by (\ref{yy})  and (\ref{idiot1}) the following matrices
 \begin{equation}   \tilde{\mathcal{V}}_{k}\left[I;J\right]=
\begin{pmatrix}  \mathcal{V}_{k}^{1}\left[I;J\right] &  \mbox{O}_{N^{p-1}}    &   \dots & \mbox{O}_{N^{p-1}}  &  \\    \mbox{O}_{N^{p-1}} &   \mathcal{V}_{k}^{2}\left[I;J\right]& \dots  & \mbox{O}_{N^{p-1}}  \\       \vdots   &  \vdots &   \ddots & \vdots\\ \mbox{O}_{N^{p-1}}  &    \mbox{O}_{N^{p-1}}  & \dots & \mathcal{V}_{k}^{N}\left[I;J\right] \end{pmatrix}\in\mathcal{M}_{N^{p}\times N^{p} }\left(\mathbb{C}\right),\quad\mbox{for all $k=1,\dots,p-1$,}  \label{Africasuc2}
\end{equation}
where   (\ref{suc2}) is respected by each of the following matrices
\begin{equation}\mathcal{V}_{k}^{1}\left[I;J\right],\mathcal{V}_{k}^{2}\left[I;J\right],\dots,\mathcal{V}_{k}^{N}\left[I;J\right]\in\mathcal{M}_{N^{p-1}\times N^{p-1} }\left(\mathbb{C}\right).\label{Congo2}
\end{equation}

Then (\ref{suc2}) induces by (\ref{yy})  and (\ref{idiot1}) analogously another matrices denoted  as follows
\begin{equation}\tilde{\mathcal{W}}_{k}\left[I;J\right]=\begin{pmatrix}  \mathcal{W}_{k}^{1}\left[I;J\right] &  \mbox{O}_{N^{p-1}}    &   \dots & \mbox{O}_{N^{p-1}}  &  \\    \mbox{O}_{N^{p-1}} &   \mathcal{W}_{k}^{2}\left[I;J\right]& \dots  & \mbox{O}_{N^{p-1}}  \\       \vdots   &  \vdots &   \ddots & \vdots\\ \mbox{O}_{N^{p-1}}  &    \mbox{O}_{N^{p-1}}  & \dots & \mathcal{W}_{k}^{N}\left[I;J\right] \end{pmatrix}\in\mathcal{M}_{N^{p}\times N^{p} }\left(\mathbb{C}\right), \quad\mbox{for all $k=2,\dots,p$,}\label{Africasuc33}\end{equation}      where  (\ref{suc33}) is respected by each of the following matrices
\begin{equation}\mathcal{W}_{k}^{1}\left[I;J\right],\mathcal{W}_{k}^{2}\left[I;J\right],\dots,\mathcal{W}_{k}^{N}\left[I;J\right]\in\mathcal{M}_{N^{p-1}\times N^{p-1} }\left(\mathbb{C}\right).\label{Congo3}
\end{equation} 

 Thus in (\ref{Africalil}),  the following matrices are by (\ref{suc1}), (\ref{suc2}), (\ref{suc33}), (\ref{Africasuc1}), (\ref{Africasuc2}) and   (\ref{Africasuc33}) naturally defined  as follows
\begin{equation} \begin{split}&\hspace{0.05 cm}
\tilde{\mathcal{O}}_{k}=\tilde{\mathcal{O}}_{k}\left[p-2,\dots,0;0,\dots,0\right]\dots\tilde{\mathcal{O}}_{k}\left[p-3,\dots,0;1,\dots,0\right]\dots \tilde{\mathcal{O}}_{k}\left[0,\dots,0;0,\dots,p-2\right], \hspace{0.15 cm}\quad\mbox{for all $k=1,\dots,p$,} \\&\hspace{0.1 cm} \tilde{\mathcal{V}}_{k}=\tilde{\mathcal{V}}_{k}\left[p-2,\dots,0;0,\dots,0\right]\dots\tilde{\mathcal{V}}_{k}\left[p-3,\dots,0;1,\dots,0\right]\dots \tilde{\mathcal{V}}_{k}\left[0,\dots,0;0,\dots,p-2\right],\quad\quad\mbox{for all $k=1,\dots,p-1$,}  \\& \tilde{\mathcal{W}}_{k}=\tilde{\mathcal{W}}_{k}\left[p-2,\dots,0;0,\dots,0\right]\dots\tilde{\mathcal{W}}_{k}\left[p-3,\dots,0;1,\dots,0\right]\dots \tilde{\mathcal{W}}_{k}\left[0,\dots,0;0,\dots,p-2\right],\quad\mbox{for all $k=2,\dots,p$.}\end{split}  \label{Africa90000}
\end{equation}
 
$\bf{The\hspace{0.1 cm}matrices:}$
$$\bf{\left\{\tilde{\mathcal{O}}_{k}^{'}\right\}_{k=2,\dots,p},\hspace{0.1 cm}  \left\{\tilde{\mathcal{V}}_{k}^{'}\right\}_{k=1,\dots,p},\hspace{0.1 cm} \left\{\tilde{\mathcal{W}}_{k}^{'}\right\}_{k=3,\dots,p}.}$$ 

Respecting the  lexicografic order, the first  matrix from (\ref{idiot2}) induces   by   (\ref{yy}),   (\ref{B1}), (\ref{idiot2se}) and (\ref{lili})    the following matrices 
 \begin{equation}      \tilde{\mathcal{O'}}_{k}\left[I;J\right]= \begin{pmatrix}  \mathcal{O'}_{k}^{1}\left[I;J\right] &  \mbox{O}_{N^{p-1}}    &   \dots & \mbox{O}_{N^{p-1}}  &  \\    \mbox{O}_{N^{p-1}} &   \mathcal{O'}_{k}^{2}\left[I;J\right]& \dots  & \mbox{O}_{N^{p-1}}  \\       \vdots   &  \vdots &   \ddots & \vdots\\ \mbox{O}_{N^{p-1}}  &    \mbox{O}_{N^{p-1}}  & \dots & \mathcal{O'}_{k}^{N}\left[I;J\right] \end{pmatrix}\in \mathcal{M}_{N^{p}\times N^{p}}\left(\mathbb{C}\right),\quad\mbox{for all $k=2,\dots,p$,}\label{AfricasucL1} 
\end{equation} 
where   (\ref{sucL1}) is respected by each of the following matrices
\begin{equation}\mathcal{O'}_{k}^{1}\left[I;J\right],\mathcal{O'}_{k}^{2}\left[I;J\right],\dots,\mathcal{O'}_{k}^{N}\left[I;J\right]\in\mathcal{M}_{N^{p-1}\times N^{p-1} }\left(\mathbb{C}\right),\label{Congo11}
\end{equation}

According to the   lexicografic order, the second matrix from (\ref{idiot2}) induces similarly by (\ref{yy})   the following matrices
  \begin{equation}     \tilde{\mathcal{V'}}_{k}\left[I;J\right]= \begin{pmatrix}  \mathcal{V'}_{k}^{1}\left[I;J\right] &  \mbox{O}_{N^{p-1}}    &   \dots & \mbox{O}_{N^{p-1}}  &  \\    \mbox{O}_{N^{p-1}} &   \mathcal{V'}_{k}^{2}\left[I;J\right]& \dots  & \mbox{O}_{N^{p-1}}  \\       \vdots   &  \vdots &   \ddots & \vdots\\ \mbox{O}_{N^{p-1}}  &    \mbox{O}_{N^{p-1}}  & \dots & \mathcal{V'}_{k}^{N}\left[I;J\right] \end{pmatrix}\in\mathcal{M}_{N^{p}\times N^{p} }\left(\mathbb{C}\right),\quad\mbox{for all $k=1,\dots,p$,}\label{AfricasucL11} 
\end{equation} 
where  (\ref{sucL11}) is respected by the following matrices
\begin{equation}\mathcal{V'}_{k}^{1}\left[I;J\right],\mathcal{V'}_{k}^{2}\left[I;J\right],\dots,\mathcal{V'}_{k}^{N}\left[I;J\right]\in\mathcal{M}_{N^{p-1}\times N^{p-1} }\left(\mathbb{C}\right).\label{Congo12}
\end{equation}
 
Then (\ref{sucL11}) induces by (\ref{yy})  and (\ref{idiot2}) another matrices denoted   as follows 
\begin{equation}\tilde{\mathcal{W'}}_{k}\left[I;J\right]=\begin{pmatrix}  \mathcal{W'}_{k}^{1}\left[I;J\right] &  \mbox{O}_{N^{p-1}}    &   \dots & \mbox{O}_{N^{p-1}}  &  \\    \mbox{O}_{N^{p-1}} &   \mathcal{W'}_{k}^{2}\left[I;J\right]& \dots  & \mbox{O}_{N^{p-1}}  \\       \vdots   &  \vdots &   \ddots & \vdots\\ \mbox{O}_{N^{p-1}}  &    \mbox{O}_{N^{p-1}}  & \dots & \mathcal{W'}_{k}^{N}\left[I;J\right] \end{pmatrix} \in \mathcal{M}_{N^{p}\times N^{p}}\left(\mathbb{C}\right),\quad\mbox{for all $k=3,\dots,p$,} \label{AfricasucL111}\end{equation}
where   (\ref{sucL111}) is respected by each of the following matrices
\begin{equation}\mathcal{W'}_{k}^{1}\left[I;J\right],\mathcal{W'}_{k}^{2}\left[I;J\right],\dots,\mathcal{W'}_{k}^{N}\left[I;J\right]\in\mathcal{M}_{N^{p-1}\times N^{p-1} }\left(\mathbb{C}\right).\label{Congo13}
\end{equation}
 
Thus in (\ref{lil}), the following matrices are by (\ref{sucL1}), (\ref{sucL11}), (\ref{sucL111}), (\ref{AfricasucL1}), (\ref{AfricasucL11}) and (\ref{AfricasucL111}) naturally defined as follows
\begin{equation} \begin{split}&\hspace{0.05 cm}
\tilde{\mathcal{O'}}_{k}=\tilde{\mathcal{O'}}_{k}\left[p-2,\dots,0;0,\dots,0\right]\dots\tilde{\mathcal{O'}}_{k}\left[p-3,\dots,0;1,\dots,0\right]\dots\tilde{\mathcal{O'}}_{k}\left[0,\dots,0;0,\dots,p-2\right], \quad\hspace{0.18 cm}\mbox{for all $k=2,\dots,p$,} \\&\hspace{0.1 cm} \tilde{\mathcal{V'}}_{k}=\tilde{\mathcal{V'}}_{k}\left[p-2,\dots,0;0,\dots,0\right]\dots\tilde{\mathcal{V'}}_{k}\left[p-3,\dots,0;1,\dots,0\right]\dots\tilde{\mathcal{V'}}_{k}\left[0,\dots,0;0,\dots,p-2\right],\quad\quad\mbox{for all $k=1,\dots,p$,}  \\& \tilde{\mathcal{W'}}_{k}=\tilde{\mathcal{W'}}_{k}\left[p-2,\dots,0;0,\dots,0\right]\dots\tilde{\mathcal{W'}}_{k}\left[p-3,\dots,0;1,\dots,0\right]\dots\tilde{\mathcal{W'}}_{k}\left[0,\dots,0;0,\dots,p-2\right],\quad\mbox{for all $k=3,\dots,p$.}\end{split}  \label{Africa90000se}
\end{equation}
 
$\bf{The\hspace{0.1 cm}matrices:}$
$$ \bf{\left\{\tilde{\mathcal{O}}_{k}^{''}\right\}_{k=1,\dots,p-1},\hspace{0.1 cm}  \left\{\tilde{\mathcal{V}}_{k}^{''}\right\}_{k=1,\dots,p-2},\hspace{0.1 cm} \left\{\tilde{\mathcal{W}}_{k}^{''}\right\}_{k=1,\dots,p}.} $$

Analysing more  (\ref{B1}) and (\ref{B111})  according to  the    lexicografic order from (\ref{ODIN}), the first  matrix from (\ref{idiot3}) induces by    (\ref{yy}), (\ref{B1}),   (\ref{idiot3se}) and  (\ref{lili}) the following matrices 
 \begin{equation}  \tilde{\mathcal{O''}}_{k}\left[I;J\right]= \begin{pmatrix}  \mathcal{O''}_{k}^{1}\left[I;J\right] &  \mbox{O}_{N^{p-1}}    &   \dots & \mbox{O}_{N^{p-1}}  &  \\    \mbox{O}_{N^{p-1}} &   \mathcal{O''}_{k}^{2}\left[I;J\right]& \dots  & \mbox{O}_{N^{p-1}}  \\       \vdots   &  \vdots &   \ddots & \vdots\\ \mbox{O}_{N^{p-1}}  &    \mbox{O}_{N^{p-1}}  & \dots & \mathcal{O''}_{k}^{N}\left[I;J\right] \end{pmatrix}\in\mathcal{M}_{N^{p-1}\times N^{p-1} }\left(\mathbb{C}\right),\quad\mbox{for all $k=1,\dots,p-1$,}  \label{AfricasucL1sese}
\end{equation}
where  (\ref{sucL1sese}) is respected by each of the following matrices
\begin{equation}\mathcal{O''}_{k}^{1}\left[I;J\right],\mathcal{O''}_{k}^{2}\left[I;J\right],\dots,\mathcal{O''}_{k}^{N}\left[I;J\right]\in\mathcal{M}_{N^{p-1}\times N^{p-1} }\left(\mathbb{C}\right),\label{Congo111}
\end{equation}
 
Similarly as previously, we consider by (\ref{yy})  the following matrices
 \begin{equation} \tilde{\mathcal{V''}}_{k}\left[I;J\right]= \begin{pmatrix}  \mathcal{V''}_{k}^{1}\left[I;J\right] &  \mbox{O}_{N^{p-1}}    &   \dots & \mbox{O}_{N^{p-1}}  &  \\    \mbox{O}_{N^{p-1}} &   \mathcal{V''}_{k}^{2}\left[I;J\right]& \dots  & \mbox{O}_{N^{p-1}}  \\       \vdots   &  \vdots &   \ddots & \vdots\\ \mbox{O}_{N^{p-1}}  &    \mbox{O}_{N^{p-1}}  & \dots & \mathcal{V''}_{k}^{N}\left[I;J\right] \end{pmatrix}\in\mathcal{M}_{N^{p}\times N^{p} }\left(\mathbb{C}\right),\quad\mbox{for all $k=1,\dots,p-2$,} \label{AfricasucL11sese}
\end{equation}
where (\ref{sucL11sese}) is respected by each of the following matrices
\begin{equation}\mathcal{V''}_{k}^{1}\left[I;J\right],\mathcal{V''}_{k}^{2}\left[I;J\right],\dots,\mathcal{V''}_{k}^{N}\left[I;J\right]\in\mathcal{M}_{N^{p-1}\times N^{p-1} }\left(\mathbb{C}\right),\label{Congo112}\end{equation}

Then (\ref{sucL11sese}) induces by (\ref{yy})   another matrices denoted   as follows 
\begin{equation}\tilde{\mathcal{W''}}_{k}\left[I;J\right]=\begin{pmatrix}  \mathcal{W''}_{k}^{1}\left[I;J\right] &  \mbox{O}_{N^{p-1}}    &   \dots & \mbox{O}_{N^{p-1}}  &  \\    \mbox{O}_{N^{p-1}} &   \mathcal{W''}_{k}^{2}\left[I;J\right]& \dots  & \mbox{O}_{N^{p-1}}  \\       \vdots   &  \vdots &   \ddots & \vdots\\ \mbox{O}_{N^{p-1}}  &    \mbox{O}_{N^{p-1}}  & \dots & \mathcal{W''}_{k}^{N}\left[I;J\right] \end{pmatrix}\in\mathcal{M}_{N^{p}\times N^{p} }\left(\mathbb{C}\right)  ,\quad\mbox{ for all $k=1,\dots,p$}.\label{AficasucL111se}\end{equation}  
where (\ref{sucL111se}) is respected by each of the following matrices
\begin{equation}\mathcal{W''}_{k}^{1}\left[I;J\right],\mathcal{W''}_{k}^{2}\left[I;J\right],\dots,\mathcal{W''}_{k}^{N}\left[I;J\right]\in\mathcal{M}_{N^{p-1}\times N^{p-1} }\left(\mathbb{C}\right).\label{Congo113}\end{equation}
    
Thus, in (\ref{Africalil}) the following matrices are by (\ref{suc1}), (\ref{suc2}),  (\ref{sucL111se}), (\ref{AfricasucL1sese}), (\ref{AfricasucL11sese}) and  (\ref{AficasucL111se})   defined as follows
\begin{equation} \begin{split}&\hspace{0.05 cm}
\tilde{\mathcal{O''}}_{k}=\tilde{\mathcal{O''}}_{k}\left[p-2,\dots,0;0,\dots,0\right]\dots\tilde{\mathcal{O''}}_{k}\left[p-3,\dots,0;1,\dots,0\right]\dots\tilde{\mathcal{O''}}_{k}\left[0,\dots,0;0,\dots,p-2\right],\quad \hspace{0.15 cm}\mbox{for all $k=1,\dots,p-1$,}  \\& \hspace{0.1 cm} \tilde{\mathcal{V''}}_{k}=\tilde{\mathcal{V''}}_{k}\left[p-2,\dots,0;0,\dots,0\right]\dots\tilde{\mathcal{V''}}_{k}\left[p-3,\dots,0;1,\dots,0\right]\dots \tilde{\mathcal{V''}}_{k}\left[0,\dots,0;0,\dots,p-2\right],\quad\quad \mbox{for all $k=1,\dots,p-2$,}  \\& \tilde{\mathcal{W''}}_{k}=\tilde{\mathcal{W''}}_{k}\left[p-2,\dots,0;0,\dots,0\right]\dots\tilde{\mathcal{W''}}_{k}\left[p-3,\dots,0;1,\dots,0\right]\dots \tilde{\mathcal{W''}}_{k}\left[0,\dots,0;0,\dots,p-2\right],\quad\mbox{for all $k=1,\dots,p$.}\end{split} \label{Afica90000se1}
\end{equation}
 
It is crucial now to observe now the    invertibility  of each of the following matrices
\begin{equation}\mathcal{L}_{k}=\mathcal{O}_{k}^{(0)}+\mathcal{O}_{k}+\mathcal{V}_{k}^{'}+\mathcal{W}_{k}^{''},\quad\mbox{for all $k= 1,\dots,p$}.\label{Africahihi}\end{equation}

We recall (\ref{lili1}), (\ref{suc1}), (\ref{90000}), (\ref{sucL11}),  (\ref{90000se}), (\ref{sucL111se}) and (\ref{90000se1}). We obtain
 \begin{equation}    \mathcal{M}_{N^{p}\times N^{p}}\left(\mathbb{C}\right)\ni\mathcal{L}_{k}\left[I;J\right]=\begin{pmatrix}  \mathcal{L}_{k}^{1}\left[I;J\right] &  \mbox{O}_{N^{p-1}}    &   \dots & \mbox{O}_{N^{p-1}}  &  \\    \mbox{O}_{N^{p-1}} &   \mathcal{L}_{k}^{2}\left[I;J\right]& \dots  & \mbox{O}_{N^{p-1}}  \\       \vdots   &  \vdots &   \ddots & \vdots\\ \mbox{O}_{N^{p-1}}  &    \mbox{O}_{N^{p-1}}  & \dots & \mathcal{L}_{k}^{N}\left[I;J\right] \end{pmatrix}\in\mathcal{M}_{N^{p}\times N^{p} }\left(\mathbb{C}\right),\quad\mbox{for all $k=1,\dots,p$,}\label{Africashobi}
\end{equation}
where   (\ref{shobi}) is respected by each of the following matrices
\begin{equation}\mathcal{L}_{k}^{1}\left[I;J\right],\mathcal{L}_{k}^{2}\left[I;J\right],\dots,\mathcal{L}_{k}^{N}\left[I;J\right]\in\mathcal{M}_{N^{p-1}\times N^{p-1} }\left(\mathbb{C}\right),\label{Congo113zambia}\end{equation}
\subsection{Systems of equations}We  simplify  the system of equations (\ref{Africabeb120}). It is  equivalent to
 \begin{equation}\begin{pmatrix} I_{N^{p}} & \frac{\tilde{\mathcal{M}}_{1,2}}{\tilde{\mathcal{M}}_{1,1}} & \frac{\tilde{\mathcal{M}}_{1,3}}{\tilde{\mathcal{M}}_{1,1}} & \mbox{O}_{N^{p}} & \mbox{O}_{N^{p}} &\dots & \mbox{O}_{N^{p}} & \mbox{O}_{N^{p}} & \mbox{O}_{N^{p}} & \mbox{O}_{N^{p}} \\ \frac{\tilde{\mathcal{M}}_{2,1}}{\tilde{\mathcal{M}}_{2,2}} & I_{N^{p}} & \frac{\tilde{\mathcal{M}}_{2,3}}{\tilde{\mathcal{M}}_{2,2}} & \frac{\tilde{\mathcal{M}}_{2,4}}{\tilde{\mathcal{M}}_{2,2}} & \mbox{O}_{N^{p}} &\dots & \mbox{O}_{N^{p}} & \mbox{O}_{N^{p}} & \mbox{O}_{N^{p}} & \mbox{O}_{N^{p}} \\ \frac{\tilde{\mathcal{M}}_{3,1}}{\tilde{\mathcal{M}}_{3,3}}  & \frac{\tilde{\mathcal{M}}_{3,2}}{\tilde{\mathcal{M}}_{3,3}}  & I_{N^{p}}  & \frac{\tilde{\mathcal{M}}_{3,4}}{\tilde{\mathcal{M}}_{3,3}}  & \frac{\tilde{\mathcal{M}}_{3,5}}{\tilde{\mathcal{M}}_{3,3}}  &\dots & \mbox{O}_{N^{p}} & \mbox{O}_{N^{p}} & \mbox{O}_{N^{p}} & \mbox{O}_{N^{p}}  \\ \mbox{O}_{N^{p}} & \frac{\tilde{\mathcal{M}}_{4,2}}{\tilde{\mathcal{M}}_{4,4}}  & \frac{\tilde{\mathcal{M}}_{4,3}}{\tilde{\mathcal{M}}_{4,4}}  & I_{N^{p}} & \frac{\tilde{\mathcal{M}}_{4,5}}{\tilde{\mathcal{M}}_{4,4}}  &\dots & \mbox{O}_{N^{p}}& \mbox{O}_{N^{p}} &\mbox{O}_{N^{p}}& \mbox{O}_{N^{p}} \\ \mbox{O}_{N^{p}}& \mbox{O}_{N^{p}} & \frac{\tilde{\mathcal{M}}_{5,3}}{\tilde{\mathcal{M}}_{5,5}} & \frac{\tilde{\mathcal{M}}_{5,4}}{\tilde{\mathcal{M}}_{5,5}}& I_{N^{p}} &\dots & \mbox{O}_{N^{p}} & \mbox{O}_{N^{p}} & \mbox{O}_{N^{p}} &\mbox{O}_{N^{p}}   \\ \vdots & \vdots & \vdots & \vdots & \vdots & \ddots & \vdots & \vdots & \vdots & \vdots \\ \mbox{O}_{N^{p}} & \mbox{O}_{N^{p}} & \mbox{O}_{N^{p}} & \mbox{O}_{N^{p}} & \mbox{O}_{N^{p}} &\dots & I_{N^{p}} & \frac{\tilde{\mathcal{M}}_{p-2,p-1}}{\tilde{\mathcal{M}}_{p-2,p-2}}  & \frac{\tilde{\mathcal{M}}_{p-2,p}}{\tilde{\mathcal{M}}_{p-2,p-2}} & \mbox{O}_{N^{p}} \\ \mbox{O}_{N^{p}} & \mbox{O}_{N^{p}} & \mbox{O}_{N^{p}} & \mbox{O}_{N^{p}} & \mbox{O}_{N^{p}} &\dots & \frac{\tilde{\mathcal{M}}_{p-1,p-2}}{\tilde{\mathcal{M}}_{p-1,p-1}}  & I_{N^{p}} & \frac{\tilde{\mathcal{M}}_{p-1,p}}{\tilde{\mathcal{M}}_{p-1,p-1}}  & \frac{\tilde{\mathcal{M}}_{p-1,p+1}}{\tilde{\mathcal{M}}_{p-1,p-1}}   \\ \mbox{O}_{N^{p}} & \mbox{O}_{N^{p}} & \mbox{O}_{N^{p}} & \mbox{O}_{N^{p}} & \mbox{O}_{N^{p}} &\dots & \frac{\tilde{\mathcal{M}}_{p,p-2}}{\tilde{\mathcal{M}}_{p,p}}  & \frac{\tilde{\mathcal{M}}_{p,p-1}}{\tilde{\mathcal{M}}_{p,p}}  & I_{N^{p}}  & \frac{\tilde{\mathcal{M}}_{p,p+1}}{\tilde{\mathcal{M}}_{p,p}}  \\ \mbox{O}_{N^{p}} & \mbox{O}_{N^{p}} & \mbox{O}_{N^{p}} & \mbox{O}_{N^{p}} & \mbox{O}_{N^{p}} &\dots & \mbox{O}_{N^{p}} & \frac{\tilde{\mathcal{M}}_{p+1,p-1}}{\mathcal{\tilde{M}}_{p+1,p+1}}  & \frac{\tilde{\mathcal{M}}_{p+1,p}}{\tilde{\mathcal{M}}_{p+1,p+1}}  & I_{N^{p}}        
\end{pmatrix}\begin{pmatrix} \tilde{Y}_{1}\\\tilde{Y}_{2}\\\tilde{Y}_{3}\\\tilde{Y}_{4}\\\tilde{Y}_{5}\\ \vdots\\ \tilde{Y}_{p-3}\\ \tilde{Y}_{p-2} \\ \tilde{Y}_{p-1}\\ \tilde{Y}_{p}\\ \tilde{Y}_{p+1}
\end{pmatrix}=\begin{pmatrix} \frac{1}{\tilde{\mathcal{M}}_{11}}W_{1}\\ W_{0} \\ W_{0} \\ W_{0}\\ W_{0}  \\   \vdots \\ W_{0} \\ W_{0}\\ W_{0}\\ W_{0}\\ W_{0}
\end{pmatrix}, \label{Africabeb120}
\end{equation}
where its  elements are matrix defined by (\ref{Africalil}), writing its  first equation   as follows
\begin{equation}\mbox{I}_{N^{p}} \tilde{Y}_{1}+ \tilde{\mathcal{N}}_{1,2}\tilde{Y}_{2}+ \tilde{\mathcal{N}}_{1,3}\tilde{Y}_{3}+\mbox{O}_{N^{p}}\tilde{Y}_{4}+\dots+\mbox{O}_{N^{p}}\tilde{Y}_{p+1}=\tilde{W}_{1},\label{Africagringo1}
\end{equation}
where we have used the following matrices
\begin{equation}\tilde{\mathcal{N}}_{1,2}=\frac{\tilde{\mathcal{M}}_{1,2}}{\tilde{\mathcal{M}}_{1,1}},\quad \tilde{\mathcal{N}}_{1,3}=\frac{\tilde{\mathcal{M}}_{1,3}}{\tilde{\mathcal{M}}_{1,1}},\quad \tilde{W}_{1}=\frac{1}{\tilde{\mathcal{M}}_{1,1}}W_{1}.\label{Africanene1} \end{equation}

 Making   substraction from the second equation in (\ref{Africabeb120}) using  (\ref{Africagringo1}), we obtain 
\begin{equation}\mbox{I}_{N^{p}} \tilde{Y}_{2}+\tilde{\mathcal{N}}_{2,3}\tilde{Y}_{3}+\tilde{\mathcal{N}}_{2,4}\tilde{Y}_{4}+\mbox{O}_{N^{p}}\tilde{Y}_{5}+\dots+\mbox{O}_{N^{p}}\tilde{Y}_{p+1}=\tilde{W}_{2},\label{Africagringo2}
\end{equation}
where we have used the following matrices
\begin{equation}\begin{split}& \tilde{\mathcal{N}}_{2,3}=\begin{pmatrix}  \tilde{\mathcal{N}}_{2,3}^{1} &  \mbox{O}_{N^{p-1}}    &   \dots & \mbox{O}_{N^{p-1}}  &  \\    \mbox{O}_{N^{p-1}} &   \tilde{\mathcal{N}}_{2,3}^{2}& \dots  & \mbox{O}_{N^{p-1}}  \\       \vdots   &  \vdots &   \ddots & \vdots\\ \mbox{O}_{N^{p-1}}  &    \mbox{O}_{N^{p-1}}  & \dots & \tilde{\mathcal{N}}_{2,3}^{N} \end{pmatrix},\hspace{0.1 cm}\tilde{\mathcal{N}}_{2,4}=\begin{pmatrix}  \tilde{\mathcal{N}}_{2,4}^{1} &  \mbox{O}_{N^{p-1}}    &   \dots & \mbox{O}_{N^{p-1}}  &  \\    \mbox{O}_{N^{p-1}} &   \tilde{\mathcal{N}}_{2,4}^{2}& \dots  & \mbox{O}_{N^{p-1}}  \\       \vdots   &  \vdots &   \ddots & \vdots\\ \mbox{O}_{N^{p-1}}  &    \mbox{O}_{N^{p-1}}  & \dots & \tilde{\mathcal{N}}_{2,4}^{N} \end{pmatrix}\in\mathcal{M}_{N^{p}\times N^{p} }\left(\mathbb{C}\right),\\& \quad\quad\quad\quad\quad\quad\quad\quad\quad\quad\quad\quad\quad\quad\quad\quad\quad\quad\quad\quad\quad\quad\quad\quad\quad\quad\quad\quad\quad\quad\quad\quad\quad\quad\quad \quad \tilde{W}_{2}=\begin{pmatrix}\tilde{V}_{2}^{1} \\ \tilde{V}_{2}^{2} \\ \vdots \\ \tilde{V}_{2}^{N}
\end{pmatrix}\in\mathcal{M}_{  N^{p}\times 1}\left(\mathbb{C}\right), \end{split}\label{Africanenee}
\end{equation} 
such that (\ref{nene2}) is respected by each of the following matrices
\begin{equation}\tilde{\mathcal{N}}_{2,3}^{1},\tilde{\mathcal{N}}_{2,3}^{2},\dots,\tilde{\mathcal{N}}_{2,3}^{N};\quad\tilde{\mathcal{N}}_{2,4}^{1},\tilde{\mathcal{N}}_{2,4}^{2},\dots,\tilde{\mathcal{N}}_{2,4}^{N}\in\mathcal{M}_{N^{p-1}\times N^{p-1} }\left(\mathbb{C}\right),\quad \tilde{V}_{2}^{1},\tilde{V}_{2}^{2},\dots,\tilde{V}_{2}^{N}\in\mathcal{M}_{1\times N^{p-1} }\left(\mathbb{C}\right).\label{Congoe}
\end{equation}

Recalling (\ref{gringo3}) and (\ref{ggrin2}), we obtain
\begin{equation}  \mbox{I}_{N^{p}} \tilde{Y}_{3}+\tilde{\mathcal{N}}_{3,4}\tilde{Y}_{4}+ \tilde{\mathcal{N}}_{3,5}\tilde{Y}_{5}+\mbox{O}_{N^{p}}\tilde{Y}_{6}+\dots+\mbox{O}_{N^{p}}\tilde{Y}_{p+1}=\tilde{W}_{3}, \label{Africagringo4} 
\end{equation} 
where we have used the following matrices
\begin{equation}\begin{split}&\tilde{\mathcal{N}}_{3,4}=\begin{pmatrix}  \tilde{\mathcal{N}}_{3,4}^{1} &  \mbox{O}_{N^{p-1}}    &   \dots & \mbox{O}_{N^{p-1}}  &  \\    \mbox{O}_{N^{p-1}} &   \tilde{\mathcal{N}}_{3,4}^{2}& \dots  & \mbox{O}_{N^{p-1}}  \\       \vdots   &  \vdots &   \ddots & \vdots\\ \mbox{O}_{N^{p-1}}  &    \mbox{O}_{N^{p-1}}  & \dots & \tilde{\mathcal{N}}_{3,4}^{N} \end{pmatrix},\hspace{0.1 cm}\tilde{\mathcal{N}}_{3,5}=\begin{pmatrix}  \tilde{\mathcal{N}}_{3,5}^{1} &  \mbox{O}_{N^{p-1}}    &   \dots & \mbox{O}_{N^{p-1}}  &  \\    \mbox{O}_{N^{p-1}} &   \tilde{\mathcal{N}}_{3,5}^{2}& \dots  & \mbox{O}_{N^{p-1}}  \\       \vdots   &  \vdots &   \ddots & \vdots\\ \mbox{O}_{N^{p-1}}  &    \mbox{O}_{N^{p-1}}  & \dots & \tilde{\mathcal{N}}_{3,5}^{N} \end{pmatrix}\in\mathcal{M}_{N^{p}\times N^{p} }\left(\mathbb{C}\right),\\& \quad\quad\quad\quad\quad\quad\quad\quad\quad\quad\quad\quad\quad\quad\quad\quad\quad\quad\quad\quad\quad\quad\quad\quad\quad\quad\quad\quad\quad\quad\quad\quad\quad\quad\quad \quad  \tilde{W}_{3}=\begin{pmatrix}\tilde{V}_{3}^{1} \\ \tilde{V}_{3}^{2} \\ \vdots \\ \tilde{V}_{3}^{N}
\end{pmatrix}\in\mathcal{M}_{  N^{p}\times 1}\left(\mathbb{C}\right), \end{split}   \label{Africanene2}
\end{equation} 
such that (\ref{nene3}) is respected by each of  the following matrices
\begin{equation}\tilde{\mathcal{N}}_{3,4}^{1},\tilde{\mathcal{N}}_{3,4}^{2},\dots,\tilde{\mathcal{N}}_{3,4}^{N};\quad\tilde{\mathcal{N}}_{3,5}^{1},\tilde{\mathcal{N}}_{3,5}^{2},\dots,\tilde{\mathcal{N}}_{3,5}^{N}\in\mathcal{M}_{N^{p-1}\times N^{p-1} }\left(\mathbb{C}\right)'\quad \tilde{V}_{3}^{1},\tilde{V}_{3}^{2},\dots,\tilde{V}_{3}^{N}\in\mathcal{M}_{1\times N^{p-1} }\left(\mathbb{C}\right).\label{Congoe1}
\end{equation}

We continue the computations imposing an induction process depending on $3 \leq l \leq p$. Assume 
\begin{equation}  \mbox{I}_{N^{p}} \tilde{Y}_{l-2}+\tilde{\mathcal{N}}_{l-2,l-1}\tilde{Y}_{l-1}+ \tilde{\mathcal{N}}_{l-2,l}\tilde{Y}_{l}+\mbox{O}_{N^{p}}\tilde{Y}_{l+1}+\dots+\mbox{O}_{N^{p}}\tilde{Y}_{p+1}=\tilde{W}_{l-2}, \label{Africagringo44se} 
\end{equation}  
and respectively 
\begin{equation}\mbox{I}_{N^{p}} \tilde{Y}_{l-1}+\tilde{\mathcal{N}}_{l-1,l}\tilde{Y}_{l}+ \tilde{\mathcal{N}}_{l-1,l+1}\tilde{Y}_{l+1}+\mbox{O}_{N^{p}}\tilde{Y}_{l+2}+\dots+\mbox{O}_{N^{p}}\tilde{Y}_{p+1}=\tilde{W}_{l-1}, \label{Africagringo44se1} 
\end{equation}  

We recall the equation $l$ from (\ref{beb120}) and also the substraction from  (\ref{gringoo1}) using (\ref{gringo44se}) concluding   (\ref{gringo55}). Recalling as well the 
substraction from   (\ref{gringo55}) using (\ref{gringo44se1}) concluding  (\ref{gringooo}), we obtain  
 \begin{equation}\mbox{I}_{N^{p}} \tilde{Y}_{1}+ \tilde{\mathcal{N}}_{l,l+1}\tilde{Y}_{l+1}+ \tilde{\mathcal{N}}_{l,l+2}\tilde{Y}_{l+2}+\mbox{O}_{N^{p}}\tilde{Y}_{4}+\dots+\mbox{O}_{N^{p}}\tilde{Y}_{p+1}=\tilde{W}_{l},\label{Africagringooo}
\end{equation} 
where we have used the following matrices
\begin{equation}\begin{split}&\tilde{\mathcal{N}}_{l,l+1}=\begin{pmatrix}  \tilde{\mathcal{N}}_{l,l+1}^{1} &  \mbox{O}_{N^{p-1}}    &   \dots & \mbox{O}_{N^{p-1}}  &  \\    \mbox{O}_{N^{p-1}} &   \tilde{\mathcal{N}}_{l,l+1}^{2}& \dots  & \mbox{O}_{N^{p-1}}  \\       \vdots   &  \vdots &   \ddots & \vdots\\ \mbox{O}_{N^{p-1}}  &    \mbox{O}_{N^{p-1}}  & \dots & \tilde{\mathcal{N}}_{l,l+1}^{N} \end{pmatrix},\hspace{0.1 cm}\tilde{\mathcal{N}}_{l,l+2}=\begin{pmatrix}  \tilde{\mathcal{N}}_{l,l+2}^{1} &  \mbox{O}_{N^{p-1}}    &   \dots & \mbox{O}_{N^{p-1}}  &  \\    \mbox{O}_{N^{p-1}} &   \tilde{\mathcal{N}}_{l,l+2}^{2}& \dots  & \mbox{O}_{N^{p-1}}  \\       \vdots   &  \vdots &   \ddots & \vdots\\ \mbox{O}_{N^{p-1}}  &    \mbox{O}_{N^{p-1}}  & \dots & \tilde{\mathcal{N}}_{l,l+2}^{N} \end{pmatrix}\in\mathcal{M}_{N^{p}\times N^{p} }\left(\mathbb{C}\right), \\& \quad\quad\quad\quad\quad\quad\quad\quad\quad\quad\quad\quad\quad\quad\quad\quad\quad\quad\quad\quad\quad\quad\quad\quad\quad\quad\quad\quad\quad\quad  \quad\quad\quad\quad\quad\quad       \tilde{W}_{l+1}=\begin{pmatrix}\tilde{V}_{l+1}^{1} \\ \tilde{V}_{l+1}^{2} \\ \vdots \\ \tilde{V}_{l+1}^{N}
\end{pmatrix}\in\mathcal{M}_{ N^{p}\times 1}\left(\mathbb{C}\right). \end{split}\label{Africanene22}
\end{equation} 
such that (\ref{nene5}) is respected by each of the following matrices
\begin{equation}\tilde{\mathcal{N}}_{l,l+1}^{1},\tilde{\mathcal{N}}_{l,l+1}^{2},\dots,\tilde{\mathcal{N}}_{l,l+1}^{N};\hspace{0.1 cm}\tilde{\mathcal{N}}_{l,l+2}^{1},\tilde{\mathcal{N}}_{l,l+2}^{2},\dots,\tilde{\mathcal{N}}_{l,l+2}^{N}\in\mathcal{M}_{N^{p-1}\times N^{p-1} }\left(\mathbb{C}\right);\quad \tilde{V}_{l+1}^{1},\tilde{V}_{l+1}^{2},\dots,\tilde{V}_{l+1}^{N}\in\mathcal{M}_{1\times N^{p-1} }\left(\mathbb{C}\right).\label{Congoe2}
\end{equation}
  
These   recurrences  (\ref{Africanene1}), (\ref{Africanene2}), (\ref{Africanene2}), (\ref{Africanene22})   are crucial for making relevant computations. Then    (\ref{Africabeb120})  is equivalent to 
 \begin{equation}\begin{pmatrix} \mbox{I}_{N^{p}} & \tilde{\mathcal{N}}_{1,2} & \tilde{\mathcal{N}}_{1,3} & \mbox{O}_{N^{p}} & \mbox{O}_{N^{p}} &\dots & \mbox{O}_{N^{p}} & \mbox{O}_{N^{p}} & \mbox{O}_{N^{p}} & \mbox{O}_{N^{p}} \\ \mbox{O}_{N^{p}} & \mbox{I}_{N^{p}} & \tilde{\mathcal{N}}_{2,3} &\tilde{\mathcal{N}}_{2,4} & \mbox{O}_{N^{p}} &\dots & \mbox{O}_{N^{p}} & \mbox{O}_{N^{p}} & \mbox{O}_{N^{p}} & \mbox{O}_{N^{p}} \\ \mbox{O}_{N^{p}} & \mbox{O}_{N^{p}} & \mbox{I}_{N^{p}} & \tilde{\mathcal{N}}_{3,4} & \tilde{\mathcal{N}}_{3,5} &\dots & \mbox{O}_{N^{p}} & \mbox{O}_{N^{p}} & \mbox{O}_{N^{p}} & \mbox{O}_{N^{p}} \\ \mbox{O}_{N^{p}} & \mbox{O}_{N^{p}} & \mbox{O}_{N^{p}} & \mbox{I}_{N^{p}} & \tilde{\mathcal{N}}_{4,5} &\dots & \mbox{O}_{N^{p}} & \mbox{O}_{N^{p}} & \mbox{O}_{N^{p}} & \mbox{O}_{N^{p}}   \\ \vdots & \vdots & \vdots & \vdots & \vdots & \ddots & \vdots & \vdots & \vdots & \vdots \\ \mbox{O}_{N^{p}} & \mbox{O}_{N^{p}} & \mbox{O}_{N^{p}} & \mbox{O}_{N^{p}} & \mbox{O}_{N^{p}} &\dots & \mbox{I}_{N^{p}} & \tilde{\mathcal{N}}_{p-2,p-1} & \tilde{\mathcal{N}}_{p-2,p} & \mbox{O}_{N^{p}} \\ \mbox{O}_{N^{p}} & \mbox{O}_{N^{p}} & \mbox{O}_{N^{p}} & \mbox{O}_{N^{p}} & \mbox{O}_{N^{p}} &\dots & \mbox{O}_{N^{p}} & \mbox{I}_{N^{p}} & \tilde{\mathcal{N}}_{p-1,p} &\tilde{\mathcal{N}}_{p-1,p+1} \\ \mbox{O}_{N^{p}} & \mbox{O}_{N^{p}} & \mbox{O}_{N^{p}} & \mbox{O}_{N^{p}} & \mbox{O}_{N^{p}} &\dots & \mbox{O}_{N^{p}} & \mbox{O}_{N^{p}} & \mbox{I}_{N^{p}} & \tilde{\mathcal{N}}_{p,p+1} \\ \mbox{O}_{N^{p}} & \mbox{O}_{N^{p}} & \mbox{O}_{N^{p}} & \mbox{O}_{N^{p}} & \mbox{O}_{N^{p}} &\dots & \mbox{O}_{N^{p}} & \mbox{O}_{N^{p}} & \mbox{O}_{N^{p}} & \mbox{I}_{N^{p}}
\end{pmatrix}\begin{pmatrix} \tilde{Y}_{1}\\ \tilde{Y}_{2}\\ \tilde{Y}_{3}\\ \tilde{Y}_{4}\\  \vdots \\ \tilde{Y}_{p-2} \\ \tilde{Y}_{p-1}\\ \tilde{Y}_{p}\\ \tilde{Y}_{p+1}
\end{pmatrix}=\begin{pmatrix} \tilde{W}_{1} \\ \tilde{W}_{2} \\ \tilde{W}_{3} \\ \tilde{W}_{4}\\    \vdots    \\ \tilde{W}_{p-2}\\ \tilde{W}_{p-1}\\ \tilde{W}_{p}\\ \tilde{W}_{p+1}
\end{pmatrix}.\label{Africabeb1se1}
\end{equation} 

Now, we  are ready to show that   (\ref{Africabeb1se1}) has unique solution. We  compute it going backwards among the equations in (\ref{Africabeb1se1}). We obtain
\begin{equation}
\tilde{Y}_{p}=\tilde{W}_{p}-\tilde{\mathcal{N}}_{p,p+1}\tilde{W}_{p+1}.\label{Africayx1}
\end{equation}

Then, we obtain
\begin{equation}
\tilde{Y}_{p-1}=\tilde{W}_{p-1}-\tilde{\mathcal{N}}_{p-1,p}\left( \tilde{W}_{p}-\tilde{\mathcal{N}}_{p,p+1}\tilde{W}_{p+1}\right)-   \tilde{\mathcal{N}}_{p-1,p+1}\tilde{W}_{p+1}.\label{Africayx2}
\end{equation}

Walking backwards among the equations in (\ref{Africabeb1se1}), we obtain
\begin{equation}
\tilde{Y}_{1}=\tilde{W}_{1}-\tilde{\mathcal{N}}_{1,2}\tilde{Y}_{2}-\tilde{\mathcal{N}}_{1,3}\tilde{Y}_{3}.\label{Africayx3}
\end{equation}

Now, the system of equations (\ref{Africabeb1se1}) is solved.
\subsection{Fischer Normalization $F$-Spaces\cite{bu2},\cite{bu3}}The Fischer Decomposition from (\ref{opsec}) gives
\begin{equation} \left(\overline{z}_{l}+2\lambda_{l}z_{l}\right) z^{J}=A_{l,J}(z,\overline{z}) Q(z,\overline{z})+C_{l,J}(z,\overline{z}),\quad  \tr \left(C_{l,J}(z,\overline{z})\right)=0,\quad\mbox{where $J\not\in \mathcal{T}_{l}$, for all  $l\in 1,\dots, N$.} \label{optsec} 
\end{equation}
 
This set of polynomials, derived from (\ref{optsec}),   gives  normalization conditions defining  certain  Spaces of Fischer Normalizations, which are constructed   iteratively from the generalized version of  the Fischer Decomposition\cite{sh}: for a given  real homogeneous polynomial of degree $p\geq 1$ in $(z,\overline{z})$ denoted by $P(z,\overline{z})$,  we  have
\begin{equation} \begin{split}& P(z,\overline{z})=P_{1}(z,\overline{z})Q(z,\overline{z})+R_{1}(z,\overline{z}),\quad\mbox{where $\tr\left(R_{1}(z,\overline{z})\right)=0$  and:}\\& \quad  R_{1}(z,\overline{z})=\displaystyle\sum_{l=1}^{N}\displaystyle\sum_{J\in\mathbb{N}^{N}\atop {\left|J\right|=p-1}} \left(a_{l,J}C_{l,J}(z,\overline{z})+b_{l,J}\overline{C_{l,J}(z,\overline{z})}\right)+R_{1,0}(z,\overline{z}),\quad\mbox{such that:}\\&  \quad\quad\quad\quad\quad\quad\quad\quad  R_{1,0}(z,\overline{z})\in \displaystyle\bigcap_{l=1}^{N}\left(\left(\displaystyle\bigcap_{J\in\mathbb{N}^{N},\hspace{0.1 cm}J\not\in \mathcal{T}_{l}\atop {\left|J\right|=p-1}} \left(\ker  C^{\star}_{l,J}  \bigcap  \ker  \overline{C}^{\star}_{l,J} \right)\right)\bigcap   \left(\displaystyle\bigcap_{J\in\mathbb{N}^{N},\hspace{0.1 cm}J\in \mathcal{T}_{l}\atop {\left|J\right|=p-1}} \left(\ker  \left( z_{l}\overline{z}^{J} \right)^{\star} \bigcap  \ker  \left( \overline{z}_{l}z^{J} \right)^{\star} \right)\right)\right),\\& P_{1}(z,\overline{z})=P_{2}(z,\overline{z})Q(z,\overline{z})+R_{2}(z,\overline{z}),\quad\mbox{where $ \tr\left(R_{2}(z,\overline{z})\right)=0$ and:}\\& \quad    R_{2}(z,\overline{z})=\displaystyle\sum_{l=1}^{N}\displaystyle\sum_{J\in\mathbb{N}^{N}\atop {\left|J\right|=p-3}} \left(a_{l,J}C_{l,J}(z,\overline{z})+b_{l,J}\overline{C_{l,J}(z,\overline{z})}\right)+R_{2,0}(z,\overline{z}),\quad\mbox{such that:} \\&\\&  \quad\quad\quad\quad\quad\quad\quad  R_{2,0}(z,\overline{z})\in \displaystyle\bigcap_{l=1}^{N}\left(\left(\displaystyle\bigcap_{J\in\mathbb{N}^{N},\hspace{0.1 cm}J\not\in \mathcal{T}_{l}\atop {\left|J\right|=p-3}}\left( \ker  C^{\star}_{l,J}  \bigcap   \ker  \overline{C}^{\star}_{l,J}\right)\right)\bigcap\left(   \displaystyle\bigcap_{J\in\mathbb{N}^{N},\hspace{0.1 cm}J\in \mathcal{T}_{l}\atop {\left|J\right|=p-3}} \ker  \left( z_{l}\overline{z}^{J} \right)^{\star} \bigcap   \ker  \left( \overline{z}_{l}z^{J} \right)^{\star} \right)\right),
\\&\quad\quad\quad\quad\vdots\quad\quad\quad\quad\quad\quad\quad\quad\vdots\quad\quad\quad\quad\quad\quad\quad\quad\vdots\quad\quad\quad\quad\quad\quad\quad\quad\vdots \quad\quad\quad\quad\quad\quad\quad\quad\quad\quad\quad\quad\quad\quad\vdots\\& P_{k}(z,\overline{z})=P_{k+1}(z,\overline{z})Q(z,\overline{z})+R_{k+1}(z,\overline{z}),\quad\mbox{where $\tr\left(R_{k+1}(z,\overline{z})\right)=0$  and:}\\&   R_{k+1}(z,\overline{z})=\displaystyle\sum_{l=1}^{N}\displaystyle\sum_{J\in\mathbb{N}^{N}\atop {\left|J\right|=p-2k-1}} \left(a_{l,J}C_{l,J}(z,\overline{z})+b_{l,J}\overline{C_{l,J}(z,\overline{z})}\right)+R_{k+1,0}(z,\overline{z}),\quad\mbox{such that:}\\&  \quad\quad\quad\quad\quad\quad\quad  R_{k+1,0}(z,\overline{z})\in \displaystyle\bigcap_{l=1}^{N}\left(\left(\displaystyle\bigcap_{J\in\mathbb{N}^{N},\hspace{0.1 cm}J\not\in \mathcal{T}_{l}\atop {\left|J\right|=p-2k-1}} \left(\ker  C^{\star}_{l,J}  \bigcap   \ker  \overline{C}^{\star}_{l,J} \right)\right)\bigcap\left( \displaystyle\bigcap_{J\in\mathbb{N}^{N},\hspace{0.1 cm}J\in \mathcal{T}_{l}\atop {\left|J\right|=p-2k-1}} \left(\ker  \left( z_{l}\overline{z}^{J} \right)^{\star} \bigcap   \ker  \left( \overline{z}_{l}z^{J} \right)^{\star}\right)\right)\right),\\&\quad\quad\quad\vdots\quad\quad\quad\quad\quad\quad\quad\quad\vdots\quad\quad\quad\quad\quad\quad\quad\quad\vdots\quad\quad\quad\quad\quad\quad\quad\quad\vdots \quad\quad\quad\quad\quad\quad\quad\quad\quad\quad\quad\quad\quad\quad\vdots\end{split}.\label{new2}
\end{equation}
where   these occurring polynomials 
\begin{equation}\left\{P_{k}(z,\overline{z})\right\}_{k=1,\dots,\left[\frac{p-1}{2}\right]},\quad \left\{R_{k}(z,\overline{z})\right\}_{k=1,\dots,\left[\frac{p-1}{2}\right]},
\label{poll2}
\end{equation}
are iteratively obtained  using the generalized version  of the Fischer Decomposition\cite{sh}. 

Recalling   strategies from \cite{bu2} and \cite{bu3}, we define 
\begin{equation}\mathcal{F}_{p},\quad p\in\mathbb{N}^{\star},\label{spartiuF}
\end{equation}
which consist in real-valued polynomials $P(z,\overline{z})$  of degree $p\geq 1$ in $(z,\overline{z})$ satisfying the normalizations:

$$ P_{k}^{\left(p\right)}(z,\overline{z})=P_{k+1}^{\left(p\right)}(z,\overline{z})Q(z,\overline{z})+R_{k+1}^{\left(p\right)}(z,\overline{z}),\quad \mbox{for all   $k=0,\dots, \left[\frac{p-1}{2}\right]$ and given $P_{0}^{\left(p\right)}(z,\overline{z})=P(z,\overline{z})$,}$$
such that
\begin{equation}
\begin{split}& R_{k+1}^{\left(p\right)}(z,\overline{z})\in \displaystyle\bigcap_{l=1}^{N}\left( \displaystyle\bigcap_{J\in\mathbb{N}^{N},\hspace{0.1 cm}J\not\in \mathcal{T}_{l}\atop {\left|J\right|=p-2k-1}} \left(\ker  C^{\star}_{l,J}  \bigcap   \ker  \overline{C}^{\star}_{l,J}\bigcap  \ker\tr\right)\right)\\&   \quad\quad\quad\quad\quad\quad\quad\quad \bigcap \left(  \displaystyle\bigcap_{J\in\mathbb{N}^{N},\hspace{0.1 cm}J\in \mathcal{T}_{l}\atop {\left|J\right|=p-2k-1}} \left(\ker  \left( z_{l}\overline{z}^{J} \right)^{\star} \bigcap   \ker  \left( \overline{z}_{l}z^{J} \right)^{\star}\right)\right),\quad \mbox{for all    $k=0,\dots, \left[\frac{p-1}{2}\right]$ .}\end{split}
\end{equation}

We  consider the Fischer Decompositions (\ref{new1}) choosing 
\begin{equation}P(z,\overline{z})=\frac{\varphi_{p}(z,\overline{z})+\overline{\varphi_{p}(z,\overline{z})}}{2},\quad\mbox{for given $p\in\mathbb{N}^{\star}$.}\label{kama2}\end{equation}

In order to show now that these Spaces of Normalizations (\ref{spartiuF}) uniquely determine  the $F$-component of the formal transformation (\ref{map}), it is required to show the linear independence, considering complex numbers, of the following set of polynomials 
\begin{equation}\left\{C_{l,J}(z,\overline{z}),\hspace{0.1 cm} \overline{C_{l,J}(z,\overline{z})}\right\}_{J\in\mathbb{N}^{N}\atop{\left|J\right|=p-1\atop{J\not\in\mathcal{T}_{l}\atop{l=1,\dots,N}}}},\quad\mbox{for all $p\in\mathbb{N}^{\star}$.}\label{Africa330}
\end{equation}

In particular, showing the linear independence of these polynomials   implies
\begin{equation}a_{l,J}=\overline{b_{l,J}},\quad\mbox{for all $J\in\mathbb{N}^{N}$ having length $p-1\geq 2$ and  $l=1,\dots,N$.}\label{Africa331}
\end{equation}

These computations are difficult  to conclude because of  the overlapping of the  homogeneous   polynomials from (\ref{Africa330}). There are recalled now the previous arguments corresponding to (\ref{330}) and (\ref{331}).  It is necessary to return in  (\ref{Africabeb120}) and (\ref{Africabeb1se1})  in order to study more carefully all their corresponding computations of its solution as follows: for any given multi-index 
\begin{equation}
  J=\left(j_{1},\dots,j_{k},\dots,j_{N}\right)\in\mathbb{N}^{N}\quad\mbox{such that $j_{1}+\dots+j_{k}+\dots+j_{N}=p-1\geq 2$,}\label{AfricaIJ1}
\end{equation}
similarly somehow to (\ref{vectori}), we consider by (\ref{AfricaIJ1}) the following vector
\begin{equation} \tilde{Z}\left[J \right]=\begin{pmatrix}a_{\left(j_{1}-2,\dots,j_{k},\dots,j_{N}\right)}^{\left(1\right)}\\ \vdots \\ a_{\left(j_{1},\dots,j_{k}-2,\dots,j_{N}\right)}^{\left(1\right)}\\ \vdots \\  a_{\left(j_{1},\dots,j_{k},\dots,j_{N}-2\right)}^{\left(1\right)}\\  \vdots\\ a_{\left(j_{1}-2,\dots,j_{k},\dots,j_{N}\right)}^{(N)}\\ \vdots \\ a_{\left(j_{1},\dots,j_{k}-2,\dots,j_{N}\right)}^{(N)}\\ \vdots \\  a_{\left(j_{1},\dots,j_{k},\dots,j_{N}-2\right)}^{(N)}.\end{pmatrix}\label{kaka}
\end{equation}

Immediately from (\ref{opt}), we obtain
\begin{equation}\left(\overline{z}_{l}+2\lambda_{l}z_{l}\right) z^{J}-A_{l,J}(z,\overline{z}) Q(z,\overline{z})=C_{l,J}(z,\overline{z}),\quad \tilde{\tr}\left(C_{l,J}(z,\overline{z})\right)=0,\quad\mbox{where $J\not\in \mathcal{T}_{l}$, for all  $l\in 1,\dots, N$.}\label{Africaopt1}
\end{equation} 

Then each very consistent  sum of terms  multiplied by $\lambda_{1}z_{1}^{2},\dots,\lambda_{N}z_{N}^{2}$ in  (\ref{optsec}), generates by  (\ref{Lambda}), (\ref{aux1}) and (\ref{AfricaIJ1}) the following terms\begin{equation} 
 \begin{pmatrix} \mbox{Aux}\Lambda  &  \mbox{O}_{N^{p-1}}& \dots & \mbox{O}_{N^{p-1}}     \\ \mbox{O}_{N^{p-1}}&    \mbox{Aux}\Lambda  & \dots & \mbox{O}_{N^{p-1}}\\ \vdots   & \vdots & \ddots & \vdots   \\ \mbox{O}_{N^{p-1}}   &  \mbox{O}_{N^{p-1}}& \dots &  \mbox{Aux}\Lambda 
\end{pmatrix}\tilde{Z}\left[J \right].\label{Africaaux1}\end{equation} 

 Similarly as previously to (\ref{Africasuc1}), (\ref{Africasuc2}), (\ref{AfricasucL1}), (\ref{AfricasucL11}), (\ref{AfricasucL1sese}) and (\ref{AfricasucL11sese}), we consider  by (\ref{AfricaIJ1}) the following matrices
 \begin{equation} \tilde{\mbox{Aux}}\left[J\right]=  \begin{pmatrix} \mbox{Aux}\left[J\right]  &  \mbox{O}_{N^{p-1}}& \dots & \mbox{O}_{N^{p-1}}     \\ \mbox{O}_{N^{p-1}}&    \mbox{Aux}\left[J\right] & \dots & \mbox{O}_{N^{p-1}}\\ \vdots   & \vdots & \ddots & \vdots   \\ \mbox{O}_{N^{p-1}}   &  \mbox{O}_{N^{p-1}}& \dots &  \mbox{Aux}\left[J\right] 
\end{pmatrix} .\label{Africacalcan1}
\end{equation}

 Defining $$\tilde{Z}^{t}=\left(\left\{a_{J}^{\left(1\right)}\right\}_{J\in\mathbb{N}^{N}\atop \left|I\right|=p-1},\left\{a_{J}^{(2)}\right\}_{J\in\mathbb{N}^{N}\atop \left|I\right|=p-1},\dots,\left\{a_{J}^{(N)}\right\}_{J\in\mathbb{N}^{N}\atop \left|I\right|=p-1}\right),$$
according to the lexicografic order related to (\ref{ODIN}), we obtain the following  system   of  equations 
\begin{equation}\left(I-\mbox{Aux}_{p}\tilde{A}\right) \tilde{Z}+\tilde{B}\overline{\tilde{Z}}=\tilde{V}\left(z_{1},z_{2},\dots,z_{N}\right), \label{Africasisise}
\end{equation} 
where $\tilde{V}\left(z_{1},z_{2},\dots,z_{N}\right)$ is a known homogeneous vector polynomial  of degree $p-2$, dealing with the following products of ,,simple'' matrices 
\begin{equation} 
\tilde{\mbox{Aux}}_{p}=\tilde{\mbox{Aux}}\left[p-2,\dots,0\right]\tilde{\mbox{Aux}}\left[p-3,1,\dots,0\right]\dots\tilde{\mbox{Aux}}\left[0,\dots,1,p-3\right] \tilde{\mbox{Aux}}\left[0,\dots,p-2\right],\quad\mbox{for all $p \geq 3$.}\label{Africacalcan2}
\end{equation}

We were using the following notations
\begin{equation}\tilde{A}=  \begin{pmatrix} A_{1} &  \mbox{O}_{N^{p-1}}& \dots & \mbox{O}_{N^{p-1}}     \\ \mbox{O}_{N^{p-1}}&    A_{2} & \dots & \mbox{O}_{N^{p-1}}\\ \vdots   & \vdots & \ddots & \vdots   \\ \mbox{O}_{N^{p-1}}   &  \mbox{O}_{N^{p-1}}& \dots &  A_{N}
\end{pmatrix} ,\quad \tilde{B}=  \begin{pmatrix} B_{1} &  \mbox{O}_{N^{p-1}}& \dots & \mbox{O}_{N^{p-1}}     \\ \mbox{O}_{N^{p-1}}&    B_{2} & \dots & \mbox{O}_{N^{p-1}}\\ \vdots   & \vdots & \ddots & \vdots   \\ \mbox{O}_{N^{p-1}}   &  \mbox{O}_{N^{p-1}}& \dots &  B_{N}
\end{pmatrix} ,\label{Africacalcan1}
\end{equation}
where   we have
$$A_{1}, A_{2},\dots, A_{N};\quad B_{1},  B_{2}, \dots,  B_{N}\in \mathcal{M}_{N^{p-1}\times N^{p-1} }\left(\mathbb{C}\right).
$$

It is   clear that (\ref{Africasisise}) has unique solution due to the existences of the following matrices
\begin{equation}\frac{1}{I_{N^{p}}-\tilde{\mbox{Aux}}_{p}\tilde{A}-\tilde{B}},\quad\frac{1}{I_{N^{p}}-\tilde{\mbox{Aux}}_{p}\tilde{A}+\tilde{B}},\label{5501}
\end{equation}
recalling the arguments related to (\ref{550se}). Moreover, we have
\begin{equation}\left\|\frac{1}{I_{N^{p}}-\tilde{\mbox{Aux}}_{p}\tilde{A}-\tilde{B}}\right\|<2,\quad\left\|\frac{1}{I_{N^{p}}-\tilde{\mbox{Aux}}_{p}\tilde{A}+\tilde{B}}\right\|<2.\label{550se11}
\end{equation}

These estimations (\ref{550se11}) may be achieved exactly as previously as (\ref{urs11}) and (\ref{urs11se}), being   important in order to establish   bounds for the solutions of the systems of equations defined as in (\ref{Africabeb1}), which give the computations of the $F$-component  of the formal equivalence (\ref{map}) as previously. There are  considered products of matrices as in    (\ref{Africa90000}), (\ref{Africa90000se})  in (\ref{Africasisise}) using  the product of simple matrices (\ref{Africacalcan2}). Then, the solution is computed, being unique in the light of the invertibility of the matrices (\ref{5501}), after there are separated the real part and the imaginary in the local defining equation as previously.  Then the norms of these matrices (\ref{5501}) naturally provide  suitable estimates   in order to conclude the convergence of the formal transformation (\ref{map})  using the rapid iterations arguments of Moser\cite{mo} and the Fischer\cite{sh} Norm:

\subsection{The Fischer Norm\cite{sh}}  We denote by $\mathbb{H}_{k}$   the space of all homogeneous polynomials of degree $k$ in   $z=\left(z_{1},z_{2},\dots,z_{N}\right)$. Then the Fischer inner product\cite{sh} is defined as follows
\begin{equation}\left<z^{\alpha};\hspace{0.1 cm}z^{\beta}\right>_{\mathcal{F}}=\left\{\substack{0,\quad\alpha\neq \beta \\ \alpha!,\quad\alpha=\beta}\right.,\quad \mbox{for $\alpha=\left(\alpha_{1},\dots, \alpha_{N}\right),\hspace{0.1 cm}\beta=\left(\beta_{1},\dots, \beta_{N}\right) \in\mathbb{N}^{N}$.}\label{ppo}\end{equation}

Then the Fischer norm \cite{sh} is defined as follows
\begin{equation}\left\|f_{k}\left(z\right)\right\|:=\left\|f_{k}\left(z\right)\right\|_{\mathcal{F}}:=\displaystyle\sum_{\left|I\right|=k\atop{I\in\mathbb{N}^{N}}}I!\left|c_{I}\right|^{2},\quad\mbox{if  $f_{k}\left(z\right):=\displaystyle\sum_{\left|I\right|=k\atop{I\in\mathbb{N}^{N}}}c_{I}z^{I}$, for $I=\left(i_{1},\dots,i_{N}\right)\in\mathbb{N}^{N}$ and $\left|I\right|=i_{1}+\dots+i_{N}$ and $k\in\mathbb{N}^{\star}$},\label{fnorm}\end{equation}
where $\left|I\right|$ and $I!$ are classically defined in (\ref{fnorm}) recalling notations from Shapiro\cite{sh}. 

We make simple computations using the Fischer inner product from (\ref{ppo}). We obtain  
\bl \label{lem}Let $f\left(z\right)$, $g\left(z\right)\in\mathbb{H}_{k}$ defining the orthogonal decomposition $f\left(z\right)=g\left(z\right)+h\left(z\right)$ with respect to the Fischer inner product. Then the following holds $$\left\|f\left(z\right)\right\|_{\mathcal{F}}=\left\|g\left(Z\right)\right\|_{\mathcal{F}}+\left\|h\left(Z\right)\right\|_{\mathcal{F}}.$$
\el

Actually, we apply the strategy from \cite{bu3}, where the Fischer norm\cite{sh}  gives convenient approximations with respect to the convergence radius of   power series defining   real-analytic submanifolds. These estimations are derived from the local transforming equations  as follows:

 Let $P(z,\overline{z})$ be defined respecting (\ref{new1})  or (\ref{new2}). Then, according to Lemma \ref{lem}, we clearly  have
\begin{equation}\left\|P_{1}(z,\overline{z})\right\|\left\|Q(z,\overline{z})\right\|\leq\left\|P_{1}(z,\overline{z})Q(z,\overline{z})\right\|\leq\left\|P(z,\overline{z})\right\|,\quad \left\|R_{1}(z,\overline{z}) \right\|\leq\left\|P(z,\overline{z})\right\|,\label{esti1}
\end{equation}
where $\left\|\cdot\right\|$ is just the Fischer norm defined in (\ref{fnorm}), being considered the following Fischer Decomposition
\begin{equation}P(z,\overline{z})=P_{1}(z,\overline{z})Q(z,\overline{z})+R_{1}(z,\overline{z}),\quad\mbox{such that $\tr\left(R_{1}(z,\overline{z})\right)=0$.}\label{kong1}
\end{equation}

Analogously, we have 
\begin{equation}\left\|P_{2}(z,\overline{z})\right\|\left\|Q(z,\overline{z})\right\|\leq\left\|P_{2}(z,\overline{z})Q(z,\overline{z})\right\|\leq\left\|P_{1}(z,\overline{z})\right\|,\quad\left\|R_{1}(z,\overline{z}) \right\|\leq\left\|P_{1}(z,\overline{z})\right\|, \label{esti1se}
\end{equation}
being considered the following Fischer Decomposition
\begin{equation}P_{1}(z,\overline{z})=P_{2}(z,\overline{z})Q(z,\overline{z})+R_{2}(z,\overline{z}),\quad\mbox{such that $\tr\left(R_{2}(z,\overline{z})\right)=0$.}\label{kong2}
\end{equation}

These two estimates (\ref{esti1}) and (\ref{esti1se}) are fundamental in order to recall computations from \cite{bu4}, especially  Remark 3.1 from \cite{bu3}, by making a suitable process  of induction. More precisely, we continue this procedure considering similar a decomposition on $P_{2}(z,\overline{z})$.

Let $k\in\mathbb{N}^{\star}$ such that 
$$k-1\in \left\{0,\dots,\left[\frac{p}{2}\right]\right\}\hspace{0.1 cm}\mbox{or}\hspace{0.1 cm}k-1\in \left\{0,\dots,\left[\frac{p-1}{2}\right]\right\}.$$

Assume that $P_{k-1}(z,\overline{z})$ is computed. Then, we have
\begin{equation}\left\|P_{k}(z,\overline{z})\right\|\left\|Q(z,\overline{z})\right\|\leq\left\|P_{k}(z,\overline{z})Q(z,\overline{z})\right\|\leq\left\|P_{k-1}(z,\overline{z})\right\|,\quad\left\|R_{k}(z,\overline{z}) \right\|\leq\left\|P_{k-1}(z,\overline{z})\right\|, \label{esti12}
\end{equation}
being considered the following Fischer Decomposition
\begin{equation}P_{k-1}(z,\overline{z})=P_{k}(z,\overline{z})Q(z,\overline{z})+R_{k}(z,\overline{z}),\quad\mbox{such that $\tr\left(R_{k}(z,\overline{z})\right)=0$.}\label{kong3}
\end{equation}

The above approach inductively describes the computations using the Fischer norm\cite{sh}.
\section{Normalizations}

Summarizing   the   computations of the formal transformation (\ref{map}), we obtain the following result:
\bp\label{proppo} Let $M\subset\mathbb{C}^{N+1}$ be a real-formal
submanifold defined near  $p=0$ as follows
\begin{equation}
w=z_{1}\overline{z}_{1}+\dots+z_{N}\overline{z}_{N}+\lambda_{1}\left(z_{1}^{2}+\overline{z}_{1}^{2}\right)+\dots+\lambda_{N}\left(z_{N}^{2}+\overline{z}_{N}^{2}\right) +\displaystyle\sum
_{k\geq 3}\varphi_{k}(z,\overline{z}), \label{ecuatie11}
\end{equation}
where   $\varphi _{k}(z,\overline{z})$ is a  
polynomial of    degree $k$ in $(z,\overline{z})$,  for all
 $k\geq 3$, such that (\ref{lambida}) is satisfied assuming
\begin{equation}\lambda_{1},\dots,\lambda_{k_{0}}\neq 0,\quad \lambda_{k_{0}+1}=\dots=\lambda_{N}=0,\quad\mbox{for some $k_{0}\in 1,\dots, N$.}\label{03}
\end{equation}  
 
Then there exists a unique formal equivalence defined as in (\ref{map}) and  normalized as follows\begin{equation}  \Re F_{1,n}(z)=0,\quad  F_{0,n+1}^{\left(k_{0}+1\right)}(z)=\dots=F_{0,n+1}^{\left(N\right)}(z)= 0,\quad\mbox{for all $n\in\mathbb{N}^{\star}$,} \label{o}
\end{equation}
which sends $M$ into the following partial normal form
$$
w'={z'}_{1}\overline{{z'}}_{1}+\dots+{z'}_{N}\overline{{z'}}_{N}+\lambda_{1}\left({z'}_{1}^{2}+\overline{{z'}}_{1}^{2}\right)+\dots+\lambda_{N}\left({z'}_{N}^{2}+\overline{{z'}}_{N}^{2}\right) +\displaystyle\sum
_{k\geq 3}\varphi'_{k}\left(z',\overline{z'}\right),  
$$where   $\varphi' _{k}(z,\overline{z})$ is a  
polynomial of    degree $k$ in $(z,\overline{z})$,  for all
 $k\geq 3$, such that there are satisfied the following normalizations
\begin{equation}\left\{\begin{split}&\frac{\varphi'_{k}\left(z',\overline{z'}\right)+\overline{\varphi'_{k}\left(z',\overline{z'}\right)}}{2}\in\mathcal{F}_{k}, \quad\mbox{for all $k\geq 3$},\\& \frac{\varphi'_{k}\left(z',\overline{z'}\right)-\overline{\varphi'_{k}\left(z',\overline{z'}\right)}}{2i}\in\mathcal{G}_{k},\quad\mbox{for all $k\geq 3$},\end{split}\right.\label{cn}
\end{equation}
and respectively, the following normalizations 
\begin{equation}\tilde{\tilde{P}}_{\frac{k}{2}}^{\left(k\right)}(z,\overline{z})\in   \displaystyle\bigcap_{i,j=1}^{N}\ker\left(\frac{\partial^{2}}{\partial z_{i}\partial\overline{z_{j}}+\lambda_{i}\partial z_{i} \partial z_{j}+\lambda_{i}\partial\overline{z_{i}}\partial \overline{z_{j}}}\right)   ,\quad\mbox{for $k$ even,} \label{coreo} 
\end{equation}
and respectively, the following normalizations
\begin{equation}  \tilde{\tilde{P}}_{\frac{k+1}{2}}^{\left(k\right)}(z,\overline{z})\in  \ker\left(\frac{\partial^{2}}{\partial z_{1} \partial \overline{z}_{1}+\dots+\partial z_{N} \partial \overline{z}_{N}}\right)\displaystyle\bigcap   \displaystyle\bigcap_{i=1}^{k_{0}}\ker \left(\frac{\partial^{3}}{\partial z_{i}^{3}}\right)  ,\quad\mbox{for $k$ odd,} \label{coreo1} 
\end{equation}
where there was used the following writing
\begin{equation}
P_{\frac{k}{2}}^{\left(k\right)}(z,\overline{z})=\tilde{P}_{\frac{k}{2}}^{\left(k\right)}(z,\overline{z})+\sqrt{-1}\tilde{\tilde{P}}_{\frac{k}{2}}^{\left(k\right)}(z,\overline{z}),\label{ui119}
\end{equation}
which is just the unique decomposition in a sum of  two polynomials with real coefficients. 
\ep 
\begin{proof} It is considered an induction process depending on $r=m+2n\in\mathbb{N}^{\star}$ in order to compute the equivalence (\ref{map}). We assume thus that  we have determined the following terms
\begin{equation}
 G_{m,n}(z),\quad\mbox{for all $m,n\in\mathbb{N}^{\star}$ with $m+2n<r$.}\label{gigel1}
\end{equation}
and respectively, the following terms
\begin{equation}
F_{m,n}(z),\quad  \quad\mbox{for all $m,n\in\mathbb{N}^{\star}$ with $m+2n<r-1$.}\label{gigel2}
\end{equation}

Moreover, we introduce by (\ref{ef}) the following notations
$$ F_{m,n}(z) =\left(F_{m,n}^{\left(1\right)}(z),\dots, F_{m,n}^{\left(N\right)}(z\right),\quad\mbox{for all $m,n\in\mathbb{N}^{\star}$ with $m+n\geq 3$}.$$
according to the following notation 
\begin{equation}F(z,w)=\left(F^{\left(1\right)}(z,w),\dots,F^{\left(N\right)}(z,w)\right).\label{ef}
\end{equation}

Moreover, we consider the following notations
$$G_{m,n}(z)=\displaystyle\sum_{i_{1}+\dots+i_{N}=m\atop I=\left(i_{1},\dots,i_{N}\right)\in\mathbb{N}^{N}} g_{I,n}z_{1}^{i_{1}}\dots z_{N}^{i_{N}}w^{n},\quad  F_{m,n}^{\left(l\right)}(z)=\displaystyle\sum_{j_{1}+\dots+j_{N}=m-1\atop J=\left(j_{1},\dots,j_{N}\right)\in\mathbb{N}^{N}}f_{J,n}^{\left(l\right)}z_{1}^{j_{1}}\dots z_{N}^{j_{N}}w^{n},\quad\mbox{for all $l=1,\dots,N$.}$$

Now, we want to determine  the following terms
$$g_{I,n},\quad\mbox{for all $I=\left(i_{1},\dots,i_{N}\right)\in\mathbb{N}^{N}$ such that $i_{1}+\dots+i_{N}=m$ and $m+2n=r$.}  
$$

The entire construction is based  on defining (\ref{spartiuG}) and  (\ref{330}). More precisely, we consider successively Fischer Decompositions\cite{sh} in order to compute the formal equivalence (\ref{map}) recalling   (\ref{spartiuG}) and  (\ref{spartiuF}), or equivalently according to the normalizations (\ref{cn}). In particular,  in order to make further computations, we study the following  
\begin{equation}\begin{split}&   
\displaystyle\sum_{m+2n=r}G_{m,n}(z)\left(z_{1}\overline{z}_{1}+\dots+z_{N}\overline{z}_{N}+\lambda_{1}\left(z_{1}^{2}+\overline{z}_{1}^{2}\right)+\dots+\lambda_{N}\left(z_{N}^{2}+\overline{z}_{N}^{2}\right)\right)^{n} \\& \quad\quad\quad\quad\quad   \quad\quad\quad\quad\quad\quad  \quad\quad\quad\quad\quad \quad\quad\quad  \quad\quad\quad\quad\quad\quad\quad\quad  \begin{tabular}{l} \rotatebox[origin=c]{270}{$=$}\end{tabular} \\& V(z,\overline{z})\\& \quad\quad\quad\quad\quad   \quad\quad\quad\quad\quad\quad  \quad\quad\quad\quad\quad \quad\quad\quad  \quad\quad\quad\quad\quad\quad\quad\quad  \begin{tabular}{l} \rotatebox[origin=c]{270}{$+$}\end{tabular} \\& 2\Re\left( \displaystyle\sum_{m+2n=r-1}\displaystyle\sum_{k=1}^{N}\left(z_{k}  \overline{ F_{m,n}^{\left(k\right)}(z)  }  +2\lambda_{k}z_{k} F_{m,n}^{\left(k\right)}(z)   +2\lambda_{k}\overline{z_{k} F_{m,n}^{\left(k\right)}(z) }\right) \right) \left(z_{1}\overline{z}_{1}+\dots+z_{N}\overline{z}_{N}+\lambda_{1}\left(z_{1}^{2}+\overline{z}_{1}^{2}\right)+\dots+\lambda_{N}\left(z_{N}^{2}+\overline{z}_{N}^{2}\right)\right)^{n} , \end{split}\label{ecuatieXYY11}
\end{equation} 
where $V(z,\overline{z})$ represents a sum of terms in $(z,\overline{z})$  depending also on (\ref{gigel1}) and (\ref{gigel2}) and on the induction hypothesis. 

Next, we recall the Fischer Decompositions (\ref{new1}) taking
$$P(z,\overline{z})=\frac{V(z,\overline{z})-\overline{V(z,\overline{z})}}{2\sqrt{-1}}.$$

Then, we have
\begin{equation}
\begin{split}&  \displaystyle\sum_{\left|I\right|=r \atop I\in\mathbb{N}^{N}}\frac{g_{I,0}z^{I}-\overline{ g_{I,0}z^{I}}}{2\sqrt{-1}}=P_{1} (z,\overline{z})Q(z,\overline{z})+R_{1,0}(z,\overline{z}),\quad\mbox{where $\tr\left(R_{1}(z,\overline{z})\right)=0$ and:}\\& \quad\quad\quad\quad\quad\quad\quad\quad \hspace{0.1 cm}   R_{1,0}(z,\overline{z})=\displaystyle\sum_{ {\left|I\right|=r}\atop I\in\mathbb{N}^{N}} \left(a_{I,0}C_{I,0}(z,\overline{z})+\overline{a}_{I,0}\overline{C_{I,0}(z,\overline{z})}\right)+R_{1,0}(z,\overline{z}),\hspace{0.1 cm}\mbox{such that:} \\& \quad\quad\quad\quad\quad \quad\quad\quad R_{1,0}(z,\overline{z})\in \left(\displaystyle\bigcap_{I\in\mathbb{N}^{N},\hspace{0.1 cm} I\not\in \mathcal{S} \atop {\left|I\right|=p}}\left( \ker  C^{\star}_{I}  \cap  \ker  \overline{C}^{\star}_{I} \right)\right) \bigcap\left( \displaystyle\bigcap_{I\in\mathbb{N}^{N},\hspace{0.1 cm} I \in \mathcal{S} \atop {\left|I\right|=p}}\left(\ker  \left(z^{I}\right)^{\star}   \cap  \ker  \left(\overline{z}^{I}\right)^{\star}\right)\right),\end{split}\label{ha1}
\end{equation}
where the corresponding  above occurring polynomials  are known   according to (\ref{new1}) in (\ref{ecuatieXYY11}). 

We obtain 
\begin{equation}G_{r,0}(z)=\displaystyle\sum_{I\in\mathbb{N}^{N},\hspace{0.1 cm} I\not\in \mathcal{S} \atop {\left|I\right|=r}}a_{I,0}z^{I}+\displaystyle\sum_{I\in\mathbb{N}^{N},\hspace{0.1 cm} I \in \mathcal{S} \atop {\left|I\right|=r}}g_{I,0}z^{I} ,\label{lele1}
\end{equation}
where the terms defining the above second sum are directly computed by identifying corresponding coefficients  in (\ref{ecuatieXYY11}).

Then,   (\ref{lele1}) provides contributions to the other terms, because we further deal with the following Fischer generalized Decompositions 
\begin{equation}\begin{split}&  \displaystyle\sum_{\left|I\right|=r-2\atop I\in\mathbb{N}^{N}}\frac{g_{I,1}z^{I}-\overline{ g_{I,1}z^{I}}}{2\sqrt{-1}} =  P_{1}(z,\overline{z})+\displaystyle\sum_{\left|I\right|=r}\frac{A_{I,0}(z)-\overline{A_{I,0}(z)}}{2\sqrt{-1}}\\& \quad \quad \quad\quad\quad\quad\quad\quad \quad\quad\hspace{0.17 cm}=P_{2}(z,\overline{z})Q(z,\overline{z})+R_{2}(z,\overline{z}),\quad\mbox{where $\tr\left(R_{2}(z,\overline{z})\right)=0$ and:}\\& \quad\quad\quad\quad\quad\quad \quad\quad \hspace{0.1 cm}  R_{2}(z,\overline{z})=\displaystyle\sum_{I\in\mathbb{N}^{N}\atop {\left|I\right|=r-2}} \left(a_{I,1}C_{I,1}(z,\overline{z})+\overline{a}_{I,1}\overline{C_{I,1}(z,\overline{z})}\right)+R_{2,0}(z,\overline{z}), \hspace{0.1 cm}\mbox{such that:}\\&\quad\quad\quad\quad   \quad \quad\quad\quad R_{2,0}(z,\overline{z})\in \left(\displaystyle\bigcap_{I\in\mathbb{N}^{N},\hspace{0.1 cm} I\not\in \mathcal{S} \atop {\left|I\right|=p-2}} \ker  C^{\star}_{I}  \cap  \ker  \overline{C}^{\star}_{I} \right) \bigcap \left(\displaystyle\bigcap_{I\in\mathbb{N}^{N},\hspace{0.1 cm} I \in \mathcal{S} \atop {\left|I\right|=p-2}}\ker  \left(z^{I}\right)^{\star}   \cap  \ker  \left(\overline{z}^{I}\right)^{\star}\right),\end{split}\label{ha2}
\end{equation}
where   all  occurring polynomials   are obtained iteratively according to (\ref{new1}). 

We obtain
\begin{equation}G_{r-2,1}(z)=\displaystyle\sum_{I\in\mathbb{N}^{N},\hspace{0.1 cm} I\not\in \mathcal{S} \atop {\left|I\right|=r-2}}a_{I,1}z^{I}+\displaystyle\sum_{I\in\mathbb{N}^{N},\hspace{0.1 cm} I \in \mathcal{S} \atop {\left|I\right|=r-2}}g_{I,1}z^{I} ,\label{lele2}
\end{equation} 
where the terms defining the above second sum are directly computed by identifying corresponding coefficients  in (\ref{ecuatieXYY11}).

These iterative computations determine inductively the $G$-component of the formal transformation using   the uniqueness of the Fischer Decomposition\cite{sh}, because  we can apply  an induction process with respect to $k=1,\dots,\left[\frac{r}{2}\right]$.
We have 
\begin{equation}\begin{split}&\displaystyle\sum_{\left|I\right|=r-2k\atop I\in\mathbb{N}^{N}}\frac{g_{I,k}z^{I}-\overline{ g_{I,k}z^{I}}}{2\sqrt{-1}}=     P_{k}(z,\overline{z})+\displaystyle\sum_{\left|I\right|=r-2k+2\atop I\in\mathbb{N}^{N}}\frac{A_{I,k-1}(z)-\overline{A_{I,k-1}(z)}}{2\sqrt{-1}}\\& \quad \quad \quad\quad\quad\quad\quad\quad \quad\quad\quad=P_{k+1}(z,\overline{z})Q(z,\overline{z})+R_{k+1}(z,\overline{z}),\quad\mbox{where $\tr\left(R_{k+1}(z,\overline{z})\right)=0$ and:}\\&  \quad\quad \quad\quad\quad\quad\quad\quad \quad \hspace{0.1 cm} R_{k+1}(z,\overline{z})=\displaystyle\sum_{I\in\mathbb{N}^{N}\atop {\left|I\right|=r-2k }} \left(a_{I,k}C_{I,k}(z,\overline{z})+\overline{a}_{I,k}\overline{C_{I,k}(z,\overline{z})}\right)+R_{k+1,0}(z,\overline{z}),\hspace{0.1 cm}\mbox{such that:}\\&\quad\quad\quad\quad  \quad\quad \quad\quad\quad R_{k+1,0}(z,\overline{z})\in \left(\displaystyle\bigcap_{I\in\mathbb{N}^{N}, \hspace{0.1 cm}I\not\in \mathcal{S} \atop {\left|I\right|=p-2k}} \left(\ker  C^{\star}_{I}  \cap  \ker  \overline{C}^{\star}_{I}\right)\right)\bigcap\left(   \displaystyle\bigcap_{I\in\mathbb{N}^{N},\hspace{0.1 cm} I \in \mathcal{S} \atop {\left|I\right|=p-2k}}\left(\ker  \left(z^{I}\right)^{\star}   \cap  \ker \left(\overline{z}^{I}\right)^{\star}\right) \right), \end{split}\end{equation}
where   all  occurring polynomials   are obtained iteratively according to (\ref{new1}). We obtain 
\begin{equation}G_{r-2k }(z)=\displaystyle\sum_{I\in\mathbb{N}^{N},\hspace{0.1 cm} I\not\in \mathcal{S} \atop {\left|I\right|=r-2k }}a_{I,k}z^{I}+\displaystyle\sum_{I\in\mathbb{N}^{N},\hspace{0.1 cm} I \in \mathcal{S} \atop {\left|I\right|=r-2k }}g_{I,k}z^{I} ,\label{lele3}
\end{equation} 
where the terms defining the above second sum are directly computed by identifying corresponding coefficients  in (\ref{ecuatieXYY11}).

 Then, according to (\ref{lele1}), (\ref{lele2}) and (\ref{lele3}), we have
\begin{equation}\displaystyle\sum_{m+2n=r}G_{m,n}(z)w^{n}=G_{r,0}(z)+G_{r-2,1}(z)w+\dots+G_{1,\frac{r-1}{2}}(z)w^{\frac{r-1}{2}}+\Im g_{0,\frac{r}{2} }(z)w^{\frac{r}{2}},\quad\mbox{for $r$ even, }\label{gee1}
\end{equation} 
and respectively, we have
\begin{equation}\displaystyle\sum_{m+2n=r}G_{m,n}(z)w^{n}=G_{r,0}(z)+G_{r-2,1}(z)w+\dots+G_{1,\frac{r-1}{2}}(z)w^{\frac{r-1}{2}},\quad\mbox{for $r$ odd}.\label{gee2}
\end{equation}
 
Now, we recall the Fischer Decompositions (\ref{new1}) taking
\begin{equation}     P(z,\overline{z})=\frac{V(z,\overline{z})+\overline{V(z,\overline{z})}}{2} +\displaystyle\sum_{m+2n=r}\displaystyle\sum_{\left|I\right|=m\atop I\in\mathbb{N}^{N}}\frac{ g_{I,n}z^{I}+\overline{ g_{I,n}z^{I}}}{2 }\left(z_{1}\overline{z}_{1}+\dots+z_{N}\overline{z}_{N}+\lambda_{1}\left(z_{1}^{2}+\overline{z}_{1}^{2}\right)+\dots+\lambda_{N}\left(z_{N}^{2}+\overline{z}_{N}^{2}\right)\right)^{n}. \end{equation}

We determine now the following terms
$$ f_{J}^{\left(l\right)},\quad\mbox{for all $J=\left(j_{1},j_{2},\dots,j_{N}\right)\in\mathbb{N}^{N}$  such that $j_{1}+\dots+j_{N}=n$ with $m+2n=r-1$, for all $l=1,\dots,N$}.
$$

Next, we consider the following Fischer Decompositions 
\begin{equation}\begin{split}&   \displaystyle\sum_{l=1}^{N}\displaystyle\sum_{\left|J\right|=r-1\atop J\in\mathbb{N}^{N}} \left(\overline{z}_{l}f_{J}^{\left(l\right)}z^{J}+z_{l}\overline{ f_{J}^{\left(l\right)}z^{I}}\right) = P_{1}(z,\overline{z})Q(z,\overline{z})+R_{1}(z,\overline{z}),\quad\mbox{where $\tr\left(R_{1}(z,\overline{z})\right)=0$ and:}\\& \quad\quad\quad\quad\quad\quad\quad  R_{1}(z,\overline{z})=\displaystyle\sum_{l=1}^{N}\displaystyle\sum_{J\in\mathbb{N}^{N}\atop {\left|J\right|=r-1}} \left(a_{l,J,0}C_{l,J,0}(z,\overline{z})+\overline{a}_{l,J,0}\overline{C_{l,J,0}(z,\overline{z})}\right)+R_{1,0}(z,\overline{z}),\quad\mbox{such that:}\\&   \quad\quad\quad \quad\quad\quad\quad  R_{1,0}(z,\overline{z})\in \displaystyle\bigcap_{l=1}^{N}\left(\left(\displaystyle\bigcap_{J\in\mathbb{N}^{N},\hspace{0.1 cm}J\not\in \mathcal{T}_{l}\atop {\left|J\right|=p-1}} \left(\ker  C^{\star}_{l,J}  \bigcap  \ker  \overline{C}^{\star}_{l,J} \right)\right)\bigcap   \left(\displaystyle\bigcap_{J\in\mathbb{N}^{N},\hspace{0.1 cm}J\in \mathcal{T}_{l}\atop {\left|J\right|=p-1}} \left(\ker  \left( z_{l}\overline{z}^{J} \right)^{\star} \bigcap  \ker  \left( \overline{z}_{l}z^{J} \right)^{\star} \right)\right)\right),\end{split}
\end{equation}
where all  occurring polynomials   are obtained iteratively according to (\ref{new2}).  

We obtain \begin{equation}F_{r-1,0}^{\left(l\right)}(z)=\displaystyle\sum_{J\in\mathbb{N}^{N},\hspace{0.1 cm} J\not\in \mathcal{T}_{l} \atop {\left|J\right|=r-1}}a_{l,J,0}z^{J}+\displaystyle\sum_{J\in\mathbb{N}^{N},\hspace{0.1 cm} J \in \mathcal{T}_{l} \atop {\left|J\right|=r-1}}f_{J,0}^{\left(l\right)}z^{J} ,\quad\mbox{for all $l=1,\dots,N$,}\label{lele1se}
\end{equation}  
where the terms defining the above second sum are directly computed by identifying corresponding coefficients in (\ref{ecuatieXYY11}). 

Then,  (\ref{lele1se}) provides contributions to the other terms, because we further deal with the following Fischer generalized Decompositions 
 \begin{equation}\begin{split}& \displaystyle\sum_{l=1}^{N}\displaystyle\sum_{\left|J\right|=r-3\atop J\in\mathbb{N}^{N}} \left(\overline{z}_{l}f_{J}^{\left(l\right)}z^{J}+z_{l}\overline{ f_{J}^{\left(l\right)}z^{J}}\right)=P_{1}(z,\overline{z})+\displaystyle\sum_{l=1}^{N}\displaystyle\sum_{\left|J\right|=r-1\atop J\in\mathbb{N}^{N}}\left(A_{J,l}(z)+\overline{A_{J,l}(z)}\right)\\& \quad\quad\quad\quad\quad\quad\quad\quad\quad\quad\quad\quad\quad\quad\quad   =P_{2}(z,\overline{z})Q(z,\overline{z})+R_{2}(z,\overline{z}),\quad\mbox{where $\tr\left(R_{2}(z,\overline{z})\right)=0$ and:}\\& \quad\quad \quad \quad \quad \quad     R_{2}(z,\overline{z})=\displaystyle\sum_{l=1}^{N}\displaystyle\sum_{J\in\mathbb{N}^{N}\atop {\left|J\right|=r-3}} \left(a_{l,J,1}C_{l,J,1}(z,\overline{z})+\overline{a}_{l,J,1}\overline{C_{l,J,1}(z,\overline{z})}\right)+R_{2,0}(z,\overline{z}),\quad\mbox{such that:} \\&    \quad\quad \quad\quad\quad\quad R_{2,0}(z,\overline{z})\in \displaystyle\bigcap_{l=1}^{N}\left(\left(\displaystyle\bigcap_{J\in\mathbb{N}^{N},\hspace{0.1 cm}J\not\in \mathcal{T}_{l}\atop {\left|J\right|=p-3}}\left( \ker  C^{\star}_{l,J}  \bigcap   \ker  \overline{C}^{\star}_{l,J}\right)\right)\bigcap\left(   \displaystyle\bigcap_{J\in\mathbb{N}^{N},\hspace{0.1 cm}J\in \mathcal{T}_{l}\atop {\left|J\right|=p-3}} \ker  \left( z_{l}\overline{z}^{J} \right)^{\star} \bigcap   \ker  \left( \overline{z}_{l}z^{J} \right)^{\star} \right)\right),   \end{split} 
 \end{equation}
where all  occurring polynomials   are obtained iteratively using  generalized versions of the Fischer Decomposition\cite{sh}.

 We obtain
\begin{equation}F_{r-3,0}^{\left(l\right)}(z)=\displaystyle\sum_{J\in\mathbb{N}^{N},\hspace{0.1 cm} J\not\in \mathcal{T}_{l} \atop {\left|J\right|=r-3}}a_{l,J,1}z^{J}+\displaystyle\sum_{J\in\mathbb{N}^{N},\hspace{0.1 cm} J \in \mathcal{T}_{l} \atop {\left|J\right|=r-3}}f_{J,1}^{\left(l\right)}z^{J} ,\quad\mbox{for all $l=1,\dots,N$,}\label{lele2se}
\end{equation}    
where the terms defining the above second sum are directly computed by identifying corresponding coefficients in (\ref{ecuatieXYY11}). 

These iterative computations determine inductively the $F$-component of the formal transformation using   the uniqueness of the Fischer Decomposition\cite{sh}, because  we can apply  an induction process with respect to $k=1,\dots,\left[\frac{r-1}{2}\right]$. We have 
\begin{equation}\begin{split}& \displaystyle\sum_{l=1}^{N}\displaystyle\sum_{\left|J\right|=r-1-2k\atop J\in\mathbb{N}^{N}} \left(\overline{z}_{l}f_{J}^{\left(l\right)}z^{J}+z_{l}\overline{ f_{J}^{\left(l\right)}z^{J}}\right)=P_{k}(z,\overline{z})+ \displaystyle\sum_{l=1}^{N}\displaystyle\sum_{\left|J\right|=r+1-2k\atop J\in\mathbb{N}^{N}}\left(A_{J,l}(z)+\overline{A_{J,l}(z)}\right)\\&\quad\quad   \quad\quad\quad \quad \quad \quad \quad\quad\quad\quad\quad\quad\quad\quad\quad =P_{k+1}(z,\overline{z})Q(z,\overline{z})+R_{k+1}(z,\overline{z}),\quad\mbox{where $\tr\left(R_{k+1}(z,\overline{z})\right)=0$ and:}\\&   \quad\quad\quad   R_{k+1}(z,\overline{z})=\displaystyle\sum_{l=1}^{N}\displaystyle\sum_{J\in\mathbb{N}^{N}\atop {\left|J\right|=r-2k-1}} \left(a_{l,J,k}C_{l,J,k}(z,\overline{z})+\overline{a}_{l,J,k}\overline{C_{l,J,k}(z,\overline{z})}\right)+R_{k+1,0}(z,\overline{z}),\quad\mbox{such that:}\\&     \quad\quad\quad R_{k+1,0}(z,\overline{z})\in \displaystyle\bigcap_{l=1}^{N}\left(\left(\displaystyle\bigcap_{J\in\mathbb{N}^{N},\hspace{0.1 cm}J\not\in \mathcal{T}_{l}\atop {\left|J\right|=p-2k-1}} \left(\ker  C^{\star}_{l,J}  \bigcap   \ker  \overline{C}^{\star}_{l,J} \right)\right)\bigcap\left( \displaystyle\bigcap_{J\in\mathbb{N}^{N},\hspace{0.1 cm}J\in \mathcal{T}_{l}\atop {\left|J\right|=p-2k-1}} \left(\ker  \left( z_{l}\overline{z}^{J} \right)^{\star} \bigcap   \ker  \left( \overline{z}_{l}z^{J} \right)^{\star}\right)\right)\right),  \end{split} 
 \end{equation}
where all  occurring  polynomials   are obtained iteratively using  generalized versions of the Fischer Decomposition\cite{sh}. 

We obtain \begin{equation}F_{r-1-2k,0}^{\left(l\right)}(z)=\displaystyle\sum_{J\in\mathbb{N}^{N},\hspace{0.1 cm} J\not\in \mathcal{T}_{l} \atop {\left|J\right|=r-1-2k}}a_{l,J,k}z^{J}+\displaystyle\sum_{J\in\mathbb{N}^{N},\hspace{0.1 cm} J \in \mathcal{T}_{l} \atop {\left|J\right|=r-1-2k}}f_{J,k}^{\left(l\right)}z^{J} ,\quad\mbox{for all $l=1,\dots,N$,}\label{lele3se}
\end{equation}    
where the terms defining the above second sum are directly computed by identifying corresponding coefficients in (\ref{ecuatieXYY11}). 

Assume  $r$   even. Now, we write   as follows
$$F_{1,n}(z):=F_{1,n}\left(z_{1},z_{2},\dots,z_{N}\right)=\begin{pmatrix}
v_{11} & v_{12} & \dots & v_{1N} \\ v_{21} & v_{22} & \dots & v_{2N} \\  \vdots  & \vdots  & \ddots & \vdots  \\ v_{N1} & v_{N2} & \dots & v_{NN}  
\end{pmatrix}\begin{pmatrix}
z_{1}\\ z_{2}\\ \vdots \\ z_{N}\end{pmatrix}.
$$ 
 
 Then, it remains to study the following sum of terms
$$2\Re\left\{\displaystyle\sum_{i,j=1}^{N}z_{i}\overline{v_{ij}z_{j}}+\displaystyle\sum_{j=1}^{N}\displaystyle\sum_{i=1}^{k_{0}}2\lambda_{i}z_{i} v_{ij}z_{j}  \right\}. 
$$

Then, $\Im V$ is determined by the first normalization condition from (\ref{coreo}). We move forward:

 Assume  $r$ is odd. Then, we write as follows
\begin{equation}F_{2,r-1}(z)=\left(F_{2,r-1}^{\left(1\right)}(z),\dots,F_{2,r-1}^{\left(N\right)}(z)\right),\quad F_{0,r}(z)=a:=\left(a_{1},\dots,a_{N}\right).\label{788a}
\end{equation}  

Now, it remains to understand the interactions of these two polynomials  in 
the local defining equation, according to the corresponding Fischer Decompositions, using  $\left<\cdot,\cdot\right>$, which is just the standard hermitian-inner product. We have
\begin{equation}\begin{split}&  \left<z,F_{2,r-1}(z) \right>+\left<a,z\right>\left<z,z\right>+\left<F_{2,r-1}(z),z\right>+\left<z,a\right>\left<z,z\right> +\left(\left<z,a\right>+\left<a,z\right>\right)\left(\displaystyle\sum_{i=1}^{k_{0}}\lambda_{i}z_{i}^{2}+\displaystyle\sum_{i=1}^{k_{0}}\lambda_{i}\overline{z}_{i}^{2}\right)\\& \\& \quad\quad\quad\quad\quad   \quad\quad\quad\quad\quad\quad  \quad\quad\quad\quad\quad \quad\quad\quad \quad\quad\quad \quad\quad\quad\quad\quad\quad\quad\quad  \begin{tabular}{l} \rotatebox[origin=c]{270}{$+$}\end{tabular} \\&\\& \displaystyle\sum_{i=1}^{k_{0}}2\lambda_{i}z_{i}F_{2,r-1}^{(i)}(z)+\displaystyle\sum_{i=1}^{k_{0}}2\lambda_{i}\overline{z_{i}F_{2,r-1}^{(i)}(z)}   +\displaystyle\sum_{i=1}^{k_{0}}2\lambda_{i}z_{i}a_{i}\left(\left<z,z\right>+\displaystyle\sum_{i=1}^{k_{0}}\lambda_{i}z_{i}^{2}+\displaystyle\sum_{i=1}^{k_{0}}\lambda_{i}\overline{z}_{i}^{2}\right)  \\& \quad\quad\quad\quad\quad   \quad\quad\quad\quad\quad\quad  \quad\quad\quad\quad\quad \quad\quad\quad  \quad\quad\quad\quad\quad\quad\quad\quad  \begin{tabular}{l} \rotatebox[origin=c]{270}{$+$}\end{tabular} \\&\displaystyle\sum_{i=1}^{k_{0}}2\lambda_{i}\overline{z}_{i}\overline{a}_{i}\left(\left<z,z\right>+\displaystyle\sum_{i=1}^{k_{0}}\lambda_{i}z_{i}^{2}+\displaystyle\sum_{i=1}^{k_{0}}\lambda_{i}\overline{z}_{i}^{2}\right)=K(z,\overline{z}),\end{split}
\end{equation} 
where $K(z,\overline{z})$ is a  determined polynomial of degree $2$ provided by corresponding Fischer Decompositions the from the local defining equation.

Writing $K(z,\overline{z})$ as a sum of homogeneous terms as follows
$$ K(z,\overline{z})=K_{3,0}(z,\overline{z})+K_{2,1}(z,\overline{z})+K_{1,2}(z,\overline{z})+K_{0,3}(z,\overline{z}),
$$
using appropriate notations, we obtain  
\begin{equation}\begin{split}& \quad\quad\quad\quad\quad\quad\quad\quad\quad\quad\quad\quad  \hspace{0.051 cm}      \displaystyle\sum_{i=1}^{k_{0}}2\lambda_{i}z_{i}F_{2,r-1}^{(i)}(z)+\left<z,a\right>\displaystyle\sum_{i=1}^{k_{0}}\lambda_{i}z_{i}^{2}+\displaystyle\sum_{i=1}^{k_{0}}2\lambda_{i}z_{i}a_{i}\displaystyle\sum_{i=1}^{k_{0}}\lambda_{i}z_{i}^{2}=K_{3,0}(z,\overline{z}),\\&  \left<F_{2,r-1}(z),z\right>+\left<z,a\right>\left<z,z\right>+\left<a,z\right>\displaystyle\sum_{i=1}^{k_{0}}\lambda_{i}z_{i}^{2}+\displaystyle\sum_{i=1}^{k_{0}}2\lambda_{i}z_{i}a_{i}  \left<z,z\right>+\displaystyle\sum_{i=1}^{k_{0}}2\lambda_{i}\overline{z}_{i}\overline{a}_{i}\displaystyle\sum_{i=1}^{k_{0}}\lambda_{i}z_{i}^{2}=K_{2,1}(z,\overline{z}),\\& \left<z,F_{2,r-1}(z)\right>+\left<a,z\right>\left<z,z\right>+\left<z,a\right>\displaystyle\sum_{i=1}^{k_{0}}\lambda_{i}\overline{z}_{i}^{2}+\displaystyle\sum_{i=1}^{k_{0}}2\lambda_{i}\overline{z}_{i}\overline{a}_{i}\left<z,z\right>+\displaystyle\sum_{i=1}^{k_{0}}2\lambda_{i}z_{i}a_{i}\displaystyle\sum_{i=1}^{k_{0}}\lambda_{i}\overline{z}_{i}^{2}=K_{1,2}(z,\overline{z}),\\& \quad\quad\quad\quad\quad\quad\quad\quad\quad\quad\quad\quad  \hspace{0.051 cm}\displaystyle\sum_{i=1}^{k_{0}}2\lambda_{i}\overline{z_{i}F_{2,r-1}^{(i)}(z)}+\left<a,z\right>\displaystyle\sum_{i=1}^{k_{0}}\lambda_{i}\overline{z}_{i}^{2}+\displaystyle\sum_{i=1}^{k_{0}}2\lambda_{i}\overline{z}_{i}\overline{a}_{i}\displaystyle\sum_{i=1}^{k_{0}}\lambda_{i}\overline{z}_{i}^{2}=K_{0,3}(z,\overline{z}).\end{split}\label{011}
\end{equation} 

It is obvious the equivalence between the first and the last equation from (\ref{011}), being conjugated, and respectively the equivalence between the first and the last equation from (\ref{011}), being conjugated as well. We obtain
$$ F_{2,r-1}(z)+\left<z,a\right>z+a\displaystyle\sum_{i=1}^{k_{0}}\lambda_{i}z_{i}^{2}+\displaystyle\sum_{i=1}^{k_{0}}2\lambda_{i}z_{i}a_{i}  z +  2\left(\lambda_{1}a_{1}\displaystyle\sum_{i=1}^{k_{0}}\lambda_{i}z_{i}^{2},\dots,\lambda_{k_{0}}a_{k_{0}}\displaystyle\sum_{i=1}^{k_{0}}\lambda_{i}z_{i}^{2},0,\dots,0\right)=\left(\frac{\partial }{\partial \overline{z}_{1}},\dots,\frac{\partial }{\partial \overline{z}_{N}}\right)\left(K_{2,1}(z,\overline{z})\right),    
$$
and equivalently, we obtain
$$
 F_{2,r-1}^{\left(j\right)}(z)+\left<z,a\right>z_{j}+a_{j}\displaystyle\sum_{i=1}^{k_{0}}\lambda_{i}z_{i}^{2}+\displaystyle\sum_{i=1}^{k_{0}}2\lambda_{i}z_{i}a_{i}  z_{j}+2\lambda_{j}a_{j}\displaystyle\sum_{i=1}^{k_{0}}\lambda_{i}z_{i}^{2}= \frac{\partial }{\partial \overline{z}_{j}}\left(K_{2,1}(z,\overline{z})\right),\quad\mbox{for all $j=1,\dots,k_{0}$}.
$$ 
 
Returning in the first equation in (\ref{011}), we have
$$ \displaystyle\sum_{j=1}^{k_{0}}2\lambda_{j}z_{j}F_{2,r-1}^{\left(j\right)}(z)+\left<z,a\right>\displaystyle\sum_{j=1}^{k_{0}}\lambda_{j}z_{j}^{2}+\displaystyle\sum_{i,j=1}^{k_{0}}2\lambda_{i}z_{i}a_{i} \lambda_{j}z_{j}^{2}=K_{3,0}(z,\overline{z}), $$
or equivalently, we have 
\begin{equation*} \begin{split}&\displaystyle\sum_{j=1}^{k_{0}}2\lambda_{j}z_{j}\left(\frac{\partial }{\partial \overline{z}_{j}}\left(K_{2,1}(z,\overline{z})\right)-\left<z,a\right>z_{j}-a_{j}\displaystyle\sum_{i=1}^{k_{0}}\lambda_{i}z_{i}^{2}-\displaystyle\sum_{i=1}^{k_{0}}2\lambda_{i}z_{i}a_{i}  z_{j}-2\lambda_{j}a_{j}\displaystyle\sum_{i=1}^{k_{0}}\lambda_{i}z_{i}^{2}\right)\\& \quad\quad\quad\quad\quad   \quad\quad\quad\quad\quad\quad  \quad\quad\quad\quad\quad \quad\quad\quad  \quad\quad\quad\quad\quad\quad\quad\quad  \begin{tabular}{l} \rotatebox[origin=c]{270}{$+$}\end{tabular} \\&\left<z,a\right>\displaystyle\sum_{j=1}^{k_{0}}\lambda_{j}z_{j}^{2}+\displaystyle\sum_{i,j=1}^{k_{0}}2\lambda_{i}z_{i}a_{i} \lambda_{j}z_{j}^{2}=K_{3,0}(z,\overline{z}). \end{split}
\end{equation*}  
 
It follows that
\begin{equation} \begin{split}&    0=\displaystyle\sum_{j=1}^{k_{0}}2\lambda_{j}z_{j}\frac{\partial }{\partial \overline{z}_{j}}\left(K_{2,1}(z,\overline{z})\right)-K_{3,0}(z,\overline{z})-\displaystyle\sum_{j=1}^{k_{0}}2\lambda_{j}z_{j} \left( \left<z,a\right>z_{j}+a_{j}\displaystyle\sum_{i=1}^{k_{0}}\lambda_{i}z_{i}^{2}+\displaystyle\sum_{i=1}^{k_{0}}2\lambda_{i}z_{i}a_{i}  z_{j}+2\lambda_{j}a_{j}\displaystyle\sum_{i=1}^{k_{0}}\lambda_{i}z_{i}^{2}\right)\\& \quad\quad\quad\quad\quad   \quad\quad\quad\quad\quad\quad  \quad\quad\quad\quad\quad \quad\quad\quad  \quad\quad\quad\quad\quad\quad\quad\quad  \begin{tabular}{l} \rotatebox[origin=c]{270}{$+$}\end{tabular} \\&\left<z,a\right>\displaystyle\sum_{j=1}^{k_{0}}\lambda_{j}z_{j}^{2}+\displaystyle\sum_{i,j=1}^{k_{0}}2\lambda_{i}z_{i}a_{i} \lambda_{j}z_{j}^{2}, \end{split}
\end{equation}
which, after few simplifications, confirms the following
\begin{equation} \begin{split}&    0=\displaystyle\sum_{j=1}^{k_{0}}2\lambda_{j}z_{j}\frac{\partial }{\partial \overline{z}_{j}}\left(K_{2,1}(z,\overline{z})\right)-K_{3,0}(z,\overline{z})-\left<z,a\right>\displaystyle\sum_{j=1}^{k_{0}}\lambda_{j}z_{j}^{2}  -\displaystyle\sum_{j=1}^{k_{0}}2\lambda_{j}z_{j} \left(  \displaystyle\sum_{i=1}^{k_{0}}2\lambda_{i}z_{i}a_{i}  z_{j}+2\lambda_{j}a_{j}\displaystyle\sum_{i=1}^{k_{0}}\lambda_{i}z_{i}^{2}\right)   . \end{split}
\end{equation}
 
Now, we focus only on the coefficients of the above terms $z_{1}^{3},\dots, z_{k_{0}}^{3}$. Then, we compute easily $a_{1},\dots, a_{k_{0}}$ according to the assumption (\ref{03}), which allows to compute the components of $2\lambda_{1}a_{1}+\overline{a_{1}},\dots, 2\lambda_{k_{0}}a_{k_{0}}+\overline{a_{k_{0}}}$. Then, we compute $F_{2,r-1}(z)$, for all $i=1,\dots, k_{0}$.

Now, $\Im V$ is determined by the first normalization condition from (\ref{coreo1}). Then, according to (\ref{lele1se}), (\ref{lele2se}) and (\ref{lele3se}), we have
\begin{equation}\left\{\begin{split}&\displaystyle\sum_{m+2n=r-1}F_{m,n}^{\left(l\right)}(z)w^{n}=F_{r-1,0}^{\left(l\right)}(z)+F_{r-3,1}^{\left(l\right)}(z)w+\dots  ,\quad\mbox{for all $l=1,\dots,k_{0}$,}\\& \displaystyle\sum_{m+2n=r-1}F_{m,n}^{\left(l\right)}(z)w^{n}=F_{r-1,0}^{\left(l\right)}(z)+F_{r-3,1}^{\left(l\right)}(z)w+\dots  ,\quad\mbox{for all $l=k_{0}+1,\dots,N$,} \end{split}\right.\label{londra1} \end{equation}  
when $r$ is even, and respectively
\begin{equation}\left\{\begin{split}&\displaystyle\sum_{m+2n=r-1}F_{m,n}^{\left(l\right)}(z)w^{n}=F_{r-1,0}^{\left(l\right)}(z)+F_{r-3,1}^{\left(l\right)}(z)w+\dots  ,\quad\mbox{for all $l=1,\dots,k_{0}$,}\\& \displaystyle\sum_{m+2n=r-1}F_{m,n}^{\left(l\right)}(z)w^{n}=F_{r-1,0}^{\left(l\right)}(z)+F_{r-3,1}^{\left(l\right)}(z)w+\dots  ,\quad\mbox{for all $l=k_{0}+1,\dots,N$,} \end{split}\right. \label{londra2}\end{equation} 
when $r$ is odd. 
\end{proof}

This is just a partial normal form, because (\ref{o}) describes an infinite number of parameters acting chaotically on the local defining equation. Their action may be further detected using careful computations in the local defining equations using some non-degeneracy conditions  as we shall see.  Recalling the approach from \cite{bu4}, we immediately obtain the following:
\bp There exists a unique $T_{0}\in \mbox{Aut}_{0}\left(M_{\infty}\right)$ such that $T_{0}\circ F$ satisfies (\ref{o}).\ep 

This   is actually based on the following fact:
\bp Assume  $\lambda_{1},\dots,\lambda_{N} $ satisfy (\ref{lambida}) such that $\lambda_{1},\dots,\lambda_{N}\neq 0 $. Then  \begin{equation}\mbox{Aut}_{0}\left(M_{\infty}\right)=\left\{(z,w)\rightarrow\left(a(w)U(z),w\right),\quad\mbox{where $a(w)=\overline{a}(w)$ is a formal mapping and $U(z)$ leaves invariant the model $M_{\infty}$}\right\}.
\end{equation}
\ep 
\begin{equation}
\end{equation}
\begin{equation}
\end{equation}
\begin{equation}
\end{equation} 
 The proofs may be easy to be concluded as exercises. Moreover, we can easily obtain similarly the formula of the group of formal automorphisms of $M_{\infty}$, when some of the Bishop invariants are vanishing according to $(3.1)$ from Huang-Yin\cite{huyi1}.
 
\section{Application of Moser's Methods\cite{mo} and Proof of Theorem \ref{t1}} 
We apply the  classical rapid iteration procedure of Moser\cite{mo} following   computations from Huang-Yin\cite{huyi1}, Coffman\cite{co3},\cite{co4} and Gong\cite{go1}.  In particular, the proof  of  Theorem  \ref{t1}  has as model the proof of the Generalization\cite{huyi1} of Huang-Yin\cite{huyi1} of the Theorem of Moser \cite{mo}. We proceed by  recalling the lines from  \cite{bu5} as follows in   order to apply   the rapid iteration procedure of Moser\cite{mo}.

\subsection{Settings from  Huang-Yin\cite{huyi1} and Moser\cite{mo}}
For $R:=\left(r ,\dots,r \right)$, in order to    use the methods based on  Moser's rapid convergence   arguments,  we introduce  the following domains
\begin{equation}\begin{split}& D_{r}=\left\{\left(z,\xi\right)\in\mathbb{C}^{N}\times\mathbb{C}^{N},\hspace{0.1 cm}\left|z_{i}\right|<r ,\hspace{0.1 cm}\left|\xi_{i}\right|<r,\hspace{0.1 cm}\mbox{for all $i=1,\dots,N$ }\right\} ,\\&        \Delta_{r}=\left\{\left(z,w\right)\in\mathbb{C}^{N+1};\hspace{0.1 cm}\left|z_{i}\right|<r ,\hspace{0.1 cm}  \left|w\right|<r^{2}\left(N+2\lambda_{1}+\dots+2\lambda_{N}\right) ,\hspace{0.1 cm}\mbox{for all $i=1,\dots,N$ } \right\},\end{split} \label{VX}\end{equation}
where we have been using the following notation
$$\xi=\left(\xi_{1},\xi_{2},\dots,\xi_{N}\right).$$

Throughout  the rest of this paper, we   use the following notations
\begin{equation}\left\|E\right\|_{r}:=\displaystyle\sup_{\left(z,\xi\right)\in D_{r}}\left|E\left(z,\xi\right)\right|,\quad \left|h\right|_{r}=\displaystyle\sup_{\left(z,w\right)\in \Delta_{r}}\left|h\left(z,w\right)\right|,\label{XYZ1}\end{equation}
where $E\left(z,\xi\right)$ is a holomorphic function defined over $\overline{D_{r}}$, and where respectively $h\left(z,w\right)$ is a holomorphic function defined over $\overline{\Delta_{r}}$.

Following Moser\cite{mo}, we  define also the following real numbers
\begin{equation}\frac{1}{2}<r'<\sigma<\rho<r\leq 1,\quad \rho =\frac{2 r'+r}{3},\quad \sigma =\frac{2r'+\rho}{3},\quad n\in\mathbb{N}^{\star}.\label{num}\end{equation}

We also recall here Lemma $4.5$ of Huang-Yin\cite{huyi1} that will be applied later:
\bl Suppose that there exists $C>0$ and a number $a>1$ such that $d_{n}\geq Ca^{n}$. Then  we have 
$$\displaystyle\sum_{n\rightarrow \infty}n^{m_{3}}d_{n}^{m_{1}}\left(1-\frac{1}{n^{m_{2}}}\right)^{d_{n}}=0,$$
for any integers $m_{1}$, $m_{2}$, $m_{3}>0$.
\el

Now, we are ready to more forward.  
\subsection{Degree Estimations from \cite{bu3}}We consider the real  submanifold $M\subset\mathbb{C}^{N+1}$ defined near $p=0$ as follows
\begin{equation}w=\Phi\left(z,\overline{z}\right)=z_{1}z_{1}+\dots+z_{N}z_{N}+\lambda_{1}\left(z_{1}^{2}+\overline{z}_{1}^{2}\right)+\dots+\lambda_{N}\left(z_{N}^{2}+\overline{z}_{N}^{2}\right)+E\left(z,\overline{z}\right),\label{r1}\end{equation}
where    $E\left(z,\overline{z}\right)=\rm{O}(3)$ is real-analytic near $z=0$ according to the hypothesis of the statement of the Theorem \ref{t1}.

We know   the   formal  equivalence, which satisfies the normalization conditions (\ref{o}), sends $M$, defined near $p=0$ by (\ref{var1}),  into the model manifold defined by (\ref{mmodel1}) up to the degree $d\geq 3$. Recallling standard notations used in \cite{bu1}, we find by (\ref{cn}) the following pair of polynomials
\begin{equation}\left(F_{\rm{nor}}^{d-1}\left(z,w\right),\hspace{0.1 cm} G_{\rm{nor}}^{\left(d\right)}\left(z,w\right)\right).\label{sol}\end{equation}

 Following the strategy of Huang-Yin \cite{huyi2},  we define the following  equivalence
\begin{equation}    \Theta\left(z,w\right):=\left(z+\widehat{F}\left(z,w\right),\hspace{0.1 cm}w+\widehat{G}\left(z,w\right)\right)=\left(z+F_{\rm{nor}}^{(d-1)}\left(z,w\right)+O_{\rm{wt}}\left(d\right),\hspace{0.1 cm}w+G_{\rm{nor}}^{\left(d\right)}\left(z,w\right)+\rm{O}_{\rm{wt}}\left(d+1\right)\right), \label{r2}\end{equation}
which sends $M$ up to the degree $d$ into the model manifold $M_{\infty}$ defined by (\ref{mmodel1}).

It remains to understand how the degree of the remaining terms changes   using the  equivalence (\ref{r2}) recalling Lemma $4.1$ from \cite{bu3} and in particular $(4.8)$, $(4,9)$, $(4.10)$ and $(4.11)$ from \cite{bu3} introducing
\begin{equation}M'=\Theta\left(M\right).
\end{equation}

Then, we obtain
\begin{equation}M': w'={z'}_{1}\overline{{z'}}_{1}+\dots+{z'}_{N}\overline{{z'}}_{N}+\lambda_{1}\left({z'}_{1}^{2}+\overline{{z'}}_{1}^{2}\right)+\dots+\lambda_{N}\left({z'}_{N}^{2}+\overline{{z'}}_{N}^{2}\right)+E'\left({z'},\overline{{z'}}\right),\label{900}\end{equation}
where we have $$\rm{Ord}\left(E'\left(z,\overline{z}\right)\right)\geq 2d-2.$$ 

Now, we are ready to make firstly suitable estimations for the equivalence (\ref{map}).     We have:
 
\subsection{Estimations for the Formal  Equivalence}  It suffices to assume $k_{0}=N$ in (\ref{03}), because this situation provides generous estimations covering the other situation. We assume also that the real-analytic  submanifold $M$ defined in (\ref{r1}) is formally equivalent to $M_{\infty}$ defined in (\ref{model}) with $E\left(z,\xi\right)$  holomorphic over $\overline{D_{r}}$ and $\rm{Ord}\left(E\left(Z,\xi\right)\right)\geq d$. Then   the following estimates hold
\begin{equation}\begin{split}&\quad\quad\quad\left\|E\left(Z,\xi\right)-J^{2d-3}\left(E\left(Z,\xi\right)\right)\right\|_{\rho}
\leq\frac{\left(2d\right)^{4N}\left\|E\right\|_{r}}{\left(r-\rho\right)^{2N}}\left(\frac{\rho}{r}\right)^{2d-2},\quad\left|\widehat{G}\left(z,w\right)\right|_{\rho}\leq  d^{2}\left(2d\right)^{4N}\left\|E\right\|_{r} \left(\frac{\rho}{r}\right)^{2d-3},
\\&\left|\nabla\widehat{G}\left(z,w\right)\right|_{\rho}\leq  \left(\frac{36}{r-\rho}+2N\right)  \frac{d^{2}\left(2d\right)^{4N}\left\|E\right\|_{r}}{N\left(r-\rho\right)}\left(\frac{\rho}{r}\right)^{\frac{2d-3}{2}},\hspace{0.1 cm}\left|\widehat{F}_{l}
\left(z,w\right)\right|_{\rho}\leq   \left(d^{2}\left(2d\right)^{4N}+d^{2}\left(2d\right)^{6N}\right)\left\|E\right\|_{r} \left(\frac{\rho}{r}\right)^{2d-2},\\&\quad\quad\quad\quad\quad\quad\quad\left|\nabla\widehat{F}_{l} \left(z,w\right)\right|_{\rho}\leq \left(\frac{36 \left(1+ d^{2}\left(2d\right)^{2N}\right)}{r-\rho}+6N \left(1+ d^{2}\left(2d\right)^{2N}\right)\right)\frac{d^{2}\left(2d\right)^{4N}\left\|E\right\|_{r}}{ \left(r-\rho\right) }\left(\frac{\rho}{r}\right)^{d-1} ,\end{split}\label{1003}\end{equation}
for all $l\in\{1,\dots,N\}$, where $J^{2d-3}\left(E\left(z,\xi\right)\right)$ is the polynomial defined by the Taylor expansion of $E\left(Z,\xi\right)$ up to the degree $2d-3$ and  $\nabla$ represents the gradient.
 
Indeed, recalling $(4.13)$ from \cite{bu3}, the first   inequality becomes clear in (\ref{1003}). Then, recalling (\ref{90000se1extra}), (\ref{90000}), (\ref{90000se}), (\ref{90000se1}), (\ref{550se}),  (\ref{esti1}), (\ref{kong1}), (\ref{esti1se}), (\ref{kong2})  , (\ref{esti12}) and (\ref{kong3}) making use of  (\ref{gee1}), (\ref{gee2}),   we obtain the second two inequalities from above according to the rationing preceding  $(4.16)$ from \cite{bu3}. Then, returning in (\ref{ecuatieXYY11}), we make use of (\ref{90000se1extra16}), (\ref{Africa90000}), (\ref{Africa90000se}),   (\ref{550se11}),  (\ref{esti1}), (\ref{kong1}), (\ref{esti1se}), (\ref{kong2})  , (\ref{esti12}), (\ref{kong3}),  (\ref{londra1}) and  (\ref{londra2}). Then, we obtain the other above inequalities according to the rationing preceding $(4.17)$ and $(4.18)$ from \cite{bu3}.
   
Now, following the lines from  \cite{bu3}, we ready to move forward:  
\subsection{The Iteration Process} Then, there exist a constant $\delta_{0}\left(d\right)>0$ depending on $n$ and independent on $E\left(Z,\xi\right)$ and $r,\sigma,\rho,r'$ defined by (\ref{num})  such that if the following inequality holds
\begin{equation}\left(\frac{36 \left(1+ d^{2}\left(2d\right)^{2N}\right)}{r-\rho}+6N \left(1+ d^{2}\left(2d\right)^{2N}\right)\right)\frac{d^{2}\left(2d\right)^{4N}\left\|E\right\|_{r}}{N\left(r-\rho\right) }\left(\frac{\rho}{r}\right)^{d-1}<\delta_{0}\left(d\right),\label{X00}\end{equation}
we have that the mapping $\Psi\left(z',w'\right):=H^{-1}\left(z',w'\right)$ is well defined in $\overline{\Delta_{\sigma}}$. Furthermore, it follows that  $\Psi\left(\Delta_{r'}\right)\subset\overline{\Delta_{\sigma}}$, $\Psi\left(\Delta_{\sigma}\right)\subset\Delta_{\rho}$, $E'\left(Z,\xi\right)$ is holomorphic in $\overline{\Delta_{\sigma}}$  and also the following inequality holds
\begin{equation}\begin{split}&\quad\quad\quad\quad\quad\quad\quad\quad\quad\quad\quad\quad\quad\quad  \left\|E'\right\|_{r'}\leq \left\|E\right\|_{r}\frac{3^{2N}\left(2d\right)^{4N}}{\left(r-r'\right)^{2N}}\left(\frac{r'}{r}\right)^{d-1}+\\&  \left\|E\right\|_{r}^{2}\left(\frac{\left(2d\right)^{4N}}{N\left(r-r'\right)}\left( \left(\frac{A\left(d\right)}{r-r'}+B\left(d\right)\right)  \left(\frac{r'}{r}\right)^{\frac{d-1}{2}}+\left(\frac{108}{r-r'}+D\left(d\right)\right)  \left(\frac{r'}{r}\right)^{\frac{2d-3}{4}}\right)+E\left(d\right)\left(\frac{r'}{r}\right)^{2d-3} \right) ,\end{split}\label{ert}\end{equation}
where we have used the following notations
\begin{equation}A\left(d\right)=324\left(1+d^{2}\left(2d\right)^{2N}\right) ,\quad B\left(d\right)=18N\left(1+d^{2}\left(2d\right)^{2N}\right) ,\quad  D\left(d\right)=6N ,\quad E\left(d\right)=100 d^{2}\left(2d\right)^{8N} .\label{AB}\end{equation}
 
The proof of these estimates follow as in \cite{bu3} using an induction argument and 
 the computations of the formal equivalence.

\subsection{Proof of Theorem \ref{t1}}The arguments are extracted and adapted  exactly as in \cite{bu4}  from Moser\cite{mo} and Huang-Yin\cite{huyi2}. Then, following Huang-Yin\cite{huyi1} and Moser\cite{mo}, we define the following sequence of real analytic submanifolds
$$M_{n}:\quad w=Q\left(z,\overline{z}\right)+E_{n}\left(z,\overline{z}\right), $$
or equivalently as follows $M_{0}:=M$, $M_{n+1}:=\Psi^{-1}_{n}\left(M_{n}\right)$, for all $n\in\mathbb{N}$.  Here $\Psi_{n}$ is the holomorphic mapping between $\Delta_{\sigma_{n}}$ and $\Delta_{\rho_{n}}$. It is clear that $d_{n}:=\rm{Ord}\left(E_{n}\left(z,\overline{z}\right)\right)\geq 2^{n}+2$, for all $n\in\mathbb{N}$. 
 
 Following Moser\cite{mo}, we define the following sequences of numbers
  $$r_{n}:=\frac{1}{2}\left(1+\frac{1}{n+1}\right),\quad \rho_{n}=\frac{r_{n+1}+2r_{n}}{3},\quad 
 \sigma_{n}=\frac{\rho_{n}+2r_{n}}{3},$$
 and we apply the estimations with $r=r_{n}$, $\rho=\rho_{n}$, $r'=r_{n+1}$, $\psi=\psi_{n}$, for all $n\in\mathbb{N}$. Following Moser\cite{mo}, we have that 
\begin{equation} \frac{r_{n+1}}{r_{n}}=1-\frac{1}{\left(n+1\right)^{2}},\quad \frac{1}{r_{n}-r_{n+1}}=\left(n+1\right)\left(n+2\right).\label{rrr}\end{equation}
We define the following sequence of real numbers  $$\epsilon_{n}:=\frac{\left\|E\right\|_{r_{n}}}{\left(r_{n}-r_{n+1}\right)^{2}},$$
and then   we obtain the following 
\begin{equation}\begin{split}&\quad\quad\quad  \epsilon_{n+1}\leq \epsilon_{n}\frac{\left(r_{n}-r_{n+1}\right)^{2}}{\left(r_{n+1}-r_{n+2}\right)^{2}}\frac{3^{2N}\left(2d_{n}\right)^{4N}}{\left(r_{n}-r_{n+1}\right)^{2N}}\left(\frac{r_{n+1}}{r_{n}}\right)^{d_{n}-1}+\epsilon_{n}^{2} \frac{\left(r_{n}-r_{n+1}\right)^{4}}{\left(r_{n+1}-r_{n+2}\right)^{2}}E\left(d_{n}\right)\left(\frac{r_{n+1}}{r_{n}}\right)^{2d_{n}-3}+\\&  \epsilon_{n}^{2} \frac{\left(r_{n}-r_{n+1}\right)^{4}}{\left(r_{n+1}-r_{n+2}\right)^{2}}\frac{\left(2d_{n}\right)^{4N}}{N\left(r_{n}-r_{n+1}\right)}\left( \left(\frac{A\left(d_{n}\right)}{r_{n}-r_{n+1}}+B\left(d_{n}\right)\right)  \left(\frac{r_{n+1}}{r_{n}}\right)^{\frac{d_{n}-1}{2}}+\left(\frac{108}{r_{n}-r_{n+1}}+D\left(d_{n}\right)\right)  \left(\frac{r_{n+1}}{r_{n}}\right)^{\frac{2d_{n}-3}{4}}\right)  ,\end{split}\label{ertQ}\end{equation}
and then we continue canonically as in \cite{bu3}. The proof is completed. The details remain as an exercise to the reader.

\section{Proofs of Theorems \ref{tA1} and \ref{tA2}}
We proceed in the light of  (\ref{var1A}),(\ref{var1B}),(\ref{var1Ase}),(\ref{var1Bse}),
(\ref{var1AX}),(\ref{var1BX}). Let's   compute:
\subsection{Transforming Equations and Models}We consider the following non-constant formal mapping
\begin{equation*}\begin{split}&w'=G(w,z),\quad z_{1}'=F_{1}(w,z),z_{2}'=F_{2}(w,z),\dots,z_{N}'=F_{N}(w,z), \\& \quad\quad \quad\quad\quad\quad\quad\hspace{0.1 cm} z_{N+1}'=F_{N+1}(w,z),z_{N+2}'=F_{N+2}(w,z),\dots,z_{N'}'=F_{N'}(w,z),\end{split}
\end{equation*}
which sends (\ref{var1A}) into (\ref{var1B})  according to the following formal expansions
\begin{equation}\begin{split}&G(w,z)=\displaystyle\sum_{m,n\geq 0}G_{m,n}(z)w^{n},\quad \left(F^{\left(1\right)}(w,z),\dots,F^{\left(N\right)}(w,z)\right)=\left(\displaystyle\sum_{m,n\geq 0}F_{m,n}^{\left(1\right)}(z)w^{n},\dots,\displaystyle\sum_{m,n\geq 0}F_{m,n}^{\left(N'\right)}(z)w^{n}\right) ,\\&\quad\quad\quad\quad\quad\quad\quad\quad\quad\quad\quad\quad\quad\quad \left(F^{\left(N+1\right)}(w,z),\dots,F^{\left(N\right)}(w,z)\right)=\left(\displaystyle\sum_{m,n\geq 0}F_{m,n}^{\left(N+1\right)}(z)w^{n},\dots,\displaystyle\sum_{m,n\geq 0}F_{m,n}^{\left(N'\right)}(z)w^{n}\right),\label{440}\end{split}
\end{equation}
where we have are homogeneous polynomials of degree $m$ in $z$, denoted as

$$G_{m,n}(z),\quad F_{m,n}^{\left(1\right)}(z),\dots,F_{m,n}^{\left(N\right)}(z),F_{m,n}^{\left(N+1\right)}(z),\dots,F_{m,n}^{\left(N'\right)}(z).$$

 We obtain  
\begin{equation}
\displaystyle\sum_{m,n\geq 0}G_{m,n}(z)w^{n}= \displaystyle\sum_{k=1}^{N'} \left(\displaystyle\sum_{m,n\geq 0}\left(F_{m,n}^{\left(k\right)}(z)w^{n}\right)\overline{\left(\displaystyle\sum_{m,n\geq 0}F_{m,n}^{\left(k\right)}(z)w^{n}\right)}+\lambda_{k}\left(\displaystyle\sum_{m,n\geq 0}F_{m,n}^{\left(k\right)}(z)w^{n}\right)^{2}+\lambda_{k}\overline{\left(\displaystyle\sum_{m,n\geq 0}F_{m,n}^{\left(k\right)}(z)w^{n}\right)^{2}}\right) . \label{ecuatieXYY}
\end{equation}

We observe the following important aspect
\begin{equation}\mbox{rank}\begin{pmatrix}\frac{\partial F_{1}}{\partial z_{1}}(0) &\dots & \frac{\partial F_{N}}{\partial z_{1}}(0)& \frac{\partial F_{N+1}}{\partial z_{1}}(0)&\dots&\frac{\partial F_{N'}}{\partial z_{1}}(0) \\ \vdots &\ddots & \vdots & \vdots &\ddots&\vdots \\ \frac{\partial F_{1}}{\partial z_{N}}(0) &\dots & \frac{\partial F_{N}}{\partial z_{N}}(0)& \frac{\partial F_{N+1}}{\partial z_{N}}(0)&\dots&\frac{\partial F_{N'}}{\partial z_{N}}(0) \\
\end{pmatrix}=N,
\end{equation}
otherwise we would obtain a contraction with (\ref{transver}). 

Indeed, (\ref{transver}) holds, because contrary we would have 
\begin{equation*}G(z,w)=\mbox{O}(2),
\end{equation*}
which implies that $G(z,w)=0$, in the light of the Fischer Decompositions defining (\ref{spartiuG}), 
and also that  $F(z,w)=0$, because of (\ref{ecuatieXYY}), obtaining again a contradiction because the formal embedding was assumed to be not constant.

Imposing the normalization procedure of Baouendi-Huang\cite{BH} similarly as in \cite{bu5}, we obtain 
\begin{equation}\begin{split}&G(z,w)=w+\mbox{O}(2),\quad F^{\left(1\right)}(w,z)=z_{1}+\mbox{O}(2),\dots, F^{\left(N\right)}(w,z)=z_{N}+\mbox{O}(2),\\&\quad\quad\quad\quad\quad\quad\quad\quad\quad\quad F^{\left(N+1\right)}(w,z)=\mbox{O}(2),\dots, F^{\left(N'\right)}(w,z)= \mbox{O}(2),\end{split}
 \end{equation} 
by eventually composing (\ref{440}) with elements of $\mbox{Aut}\left(M_{\infty}\right)$ and   $\mbox{Aut}\left(M_{\infty}'\right)$, giving sense to the assumption (\ref{maro}). 

We study the equivalence defined by (\ref{440}) recalling the strategy from \cite{bu5}. Eventually composing again (\ref{440}) with elements of $\mbox{Aut}\left(M_{\infty}\right)$ and   $\mbox{Aut}\left(M_{\infty}'\right)$, we obtain
\begin{equation}\Re F_{1,n+1}^{\left(1\right)}(z)=\dots=\Re F_{1,n+1}^{\left(N\right)}(z)=  0,\quad\mbox{for all $n\in\mathbb{N}^{\star}$.}\label{bev}
  \end{equation}

Separating the imaginary side from the real part in (\ref{ecuatieXYY}), we obtain
\begin{equation}G(z,w)=w,
\end{equation}
in the light of the Fischer Decompositions defining (\ref{spartiuG}).

Now, in order to move forward with these computations, we write as follows
\begin{equation} \left(F^{\left(1\right)},F^{\left(2\right)},\dots,F^{\left(N\right)}\right)(w,z)=\left(z_{1},z_{2},\dots,z_{N}\right)+A\left(z_{1},z_{2},\dots,z_{N}\right)w+B_{0}\left(z_{1},z_{2},\dots,z_{N}\right)+B_{1}\left(z_{1},z_{2},\dots,z_{N}\right)+\mbox{O}(4),\label{882}
\end{equation}
where we deal with
\begin{itemize}
\item $B_{0}\left(z_{1},z_{2},\dots,z_{N}\right)$ is a vector   polynomial of degree $2$ in  $\left(z_{1},z_{2},\dots,z_{N}\right)$,
\item $B_{1}\left(z_{1},z_{2},\dots,z_{N}\right)$ is a vector   polynomial of degree $3$ in  $\left(z_{1},z_{2},\dots,z_{N}\right)$,
\item  $A\left(z_{1},z_{2},\dots,z_{N}\right)$ is a linear form in $\left(z_{1},z_{2},\dots,z_{N}\right)$,   written   as follows
\begin{equation} A\left(z_{1},z_{2},\dots,z_{N}\right)=\begin{pmatrix}
a_{11} & a_{12} &\dots& a_{1N} \\ a_{21} & a_{22} &\dots& a_{2N} \\  \vdots & \vdots &\ddots& \vdots \\ a_{N1} & a_{N2} &\dots& a_{NN} \\ 
\end{pmatrix}\begin{pmatrix}
z_{1} \\ z_{2}\\ \vdots \\ z_{N}
\end{pmatrix},\label{lalala}
\end{equation}\end{itemize}
having in mind (\ref{bev}), according to the following notations
\begin{equation}\begin{split}& B_{0}\left(z_{1},z_{2},\dots,z_{N}\right)=\left(B_{0}^{(1)},B_{0}^{(2)},\dots,B_{0}^{\left(N\right)}\right)\left(z_{1},z_{2},\dots,z_{N}\right),\\& B_{1}\left(z_{1},\dots,z_{N}\right)=\left(B_{1}^{(1)},B_{1}^{(2)},\dots,B_{1}^{\left(N\right)}\right)\left(z_{1},z_{2},\dots,z_{N}\right),\end{split}\label{lalala1}
\end{equation}
 
Focusing on the real part  in (\ref{ecuatieXYY}), we obtain 
\begin{equation}\displaystyle\sum_{l=1}^{N}B_{0}^{\left(l\right)} \overline{z}_{l}=0,
\end{equation}
because there are no homogeneous terms of degree $3$ in the right hand side of (\ref{ecuatieXYY}), concluding that
$$ B_{0}\left(z_{1},z_{2}\dots,z_{N}\right)=0.
$$

Moreover, we have
\begin{equation}\begin{split}& 2Re\left(\displaystyle\sum_{k,l=1}^{N}a_{kl}z_{l}\overline{z}_{k} +\displaystyle\sum_{l =1}^{N} 2\lambda_{l}a_{lk}z_{l}z_{k}\right)\left(z_{1}\overline{z}_{1}+\dots+z_{N}\overline{z}_{N}+\lambda_{1}\left(z_{1}^{2}+\overline{z}_{1}^{2}\right)+\dots+\lambda_{N}\left(z_{N}^{2}+\overline{z}_{N}^{2}\right)\right)\\&+  2\Re\left(\displaystyle\sum_{l=1}^{N}B^{\left(l\right)}_{1} \overline{z}_{l}+\sum_{l=1}^{N}B^{\left(l\right)}_{1} 2\lambda_{l} z_{l}\right)\\& \quad\quad\quad\quad\quad   \quad\quad\quad\quad\quad\quad  \quad\quad\quad\quad\quad \quad\quad\quad  \quad\quad\quad\quad\quad\quad\quad\quad  \begin{tabular}{l} \rotatebox[origin=c]{270}{$=$}\end{tabular} \\&-\left(\displaystyle\sum_{l=N+1}^{N'}F_{2,0}^{\left(l\right)}(z)\right)\overline{\left(\displaystyle\sum_{l=N+1}^{N'}F_{2,0}^{\left(l\right)}(z)\right)}                
-\displaystyle\sum_{l=N+1}^{N'}\lambda'_{l}\left( F_{2,0}^{\left(l\right)}(z)\right)^{2}-\displaystyle\sum_{l=N+1}^{N'}\lambda'_{l}  \left(\overline{ F_{2,0}^{\left(l\right)}(z)}\right)^{2} .\end{split}\label{ect}
\end{equation}

Recalling Lemma $3.2$ from Huang\cite{huang1}, we obtain
  $$F_{2,0}^{\left(l\right)}(z)=0,\hspace{0.1 cm}\mbox{for all $l=N+1,\dots, N'$.}
$$

According to the Fischer Decompositions defined by (\ref{spartiuF}), we obtain
$$ A\left(z_{1},z_{2},\dots,z_{N}\right)=0,\quad B_{1}\left(z_{1},\dots,z_{N}\right)=0.$$

We move forward:
\subsection{Proof of Theorem \ref{tA1}}  It is continued the previous analysis in (\ref{ecuatieXYY}). We have
\begin{equation} \left(F^{\left(1\right)},F^{\left(2\right)},\dots,F^{\left(N\right)}\right)(w,z)=\left(z_{1},z_{2},\dots,z_{N}\right)+A_{0}\left(z_{1},z_{2},\dots,z_{N}\right)w+B_{0}\left(z_{1},z_{2},\dots,z_{N}\right)+vw^{2}+\mbox{O}(5),\label{801}
\end{equation}
where we deal with
\begin{itemize}
\item $B_{0}\left(z_{1},z_{2},\dots,z_{N}\right)$ is a vector homogeneous polynomial of degree $4$ in  $\left(z_{1},z_{2},\dots,z_{N}\right)$,
\item $A_{0}\left(z_{1},z_{2},\dots,z_{N}\right)$ is a vector homogeneous polynomial of degree $2$ in  $\left(z_{1},z_{2},\dots,z_{N}\right)$, 
\item $v\in\mathbb{C}^{N}$,
\end{itemize}
according to the following notations
$$B_{0}\left(z_{1},z_{2},\dots,z_{N}\right)=\left(B_{0}^{(1)},B_{0}^{(2)},\dots,B_{0}^{\left(N\right)}\right)\left(z_{1},z_{2},\dots,z_{N}\right),$$  $$ A_{0}\left(z_{1},z_{2},\dots,z_{N}\right)=\left(A_{0}^{(1)},A_{0}^{(2)},\dots,A_{0}^{\left(N\right)}\right)\left(z_{1},z_{2},\dots,z_{N}\right),$$ $$ v=\left(v_{1},v_{2},\dots,v_{N}\right)\in\mathbb{C}^{N}.$$

Extracting  terms of degree $5$ in $(z,\overline{z})$ from (\ref{ecuatieXYY}), we obtain
  \begin{equation}\begin{split}&\quad\quad\quad\quad\quad\quad\quad\quad\quad\quad\quad\quad\quad\quad      \Re \left( \displaystyle\sum_{l=1}^{N} \overline{z}_{l}A_{0}^{\left(l\right)} +  \displaystyle\sum_{l=1}^{N} 2\lambda_{l}z_{l}A_{0}^{\left(l\right)}\right) \left(z_{1}\overline{z}_{1}+\dots+z_{N}\overline{z}_{N}+\lambda_{1}\left(z_{1}^{2}+\overline{z}_{1}^{2}\right)+\dots+\lambda_{N}\left(z_{N}^{2}+\overline{z}_{N}^{2}\right)  \right)\\& \quad\quad\quad\quad\quad   \quad\quad\quad\quad\quad\quad  \quad\quad\quad\quad\quad \quad\quad\quad  \quad\quad\quad\quad\quad\quad\quad\quad  \begin{tabular}{l} \rotatebox[origin=c]{270}{$+$}\end{tabular} \\&   \Re \left( \displaystyle\sum_{l=1}^{N} \overline{z}_{l}B_{0}^{\left(l\right)} +\displaystyle\sum_{l=1}^{N} 2\lambda_{l}z_{l}B_{0}^{\left(l\right)}\right)  \\& \quad\quad\quad\quad\quad   \quad\quad\quad\quad\quad\quad  \quad\quad\quad\quad\quad \quad\quad\quad  \quad\quad\quad\quad\quad\quad\quad\quad  \begin{tabular}{l} \rotatebox[origin=c]{270}{$+$}\end{tabular} \\&  2\Re \left( \displaystyle\sum_{l=1}^{N} \overline{z}_{l}\overline{v_{l}} +\displaystyle\sum_{l=1}^{N}2\lambda_{l}z_{l} v_{l}\right)  \left(z_{1}\overline{z}_{1}+\dots+z_{N}\overline{z}_{N}+\lambda_{1}\left(z_{1}^{2}+\overline{z}_{1}^{2}\right)+\dots+\lambda_{N}\left(z_{N}^{2}+\overline{z}_{N}^{2}\right) \right)^{2}=0.
\end{split}\end{equation}

According to the Fischer Decompositions defined by (\ref{spartiuF}), it follows that
$$B_{0}\left(z_{1},z_{2},\dots,z_{N}\right)=0,\quad A_{0}\left(z_{1},z_{2},\dots,z_{N}\right)=0,\quad v=0.$$ 

We move forward to the next degree in (\ref{ecuatieXYY}). We have
\begin{equation}\left(F^{\left(1\right)},F^{\left(2\right)},\dots,F^{\left(N\right)}\right)(w,z)=\left(z_{1},z_{2},\dots,z_{N}\right)+A\left(z_{1},z_{2},\dots,z_{N}\right)w^{2}+B_{1}\left(z_{1},\dots,z_{N}\right)w+B_{0}\left(z_{1},z_{2},\dots,z_{N}\right)+ \mbox{O}(5),\label{901}
\end{equation}
where we deal with
\begin{itemize}
\item $B_{1}\left(z_{1},z_{2},\dots,z_{N}\right)$ is a vector homogeneous polynomial of degree $5$ in  $\left(z_{1},z_{2},\dots,z_{N}\right)$,
\item $B_{2}\left(z_{1},z_{2},\dots,z_{N}\right)$ is a vector homogeneous polynomial of degree $3$ in  $\left(z_{1},z_{2},\dots,z_{N}\right)$,
\item $A_{0}\left(z_{1},z_{2},\dots,z_{N}\right)$ is a linear form in  $\left(z_{1},z_{2},\dots,z_{N}\right)$  defined as in (\ref{lalala}),
\end{itemize}
according to the following notations 
\begin{equation*}
\begin{split}&B_{0}\left(z_{1},z_{2},\dots,z_{N}\right)=\left(B_{0}^{(1)},B_{0}^{(2)},\dots,B_{0}^{\left(N\right)}\right)\left(z_{1},z_{2},\dots,z_{N}\right) ,\\& B_{1}\left(z_{1},z_{2},\dots,z_{N}\right)=\left(B_{1}^{(1)},B_{1}^{(2)},\dots,B_{1}^{\left(N\right)}\right)\left(z_{1},z_{2},\dots,z_{N}\right) .\end{split}
\end{equation*}

We obtain
 \begin{equation}\begin{split}&\quad\quad\quad\quad\quad\quad\quad\quad\quad\quad\quad\quad    2 Re\left(\displaystyle\sum_{k,l=1}^{N}a_{kl}z_{l}\overline{z}_{k} +\displaystyle\sum_{k =1}^{N} 2\lambda_{k}a_{lk}z_{l}z_{k}\right) \left(z_{1}\overline{z}_{1}+\dots+z_{N}\overline{z}_{N}+\lambda_{1}\left(z_{1}^{2}+\overline{z}_{1}^{2}\right)+\dots+\lambda_{N}\left(z_{N}^{2}+\overline{z}_{N}^{2}\right) \right)^{2}\\& \quad\quad\quad\quad\quad   \quad\quad\quad\quad\quad\quad  \quad\quad\quad\quad\quad \quad\quad\quad  \quad\quad\quad\quad\quad\quad\quad\quad  \begin{tabular}{l} \rotatebox[origin=c]{270}{$+$}\end{tabular} \\& 2Re\left(\displaystyle\sum_{l=1}^{N}B_{0}^{\left(l\right)}\overline{z}_{k} +  \displaystyle\sum_{l=1}^{N}2\lambda_{l} z_{l}B_{0}^{\left(l\right)} \right) \left(z_{1}\overline{z}_{1}+\dots+z_{N}\overline{z}_{N}+\lambda_{1}\left(z_{1}^{2}+\overline{z}_{1}^{2}\right)+\dots+\lambda_{N}\left(z_{N}^{2}+\overline{z}_{N}^{2}\right)\right) \\& \quad\quad\quad\quad\quad   \quad\quad\quad\quad\quad\quad  \quad\quad\quad\quad\quad \quad\quad\quad  \quad\quad\quad\quad\quad\quad\quad\quad  \begin{tabular}{l} \rotatebox[origin=c]{270}{$+$}\end{tabular} \\&2\Re\left(\displaystyle\sum_{l=1}^{N}B_{1}^{\left(l\right)}\overline{z}_{l} +  \displaystyle\sum_{l=1}^{N}2\lambda_{l} z_{l}B_{1}^{\left(l\right)} \right)  \\& \quad\quad\quad\quad\quad   \quad\quad\quad\quad\quad\quad  \quad\quad\quad\quad\quad \quad\quad\quad  \quad\quad\quad\quad\quad\quad\quad\quad  \begin{tabular}{l} \rotatebox[origin=c]{270}{$=$}\end{tabular} \\&-\left(\displaystyle\sum_{l=N+1}^{N'}F_{3,0}^{\left(l\right)}(z)\right)\overline{\left(\displaystyle\sum_{l=N+1}^{N'}F_{3,0}^{\left(l\right)}(z)\right)}      
-\displaystyle\sum_{l=N+1}^{N'}\lambda'_{l}\left( F_{3,0}^{\left(l\right)}(z)\right)^{2}-\displaystyle\sum_{l=N+1}^{N'}\lambda'_{l}  \left(\overline{ F_{3,0}^{\left(l\right)}(z)}\right)^{2}.\end{split}\label{ect1}
\end{equation}

Recalling Lemma $3.2$ from Huang\cite{huang1}, we obtain
 $$F_{3,0}^{\left(l\right)}(z)=0,\hspace{0.1 cm}\mbox{for all $l=N+1,\dots, N'$.}
$$

According to the Fischer Decompositions defined by (\ref{spartiuF}), it follows that$$B_{0}\left(z_{1},z_{2},\dots,z_{N}\right)=0,\quad B_{1}\left(z_{1},z_{2},\dots,z_{N}\right)=0.$$

It becomes clear now the induction process administrated straightforwardly in (\ref{ecuatieXYY}) providing obviously the equivalence class from (\ref{bebe1}).  
 
\subsection{Proof of Theorem \ref{tA2}}Let $M\subset\mathbb{C}^{N+1}$ be a real-formal
submanifold defined near  $p=0$ as in (\ref{ecuatie}), and respectively, consider 
 another real-formal submanifold $M'\subset\mathbb{C}^{N+1}$ defined near $p=0$ as in (\ref{ecuatie1}) such that (\ref{diag1}) holds.

It is important to observe that any formal embedding of $M$ into the corresponding model, can be written as follows
\begin{equation}
(w,z)\rightarrow \left(w+\mbox{O}(2),z+\mbox{O}(2),\mbox{O}(2)\right).
\end{equation}

More precisely, in the light of the previous computations, we have
\begin{equation}
(w,z)\rightarrow \left(w+\mbox{O}(2),z+\mbox{O}(2),0\right).
\end{equation}

We take an embedding of similar type for $M'$. It follows that
\begin{equation}\varphi_{k}(z,\overline{z})=\varphi' _{k}\left(z',\overline{z'}\right)|_{z'=(z,0)},\quad\mbox{for all
 $k\geq 3$. } \label{yuk}\end{equation}
 
It is an exercise to show that (\ref{yuk}) does not depend on the chosen embeddings. Moreover, we obtain easily the last part of the statement in the light of the previous computations by easy and straightforward computations.

 \section{Ingredients and Computations using Formal Holomorphic Segre Preserving Mappings}
 Similarly as before, we   consider  indications from my professor\cite{za}. More precisely, we have to study how a Formal (Holomorphic) Segre Preserving Equivalence occurs in the local defining equations. Recalling  the terminology from \cite{V77}, we recall  the local defining equations  (\ref{ecuatie}) and (\ref{ecuatie1})   as follows 
\begin{equation}\begin{split}&
\frac{w+\overline{w}}{2}=z_{1}\overline{z}_{1}+\dots+z_{N}\overline{z}_{N}+\lambda_{1}\left(z_{1}^{2}+\overline{z}_{1}^{2}\right)+\dots+\lambda_{N}\left(z_{N}^{2}+\overline{z}_{N}^{2}\right)+\displaystyle\sum
_{k\geq 3}\frac{\varphi_{k}(z,\overline{z})+\overline{\varphi_{k}(z,\overline{z})}}{2},\\&\frac{w-\overline{w}}{2\sqrt{-1}}= \displaystyle\sum
_{k\geq 3}\frac{\varphi_{k}(z,\overline{z})-\overline{\varphi_{k}(z,\overline{z})}}{2\sqrt{-1}}, \end{split} \label{ecuatielala1}
\end{equation}
and respectively 
\begin{equation}\begin{split}&\frac{w'+\overline{w'}}{2}= {z'}_{1}\overline{{z'}}_{1}+\dots+{z'}_{N}\overline{{z'}}_{N}+\lambda_{1}\left({z'}_{1}^{2}+\overline{{z'}}_{1}^{2}\right)+\dots+\lambda_{N}\left({z'}_{N}^{2}+\overline{{z'}}_{N}^{2}\right)+\displaystyle\sum _{k\geq
3}\frac{\varphi'_{k}\left(z',\overline{z'}\right)+\overline{\varphi'_{k}\left(z',\overline{z'}\right)}}{2 },\\& \frac{w'-\overline{w'}}{2\sqrt{-1}}=  \displaystyle\sum _{k\geq
3}\frac{\varphi'_{k}\left(z',\overline{z'}\right)-\overline{\varphi'_{k}\left(z',\overline{z'}\right)}}{2\sqrt{-1}}.\end{split} \label{ecuatie1lala1}
\end{equation}

Now, we proceed as follows:
\subsection{Changes of Coordinates}We use the following notations:
 \begin{equation*} \overline{w}=\nu,\hspace{0.1 cm}\overline{w'}=\nu',\quad \overline{z}_{1}=\xi_{1},\hspace{0.1 cm}\overline{z}_{2}=\xi_{2},\dots \overline{z}_{N}=\xi_{N} ,\quad \overline{z'}_{1}={\xi'}_{1},\hspace{0.1 cm}\overline{z'}_{2}={\xi'}_{2},\dots \overline{z'}_{N}={\xi'}_{N}.
    \end{equation*}

Now, let $N\subset\mathbb{C}^{2N+2}$   defined near  $p=0$ as follows
\begin{equation}\begin{split}&
\frac{w+\nu}{2}=z_{1}\xi_{1}+\dots+z_{N}\xi_{N}+\lambda_{1}\left(z_{1}^{2}+\xi_{1}^{2}\right)+\dots+\lambda_{N}\left(z_{N}^{2}+\xi_{N}^{2}\right)+\displaystyle\sum
_{k\geq 3}\frac{\varphi_{k}\left(z,\xi\right)+\overline{\varphi}_{k}\left(\xi,z\right)}{2},\\&\frac{w-\nu}{2\sqrt{-1}}= \displaystyle\sum
_{k\geq 3}\frac{\varphi_{k}\left(z,\xi\right)-\overline{\varphi}_{k}\left(\xi,z\right)}{2\sqrt{-1}}, \end{split} \label{ecuatielala}
 \end{equation} 
and respectively, let   $N'\subset\mathbb{C}^{2N+2}$ defined near $p=0$ as follows
\begin{equation}\begin{split}&\frac{w'+\nu'}{2}= {z'}_{1}{\xi'}_{1}+\dots+{z'}_{N}{\xi'}_{N}+\lambda_{1}\left({z'}_{1}^{2}+{\xi'}_{1}^{2}\right)+\dots+\lambda_{N}\left({z'}_{N}^{2}+{\xi'}_{N}^{2}\right)+\displaystyle\sum _{k\geq
3}\frac{{\varphi'}_{k}\left(z',\xi'\right)+\overline{\varphi'}_{k}\left(\xi',z'\right)}{2},\\& \frac{w'-\nu'}{2\sqrt{-1}}=  \displaystyle\sum _{k\geq
3}\frac{{\varphi'}_{k}\left(z',\xi'\right)-\overline{\varphi'}_{k}\left(\xi',z'\right)}{2\sqrt{-1}},\end{split} \label{ecuatie1lala}
\end{equation}
in the light of    (\ref{ecuatielala1})  and (\ref{ecuatie1lala1}).

Rewriting now (\ref{ecuatielala}), we obtain 
\begin{equation}\begin{split}&
w=z_{1}\xi_{1}+\dots+z_{N}\xi_{N}+\lambda_{1}\left(z_{1}^{2}+\xi_{1}^{2}\right)+\dots+\lambda_{N}\left(z_{N}^{2}+\xi_{N}^{2}\right)+\displaystyle\sum
_{k\geq 3} \varphi_{k}\left(z,\xi\right) ,\\& \nu=z_{1}\xi_{1}+\dots+z_{N}\xi_{N}+\lambda_{1}\left(z_{1}^{2}+\xi_{1}^{2}\right)+\dots+\lambda_{N}\left(z_{N}^{2}+\xi_{N}^{2}\right)+\displaystyle\sum
_{k\geq 3} \overline{\varphi}_{k}\left( \xi,z\right), \end{split} \label{ecuatielalase}
 \end{equation} 
but we continue to use (\ref{ecuatie1lala}) in its current shape in order to work with Formal Holomorphic Segre Preserving Mappings defined as follows:

\subsection{Transforming Equations} Let a  Formal (Holomorphic) Segre Preserving Equivalence between  $N$ and $N'$, which fixes the point
$0\in\mathbb{C}^{2N+2}$, denoted as follows
$$\left(z',\xi';w',\nu'\right)=\left(F(z,w), \tilde{F}\left(\xi,\nu\right);G(z,w),\tilde{G}\left(\xi,\nu\right)\right).$$  

Then, for  $w$ and  $\nu$ defined by (\ref{ecuatielalase}), we have 
\begin{equation}\begin{split}&\frac{G(z,w)+\tilde{G}\left(\xi,\nu\right)}{2}= Q\left(
F(z,w),\tilde{F}\left(\xi,\nu\right)\right)+\displaystyle\sum _{k\geq
3}\frac{{\varphi'}_{k}\left(
F(z,w),\tilde{F}\left(\xi,\nu\right)\right)+\overline{\varphi'}_{k}\left(
 \tilde{F}\left(\xi,\nu\right),F(z,w)\right)}{2},\\& \frac{G(z,w)-\tilde{G}\left(\xi,\nu\right)}{2\sqrt{-1}}=  \displaystyle\sum _{k\geq
3}\frac{{\varphi'}_{k}\left(
F(z,w),\tilde{F}\left(\xi,\nu\right)\right)-\overline{\varphi'}_{k}\left(
 \tilde{F}\left(\xi,\nu\right),F(z,w)\right)}{2\sqrt{-1}}.\end{split} \label{eclala}
\end{equation}

In order to understand better the interactions of terms in (\ref{eclala}), we write this (formal) equivalence as follows
\begin{equation}\left(F(z,w),\tilde{F}\left(\xi,\nu\right);G(z,w),\tilde{G}\left(\xi,\nu\right)\right)
=\left(\displaystyle\sum_{m,n\geq 0}F_{m,n}(z)w^{n},\displaystyle\sum_{m,n\geq 0}\tilde{F}_{m,n}\left(\xi\right)\nu^{n};\displaystyle\sum_{m,n\geq 0}G_{m,n}(z)w^{n},\displaystyle\sum_{m,n\geq 0}\tilde{G}_{m,n}\left(\xi\right)\nu^{n}\right), 
\label{maplala}\end{equation}
where $G_{m,n}(z)$, $F_{m,n}(z)$ are homogeneous polynomials of degree
$m$ in $z$ for all  $ n\in\mathbb{N}$, and respectively $\tilde{G}_{m,n}\left(\xi\right)$, $\tilde{F}_{m,n}\left(\xi\right)$ are homogeneous polynomials of degree
$m$ in $z$ for all  $ n\in\mathbb{N}$. 

We conclude then by (\ref{eclala}) and (\ref{maplala}) that
\begin{equation}\left.\begin{split}\displaystyle & \quad    \sum
_{m,n\geq0}\left( G_{m,n}(z)\left(Q\left(z,\xi\right) +\displaystyle\sum
_{k\geq 3}\varphi_{k}\left(z,\xi\right)\right)^{n}+ \tilde{G}_{m,n}\left(\xi\right)\left(Q\left(z,\xi\right) +\displaystyle\sum
_{k\geq 3}\overline{\varphi}_{k}(\xi,z)\right)^{n}\right) \\& \quad\quad\quad\quad\quad   \quad\quad\quad\quad\quad\quad\quad \quad\quad\quad\quad\quad   \quad\quad\quad\quad\quad   \begin{tabular}{l} \rotatebox[origin=c]{270}{$=$}\end{tabular} \\& \quad\hspace{0.3 cm}  2 Q\left(\displaystyle\sum _{m,n \geq 0}
F_{m,n}(z)\left(Q\left(z,\xi\right) +\displaystyle\sum
_{k\geq 3}\varphi_{k}\left(z,\xi\right)\right)^{n},\tilde{F}_{m,n}\left(\xi\right)\left(Q\left(z,\xi\right) +\displaystyle\sum
_{k\geq 3}\overline{\varphi}_{k}(\xi,z)\right)^{n}\right)  \\& \quad\quad\quad\quad\quad   \quad\quad\quad\quad\quad\quad\quad \quad\quad\quad\quad\quad   \quad\quad\quad\quad\quad   \begin{tabular}{l} \rotatebox[origin=c]{270}{$+$}\end{tabular} \\& \displaystyle\sum
_{k\geq 3}\varphi_{k}'
\left(\displaystyle\sum _{m,n \geq
0}F_{m,n}(z)\left(\left(z,\xi\right) +\displaystyle\sum
_{k\geq 3}\varphi_{k}\left(z,\xi\right)\right)^{n},\tilde{F}_{m,n}\left(\xi\right)\left(Q\left(z,\xi\right) +\displaystyle\sum
_{k\geq 3}\overline{\varphi}_{k}(\xi,z)\right)^{n}\right)\\& \quad\quad\quad\quad\quad   \quad\quad\quad\quad\quad\quad\quad \quad\quad\quad\quad\quad   \quad\quad\quad\quad\quad   \begin{tabular}{l} \rotatebox[origin=c]{270}{$+$}\end{tabular} \\& \displaystyle\sum
_{k\geq 3}\overline{{\varphi'}}_{k} 
\left(\tilde{F}_{m,n}\left(\xi\right)\left(Q\left(z,\xi\right) +\displaystyle\sum
_{k\geq 3}\overline{\varphi}_{k}(\xi,z)\right)^{n},\displaystyle\sum _{m,n \geq
0}F_{m,n}(z)\left(\left(z,\xi\right) +\displaystyle\sum
_{k\geq 3}\varphi_{k}\left(z,\xi\right)\right)^{n} \right),
\end{split}\right.
\label{eculala1}\end{equation}
and respectively, the following
\begin{equation}\left.\begin{split}    \displaystyle & \quad\hspace{0.1 cm}    \sum
_{m,n\geq0}\left( G_{m,n}(z)\left(Q\left(z,\xi\right) +\displaystyle\sum
_{k\geq 3}\varphi_{k}\left(z,\xi\right)\right)^{n}+ \tilde{G}_{m,n}\left(\xi\right)\left(Q\left(z,\xi\right) +\displaystyle\sum
_{k\geq 3}\overline{\varphi}_{k}(\xi,z)\right)^{n}\right) \\& \quad\quad\quad\quad\quad   \quad\quad\quad\quad\quad\quad\quad \quad\quad\quad\quad\quad   \quad\quad\quad\quad\quad   \begin{tabular}{l} \rotatebox[origin=c]{270}{$=$}\end{tabular} \\&    \displaystyle\sum
_{k\geq 3}\varphi_{k}'
\left(\displaystyle\sum _{m,n \geq
0}F_{m,n}(z)\left(\left(z,\xi\right) +\displaystyle\sum
_{k\geq 3}\varphi_{k}\left(z,\xi\right)\right)^{n},\tilde{F}_{m,n}\left(\xi\right)\left(Q\left(z,\xi\right) +\displaystyle\sum
_{k\geq 3}\overline{\varphi}_{k}(\xi,z)\right)^{n}\right)\\& \quad\quad\quad\quad\quad   \quad\quad\quad\quad\quad\quad\quad \quad\quad\quad\quad\quad   \quad\quad\quad\quad\quad   \begin{tabular}{l} \rotatebox[origin=c]{270}{$+$}\end{tabular} \\& \displaystyle\sum
_{k\geq 3}\overline{{\varphi'}}_{k} 
\left(\tilde{F}_{m,n}\left(\xi\right)\left(Q\left(z,\xi\right) +\displaystyle\sum
_{k\geq 3}\overline{\varphi}_{k}(\xi,z)\right)^{n},\displaystyle\sum _{m,n \geq
0}F_{m,n}(z)\left(\left(z,\xi\right) +\displaystyle\sum
_{k\geq 3}\varphi_{k}\left(z,\xi\right)\right)^{n} \right),
\end{split}\right.
\label{eculala2}\end{equation}

 We obtain $G_{0,0}(z)=0$, $\tilde{G}_{0,0}(z)=0$ and
 $F_{0,0}(z)=0$, $\tilde{F}_{0,0}(z)=0$, because the equivalence (\ref{maplala})   fixes  the
point $0\in\mathbb{C}^{2N+2}$. Moreover, collecting the terms of bidegree $(1,0)$ and $(1,1)$ in
$\left(z,\xi\right)$ from (\ref{eculala1}) and (\ref{eculala2}), we obtain   
$G_{1,0}(z)=0$  and   $\tilde{G}_{1,0}(z)=0$ , because
  \begin{equation}G_{0,1}Q\left(z,\xi\right)=Q\left(
F_{1,0}(z),\tilde{F}_{1,0}\left(\xi\right)\right).\label{gelala}\end{equation} 

We can thus assume $G_{0,1}(z)=1$, $\tilde{G}_{0,1}\left(\xi\right)=1$ and $F_{1,0}(z)=z$, $\tilde{F}_{1,0}\left(\xi\right)=\xi$  by   composing (\ref{map})
with a   linear holomorphic automorphism of the following  model 
\begin{equation}
w=z_{1}\xi_{1}+\dots+z_{N}\xi_{N}+\lambda_{1}\left(z_{1}^{2}+\xi_{1}^{2}\right)+\dots+\lambda_{N}\left(z_{N}^{2}+\xi_{N}^{2}\right) ,\quad w=\nu, \label{model1lala}
\end{equation}
where $\lambda_{1},\dots,\lambda_{N}$ satisfy (\ref{lambidalala}).

Following the previous computations lines, we study  interactions of   homogeneous terms of the formal    equivalence (\ref{maplala})  in  (\ref{eculala1}) and (\ref{eculala2}). The  computational obstacles are eliminated using iterative Fischer Decompositions\cite{sh} as follows:

\subsection{Fischer Decompositions\cite{sh}}Recalling the strategy from \cite{bu2} according to the assumption of complexification in the local defining equations considered previously, we define by abuse of notation, according to (\ref{model}), the following differential operator
\begin{equation} \tr=\frac{\partial^{2}}{\partial z_{1}\partial\xi_{1}}+\dots+\frac{\partial^{2}}{\partial z_{N}\partial\xi_{N}}+\lambda_{1}\left( \frac{\partial^{2} }{\partial z_{1}^{2}}+\frac{\partial^{2} }{ \partial\xi_{1}^{2}} \right)+\dots+\lambda_{N}\left(\frac{\partial^{2} }{\partial z_{N}^{2}}+\frac{\partial^{2} }{\partial\xi_{N}^{2}}\right),\label{tracelala}
\end{equation}
where $\lambda_{1},\dots,\lambda_{N}$ satisfy (\ref{lambidalala}).

This differential operator (\ref{tracelala}) is just the Fischer differential operator associated to the complexification of the polynomial (\ref{model}).  Then, recalling the Fischer Decomposition from Shapiro\cite{sh}, we   write by (\ref{model1lala})  and (\ref{yy}) uniquely, by abuse of notation, as follows
\begin{equation} \begin{split}& \quad\quad\quad\quad\quad\hspace{0.15 cm}  z^{I} =A\left(z,\xi\right)Q\left(z,\xi\right)+C\left(z,\xi\right),\quad\hspace{0.2 cm}\mbox{where $\tilde{\tr}\left(C\left(z,\xi\right)\right)=0$, }\\&\left(  \overline{z}_{l}+2\lambda_{l}z_{l}\right)   z^{J} =C_{l}\left(z,\xi\right)Q\left(z,\xi\right)+D_{l}\left(z,\xi\right),\quad\mbox{where $\tilde{\tr}\left(D_{l}\left(z,\xi\right)\right)=0$, for all $l=1,\dots,N$,}\end{split} 
\label{V2lala1}
\end{equation}  
and respectively, we write as follows 
\begin{equation} \begin{split}& \quad\quad\quad\quad\quad\hspace{0.15 cm}  \xi^{I} =\tilde{A}\left(z,\xi\right)Q\left(z,\xi\right)+\tilde{C}\left(z,\xi\right),\quad\mbox{where $\tilde{\tr}\left(\tilde{C}\left(z,\xi\right)\right)=0$, }\\&\left(  z_{l}+2\lambda_{l}\xi_{l}\right)   \xi^{J} =\tilde{C}_{l}\left(z,\xi\right)Q\left(z,\xi\right)+\tilde{D}_{l}\left(z,\xi\right), \hspace{0.16 cm}\mbox{where $\tilde{\tr}\left(\tilde{D}_{l}\left(z,\xi\right)\right)=0$, for all $l=1,\dots,N$,}\end{split} 
\label{V2lala2}
\end{equation} 

Now, similarly to (\ref{lap}),  it is naturally required to define,  abuse of notation, the following  sets of multi-indexes
\begin{equation} \begin{split}&\hspace{0.1 cm}\mathcal{S}=\left\{ \mbox{$I=\left(i_{1},i_{2},\dots,i_{N}\right) \in\mathbb{N}^{N} $ such that $\tilde{\tr} \left(z^{I}\right)=0$ and $i_{1}+i_{2}+\dots+i_{N}\geq 3$ }\right\},\\&\mathcal{T}_{k}=\left\{ \mbox{$J=\left(j_{1},j_{2},\dots,j_{N}\right) \in\mathbb{N}^{N} $ such that $\tilde{\tr}\left(\left(  \xi_{k}+2\lambda_{k}z_{k}\right) z^{J}\right)=0$ and $j_{1}+j_{2}+\dots+j_{N}\geq 2$ }\right\},\quad\mbox{for all $k=1,\dots,N$.}\end{split} \label{788lala}
\end{equation}

Respectively, we define the 
following  sets of multi-indexes
\begin{equation} \begin{split}&\hspace{0.1 cm}\tilde{\mathcal{S}}=\left\{ \mbox{$I=\left(i_{1},i_{2},\dots,i_{N}\right) \in\mathbb{N}^{N} $ such that $\tilde{\tr} \left(\xi^{I}\right)=0$ and $i_{1}+i_{2}+\dots+i_{N}\geq 3$ }\right\},\\&\tilde{\mathcal{T}}_{k}=\left\{ \mbox{$J=\left(j_{1},j_{2},\dots,j_{N}\right) \in\mathbb{N}^{N} $ such that $\tilde{\tr}\left(\left(  z_{k}+2\lambda_{k}\xi_{k}\right) \xi^{J}\right)=0$ and $j_{1}+j_{2}+\dots+j_{N}\geq 2$ }\right\},\quad\mbox{for all $k=1,\dots,N$.}\end{split} \label{788lala1}
\end{equation}

On the other hand, it is important to observe that 
$$\tilde{\mathcal{S}}= \mathcal{S},\hspace{0.1 cm}\tilde{\mathcal{T}}_{k}=\mathcal{T}_{k} ,\quad \mbox{for all $k=1,\dots, N$.}$$

These sets  (\ref{788lala}) are important  in order to   define  Fischer Normalization Spaces\cite{bu2} as previously, shaping suitably  normalizations in the local defining equations.   We proceed as follows:
 \section{Fischer $\left(G,\tilde{G}\right)$-Decompositions\cite{bu2},\cite{bu3}}For any multi-index $\tilde{I}\not\in \mathcal{S}$ of length $p$, we consider by (\ref{model}) the following Fischer Decomposition 
\begin{equation}\begin{split}& z^{\tilde{I}}=A\left(z,\xi\right) Q\left(z,\xi\right)+C\left(z,\xi\right),\hspace{0.1 cm} \tilde{\tr}\left(C\left(z,\xi\right)\right)=0,\quad\mbox{where  $\tilde{I}=\left(\tilde{i}_{1},\tilde{i}_{2},\dots,\tilde{i}_{N}\right)\in\mathbb{N}^{N}$,}\\& \xi^{\tilde{I}}=\tilde{A}\left(z,\xi\right) Q\left(z,\xi\right)+\tilde{C}\left(z,\xi\right),\hspace{0.1 cm} \tilde{\tr}\left(\tilde{C}\left(z,\xi\right)\right)=0,\quad\mbox{where  $\tilde{I}=\left(\tilde{i}_{1},\tilde{i}_{2},\dots,\tilde{i}_{N}\right)\in\mathbb{N}^{N}$.}\end{split}\label{oplala}
\end{equation}

These homogeneous polynomials $A\left(z,\xi\right)$, $\tilde{A}\left(z,\xi\right)$ and $C\left(z,\xi\right)$, $\tilde{C}\left(z,\xi\right)$ are uniquely determined according to Shapiro\cite{sh}. We compute them straightforwardly  in (\ref{op}),  writing   as follows
\begin{equation} A\left(z,\xi\right)=\displaystyle\sum_{I,J\in\mathbb{N}^{N}\atop\left|I\right|+\left|J\right|=p-2} a_{I;J}z^{I}\overline{z}^{J}, \quad C\left(z,\xi\right)=\displaystyle\sum_{I,J\in\mathbb{N}^{N}\atop\left|I\right|+\left|J\right|=p} c_{I;J}z^{I}\overline{z}^{J},\quad\quad\tilde{A}\left(z,\xi\right)=\displaystyle\sum_{I,J\in\mathbb{N}^{N}\atop\left|I\right|+\left|J\right|=p-2} \tilde{a}_{I;J}z^{I}\overline{z}^{J}, \quad \tilde{C}\left(z,\xi\right)=\displaystyle\sum_{I,J\in\mathbb{N}^{N}\atop\left|I\right|+\left|J\right|=p} \tilde{c}_{I;J}z^{I}\overline{z}^{J}.\label{oppvlala}
\end{equation}

Now, taking in consideration   the complexifications of (\ref{909se}), (\ref{B11}),  (\ref{909se11})    and (\ref{B111}) by applying (\ref{tracelala}) in (\ref{oplala}) in order to determine (\ref{oppvlala}) according to (\ref{yy}), we  obtain
\begin{equation}\begin{split}& \quad   \lambda_{1}\tilde{i}_{1}\left(\tilde{i}_{1}-1\right)z_{1}^{\tilde{i}_{1}-2}\dots z_{N}^{\tilde{i}_{N}}+\dots+\lambda_{k}\tilde{i}_{k}\left(\tilde{i}_{k}-1\right)z_{1}^{\tilde{i}_{1}}\dots z_{k}^{\tilde{i}_{k}-2}\dots z_{N}^{\tilde{i}_{N}}+\dots+\lambda_{N}\tilde{i}_{N}\left(\tilde{i}_{N}-1\right)z_{1}^{\tilde{i}_{1}}\dots z_{N}^{\tilde{i}_{N}-2} \\& \quad\quad\quad\quad\quad   \quad\quad\quad\quad\quad\quad\quad \quad\quad\quad\quad\quad   \quad\quad\quad\quad\quad   \begin{tabular}{l} \rotatebox[origin=c]{270}{$=$}\end{tabular} \\& \tilde{\tr}\left(A\left(z,\xi\right) \right)Q\left(z,\xi\right)+\displaystyle\sum_{k=1}^{N}\frac{\partial  }{\partial z_{k} }\left(A\left(z,\xi\right)\right)\left(z_{k}\left(1+4\lambda_{k}^{2}\right)+4\lambda_{k}\xi_{k} \right)+\displaystyle\sum_{k=1}^{N}\frac{\partial  }{\partial \xi_{k} }\left(A \left(z,\xi\right)\right)\left(\xi_{k}\left(1+4\lambda_{k}^{2}\right)+4\lambda_{k}z_{k}\right)\\& \quad\quad\quad\quad\quad   \quad\quad\quad\quad\quad\quad\quad \quad\quad\quad\quad\quad   \quad\quad\quad\quad\quad   \begin{tabular}{l} \rotatebox[origin=c]{270}{$+$}\end{tabular} \\&  \quad\quad\quad\quad\quad\quad\quad\quad\quad\quad\quad\quad\quad\quad\quad\quad A\left(z,\xi\right)\left(N+4\lambda_{1}^{2}+\dots+4\lambda_{N}^{2}\right),\end{split}\label{VVV1la1}
\end{equation} 
and respectively, we obtain
\begin{equation}\begin{split}& \quad\quad  \lambda_{1}\tilde{i}_{1}\left(\tilde{i}_{1}-1\right)\xi_{1}^{\tilde{i}_{1}-2}\dots \xi_{N}^{\tilde{i}_{N}}+\dots+\lambda_{k}\tilde{i}_{k}\left(\tilde{i}_{k}-1\right)\xi_{1}^{\tilde{i}_{1}}\dots \xi_{k}^{\tilde{i}_{k}-2}\dots \xi_{N}^{\tilde{i}_{N}}+\dots+\lambda_{N}\tilde{i}_{N}\left(\tilde{i}_{N}-1\right)\xi_{1}^{\tilde{i}_{1}}\dots \xi_{N}^{\tilde{i}_{N}-2} \\& \quad\quad\quad\quad\quad   \quad\quad\quad\quad\quad\quad\quad \quad\quad\quad\quad\quad   \quad\quad\quad\quad\quad   \begin{tabular}{l} \rotatebox[origin=c]{270}{$=$}\end{tabular} \\& \tilde{\tr}\left(\tilde{A}\left(z,\xi\right) \right)Q\left(z,\xi\right)+\displaystyle\sum_{k=1}^{N}\frac{\partial  }{\partial z_{k} }\left(\tilde{A}\left(z,\xi\right)\right)\left(z_{k}\left(1+4\lambda_{k}^{2}\right)+4\lambda_{k}\xi_{k} \right)+\displaystyle\sum_{k=1}^{N}\frac{\partial  }{\partial \xi_{k} }\left(\tilde{A}\left(z,\xi\right)\right)\left(\xi_{k}\left(1+4\lambda_{k}^{2}\right)+4\lambda_{k}z_{k}\right)\\& \quad\quad\quad\quad\quad   \quad\quad\quad\quad\quad\quad  \quad\quad\quad\quad\quad \quad\quad\quad  \quad\quad\quad  \begin{tabular}{l} \rotatebox[origin=c]{270}{$+$}\end{tabular} \\& \quad\quad\quad\quad\quad \quad\quad\quad\quad\quad\quad\quad\quad\quad\quad\quad  \tilde{A}\left(z,\xi\right)\left(N+4\lambda_{1}^{2}+\dots+4\lambda_{N}^{2}\right).\end{split}\label{VVV1la2}
\end{equation} 

On the other hand, we   have
 \begin{equation}\begin{split}& \quad\quad\quad\quad\quad\quad\quad\quad\quad\quad\quad\quad\quad\quad\quad\quad\quad\quad\quad\quad\quad\quad\quad\quad\quad \tilde{\tr}\left(A\left(z,\xi\right)\right)Q\left(z,\xi\right)\\& \quad\quad\quad\quad\quad   \quad\quad\quad\quad\quad\quad  \quad\quad\quad\quad\quad \quad\quad\quad  \quad\quad\quad\quad\quad\quad\quad\quad  \begin{tabular}{l} \rotatebox[origin=c]{270}{$=$}\end{tabular} \\&  \left(z_{1}\xi_{1}+\dots+z_{N}\xi_{N}+\lambda_{1}\left(z_{1}^{2}+\xi_{1}^{2}\right)+\dots+\lambda_{N}\left(z_{N}^{2}+\xi_{N}^{2}\right)\right)\left(\displaystyle\sum_{k=1}^{N}\displaystyle\sum_{I,J\in\mathbb{N}^{N}\atop\left|I\right|+\left|J\right|=p-2} a_{I;J}i_{k}j_{k} z_{1}^{i_{1}}\dots z_{k}^{i_{k}-1} \dots z_{N}^{i_{N}} \xi_{1}^{j_{1}}\dots  \xi_{k}^{j_{k}-1}\dots  \xi_{N}^{j_{N}}\right.\\& \quad\quad\quad\quad\quad   \quad\quad\quad\quad\quad\quad  \quad\quad\quad\quad\quad \quad\quad\quad  \quad\quad\quad\quad\quad\quad\quad\quad  \begin{tabular}{l} \rotatebox[origin=c]{270}{$+$}\end{tabular} \\&    \left.\displaystyle\sum_{k=1}^{N}\displaystyle\sum_{I,J\in\mathbb{N}^{N}\atop\left|I\right|+\left|J\right|=p-2} a_{I;J}\lambda_{k} i_{k}\left(i_{k}-1\right) z_{1}^{i_{1}}\dots z_{k}^{i_{k}-2}\dots z_{N}^{i_{N}} \xi_{1}^{j_{1}}\dots   \xi_{N}^{j_{N}}+\displaystyle\sum_{k=1}^{N}\displaystyle\sum_{I,J\in\mathbb{N}^{N}\atop\left|I\right|+\left|J\right|=p-2} a_{I;J}\lambda_{k} j_{k}\left(j_{k}-1\right) z_{1}^{i_{1}}\dots  z_{N}^{i_{N}}  \xi_{1}^{j_{1}}\dots  \xi_{k}^{j_{k}-2}\dots  \xi_{N}^{j_{N}}\right) .\end{split}
\label{B1la1}\end{equation}
and respectively, we obtain
 \begin{equation}\begin{split}& \quad\quad\quad\quad\quad\quad\quad\quad\quad\quad\quad\quad\quad\quad\quad\quad\quad\quad\quad\quad\quad\quad\quad\quad\quad \tilde{\tr}\left(\tilde{A}\left(z,\xi\right)\right)Q\left(z,\xi\right)\\& \quad\quad\quad\quad\quad   \quad\quad\quad\quad\quad\quad  \quad\quad\quad\quad\quad \quad\quad\quad  \quad\quad\quad\quad\quad\quad\quad\quad  \begin{tabular}{l} \rotatebox[origin=c]{270}{$=$}\end{tabular} \\&  \left(z_{1}\overline{z}_{1}+\dots+z_{N}\overline{z}_{N}+\lambda_{1}\left(z_{1}^{2}+\xi_{1}^{2}\right)+\dots+\lambda_{N}\left(z_{N}^{2}+\xi_{N}^{2}\right)\right)\left(\displaystyle\sum_{k=1}^{N}\displaystyle\sum_{I,J\in\mathbb{N}^{N}\atop\left|I\right|+\left|J\right|=p-2} \tilde{a}_{I;J}i_{k}j_{k} z_{1}^{i_{1}}\dots z_{k}^{i_{k}-1} \dots z_{N}^{i_{N}} \xi_{1}^{j_{1}}\dots  \xi_{k}^{j_{k}-1}\dots  \xi_{N}^{j_{N}} \right.\\& \quad\quad\quad\quad\quad   \quad\quad\quad\quad\quad\quad  \quad\quad\quad\quad\quad \quad\quad\quad  \quad\quad\quad\quad\quad\quad\quad\quad  \begin{tabular}{l} \rotatebox[origin=c]{270}{$+$}\end{tabular} \\&\left.\displaystyle\sum_{k=1}^{N}\displaystyle\sum_{I,J\in\mathbb{N}^{N}\atop\left|I\right|+\left|J\right|=p-2} \tilde{a}_{I;J}\lambda_{k} i_{k}\left(i_{k}-1\right) z_{1}^{i_{1}}\dots z_{k}^{i_{k}-2}\dots z_{N}^{i_{N}} \xi_{1}^{j_{1}}\dots   \xi_{N}^{j_{N}}+\displaystyle\sum_{k=1}^{N}\displaystyle\sum_{I,J\in\mathbb{N}^{N}\atop\left|I\right|+\left|J\right|=p-2} \tilde{a}_{I;J}\lambda_{k} j_{k}\left(j_{k}-1\right) z_{1}^{i_{1}}\dots  z_{N}^{i_{N}}  \xi_{1}^{j_{1}}\dots  \xi_{k}^{j_{k}-2}\dots  \xi_{N}^{j_{N}}\right) .\end{split}
\label{B1la2}\end{equation}
 
The interactions of homogeneous terms  are very complicated in the last two equations. We have thus to better organize all these terms  depending on their contributions at each  degree level in the last two equations, observing by (\ref{B111}) the following 

\begin{equation}\begin{split}&\quad\quad\quad\quad\quad\quad\quad\quad\quad\quad\quad\quad\quad\quad\quad\quad\quad\quad\quad\quad\quad\quad\hspace{0.1 cm}\tilde{\tr}\left(A\left(z,\xi\right)\right)Q\left(z,\xi\right)\\& \quad\quad\quad\quad\quad   \quad\quad\quad\quad\quad\quad  \quad\quad\quad\quad\quad \quad\quad\quad  \quad\quad\quad\quad\quad\quad\quad\quad  \begin{tabular}{l} \rotatebox[origin=c]{270}{$=$}\end{tabular} \\&  \left(z_{1}\xi_{1}+\dots+z_{N}\xi_{N}+\lambda_{1}\left(z_{1}^{2}+\xi_{1}^{2}\right)+\dots+\lambda_{N}\left(z_{N}^{2}+\xi_{N}^{2}\right)\right)\left( \displaystyle\sum_{I,J\in\mathbb{N}^{N}\atop\left|I\right|+\left|J\right|=p-2} b_{I;J} z_{1}^{i_{1}}\dots z_{k}^{i_{k}-1} \dots z_{N}^{i_{N}} \xi_{1}^{j_{1}}\dots \xi_{k}^{j_{k}-1}\dots \xi_{N}^{j_{N}} \right.\\& \quad\quad\quad\quad\quad   \quad\quad\quad\quad\quad\quad  \quad\quad\quad\quad\quad \quad\quad\quad  \quad\quad\quad\quad\quad\quad\quad\quad  \begin{tabular}{l} \rotatebox[origin=c]{270}{$+$}\end{tabular} \\& \quad\quad\quad\quad\quad  \left.  \displaystyle\sum_{I,J\in\mathbb{N}^{N}\atop\left|I\right|+\left|J\right|=p-2} {b'}_{I;J} z_{1}^{i_{1}}\dots z_{k}^{i_{k}-2}\dots z_{N}^{i_{N}}\xi_{1}^{j_{1}}\dots  \xi_{N}^{j_{N}}+ \displaystyle\sum_{I,J\in\mathbb{N}^{N}\atop\left|I\right|+\left|J\right|=p-2} {b''}_{I;J} z_{1}^{i_{1}}\dots  z_{N}^{i_{N}} \xi_{1}^{j_{1}}\dots \xi_{k}^{j_{k}-2}\dots \xi_{N}^{j_{N}}\right) ,\end{split}
\label{B111la1}\end{equation}
and respectively, the following 
 \begin{equation}\begin{split}&\quad\quad\quad\quad\quad\quad\quad\quad\quad\quad\quad\quad\quad\quad\quad\quad\quad\quad\quad\quad\quad\quad\hspace{0.1 cm}\tilde{\tr}\left(\tilde{A}\left(z,\xi\right)\right)Q\left(z,\xi\right)\\& \quad\quad\quad\quad\quad   \quad\quad\quad\quad\quad\quad  \quad\quad\quad\quad\quad \quad\quad\quad  \quad\quad\quad\quad\quad\quad\quad\quad  \begin{tabular}{l} \rotatebox[origin=c]{270}{$=$}\end{tabular} \\&  \left(z_{1}\xi_{1}+\dots+z_{N}\xi_{N}+\lambda_{1}\left(z_{1}^{2}+\xi_{1}^{2}\right)+\dots+\lambda_{N}\left(z_{N}^{2}+\xi_{N}^{2}\right)\right)\left( \displaystyle\sum_{I,J\in\mathbb{N}^{N}\atop\left|I\right|+\left|J\right|=p-2} \tilde{b}_{I;J} z_{1}^{i_{1}}\dots z_{k}^{i_{k}-1} \dots z_{N}^{i_{N}} \xi_{1}^{j_{1}}\dots \xi_{k}^{j_{k}-1}\dots \xi_{N}^{j_{N}} \right.\\& \quad\quad\quad\quad\quad   \quad\quad\quad\quad\quad\quad  \quad\quad\quad\quad\quad \quad\quad\quad  \quad\quad\quad\quad\quad\quad\quad\quad  \begin{tabular}{l} \rotatebox[origin=c]{270}{$+$}\end{tabular} \\& \quad\quad\quad\quad\quad  \left.  \displaystyle\sum_{I,J\in\mathbb{N}^{N}\atop\left|I\right|+\left|J\right|=p-2} {\tilde{b}'}_{I;J} z_{1}^{i_{1}}\dots z_{k}^{i_{k}-2}\dots z_{N}^{i_{N}}\xi_{1}^{j_{1}}\dots  \xi_{N}^{j_{N}}+ \displaystyle\sum_{I,J\in\mathbb{N}^{N}\atop\left|I\right|+\left|J\right|=p-2} {\tilde{b}''}_{I;J} z_{1}^{i_{1}}\dots  z_{N}^{i_{N}} \xi_{1}^{j_{1}}\dots \xi_{k}^{j_{k}-2}\dots \xi_{N}^{j_{N}}\right) ,\end{split}
\label{B111la2}\end{equation}
where we have used by (\ref{yy}) the following notations
\begin{equation} \begin{split}& {b}_{I;J}=\left(i_{1}+1\right)\left(j_{1}+1\right)a_{\left(i_{1}+1,\dots,i_{k}-1,\dots,i_{N};j_{1}+1,\dots,j_{k}-1,\dots,j_{N}\right)}+\dots +i_{k}j_{k}a_{\left(i_{1},\dots,i_{k},\dots,i_{N};j_{1},\dots,j_{k},\dots,j_{N}\right)}\\&\quad\quad\quad\quad\quad\quad\quad\quad\quad\quad\quad\quad\quad\quad\quad\quad\quad\quad\quad\quad\quad\quad\quad\quad\quad\quad\quad\quad\quad\quad\quad\quad\quad\quad\quad\quad\quad\quad\hspace{0.1 cm} +\\& \quad\quad\quad\quad\quad\quad\quad\quad\quad\quad\quad\quad\quad\quad\quad\quad\quad\quad\quad\quad\quad\quad\quad\quad\quad\quad\quad\quad\quad\quad\quad\quad\quad\quad\quad\quad\quad\quad\hspace{0.1 cm}\begin{tabular}{l} \rotatebox[origin=c]{270}{$\dots$}\end{tabular}\\&\quad\quad\quad\quad\quad\quad\quad\quad\quad\quad\quad\quad\quad\quad\quad\quad\quad\quad\quad\quad\quad\quad\quad\quad\quad\quad\quad\quad\quad\quad\quad\quad\quad\quad\quad\quad\quad\quad\hspace{0.1 cm} +\\&\quad\quad\quad\quad\quad\quad\quad\quad\quad\quad\quad\quad\quad\quad\quad\quad\quad\quad\quad\quad\quad\quad\quad\left(i_{N}+1\right)\left(j_{N}+1\right)a_{\left(i_{1},\dots,i_{k}-1,\dots,i_{N}+1;j_{1},\dots,j_{k}-1,\dots,j_{N}+1\right)},  \\&{b'}_{I;J}=\lambda_{1}  \left(i_{1}+1\right)\left(i_{1}+2\right) a_{\left(i_{1}+2,\dots,i_{k}-2,\dots,i_{N};j_{1},\dots,j_{N}\right)}+\dots
+\lambda_{k} i_{k}\left(i_{k}-1\right)a_{\left(i_{1},\dots,i_{k},\dots,i_{N};j_{1},\dots,j_{N}\right)}\\&\quad\quad\quad\quad\quad\quad\quad\quad\quad\quad\quad\quad\quad\quad\quad\quad\quad\quad\quad\quad\quad\quad\quad\quad\quad\quad\quad\quad\quad\quad\quad\quad\quad\quad\quad\quad\quad\quad\hspace{0.1 cm} +\\& \quad\quad\quad\quad\quad\quad\quad\quad\quad\quad\quad\quad\quad\quad\quad\quad\quad\quad\quad\quad\quad\quad\quad\quad\quad\quad\quad\quad\quad\quad\quad\quad\quad\quad\quad\quad\quad\quad\hspace{0.1 cm}\begin{tabular}{l} \rotatebox[origin=c]{270}{$\dots$}\end{tabular}\\&\quad\quad\quad\quad\quad\quad\quad\quad\quad\quad\quad\quad\quad\quad\quad\quad\quad\quad\quad\quad\quad\quad\quad\quad\quad\quad\quad\quad\quad\quad\quad\quad\quad\quad\quad\quad\quad\quad\hspace{0.1 cm} +\\&\quad\quad\quad\quad\quad\quad\quad\quad\quad\quad\quad\quad\quad\quad\quad\quad\quad\quad\quad\quad\quad\quad\quad\lambda_{N} \left(i_{N}+1\right)\left(i_{N}+2\right)a_{\left(i_{1},\dots,i_{k}-2,\dots,i_{N}+2;j_{1},\dots,j_{N}\right)},   \\&{b''}_{I;J}=\lambda_{1}\left(j_{1}+1\right)\left(j_{1}+2\right) a_{\left(i_{1},\dots,i_{N};j_{1}+2,\dots,j_{k}-2,\dots,j_{N}\right)}+ \dots+\lambda_{k}j_{k}\left(j_{k}-1\right) a_{\left(i_{1},\dots,i_{N};j_{1},\dots,j_{k}\dots,j_{N}\right)}\\&\quad\quad\quad\quad\quad\quad\quad\quad\quad\quad\quad\quad\quad\quad\quad\quad\quad\quad\quad\quad\quad\quad\quad\quad\quad\quad\quad\quad\quad\quad\quad\quad\quad\quad\quad\quad\quad\quad\hspace{0.1 cm} +\\& \quad\quad\quad\quad\quad\quad\quad\quad\quad\quad\quad\quad\quad\quad\quad\quad\quad\quad\quad\quad\quad\quad\quad\quad\quad\quad\quad\quad\quad\quad\quad\quad\quad\quad\quad\quad\quad\quad\hspace{0.1 cm}\begin{tabular}{l} \rotatebox[origin=c]{270}{$\dots$}\end{tabular}\\&\quad\quad\quad\quad\quad\quad\quad\quad\quad\quad\quad\quad\quad\quad\quad\quad\quad\quad\quad\quad\quad\quad\quad\quad\quad\quad\quad\quad\quad\quad\quad\quad\quad\quad\quad\quad\quad\quad\hspace{0.1 cm} +\\&\quad\quad\quad\quad\quad\quad\quad\quad\quad\quad\quad\quad\quad\quad\quad\quad\quad\quad\quad\quad\quad\quad\quad\lambda_{N}\left(j_{N}+1\right)\left(j_{N}+2\right)a_{\left(i_{1},\dots,i_{N};j_{1},\dots,j_{k}-2\dots,j_{N}+2\right)}.   \end{split} \label{99la1}
\end{equation}
and respectively, the following notations
\begin{equation} \begin{split}& {\tilde{b}}_{I;J}=\left(i_{1}+1\right)\left(j_{1}+1\right)\tilde{a}_{\left(i_{1}+1,\dots,i_{k}-1,\dots,i_{N};j_{1}+1,\dots,j_{k}-1,\dots,j_{N}\right)}+ \dots+i_{k}j_{k}\tilde{a}_{\left(i_{1},\dots,i_{k},\dots,i_{N};j_{1},\dots,j_{k},\dots,j_{N}\right)}\\&\quad\quad\quad\quad\quad\quad\quad\quad\quad\quad\quad\quad\quad\quad\quad\quad\quad\quad\quad\quad\quad\quad\quad\quad\quad\quad\quad\quad\quad\quad\quad\quad\quad\quad\quad\quad\quad\quad\hspace{0.1 cm} +\\& \quad\quad\quad\quad\quad\quad\quad\quad\quad\quad\quad\quad\quad\quad\quad\quad\quad\quad\quad\quad\quad\quad\quad\quad\quad\quad\quad\quad\quad\quad\quad\quad\quad\quad\quad\quad\quad\quad\hspace{0.1 cm}\begin{tabular}{l} \rotatebox[origin=c]{270}{$\dots$}\end{tabular}\\&\quad\quad\quad\quad\quad\quad\quad\quad\quad\quad\quad\quad\quad\quad\quad\quad\quad\quad\quad\quad\quad\quad\quad\quad\quad\quad\quad\quad\quad\quad\quad\quad\quad\quad\quad\quad\quad\quad\hspace{0.1 cm} +\\&\quad\quad\quad\quad\quad\quad\quad\quad\quad\quad\quad\quad\quad\quad\quad\quad\quad\quad\quad\quad\quad\quad\quad\left(i_{N}+1\right)\left(j_{N}+1\right)\tilde{a}_{\left(i_{1},\dots,i_{k}-1,\dots,i_{N}+1;j_{1},\dots,j_{k}-1,\dots,j_{N}+1\right)},  \\&{\tilde{b}'}_{I;J}=\lambda_{1}  \left(i_{1}+1\right)\left(i_{1}+2\right) \tilde{a}_{\left(i_{1}+2,\dots,i_{k}-2,\dots,i_{N};j_{1},\dots,j_{N}\right)}+\dots+\lambda_{k} i_{k}\left(i_{k}-1\right)\tilde{a}_{\left(i_{1},\dots,i_{k},\dots,i_{N};j_{1},\dots,j_{N}\right)}\\&\quad\quad\quad\quad\quad\quad\quad\quad\quad\quad\quad\quad\quad\quad\quad\quad\quad\quad\quad\quad\quad\quad\quad\quad\quad\quad\quad\quad\quad\quad\quad\quad\quad\quad\quad\quad\quad\quad\hspace{0.1 cm} +\\& \quad\quad\quad\quad\quad\quad\quad\quad\quad\quad\quad\quad\quad\quad\quad\quad\quad\quad\quad\quad\quad\quad\quad\quad\quad\quad\quad\quad\quad\quad\quad\quad\quad\quad\quad\quad\quad\quad\hspace{0.1 cm}\begin{tabular}{l} \rotatebox[origin=c]{270}{$\dots$}\end{tabular}\\&\quad\quad\quad\quad\quad\quad\quad\quad\quad\quad\quad\quad\quad\quad\quad\quad\quad\quad\quad\quad\quad\quad\quad\quad\quad\quad\quad\quad\quad\quad\quad\quad\quad\quad\quad\quad\quad\quad\hspace{0.1 cm} +\\&\quad\quad\quad\quad\quad\quad\quad\quad\quad\quad\quad\quad\quad\quad\quad\quad\quad\quad\quad\quad\quad\quad\quad\lambda_{N} \left(i_{N}+1\right)\left(i_{N}+2\right)\tilde{a}_{\left(i_{1},\dots,i_{k}-2,\dots,i_{N}+2;j_{1},\dots,j_{N}\right)},   \\&{\tilde{b}''}_{I;J}=\lambda_{1}\left(j_{1}+1\right)\left(j_{1}+2\right) \tilde{a}_{\left(i_{1},\dots,i_{N};j_{1}+2,\dots,j_{k}-2,\dots,j_{N}\right)}+\dots +\lambda_{k}j_{k}\left(j_{k}-1\right) \tilde{a}_{\left(i_{1},\dots,i_{N};j_{1},\dots,j_{k}\dots,j_{N}\right)}\\&\quad\quad\quad\quad\quad\quad\quad\quad\quad\quad\quad\quad\quad\quad\quad\quad\quad\quad\quad\quad\quad\quad\quad\quad\quad\quad\quad\quad\quad\quad\quad\quad\quad\quad\quad\quad\quad\quad\hspace{0.1 cm} +\\& \quad\quad\quad\quad\quad\quad\quad\quad\quad\quad\quad\quad\quad\quad\quad\quad\quad\quad\quad\quad\quad\quad\quad\quad\quad\quad\quad\quad\quad\quad\quad\quad\quad\quad\quad\quad\quad\quad\hspace{0.1 cm}\begin{tabular}{l} \rotatebox[origin=c]{270}{$\dots$}\end{tabular}\\&\quad\quad\quad\quad\quad\quad\quad\quad\quad\quad\quad\quad\quad\quad\quad\quad\quad\quad\quad\quad\quad\quad\quad\quad\quad\quad\quad\quad\quad\quad\quad\quad\quad\quad\quad\quad\quad\quad\hspace{0.1 cm} +\\&\quad\quad\quad\quad\quad\quad\quad\quad\quad\quad\quad\quad\quad\quad\quad\quad\quad\quad\quad\quad\quad\quad\quad\lambda_{N}\left(j_{N}+1\right)\left(j_{N}+2\right)\tilde{a}_{\left(i_{1},\dots,i_{N};j_{1},\dots,j_{k}-2\dots,j_{N}+2\right)}.   \end{split} \label{99la2}
\end{equation}
  
In order to simplify the computations in (\ref{VVV1la1}) and (\ref{VVV1la2}), we introduce by (\ref{yy}) and (\ref{99la1}) and (\ref{99la2}),  the following vectors

\begin{equation} \begin{split}&
X_{1}=X_{1}\left[I;J\right]=\begin{pmatrix}a_{\left(i_{1}+1,\dots,i_{k}-1,\dots,i_{N};j_{1}+1,\dots,j_{k}-1,\dots,j_{N}\right)} \\ \vdots \\ a_{\left(i_{1},\dots,i_{k},\dots,i_{N};j_{1},\dots,j_{k},\dots,j_{N}\right)} \\ \vdots \\  a_{\left(i_{1},\dots,i_{k}-1,\dots,i_{N}+1;j_{1},\dots,j_{k}-1,\dots,j_{N}+1\right)},\end{pmatrix},\quad \tilde{X}_{1}=\tilde{X}_{1}\left[I;J\right]=\begin{pmatrix}a_{\left(i_{1}+1,\dots,i_{k}-1,\dots,i_{N};j_{1}+1,\dots,j_{k}-1,\dots,j_{N}\right)} \\ \vdots \\ a_{\left(i_{1},\dots,i_{k},\dots,i_{N};j_{1},\dots,j_{k},\dots,j_{N}\right)} \\ \vdots \\  a_{\left(i_{1},\dots,i_{k}-1,\dots,i_{N}+1;j_{1},\dots,j_{k}-1,\dots,j_{N}+1\right)},\end{pmatrix}\\& \quad\quad\quad\quad\hspace{0.1 cm}    X_{2}=X_{2}\left[I;J\right]=\begin{pmatrix}a_{\left(i_{1}+2,\dots,i_{k}-2,\dots,i_{N};j_{1},\dots,j_{N}\right)}\\ \vdots \\ a_{\left(i_{1},\dots,i_{k},\dots,i_{N};j_{1},\dots,j_{N}\right)} \\ \vdots \\  a_{\left(i_{1},\dots,i_{k}-2,\dots,i_{N}+2;j_{1},\dots,j_{N}\right)}\end{pmatrix},\quad  X_{3}=X_{3}\left[I;J\right]=\begin{pmatrix}a_{\left(i_{1},\dots,i_{N};j_{1}+2,\dots,j_{k}-2,\dots,j_{N}\right)} \\ \vdots \\ a_{\left(i_{1},\dots,i_{N};j_{1},\dots,j_{k}\dots,j_{N}\right)} \\ \vdots \\  a_{\left(i_{1},\dots,i_{N};j_{1},\dots,j_{k}-2\dots,j_{N}+2\right)}\end{pmatrix},\\& \quad\quad\quad\quad\hspace{0.1 cm}    \tilde{X}_{2}=\tilde{X}_{2}\left[I;J\right]=\begin{pmatrix}a_{\left(i_{1}+2,\dots,i_{k}-2,\dots,i_{N};j_{1},\dots,j_{N}\right)}\\ \vdots \\ a_{\left(i_{1},\dots,i_{k},\dots,i_{N};j_{1},\dots,j_{N}\right)} \\ \vdots \\  a_{\left(i_{1},\dots,i_{k}-2,\dots,i_{N}+2;j_{1},\dots,j_{N}\right)}\end{pmatrix},\quad  \tilde{X}_{3}=\tilde{X}_{3}\left[I;J\right]=\begin{pmatrix}a_{\left(i_{1},\dots,i_{N};j_{1}+2,\dots,j_{k}-2,\dots,j_{N}\right)} \\ \vdots \\ a_{\left(i_{1},\dots,i_{N};j_{1},\dots,j_{k}\dots,j_{N}\right)} \\ \vdots \\  a_{\left(i_{1},\dots,i_{N};j_{1},\dots,j_{k}-2\dots,j_{N}+2\right)}\end{pmatrix}. \label{vectorila} 
\end{split} \end{equation}

The first very consistent sum of homogeneous terms  is obviously multiplied in  (\ref{B111la1}) and (\ref{B111la2}) by $$\mbox{$\lambda_{1}z_{1}^{2},\dots,\lambda_{N}z_{N}^{2}$, $z_{1}\xi_{1},\dots,z_{N}\xi_{N}$ and $\lambda_{1}\xi_{1}^{2},\dots,\lambda_{N}\xi_{N}^{2}$.}$$ 

This sum    generates obviously  by  (\ref{Lambda})  and (\ref{vectorila}) in (\ref{VVV1la1}) and (\ref{VVV1la2})  the following terms
\begin{equation}\mathcal{A}X_{1},\quad  \left(\Lambda\mathcal{A}\right)X_{1},\quad \left(\Lambda\mathcal{A} A\right)X_{1}, \quad\quad \mathcal{A}\tilde{X}_{1},\quad  \left(\Lambda\mathcal{A}\right)\tilde{X}_{1},\quad \left(\Lambda\mathcal{A} A\right)\tilde{X}_{1}.
  \label{idiot1lala}
\end{equation}
 
The second very consistent sum of homogeneous terms  is obviously multiplied in  (\ref{B111la1}) and (\ref{B111la2}) by $$\mbox{$\lambda_{1}z_{1}^{2},\dots,\lambda_{N}z_{N}^{2}$, $z_{1}\xi_{1},\dots,z_{N}\xi_{N}$ and $\lambda_{1}\xi_{1}^{2},\dots,\lambda_{N}\xi_{N}^{2}$.}$$

This sum    generates obviously  by  (\ref{Lambda})  and (\ref{vectorila}) in (\ref{VVV1la1}) and (\ref{VVV1la2})  the following terms
\begin{equation}
\left( \mathcal{A'}\Lambda\right)X_{2} ,\quad \left(\Lambda\mathcal{A'}\Lambda\right)X_{2} ,\quad \left(\Lambda\mathcal{A'}\Lambda\right)X_{2},\quad\quad \left( \mathcal{A'}\Lambda\right)\tilde{X}_{2} ,\quad \left(\Lambda\mathcal{A'}\Lambda\right)\tilde{X}_{2} ,\quad \left(\Lambda\mathcal{A'}\Lambda\right)\tilde{X}_{2}. \label{idiot2lala}
\end{equation}

The third very consistent sum of homogeneous terms  is obviously multiplied in  (\ref{B111la1}) and (\ref{B111la2}) by $$\mbox{$\lambda_{1}z_{1}^{2},\dots,\lambda_{N}z_{N}^{2}$, $z_{1}\xi_{1},\dots,z_{N}\xi_{N}$ and $\lambda_{1}\xi_{1}^{2},\dots,\lambda_{N}\xi_{N}^{2}$.}$$

This sum    generates obviously  by  (\ref{Lambda})  and (\ref{vectorila}) in (\ref{VVV1la1}) and (\ref{VVV1la2})  the following terms
\begin{equation}
\left(\ \mathcal{A''}\Lambda\right)X_{3} ,\quad \left(\Lambda\mathcal{A''}\Lambda\right)X_{3} ,\quad \left(\Lambda\mathcal{A''}\Lambda\right)X_{3},\quad\quad \left(\ \mathcal{A''}\Lambda\right)\tilde{X}_{3} ,\quad \left(\Lambda\mathcal{A''}\Lambda\right)\tilde{X}_{3} ,\quad \left(\Lambda\mathcal{A''}\Lambda\right)\tilde{X}_{3}. \label{idiot3lala}
\end{equation}

 The role of the  matrix $\Lambda$  is crucial   in order to  understand  the very complicated interactions of of homogeneous terms in (\ref{VVV1la1})  and (\ref{VVV1la2})  by solving non-trivial systems of equations: 

   \subsection{Systems of Equations}It is imposed  the standard lexicografic order corresponding to \begin{equation}\left(z_{1},z_{2},\dots,z_{N}, \xi_{1},\xi_{2},\dots,\xi_{N}\right)\in\mathbb{C}^{N}\times\mathbb{C}^{N},  
\label{ODINla}
\end{equation}   considering     the following vectors   
\begin{equation} Y_{1}^{t}=\left\{\left(a_{I;0}\right)\right\}_{I\in\mathbb{N}^{N}\atop \left|I\right|=p},\hspace{0.1 cm} Y_{2}^{t}=\left\{\left(a_{I;J}\right)\right\}_{I,J\in\mathbb{N}^{N}\atop \left|I\right|=p-1, \left|J\right|=1},\dots,  Y_{k}^{t}=\left\{\left(a_{I;J}\right)\right\}_{I,J\in\mathbb{N}^{N}\atop \left|I\right|=p-k+1, \left|J\right|=k-1},\dots, Y_{p+1}^{t}=\left\{\left(a_{0;J}\right)\right\}_{J\in\mathbb{N}^{N}\atop   \left|J\right|=p},\label{91Ala1}
\end{equation}
and respectively, the following vectors
\begin{equation} \tilde{Y}_{1}^{t}=\left\{\left(\tilde{a}_{I;0}\right)\right\}_{I\in\mathbb{N}^{N}\atop \left|I\right|=p},\hspace{0.1 cm} \tilde{Y}_{2}^{t}=\left\{\left(\tilde{a}_{I;J}\right)\right\}_{I,J\in\mathbb{N}^{N}\atop \left|I\right|=p-1, \left|J\right|=1},\dots,  \tilde{Y}_{k}^{t}=\left\{\left(\tilde{a}_{I;J}\right)\right\}_{I,J\in\mathbb{N}^{N}\atop \left|I\right|=p-k+1, \left|J\right|=k-1},\dots, \tilde{Y}_{p+1}^{t}=\left\{\left(\tilde{a}_{0;J}\right)\right\}_{J\in\mathbb{N}^{N}\atop   \left|J\right|=p},\label{91Ala2}
\end{equation}
where  $p=\tilde{i}_{1}+\dots+\tilde{i}_{N}$. 

More exactly, we construct   systems of equations by  extracting homogeneous terms in (\ref{VVV1la1}) and  (\ref{VVV1la2}) in order to determine the polynomials (\ref{oppvlala}) using (\ref{oplala}).   More precisely,   we construct systems of equations similar to (\ref{beb1})   by replacing (\ref{91A}) with (\ref{91Ala1}) and (\ref{91Ala2}). The  further computations are obvious, because there are recalled (\ref{lil}),(\ref{lili1}),(\ref{lili2}),(\ref{lili3}),(\ref{90000se1extra}),(\ref{suc1}),(\ref{suc2}), (\ref{suc33}),(\ref{90000})),(\ref{sucL1}),(\ref{sucL11}),(\ref{90000se}),(\ref{sucL1sese}),(\ref{sucL11sese}), (\ref{sucL111se}),(\ref{90000se1}) in the light of (\ref{shobi}),(\ref{hihi1}),(\ref{shobix1}),(\ref{hihi2}),(\ref{shobix2}),(\ref{hihi1x}),(\ref{shobiy1}),(\ref{hihi2x}),(\ref{shobiy2}). Then, we conclude   systems of equations similar to (\ref{beb120}), which are solved along    (\ref{gringo1}),(\ref{nene1}),(\ref{ggrin1}),(\ref{gringo2}),(\ref{nene2}),(\ref{gringo3}),(\ref{ggrin2}),(\ref{gringo4}),(\ref{nene3}),(\ref{gringo5}),(\ref{gringo4se}),(\ref{nene4}),(\ref{gringo44se}),(\ref{gringo44se1}),(\ref{gringoo1}),(\ref{gringo55}),(\ref{gringooo}),(\ref{nene5}),  (\ref{gringooo1}),(\ref{nene4441}),(\ref{gringooo2}),(\ref{nene4442}),(\ref{gringooo3}),(\ref{gringooo31}),(\ref{beb1se1}),(\ref{yx1}),(\ref{yx2}),(\ref{yx3}) according to (\ref{vulep1}),(\ref{tigru}),(\ref{girafa}),(\ref{2211}),(\ref{urs1}),(\ref{vulep2}),(\ref{asia}),(\ref{bib}),(\ref{liliwww1}),(\ref{indexes}),  (\ref{shobiV1}),(\ref{shobiV2}),\ref{shobiV11}),(\ref{shobiV21}),(\ref{detalii}),(\ref{viena1}),(\ref{cazuri}),(\ref{viena2}),(\ref{cazuri1}),(\ref{viena3}),(\ref{congo}),(\ref{bobo}),(\ref{yor}),(\ref{urs11}),(\ref{york}),(\ref{ioc}),\ref{expan}),(\ref{iocc}),(\ref{yor2}), (\ref{urs11se}),  (\ref{urs11sese}),(\ref{vulpe}),(\ref{9909}),(\ref{urs11sesese}). Therefore, such systems of equations have unique solutions according to the previous arguments. 

\subsection{Fischer Normalization $\left(G,\tilde{G}\right)$-Spaces\cite{bu4},\cite{bu5}}The Fischer Decomposition from (\ref{op}) gives
\begin{equation}\begin{split}& z^{I}=A_{I}\left(z,\xi\right) Q\left(z,\xi\right)+C_{I}\left(z,\xi\right),\quad\mbox{where}\hspace{0.1 cm} \tilde{\tr}\left(C_{I}\left(z,\xi\right))\right)=0,\\& \xi^{I}=\tilde{A}_{I}\left(z,\xi\right) Q\left(z,\xi\right)+\tilde{C}_{I}\left(z,\xi\right),\quad\mbox{where}\hspace{0.1 cm} \tilde{\tr}\left(\tilde{C}_{I}\left(z,\xi\right))\right)=0, \end{split}\label{optlala}
\end{equation}
where $I\not\in \mathcal{\tilde{S}}$ having   length $p$.
 
This set of homogeneous polynomials, derived from (\ref{optlala}),   gives   normalizations based on  certain   Normalization Double-Spaces, which is defined using the generalized version of  the Fischer Decomposition\cite{sh}: for a given  homogeneous polynomial of degree $p\geq 1$ in $\left(z,\xi\right)$ denoted by $P\left(z,\xi\right)$, we  have:

\begin{equation} \begin{split}& P\left(z,\xi\right)=P_{1}\left(z,\xi\right)Q\left(z,\xi\right)+R_{1}\left(z,\xi\right),\quad\mbox{where $\tr\left(R_{1}\left(z,\xi\right)\right)=0$ and:}\\& \quad\quad\quad\quad\hspace{0.1 cm}  R_{1}\left(z,\xi\right)=\displaystyle\sum_{I\in\mathbb{N}^{N}\atop {\left|I\right|=p}} \left(a_{I}C_{I}\left(z,\xi\right)+b_{I}\tilde{C}_{I}\left(z,\xi\right)\right)+R_{1,0}\left(z,\xi\right), \\&\quad\quad\quad\quad\quad\quad \quad\quad\quad R_{1,0}\left(z,\xi\right)\in \left(\displaystyle\bigcap_{I\in\mathbb{N}^{N},\hspace{0.1 cm} I\not\in \mathcal{S} \atop {\left|I\right|=p}}\left( \ker  C^{\star}_{I}  \cap  \ker  \tilde{C}^{\star}_{I} \right)\right) \bigcap\left( \displaystyle\bigcap_{I\in\mathbb{N}^{N},\hspace{0.1 cm} I \in \mathcal{S} \atop {\left|I\right|=p}}\left(\ker  \left(z^{I}\right)^{\star}   \cap  \ker  \left(\xi^{I}\right)^{\star}\right)\right) ,\quad\mbox{such that:}
\\& P_{1}\left(z,\xi\right)=P_{2}\left(z,\xi\right)Q\left(z,\xi\right)+R_{2}\left(z,\xi\right),\quad\mbox{where $\tr\left(R_{2}\left(z,\xi\right)\right)=0$ and:}\\& \quad\quad\quad\quad\hspace{0.1 cm}  R_{2}\left(z,\xi\right)=\displaystyle\sum_{I\in\mathbb{N}^{N}\atop {\left|I\right|=p-2}} \left(a_{I}C_{I}\left(z,\xi\right) +b_{I}\tilde{C}_{I}\left(z,\xi\right) \right)+R_{2,0}\left(z,\xi\right),\quad\mbox{such that:}\\&\quad\quad\quad\quad \quad\quad \quad\quad\quad R_{2,0}\left(z,\xi\right)\in \left(\displaystyle\bigcap_{I\in\mathbb{N}^{N},\hspace{0.1 cm} I\not\in \mathcal{S} \atop {\left|I\right|=p-2}}\left( \ker  C^{\star}_{I}  \cap  \ker  \tilde{C}^{\star}_{I}\right) \right) \bigcap \left(\displaystyle\bigcap_{I\in\mathbb{N}^{N},\hspace{0.1 cm} I \in \mathcal{S} \atop {\left|I\right|=p-2}}\left(\ker  \left(z^{I}\right)^{\star}   \cap  \ker  \left(\xi^{I}\right)^{\star}\right)\right),\\&\quad\quad\quad\quad\vdots\quad\quad\quad\quad\quad\quad\quad\quad\vdots\quad\quad\quad\quad\quad\quad\quad\quad\vdots\quad\quad\quad\quad\quad\quad\quad\quad\vdots \quad\quad\quad\quad\quad\quad\quad\quad\quad\quad\quad\quad\quad\quad\vdots
\\& P_{k}\left(z,\xi\right)=P_{k+1}\left(z,\xi\right)Q\left(z,\xi\right)+R_{k+1}\left(z,\xi\right),\quad\mbox{where $\tr\left(R_{k+1}\left(z,\xi\right)\right)=0$ and:}\\& \quad\quad\quad \quad \hspace{0.1 cm} R_{k+1}\left(z,\xi\right)=\displaystyle\sum_{I\in\mathbb{N}^{N}\atop {\left|I\right|=p-2k}} \left(a_{I}C_{I}\left(z,\xi\right)+b_{I}\tilde{C}_{I}\left(z,\xi\right) \right)+R_{k+1,0} \left(z,\xi\right),\quad\mbox{such that:}\\&\quad\quad\quad\quad \quad\quad \quad\quad\quad R_{k+1,0}\left(z,\xi\right)\in \left(\displaystyle\bigcap_{I\in\mathbb{N}^{N}, \hspace{0.1 cm}I\not\in \mathcal{S} \atop {\left|I\right|=p-2k}} \left(\ker  C^{\star}_{I}  \cap  \ker  \tilde{C}^{\star}_{I}\right)\right)\bigcap\left(   \displaystyle\bigcap_{I\in\mathbb{N}^{N},\hspace{0.1 cm} I \in \mathcal{S} \atop {\left|I\right|=p-2k}}\left(\ker  \left(z^{I}\right)^{\star}   \cap  \ker \left(\xi^{I}\right)^{\star}\right) \right),\\&\quad\quad\quad\quad\vdots\quad\quad\quad\quad\quad\quad\quad\quad\vdots\quad\quad\quad\quad\quad\quad\quad\quad\vdots\quad\quad\quad\quad\quad\quad\quad\quad\vdots \quad\quad\quad\quad\quad\quad\quad\quad\quad\quad\quad\quad\quad\quad\vdots\end{split} \label{new1lala1}
\end{equation}
where   these occurring polynomials 
\begin{equation}\left\{P_{k}\left(z,\xi\right)\right\}_{k=1,\dots,\left[\frac{p}{2}\right]},\quad \left\{R_{k}\left(z,\xi\right)\right\}_{k=1,\dots,\left[\frac{p}{2}\right]},
\label{poll1lala1}
\end{equation}
are iteratively obtained  using the generalized version  of the Fischer Decomposition\cite{sh}.  

Recalling  the previous strategies, we define the analogue of the Normalization Space \ref{spartiuG} as follows
\begin{equation}\mathcal{G}_{p},\quad p\in\mathbb{N}^{\star},\label{spartiuGlala1}
\end{equation} 
which consist in real-valued polynomials $P(z,\xi)$  of degree $p\geq 1$ in $\left(z,\xi\right))$ satisfying the normalizations:
$$ P_{k}^{\left(p\right)}\left(z,\xi\right)=P_{k+1}^{\left(p\right)}\left(z,\xi\right)Q\left(z,\xi\right)+R_{k+1}^{\left(p\right)}\left(z,\xi\right),\hspace{0.1 cm}\mbox{for all  $k=0,\dots, \left[\frac{p}{2}\right]$ and given $P_{0}^{\left(p\right)}(z,\overline{z})=P(z,\overline{z})$,}  $$ 
such that 
$$ R_{k+1,0}^{\left(p\right)}\left(z,\xi\right)\in \left(\displaystyle\bigcap_{I\in\mathbb{N}^{N},\hspace{0.1 cm}I\not\in\mathcal{S}\atop {\left|I\right|=p-2k}}\left( \ker  C^{\star}_{I}  \cap  \ker  \tilde{C}^{\star}_{I}\cap \ker  \tr\right)\right)\bigcap\left(\displaystyle\bigcap_{I\in\mathbb{N}^{N},\hspace{0.1 cm} I \in \mathcal{S} \atop {\left|I\right|=p-2k}}\left(\ker  \left(z^{I}\right)^{\star}   \bigcap   \ker \left(\xi^{I}\right)^{\star}\right)\right),\hspace{0.1 cm}\mbox{for all  $k=0,\dots, \left[\frac{p}{2}\right]$.}$$

We consider the Fischer Decompositions (\ref{new1lala1}) choosing 
\begin{equation}P\left(z,\xi\right)=\frac{\varphi_{p}\left(z,\xi\right)-\overline{\varphi}_{p} \left( \xi,z\right) }{2\sqrt{-1}},\quad\mbox{for given $p\in\mathbb{N}^{\star}$.}\label{kama1lala1}\end{equation}

The   $\left(G,\tilde{G}\right)$-components of the formal transformation (\ref{maplala}) are  iteratively computed and uniquely determined by (\ref{spartiuGlala1}) according to (\ref{new1lala1}). 

It remains to prove   the linear independence, considering complex numbers, of the following set of polynomials 
\begin{equation}\left\{C_{I}\left(z,\xi\right),\hspace{0.1 cm} \tilde{C}_{I}\left(z,\xi\right) \right\}_{I\in\mathbb{N}^{N}\atop {{\left|I\right|=p}\atop I\not\in\mathcal{S}}},\quad\mbox{for all $p\in\mathbb{N}^{\star}$.}\label{330lala1}
\end{equation}

Again, the  computations are non-trivial because of  the overlapping of the above  homogeneous   polynomials from (\ref{330lala1}). It is thus necessary to return to (\ref{beb1se1}) in order to study more carefully the computations of their solution recalling (\ref{yx1}),
(\ref{yx2}) and (\ref{yx3}) together with  other relevant conclusions related to (\ref{beb120}) exactly as previously.  The non-triviality of the Fischer Decomposition\cite{sh} forces to consider similar approaches. In particular, we consider by (\ref{IJ1}) the following vectors
\begin{equation*} Z\left[I \right]=\begin{pmatrix}a_{\left(i_{1}-2,\dots,i_{k},\dots,i_{N}\right) }\\ \vdots \\ a_{\left(i_{1},\dots,i_{k}-2,\dots,i_{N}\right) }\\ \vdots \\  a_{\left(i_{1},\dots,i_{k},\dots,i_{N}-2\right) }\end{pmatrix},\quad \xi\left[I \right]=\begin{pmatrix}\tilde{a}_{\left(i_{1}-2,\dots,i_{k},\dots,i_{N}\right) }\\ \vdots \\ \tilde{a}_{\left(i_{1},\dots,i_{k}-2,\dots,i_{N}\right) }\\ \vdots \\  \tilde{a}_{\left(i_{1},\dots,i_{k},\dots,i_{N}-2\right) }\end{pmatrix}.
\end{equation*}

Immediately from (\ref{optlala}), we obtain
\begin{equation}\begin{split}& z^{I}-A_{I}\left(z,\xi\right)  Q\left(z,\xi\right)=C_{I}\left(z,\xi\right),\quad\mbox{for all $I\not\in \mathcal{S}$ having   length $p$,}\\& \xi^{I}-\tilde{A}_{I}\left(z,\xi\right)  Q\left(z,\xi\right)=\tilde{C}_{I}\left(z,\xi\right),\quad\mbox{for all $I\not\in \mathcal{S}$ having   length $p$.}\end{split}\label{opt1lala}
\end{equation} 

Now, in order to uniquely solve (\ref{opt1lala}), we define
  $$Z^{t}=\left\{a_{I}\right\}_{I\in\mathbb{N}^{N}\atop \left|I\right|=p},\quad \xi=\left\{\tilde{a}_{I}\right\}_{I\in\mathbb{N}^{N}\atop \left|I\right|=p},$$
according to the lexicografic order related to (\ref{ODINla}) recalling (\ref{aux1}),(\ref{calcan1}) and (\ref{calcan2}).

We obtain the following system    of  equations 
\begin{equation}\begin{split}&\left(I_{N^{p}}-\mbox{Aux}_{p}A\right) Z+B\xi=W\left(z_{1},z_{2},\dots,z_{N}\right),\\& \left(I_{N^{p}}-\mbox{Aux}_{p}A\right) \xi+BZ=\tilde{W}\left(\xi_{1},\xi_{2},\dots,\xi_{N}\right),\end{split}\label{sisilala}
\end{equation} 
where $A$, $B\in \mathcal{M}_{N^{p}\times N^{p}}\left(\mathbb{C}\right)$ and $W\left(z_{1},z_{2},\dots,z_{N}\right)$ and $\tilde{W}\left(\xi_{1},\xi_{2},\dots,\xi_{N}\right)$ are known vector homogeneous polynomials  of degree $p-2$, which obviously has unique solution in the light of existences of the matrices defined in (\ref{550}) according to (\ref{iok1}) and (\ref{iok2}).  

The   computations of the $\left(G,\tilde{G}\right)$-components of the formal equivalence (\ref{maplala}) are concluded therefore by considering products of matrices as in  (\ref{90000se1extra}), (\ref{90000}), (\ref{90000se}), (\ref{90000se1}) in (\ref{sisilala}) using  the product of simple matrices (\ref{calcan2}) as solutions of the systems of equations like (\ref{sisilala}).  

We move forward:

\section{Fischer  $\left(F,\tilde{F}\right)$-Decompositions\cite{bu4},\cite{bu5}}
For any multi-index $$\mbox{$\tilde{J}=\left(\tilde{j}_{1},\tilde{j}_{2},\dots,\tilde{j}_{N}\right)\not\in \tilde{\mathcal{T}}_{l}$, for all   $l\in 1,\dots, N$}$$   we consider by (\ref{788}) the following Fischer Decompositions 
\begin{equation}\begin{split}& \left(\xi_{l}+2\lambda_{l}z_{l}\right) z^{\tilde{J}}=A_{l}\left(z,\xi\right) Q\left(z,\xi\right)+C_{l}\left(z,\xi\right),\quad \tilde{\tr}\left(C_{l}\left(z,\xi\right)\right)=0,\quad\mbox{for all   $l\in 1,\dots, N$.}\\& \left(z_{l}+2\lambda_{l}\xi_{l}\right) z^{\tilde{J}}=\tilde{A}_{l}\left(z,\xi\right) Q\left(z,\xi\right)+\tilde{C}_{l}\left(z,\xi\right),\quad \tilde{\tr}\left(\tilde{C}_{l}\left(z,\xi\right)\right)=0,\quad\mbox{for all   $l\in 1,\dots, N$.}\end{split}  \label{opseclala} 
\end{equation}

These homogeneous polynomials
\begin{equation*}
\begin{split}& A_{1}\left(z,\xi\right),A_{2}\left(z,\xi\right),\dots,A_{N}\left(z,\xi\right), \tilde{A}_{1}\left(z,\xi\right),\tilde{A}_{2}\left(z,\xi\right),\dots,\tilde{A}_{N}\left(z,\xi\right);\\& C_{1}\left(z,\xi\right),C_{2}\left(z,\xi\right),\dots,C_{N}\left(z,\xi\right), \tilde{C}_{1}\left(z,\xi\right),\tilde{C}_{2}\left(z,\xi\right),\dots,\tilde{C}_{N}\left(z,\xi\right), \end{split}
\end{equation*}
are uniquely determined according to Shapiro\cite{sh}, being directly  computed  from (\ref{opsec}) writing   as follows
\begin{equation}\begin{split}& A_{l}\left(z,\xi\right)=\displaystyle\sum_{I,J\in\mathbb{N}^{N}\atop\left|I\right|+\left|J\right|=p-2} a_{I,J}^{\left(l\right)}z^{I}\xi^{J}, \quad \tilde{A}_{l}\left(z,\xi\right)=\displaystyle\sum_{I,J\in\mathbb{N}^{N}\atop\left|I\right|+\left|J\right|=p-2} \tilde{a}_{I,J}^{\left(l\right)}z^{I}\xi^{J},\quad\mbox{for all   $l\in 1,\dots, N$,} \\& C_{l}\left(z,\xi\right)=\displaystyle\sum_{I,J\in\mathbb{N}^{N}\atop\left|I\right|+\left|J\right|=p} c_{I,J}^{\left(l\right)}z^{I}\xi^{J},\quad \tilde{C}_{l}\left(z,\xi\right)=\displaystyle\sum_{I,J\in\mathbb{N}^{N}\atop\left|I\right|+\left|J\right|=p} \tilde{c}_{I,J}^{\left(l\right)}z^{I}\xi^{J} ,\quad\mbox{for all   $l\in 1,\dots, N$.} \end{split}\label{oppvseclala}
\end{equation}

We apply now the operator $\tilde{tr}$ in (\ref{opseclala}) using (\ref{oppvseclala}) and recalling (\ref{909se11}). We obtain 

\begin{equation}\begin{split}& \lambda_{1}^{2}\tilde{j}_{1}\left(\tilde{j}_{1}-1\right)z_{1}^{\tilde{j}_{1}-2}\dots z_{l}^{\tilde{j}_{l}+1} \dots z_{N}^{\tilde{j}_{N}}+\dots+\lambda_{l}\lambda_{1}\tilde{j}_{l}\left(\tilde{j}_{l}+1\right)z_{1}^{\tilde{j}_{1}}\dots z_{l}^{\tilde{j}_{l}-1} \dots z_{N}^{\tilde{j}_{N}}+\dots+\lambda_{N}\lambda_{1}\tilde{j}_{N}\left(\tilde{j}_{N}-1\right)z_{1}^{\tilde{j}_{1}}\dots z_{l}^{\tilde{j}_{l}+1} \dots z_{N}^{\tilde{j}_{N}-2} \\&      \xi_{l}\left(\lambda_{1} \tilde{j}_{1}\left(\tilde{j}_{1}-1\right)z_{1}^{\tilde{j}_{1}-2}  \dots z_{N}^{\tilde{j}_{N}}+\dots+ \lambda_{N}\tilde{j}_{N}\left(\tilde{j}_{N}-1\right)z_{1}^{\tilde{j}_{1}}  \dots z_{N}^{\tilde{j}_{N}-2}\right)+j_{l} z_{1}^{\tilde{j}_{1}}\dots  z_{l}^{\tilde{j}_{l}-1} \dots z_{N}^{\tilde{j}_{N}}  \\& \quad\quad\quad\quad\quad   \quad\quad\quad\quad\quad\quad  \quad\quad\quad\quad\quad \quad\quad\quad  \quad\quad\quad\quad\quad\quad\quad\quad  \begin{tabular}{l} \rotatebox[origin=c]{270}{$=$}\end{tabular} \\&\tilde{\tr}\left(A_{l}\left(z,\xi\right) \right)Q\left(z,\xi\right)+\displaystyle\sum_{k=1}^{N}\left(A_{l}\right)_{z_{k}}\left(z,\xi\right)\left(z_{k}\left(1+4\lambda_{k}^{2}\right)+4\lambda_{k}\overline{z}_{k} \right)+\displaystyle\sum_{k=1}^{N}\left(A_{l}\right)_{\overline{z}_{k}}\left(z,\xi\right)\left(\overline{z}_{k}\left(1+4\lambda_{k}^{2}\right)+4\lambda_{k}z_{k}\right) +A_{l}\left(z,\xi\right)\left(N+4\displaystyle\sum_{k=1}^{N}\lambda_{k}^{2}\right),\end{split}\label{VVV1selala1}
\end{equation} 
and respectively 
\begin{equation}\begin{split}& \lambda_{1}^{2}\tilde{j}_{1}\left(\tilde{j}_{1}-1\right)\xi_{1}^{\tilde{j}_{1}-2}\dots z_{l}^{\tilde{j}_{l}+1} \dots \xi_{N}^{\tilde{j}_{N}}+\dots+\lambda_{l}\lambda_{1}\tilde{j}_{l}\left(\tilde{j}_{l}+1\right)\xi_{1}^{\tilde{j}_{1}}\dots \xi_{l}^{\tilde{j}_{l}-1} \dots \xi_{N}^{\tilde{j}_{N}}+\dots+\lambda_{N}\lambda_{1}\tilde{j}_{N}\left(\tilde{j}_{N}-1\right)\xi_{1}^{\tilde{j}_{1}}\dots z_{l}^{\tilde{j}_{l}+1} \dots \xi_{N}^{\tilde{j}_{N}-2} \\&      z_{l}\left(\lambda_{1} \tilde{j}_{1}\left(\tilde{j}_{1}-1\right)\xi_{1}^{\tilde{j}_{1}-2}  \dots \xi_{N}^{\tilde{j}_{N}}+\dots+ \lambda_{N}\tilde{j}_{N}\left(\tilde{j}_{N}-1\right)\xi_{1}^{\tilde{j}_{1}}  \dots \xi_{N}^{\tilde{j}_{N}-2}\right)+j_{l} \xi_{1}^{\tilde{j}_{1}}\dots  \xi_{l}^{\tilde{j}_{l}-1} \dots \xi_{N}^{\tilde{j}_{N}}  \\& \quad\quad\quad\quad\quad   \quad\quad\quad\quad\quad\quad  \quad\quad\quad\quad\quad \quad\quad\quad  \quad\quad\quad\quad\quad\quad\quad\quad  \begin{tabular}{l} \rotatebox[origin=c]{270}{$=$}\end{tabular} \\&\tilde{\tr}\left(A_{l}\left(z,\xi\right) \right)Q\left(z,\xi\right)+\displaystyle\sum_{k=1}^{N}\left(A_{l}\right)_{z_{k}}\left(z,\xi\right)\left(z_{k}\left(1+4\lambda_{k}^{2}\right)+4\lambda_{k}\overline{z}_{k} \right)+\displaystyle\sum_{k=1}^{N}\left(A_{l}\right)_{\overline{z}_{k}}\left(z,\xi\right)\left(\overline{z}_{k}\left(1+4\lambda_{k}^{2}\right)+4\lambda_{k}z_{k}\right) +A_{l}\left(z,\xi\right)\left(N+4\displaystyle\sum_{k=1}^{N}\lambda_{k}^{2}\right),\end{split}\label{VVV1selala2}
\end{equation}  
for all $l=1,\dots,N$. 

On the other hand, we have
 \begin{equation}\begin{split}& \quad\quad\quad\quad\quad\quad\quad\quad\quad\quad\quad\quad\quad\quad\quad\quad\quad\quad\quad\quad\quad\quad\quad\quad\quad \tilde{\tr}\left(A_{l}\left(z,\xi\right)\right)Q\left(z,\xi\right)\\& \quad\quad\quad\quad\quad   \quad\quad\quad\quad\quad\quad  \quad\quad\quad\quad\quad \quad\quad\quad  \quad\quad\quad\quad\quad\quad\quad\quad  \begin{tabular}{l} \rotatebox[origin=c]{270}{$=$}\end{tabular} \\& \left(z_{1}\xi_{1}+\dots+z_{N}\xi_{N}+\lambda_{1}\left(z_{1}^{2}+\xi_{1}^{2}\right)+\dots+\lambda_{N}\left(z_{N}^{2}+\xi_{N}^{2}\right)\right)\left(\displaystyle\sum_{k=1}^{N}\displaystyle\sum_{I,J\in\mathbb{N}^{N}\atop\left|I\right|+\left|J\right|=p-2} a_{I,J}^{\left(l\right)}i_{k}j_{k} z_{1}^{i_{1}}\dots z_{k}^{i_{k}-1} \dots z_{N}^{i_{N}} \xi_{1}^{j_{1}}\dots \xi_{k}^{j_{k}-1}\dots \xi_{N}^{j_{N}} \right.\\& \quad\quad\quad\quad\quad   \quad\quad\quad\quad\quad\quad  \quad\quad\quad\quad\quad \quad\quad\quad  \quad\quad\quad\quad\quad\quad\quad\quad  \begin{tabular}{l} \rotatebox[origin=c]{270}{$+$}\end{tabular} \\&\left.\displaystyle\sum_{k=1}^{N}\displaystyle\sum_{I,J\in\mathbb{N}^{N}\atop\left|I\right|+\left|J\right|=p-2} a_{I,J}^{\left(l\right)}\lambda_{k} i_{k}\left(i_{k}-1\right) z_{1}^{i_{1}}\dots z_{k}^{i_{k}-2}\dots z_{N}^{i_{N}}\xi_{1}^{j_{1}}\dots  \xi_{N}^{j_{N}}+\displaystyle\sum_{k=1}^{N}\displaystyle\sum_{I,J\in\mathbb{N}^{N}\atop\left|I\right|+\left|J\right|=p-2} a_{I,J}^{\left(l\right)}\lambda_{k} j_{k}\left(j_{k}-1\right) z_{1}^{i_{1}}\dots  z_{N}^{i_{N}} \xi_{1}^{j_{1}}\dots \xi_{k}^{j_{k}-2}\dots \xi_{N}^{j_{N}}\right) ,\end{split}
\label{B1selala1}\end{equation}
and respectively
\begin{equation}\begin{split}& \quad\quad\quad\quad\quad\quad\quad\quad\quad\quad\quad\quad\quad\quad\quad\quad\quad\quad\quad\quad\quad\quad\quad\quad\quad \tilde{\tr}\left(\tilde{A}_{l}\left(z,\xi\right)\right)Q\left(z,\xi\right)\\& \quad\quad\quad\quad\quad   \quad\quad\quad\quad\quad\quad  \quad\quad\quad\quad\quad \quad\quad\quad  \quad\quad\quad\quad\quad\quad\quad\quad  \begin{tabular}{l} \rotatebox[origin=c]{270}{$=$}\end{tabular} \\&  \left(z_{1}\xi_{1}+\dots+z_{N}\xi_{N}+\lambda_{1}\left(z_{1}^{2}+\xi_{1}^{2}\right)+\dots+\lambda_{N}\left(z_{N}^{2}+\xi_{N}^{2}\right)\right)\left(\displaystyle\sum_{k=1}^{N}\displaystyle\sum_{I,J\in\mathbb{N}^{N}\atop\left|I\right|+\left|J\right|=p-2} \tilde{a}_{I,J}^{\left(l\right)}i_{k}j_{k} z_{1}^{i_{1}}\dots z_{k}^{i_{k}-1} \dots z_{N}^{i_{N}} \xi_{1}^{j_{1}}\dots \xi_{k}^{j_{k}-1}\dots \xi_{N}^{j_{N}} \right.\\& \quad\quad\quad\quad\quad   \quad\quad\quad\quad\quad\quad  \quad\quad\quad\quad\quad \quad\quad\quad  \quad\quad\quad\quad\quad\quad\quad\quad  \begin{tabular}{l} \rotatebox[origin=c]{270}{$+$}\end{tabular} \\&\left.\displaystyle\sum_{k=1}^{N}\displaystyle\sum_{I,J\in\mathbb{N}^{N}\atop\left|I\right|+\left|J\right|=p-2} \tilde{a}_{I,J}^{\left(l\right)}\lambda_{k} i_{k}\left(i_{k}-1\right) z_{1}^{i_{1}}\dots z_{k}^{i_{k}-2}\dots z_{N}^{i_{N}}\xi_{1}^{j_{1}}\dots  \xi_{N}^{j_{N}}+\displaystyle\sum_{k=1}^{N}\displaystyle\sum_{I,J\in\mathbb{N}^{N}\atop\left|I\right|+\left|J\right|=p-2} \tilde{a}_{I,J}^{\left(l\right)}\lambda_{k} j_{k}\left(j_{k}-1\right) z_{1}^{i_{1}}\dots  z_{N}^{i_{N}} \xi_{1}^{j_{1}}\dots \xi_{k}^{j_{k}-2}\dots \xi_{N}^{j_{N}}\right) ,\end{split}
\label{B1selala2}\end{equation}
for all $l=1,\dots,N$. 

Clearly, the above interactions of terms are very complicated, but we can use  the previous strategy according to the lexicografic order related to (\ref{ODINla}) using the $\left(F,\tilde{F}\right)$-components of the formal equivalence (\ref{maplala}). More precisely, we  extract homogeneous terms in (\ref{VVV1selala1}) and (\ref{VVV1selala2})  considering   the following vectors
\begin{equation}\begin{split}&  W_{1}^{t}=\left\{\left(a_{I,0}^{\left(1\right)}\right),\left(a_{I,0}^{(2)}\right),\dots,\left(a_{I,0}^{\left(N\right)}\right)\right\}_{I\in\mathbb{N}^{N}\atop \left|I\right|=p},\\& W_{2}^{t}=\left\{\left(a_{I,J}^{\left(1\right)}\right),\left(a_{I,J}^{(2)}\right),\dots,\left(a_{I,J}^{\left(N\right)}\right)\right\}_{I,J\in\mathbb{N}^{N}\atop \left|I\right|=p-1, \left|J\right|=1}\\& \quad\quad\quad\quad\quad\quad\quad\vdots \\& W_{k}^{t}=\left\{\left(a_{I,J}^{\left(1\right)}\right),\left(a_{I,J}^{(2)}\right),\dots,\left(a_{I,J}^{\left(N\right)}\right)\right\}_{I,J\in\mathbb{N}^{N}\atop \left|I\right|=p-k+1, \left|J\right|=k-1}\\& \quad\quad\quad\quad\quad\quad\quad\vdots \\& W_{p+1}^{t}=\left\{\left(a_{0,J}^{\left(1\right)}\right),\left(a_{0,J}^{(2)}\right),\dots,\left(a_{0,J}^{\left(N\right)}\right)\right\}_{J\in\mathbb{N}^{N}\atop   \left|J\right|=p},\end{split}\quad\quad \begin{split}&  \tilde{W}_{1}^{t}=\left\{\left(\tilde{a}_{I,0}^{\left(1\right)}\right),\left(\tilde{a}_{I,0}^{(2)}\right),\dots,\left(\tilde{a}_{I,0}^{\left(N\right)}\right)\right\}_{I\in\mathbb{N}^{N}\atop \left|I\right|=p},\\&\tilde{W}_{2}^{t}=\left\{\left(a_{I,J}^{\left(1\right)}\right),\left(\tilde{a}_{I,J}^{(2)}\right),\dots,\left(\tilde{a}_{I,J}^{\left(N\right)}\right)\right\}_{I,J\in\mathbb{N}^{N}\atop \left|I\right|=p-1, \left|J\right|=1}\\& \quad\quad\quad\quad\quad\quad\quad\vdots \\& \tilde{W}_{k}^{t}=\left\{\left(\tilde{a}_{I,J}^{\left(1\right)}\right),\left(\tilde{a}_{I,J}^{(2)}\right),\dots,\left(\tilde{a}_{I,J}^{(N)}\right)\right\}_{I,J\in\mathbb{N}^{N}\atop \left|I\right|=p-k+1, \left|J\right|=k-1}\\& \quad\quad\quad\quad\quad\quad\quad\vdots \\& \tilde{W}_{p+1}^{t}=\left\{\left(\tilde{a}_{0,J}^{\left(1\right)}\right),\left(\tilde{a}_{0,J}^{(2)}\right),\dots,\left(\tilde{a}_{0,J}^{\left(N\right)}\right)\right\}_{J\in\mathbb{N}^{N}\atop   \left|J\right|=p},\end{split}
\label{91ASeelala1}
\end{equation}
where $p-1=\tilde{j}_{1}+\dots+\tilde{j}_{N}$.

\subsection{Fischer Normalization $\left(F,\tilde{F}\right)$-Spaces\cite{bu2},\cite{bu3}}The Fischer Decomposition from (\ref{opseclala}) gives
\begin{equation}\begin{split}& \left(\overline{z}_{l}+2\lambda_{l}z_{l}\right) z^{J}=A_{l,J}\left(z,\xi\right) Q(\left(z,\xi\right)+C_{l,J}\left(z,\xi\right),\quad \tilde{\tr}\left(C_{l,J}\left(z,\xi\right)\right)=0,\quad\mbox{where $J\not\in \tilde{\mathcal{T}}_{l}$},\quad\mbox{for all   $l\in 1,\dots, N$.}\\& \left(\xi_{l}+2\lambda_{l}\xi_{l}\right) z^{J}=\tilde{A}_{l,J}\left(z,\xi\right) Q\left(z,\xi\right)+\tilde{C}_{l,J}\left(z,\xi\right),\quad \tilde{\tr}\left(\tilde{C}_{l,J}\left(z,\xi\right)\right)=0,\quad\mbox{where $J\not\in \tilde{\mathcal{T}}_{l}$,} \quad\mbox{for all   $l\in 1,\dots, N$.}\end{split} \label{optseclala} 
\end{equation}

This set of polynomials, derived from (\ref{optseclala}),   gives   normalization conditions defining  certain  Double-Spaces of Fischer Normalizations, which are constructed   iteratively from the generalized version of  the Fischer Decomposition\cite{sh}: for a given   homogeneous polynomial of degree $p\geq 1$ in $\left(z,\xi\right)$ denoted by $P\left(z,\xi\right)$,  we  have
\begin{equation} \begin{split}& P\left(z,\xi\right)=P_{1}\left(z,\xi\right)Q\left(z,\xi\right)+R_{1} \left(z,\xi\right)=\displaystyle\sum_{l=1}^{N}\displaystyle\sum_{J\in\mathbb{N}^{N}\atop {\left|J\right|=p-1}} \left(a_{l,J}C_{l,J}\left(z,\xi\right)+b_{l,J}\tilde{C}_{l,J}\left(z,\xi\right)\right)+R_{1,0}\left(z,\xi\right),\quad\mbox{such that:}\\&  \quad\quad\quad\quad\quad\quad\quad\quad R_{1,0}\left(z,\xi\right)\in \displaystyle\bigcap_{l=1}^{N}\left(\left(\displaystyle\bigcap_{J\in\mathbb{N}^{N},\hspace{0.1 cm}J\not\in \mathcal{T}_{l}\atop {\left|J\right|=p-1}} \left(\ker  C^{\star}_{l,J}  \bigcap  \ker  \tilde{C}^{\star}_{l,J} \right)\right)\bigcap   \left(\displaystyle\bigcap_{J\in\mathbb{N}^{N},\hspace{0.1 cm}J\in \mathcal{T}_{l}\atop {\left|J\right|=p-1}} \left(\ker  \left( z_{l}\xi^{J} \right)^{\star} \bigcap  \ker  \left( \xi_{l}z^{J} \right)^{\star} \right)\right)\right),\\& P_{1}\left(z,\xi\right)=P_{2}\left(z,\xi\right)Q\left(z,\xi\right)+R_{2}\left(z,\xi\right),\quad\mbox{where $\tr\left(R_{2}\left(z,\xi\right)\right)=0$ and:}\\& \quad    R_{2}\left(z,\xi\right)=\displaystyle\sum_{l=1}^{N}\displaystyle\sum_{J\in\mathbb{N}^{N}\atop {\left|J\right|=p-3}} \left(a_{l,J}C_{l,J}(z,\overline{z})+b_{l,J}\tilde{C}_{l,J}\left(z,\xi\right)\right)+R_{2,0}\left(z,\xi\right),\quad\mbox{such that:} \\&\\&  \quad\quad\quad\quad\quad\quad\quad  R_{2,0}\left(z,\xi\right)\in \displaystyle\bigcap_{l=1}^{N}\left(\left(\displaystyle\bigcap_{J\in\mathbb{N}^{N},\hspace{0.1 cm}J\not\in \mathcal{T}_{l}\atop {\left|J\right|=p-3}}\left( \ker  C^{\star}_{l,J}  \bigcap   \ker  \tilde{C}^{\star}_{l,J}\right)\right)\bigcap\left(   \displaystyle\bigcap_{J\in\mathbb{N}^{N},\hspace{0.1 cm}J\in \mathcal{T}_{l}\atop {\left|J\right|=p-3}} \ker  \left( z_{l}\xi^{J} \right)^{\star} \bigcap   \ker  \left( \xi_{l}z^{J} \right)^{\star} \right)\right),
\\&\quad\quad\quad\quad\vdots\quad\quad\quad\quad\quad\quad\quad\quad\vdots\quad\quad\quad\quad\quad\quad\quad\quad\vdots\quad\quad\quad\quad\quad\quad\quad\quad\vdots \quad\quad\quad\quad\quad\quad\quad\quad\quad\quad\quad\quad\quad\quad\vdots\\& P_{k}\left(z,\xi\right)=P_{k+1}\left(z,\xi\right)Q\left(z,\xi\right)+R_{k+1}(z,\overline{z}),\quad\mbox{where $\tr\left(R_{k+1}\left(z,\xi\right)\right)=0$ and:}\\&   R_{k+1}\left(z,\xi\right)=\displaystyle\sum_{l=1}^{N}\displaystyle\sum_{J\in\mathbb{N}^{N}\atop {\left|J\right|=p-2k-1}} \left(a_{l,J}C_{l,J}\left(z,\xi\right)+b_{l,J}\tilde{C}_{l,J}\left(z,\xi\right)\right)+R_{k+1,0}\left(z,\xi\right),\quad\mbox{such that:}\\&  \quad\quad\quad\quad\quad\quad\quad R_{k+1,0}\left(z,\xi\right)\in \displaystyle\bigcap_{l=1}^{N}\left(\left(\displaystyle\bigcap_{J\in\mathbb{N}^{N},\hspace{0.1 cm}J\not\in \mathcal{T}_{l}\atop {\left|J\right|=p-2k-1}} \left(\ker  C^{\star}_{l,J}  \bigcap   \ker  \tilde{C}^{\star}_{l,J} \right)\right)\bigcap\left( \displaystyle\bigcap_{J\in\mathbb{N}^{N},\hspace{0.1 cm}J\in \mathcal{T}_{l}\atop {\left|J\right|=p-2k-1}} \left(\ker  \left( z_{l}\xi^{J} \right)^{\star} \bigcap   \ker  \left( \xi_{l}z^{J} \right)^{\star}\right)\right)\right),\\&\quad\quad\quad\vdots\quad\quad\quad\quad\quad\quad\quad\quad\vdots\quad\quad\quad\quad\quad\quad\quad\quad\vdots\quad\quad\quad\quad\quad\quad\quad\quad\vdots \quad\quad\quad\quad\quad\quad\quad\quad\quad\quad\quad\quad\quad\quad\vdots\end{split}.\label{new2lala}
\end{equation}
where   these occurring polynomials 
\begin{equation}\left\{P_{k}\left(z,\xi\right)\right\}_{k=1,\dots,\left[\frac{p-1}{2}\right]},\quad \left\{R_{k}\left(z,\xi\right)\right\}_{k=1,\dots,\left[\frac{p-1}{2}\right]},
\label{poll2lala}
\end{equation}
are iteratively obtained  using the generalized version  of the Fischer Decomposition\cite{sh}. 

Recalling   strategies from \cite{bu2} and \cite{bu3}, we define 
\begin{equation}\tilde{\mathcal{F}}_{p},\quad p\in\mathbb{N}^{\star},\label{spartiuFlala}
\end{equation}
which consist in real-valued polynomials $P(z,\overline{z})$  of degree $p\geq 1$ in $(z,\overline{z})$ satisfying the normalizations:
$$ P_{k}^{\left(p\right)}\left(z,\xi\right)=P_{k+1}^{\left(p\right)}\left(z,\xi\right)Q\left(z,\xi\right)+R_{k+1}^{\left(p\right)}\left(z,\xi\right),\quad \mbox{for all   $k=0,\dots, \left[\frac{p-1}{2}\right]$ and given $P_{0}^{\left(p\right)}\left(z,\xi\right)=P\left(z,\xi\right)$,}$$
such that
\begin{equation}
\begin{split}& R_{k+1}^{\left(p\right)}\left(z,\xi\right)\in \displaystyle\bigcap_{l=1}^{N}\left( \displaystyle\bigcap_{J\in\mathbb{N}^{N},\hspace{0.1 cm}J\not\in \mathcal{T}_{l}\atop {\left|J\right|=p-2k-1}} \left(\ker  C^{\star}_{l,J}  \bigcap   \ker  \tilde{C}^{\star}_{l,J}\bigcap  \ker\tr\right)\right)\\& \quad\quad\quad\quad\quad\quad\quad\quad\bigcap\left(  \displaystyle\bigcap_{J\in\mathbb{N}^{N},\hspace{0.1 cm}J\in \mathcal{T}_{l}\atop {\left|J\right|=p-2k-1}} \left(\ker  \left( z_{l}\xi^{J} \right)^{\star} \bigcap   \ker  \left( \xi_{l}z^{J} \right)^{\star}\right)\right),\quad \mbox{for all    $k=0,\dots, \left[\frac{p-1}{2}\right]$ .} \end{split}
\end{equation}

We  consider the Fischer Decompositions (\ref{new1}) choosing 
\begin{equation}P\left(z,\xi\right))=\frac{\varphi_{p}\left(z,\xi\right)+\overline{\varphi}_{p}\left(\xi,z\right) }{2},\quad\mbox{for given $p\in\mathbb{N}^{\star}$.}\label{kama2lala}\end{equation}  
  
It remains to prove that these Double-Spaces of Fischer-Normalizations (\ref{spartiuFlala}) uniquely determine  the $\left(F,\tilde{F}\right)$-component of the formal equivalence (\ref{maplala}). It suffices to prove the linear independence, considering complex numbers, of the following set of polynomials 
\begin{equation}\left\{C_{l,J}\left(z,\xi\right),\hspace{0.1 cm} \tilde{C}_{l,J}\left(z,\xi\right)\right\}_{J\in\mathbb{N}^{N}\atop{\left|J\right|=p-1\atop{J\not\in\mathcal{T}_{l}\atop{l=1,\dots,N}}}},\quad\mbox{for all $p\in\mathbb{N}^{\star}$.}\label{Africa330lala}
\end{equation}

These computations are difficult  to conclude because of  the overlapping of the  homogeneous   polynomials from (\ref{Africa330lala}). There are recalled as previously the   arguments corresponding to (\ref{330}) and (\ref{331}). Then,  we consider by (\ref{vectori}) and (\ref{AfricaIJ1}) the following vectors 
  
\begin{equation} \tilde{Z}\left[J \right]=\begin{pmatrix}a_{\left(j_{1}-2,\dots,j_{k},\dots,j_{N}\right)}^{\left(1\right)}\\ \vdots \\ a_{\left(j_{1},\dots,j_{k}-2,\dots,j_{N}\right)}^{\left(1\right)}\\ \vdots \\  a_{\left(j_{1},\dots,j_{k},\dots,j_{N}-2\right)}^{\left(1\right)}\\  \vdots\\ a_{\left(j_{1}-2,\dots,j_{k},\dots,j_{N}\right)}^{(N)}\\ \vdots \\ a_{\left(j_{1},\dots,j_{k}-2,\dots,j_{N}\right)}^{(N)}\\ \vdots \\  a_{\left(j_{1},\dots,j_{k},\dots,j_{N}-2\right)}^{(N)}\end{pmatrix},\quad \xi\left[J \right]=\begin{pmatrix}\tilde{a}_{\left(j_{1}-2,\dots,j_{k},\dots,j_{N}\right)}^{\left(1\right)}\\ \vdots \\ \tilde{a}_{\left(j_{1},\dots,j_{k}-2,\dots,j_{N}\right)}^{\left(1\right)}\\ \vdots \\  \tilde{a}_{\left(j_{1},\dots,j_{k},\dots,j_{N}-2\right)}^{\left(1\right)}\\  \vdots\\ \tilde{a}_{\left(j_{1}-2,\dots,j_{k},\dots,j_{N}\right)}^{(N)}\\ \vdots \\ \tilde{a}_{\left(j_{1},\dots,j_{k}-2,\dots,j_{N}\right)}^{(N)}\\ \vdots \\  \tilde{a}_{\left(j_{1},\dots,j_{k},\dots,j_{N}-2\right)}^{(N)}\end{pmatrix}.\label{kakalala}
\end{equation}

Immediately from (\ref{opt}), we obtain
\begin{equation}\begin{split}&\left(\overline{z}_{l}+2\lambda_{l}z_{l}\right) z^{J}-A_{l,J}\left(z,\xi\right) Q\left(z,\xi\right)=C_{l,J}\left(z,\xi\right),\quad \tilde{\tr}\left(C_{l,J}\left(z,\xi\right)\right)=0,\quad\mbox{where $J\not\in \tilde{\mathcal{T}}_{l}$, }\\&\left(z_{l}+2\lambda_{l}\xi_{l}\right) \xi^{J}-\tilde{A}_{l,J}\left(z,\xi\right) Q\left(z,\xi\right)=\tilde{C}_{l,J}\left(z,\xi\right),\quad \tilde{\tr}\left(\tilde{C}_{l,J}\left(z,\xi\right)\right)=0,\quad\mbox{where $J\not\in \tilde{\mathcal{T}}_{l}$, } \end{split}\label{Africaopt1lala}
\end{equation} 
for all  $l\in 1,\dots, N$. 
  
As previously, we construct systems of equations according to the lexicografic order related to (\ref{ODINla}), using
\begin{equation*}
\begin{split}& \tilde{Z}^{t}=\left(\left\{a_{J}^{\left(1\right)}\right\}_{J\in\mathbb{N}^{N}\atop \left|I\right|=p-1},\left\{a_{J}^{(2)}\right\}_{J\in\mathbb{N}^{N}\atop \left|I\right|=p-1},\dots,\left\{a_{J}^{(N)}\right\}_{J\in\mathbb{N}^{N}\atop \left|I\right|=p-1}\right),\\& \tilde{\xi}^{t}=\left(\left\{\tilde{a}_{J}^{\left(1\right)}\right\}_{J\in\mathbb{N}^{N}\atop \left|I\right|=p-1},\left\{\tilde{a}_{J}^{(2)}\right\}_{J\in\mathbb{N}^{N}\atop \left|I\right|=p-1},\dots,\left\{\tilde{a}_{J}^{(N)}\right\}_{J\in\mathbb{N}^{N}\atop \left|I\right|=p-1}\right).\end{split}
\end{equation*}
 
 We obtain  the following  system   of  equations 
\begin{equation}\begin{split}&\left(I-\mbox{Aux}_{p}\tilde{A}\right) \tilde{Z}+\tilde{B} \tilde{\xi} =\tilde{V}\left(z_{1},z_{2},\dots,z_{N}\right),\\& \left(I-\mbox{Aux}_{p}\tilde{A}\right) \tilde{\xi}+\tilde{B}\tilde{Z}=\tilde{W}\left(\xi_{1},\xi_{2},\dots,\xi_{N}\right),\end{split} \label{Afrla}
\end{equation} 
where $\tilde{V}\left(z_{1},z_{2},\dots,z_{N}\right)$ and $\tilde{W}\left(\xi_{1},\xi_{2},\dots,\xi_{N}\right)$ are   known homogeneous vector polynomials  of degree $p-2$,  according with the following notations
\begin{equation}\tilde{A}=  \begin{pmatrix} A_{1} &  \mbox{O}_{N^{p-1}}& \dots & \mbox{O}_{N^{p-1}}     \\ \mbox{O}_{N^{p-1}}&    A_{2} & \dots & \mbox{O}_{N^{p-1}}\\ \vdots   & \vdots & \ddots & \vdots   \\ \mbox{O}_{N^{p-1}}   &  \mbox{O}_{N^{p-1}}& \dots &  A_{N}
\end{pmatrix} ,\quad \tilde{B}=  \begin{pmatrix} B_{1} &  \mbox{O}_{N^{p-1}}& \dots & \mbox{O}_{N^{p-1}}     \\ \mbox{O}_{N^{p-1}}&    B_{2} & \dots & \mbox{O}_{N^{p-1}}\\ \vdots   & \vdots & \ddots & \vdots   \\ \mbox{O}_{N^{p-1}}   &  \mbox{O}_{N^{p-1}}& \dots &  B_{N}
\end{pmatrix} ,\label{Africacalcan1}
\end{equation}
where   we have
$$A_{1}, A_{2},\dots, A_{N};\quad B_{1},  B_{2}, \dots,  B_{N}\in \mathcal{M}_{N^{p-1}\times N^{p-1} }\left(\mathbb{C}\right).
$$

 There are  considered products of matrices as in    (\ref{Africa90000}), (\ref{Africa90000se})  in (\ref{Afrla}) using   (\ref{Africacalcan2}). Then, the solution is computed, being unique in the light of the invertibility of the following matrices
\begin{equation}\frac{1}{I_{N^{p}}-\tilde{\mbox{Aux}}_{p}\tilde{A}-\tilde{B}},\quad\frac{1}{I_{N^{p}}-\tilde{\mbox{Aux}}_{p}\tilde{A}+\tilde{B}},\label{5501}
\end{equation}
recalling the arguments related to (\ref{550se}) in order to solve (\ref{Afrla}). 

Now, we are ready to move forward:
 \section{Normalizations using Formal Segre-Holomorphic   Equivalences}
Summarizing   the   computations of the formal transformation (\ref{maplala}), we obtain the following result:

\bp\label{proppp} Let $\mathcal{M}\subset\mathbb{C}^{2N+2}$ be a real-formal
submanifold defined near  $p=0$ as follows
\begin{equation}
w=z_{1}\xi_{1}+\dots+z_{N}\xi_{N}+\lambda_{1}\left(z_{1}^{2}+\xi_{1}^{2}\right)+\dots+\lambda_{N}\left(z_{N}^{2}+\xi_{N}^{2}\right) +\displaystyle\sum
_{k\geq 3}\varphi_{k}\left(z,\xi\right), \label{ecuatie11lala}
\end{equation}
where   $\varphi _{k}\left(z,\xi\right)$ is a  
polynomial of    degree $k$ in $\left(z,\xi\right)$,  for all
 $k\geq 3$, such that (\ref{lambida}) is satisfied assuming
\begin{equation}\lambda_{1},\dots,\lambda_{k_{0}}\neq 0,\quad \lambda_{k_{0}+1}=\dots=\lambda_{N}=0,\quad\mbox{for some $k_{0}\in 1,\dots, N$.}\label{03lala}
\end{equation}  
 
Then there exists a unique formal equivalence defined as in (\ref{maplala}) and  normalized as follows\begin{equation}  \Re \left(F_{1,n}+\tilde{F}_{1,n}\right)(z)=0,\quad  F_{0,n+1}^{\left(k_{0}+1\right)}(z)=\dots=F_{0,n+1}^{\left(N\right)}(z)= 0,\quad \tilde{F}_{0,n+1}^{\left(k_{0}+1\right)}(z)=\dots=\tilde{F}_{0,n+1}^{\left(N\right)}(z)= 0 \quad\mbox{for all $n\in\mathbb{N}^{\star}$,} \label{olala}
\end{equation}
that sends $\mathcal{M}$ into the following partial normal form
\begin{equation}
w'={z'}_{1}{\xi'}_{1}+\dots+{z'}_{N}{\xi'}_{N}+\lambda_{1}\left({z'}_{1}^{2}+{\xi'}_{1}^{2}\right)+\dots+\lambda_{N}\left({z'}_{N}^{2}+{\xi'}_{N}^{2}\right) +\displaystyle\sum
_{k\geq 3}\varphi'_{k}\left(z',{\xi'}\right), \label{ecuatie111lala}
\end{equation}
where   $\varphi' _{k}\left(z,\xi\right)$ is a  
polynomial of    degree $k$ in $\left(z,\xi\right)$,  for all
 $k\geq 3$, such that there are satisfied the following normalizations
\begin{equation}\left\{\begin{split}&\frac{\varphi'_{k}\left(z',{\xi'}\right)+\overline{\varphi'}_{k}\left( {\xi'},z'\right)}{2}\in\tilde{\mathcal{F}}_{k}, \quad\mbox{for all $k\geq 3$},\\& \frac{\varphi'_{k}\left(z',{\xi'}\right)-\overline{\varphi'}_{k}\left( {\xi'},z'\right)}{2i}\in\tilde{\mathcal{G}}_{k},\quad\mbox{for all $k\geq 3$},\end{split}\right.\label{cnlala}
\end{equation}
and respectively, the following normalizations 
\begin{equation}\tilde{\tilde{P}}_{\frac{k}{2}}^{\left(k\right)}(z,\overline{z})\in   \displaystyle\bigcap_{i,j=1}^{N}\ker\left(\frac{\partial^{2}}{\partial z_{i}\partial\overline{z_{j}}+\lambda_{i}\partial z_{i} \partial z_{j}+\lambda_{i}\partial\overline{z_{i}}\partial \overline{z_{j}}}\right)   ,\quad\mbox{for $k$ even,} \label{cnlala1} 
\end{equation}
and respectively, the following normalizations
\begin{equation}  \tilde{\tilde{P}}_{\frac{k+1}{2}}^{\left(k\right)}(z,\overline{z})\in  \ker\left(\frac{\partial^{2}}{\partial z_{1} \partial \overline{z}_{1}+\dots+\partial z_{N} \partial \overline{z}_{N}}\right)\displaystyle\bigcap   \displaystyle\bigcap_{i=1}^{k_{0}}\ker \left(\frac{\partial^{3}}{\partial z_{i}^{3}}\right)  ,\quad\mbox{for $k$ odd,} \label{coreo1lap} 
\end{equation}
according to the writing (\ref{ui119}). 
 \ep   
 \begin{proof}
 It is considered an induction process depending on $r=m+2n\in\mathbb{N}^{\star}$ in order to compute the equivalence (\ref{maplala}). We assume thus that  we have determined the following terms
\begin{equation}
 G_{m,n}(z),\quad \tilde{G}_{m,n}\left(\xi\right),\quad \mbox{for all $m,n\in\mathbb{N}^{\star}$ with $m+2n<r$.}\label{gigel1lala}
\end{equation}
and respectively, the following terms
\begin{equation}
F_{m,n}(z),\quad  \tilde{F}_{m,n}\left(\xi\right),\quad  \mbox{for all $m,n\in\mathbb{N}^{\star}$ with $m+2n<r-1$.}\label{gigel2lala}
\end{equation}

Moreover, we introduce by (\ref{ef}) the following notations
$$ F_{m,n}(z) =\left(F_{m,n}^{\left(1\right)}(z),\dots, F_{m,n}^{\left(N\right)}(z\right),\quad \tilde{F}_{m,n}\left(\xi\right) =\left(\tilde{F}_{m,n}^{\left(1\right)}\left(\xi\right),\dots, \tilde{F}_{m,n}^{\left(N\right)}\left(\xi\right)\right),\quad\mbox{for all $m,n\in\mathbb{N}^{\star}$ with $m+n\geq 3$},$$
according to the following notation
\begin{equation}F(z,w)=\left(F^{\left(1\right)}(z,w),\dots,F^{\left(N\right)}(z,w)\right),\quad \tilde{F}\left(\xi,\nu\right)=\left(\tilde{F}^{\left(1\right)}\left(\xi,\nu\right),\dots,\tilde{F}^{\left(N\right)}\left(\xi,\nu\right)\right).\label{eflala}
\end{equation}

Moreover, we consider the following notations
\begin{equation}
\begin{split}& G_{m,n}(z)=\displaystyle\sum_{i_{1}+\dots+i_{N}=m\atop I=\left(i_{1},\dots,i_{N}\right)\in\mathbb{N}^{N}} g_{I,n}z_{1}^{i_{1}}\dots z_{N}^{i_{N}}w^{n},\quad  F_{m,n}^{\left(l\right)}(z)=\displaystyle\sum_{j_{1}+\dots+j_{N}=m-1\atop J=\left(j_{1},\dots,j_{N}\right)\in\mathbb{N}^{N}}f_{J,n}^{\left(l\right)}z_{1}^{j_{1}}\dots z_{N}^{j_{N}}w^{n},  \\& \tilde{G}_{m,n}(z)=\displaystyle\sum_{i_{1}+\dots+i_{N}=m\atop I=\left(i_{1},\dots,i_{N}\right)\in\mathbb{N}^{N}} \tilde{g}_{I,n}z_{1}^{i_{1}}\dots z_{N}^{i_{N}}w^{n},\quad  \tilde{F}_{m,n}^{\left(l\right)}(z)=\displaystyle\sum_{j_{1}+\dots+j_{N}=m-1\atop J=\left(j_{1},\dots,j_{N}\right)\in\mathbb{N}^{N}}\tilde{f}_{J,n}^{\left(l\right)}z_{1}^{j_{1}}\dots z_{N}^{j_{N}}w^{n}, \end{split}
\end{equation}
for all $l=1,\dots,N$. 

We want to determine  the following terms
$$g_{I,n},\quad \tilde{g}_{I,n},\quad\mbox{for all $I=\left(i_{1},\dots,i_{N}\right)\in\mathbb{N}^{N}$ such that $i_{1}+\dots+i_{N}=m$ and $m+2n=r$.}  
$$

 The entire construction is based  on defining (\ref{spartiuGlala1}) and (\ref{spartiuFlala})  by considering as previously    successively Fischer Decompositions\cite{sh} in order to compute the formal equivalence (\ref{maplala}) according to    (\ref{new1lala1}) and  (\ref{new2lala}). In particular, we have the following
 \begin{equation}  
\displaystyle\sum_{m+2n=r}\frac{ G_{m,n}(z)- \tilde{G}_{m,n}\left(\xi\right)}{2\sqrt{-1}}\left(z_{1}\xi_{1}+\dots+z_{N}\xi_{N}+\lambda_{1}\left(z_{1}^{2}+\xi_{1}^{2}\right)+\dots+\lambda_{N}\left(z_{N}^{2}+\xi_{N}^{2}\right)\right)^{n}   =  \frac{\varphi_{r}\left(z,\xi\right)-\overline{\varphi_{r}}\left(\xi,z\right)}{2\sqrt{-1}}+V\left(z,\xi\right), \label{ecuatieXYY11lala}
\end{equation} 
where $V\left(z,\xi\right)$, which  represents a sum of known terms depending on $\left(z,\xi\right)$  depending also on (\ref{gigel1lala}) and (\ref{gigel2lala}) according to the induction hypothesis. Then the above right-hand side may be decomposed as in  (\ref{new1lala1}).  We obtain
\begin{equation}
\begin{split}&  \displaystyle\sum_{\left|I\right|=r \atop I\in\mathbb{N}^{N}}\frac{g_{I,0}z^{I}-
\tilde{g}_{I,0}\xi^{I}}{2\sqrt{-1}}=P_{1} \left(z,\xi\right)Q\left(z,\xi\right)+R_{1,0}\left(z,\xi\right),\quad\mbox{where $\tr\left(R_{1}\left(z,\xi\right)\right)=0$ and:}\\& \quad\quad\quad\quad\quad\quad\quad\quad \hspace{0.1 cm}   R_{1,0}\left(z,\xi\right)=\displaystyle\sum_{ {\left|I\right|=r}\atop I\in\mathbb{N}^{N}} \left(a_{I,0}C_{I,0}\left(z,\xi\right)+\tilde{a}_{I,0}\tilde{C}_{I,0}\left(z,\xi\right) \right)+R_{1,0}\left(z,\xi\right),\hspace{0.1 cm}\mbox{such that:} \\& \quad\quad\quad\quad\quad \quad\quad\quad R_{1,0}\left(z,\xi\right)\in \left(\displaystyle\bigcap_{I\in\mathbb{N}^{N},\hspace{0.1 cm} I\not\in \mathcal{S} \atop {\left|I\right|=p}}\left( \ker  C^{\star}_{I}  \cap  \ker  \tilde{C}^{\star}_{I} \right)\right) \bigcap\left( \displaystyle\bigcap_{I\in\mathbb{N}^{N},\hspace{0.1 cm} I \in \mathcal{S} \atop {\left|I\right|=p}}\left(\ker  \left(z^{I}\right)^{\star}   \cap  \ker  \left(\xi^{I}\right)^{\star}\right)\right) ,\end{split}\label{ha1lala}
\end{equation}
where the corresponding  above occurring polynomials  are known   according to (\ref{new1lala1}) in (\ref{ecuatieXYY11}).  We obtain 
\begin{equation}G_{r,0}(z)=\displaystyle\sum_{I\in\mathbb{N}^{N},\hspace{0.1 cm} I\not\in \mathcal{S} \atop {\left|I\right|=r}}a_{I,0}z^{I}+\displaystyle\sum_{I\in\mathbb{N}^{N},\hspace{0.1 cm} I \in \mathcal{S} \atop {\left|I\right|=r}}g_{I,0}z^{I} ,\quad \tilde{G}_{r,0}\left(\xi\right)=-\displaystyle\sum_{I\in\mathbb{N}^{N},\hspace{0.1 cm} I\not\in \mathcal{S} \atop {\left|I\right|=r}} \tilde{a}_{I,0}\xi^{I}+\displaystyle\sum_{I\in\mathbb{N}^{N},\hspace{0.1 cm} I \in \mathcal{S} \atop {\left|I\right|=r}}\tilde{g}_{I,0}\xi^{I},\label{lele1lala}
\end{equation}
where the terms defining the above second sum are directly computed by identifying corresponding coefficients  in (\ref{ecuatieXYY11lala}).

Then,   (\ref{lele1}) provides contributions to the other terms, because we further deal with the following Fischer generalized Decompositions 
\begin{equation}\begin{split}&  \displaystyle\sum_{\left|I\right|=r-2\atop I\in\mathbb{N}^{N}}\frac{g_{I,1}z^{I}-\tilde{g}_{I,1}\xi^{I}}{2\sqrt{-1}} =  P_{1}\left(z,\xi\right)+\displaystyle\sum_{\left|I\right|=r}\frac{A_{I,0}(z)-\tilde{A}_{I,0}\left(\xi\right)}{2\sqrt{-1}}\\& \quad \quad \quad\quad\quad\quad\quad\quad \quad\quad\hspace{0.17 cm}=P_{2}\left(z,\xi\right)Q\left(z,\xi\right)+R_{2}\left(z,\xi\right),\quad\mbox{where $\tr\left(R_{2}\left(z,\xi\right)\right)=0$ and:}\\& \quad\quad\quad\quad\quad\quad \quad\quad \hspace{0.1 cm}  R_{2}\left(z,\xi\right)=\displaystyle\sum_{I\in\mathbb{N}^{N}\atop {\left|I\right|=r-2}} \left(a_{I,1}C_{I,1}\left(z,\xi\right)+\tilde{a}_{I,1}\tilde{C}_{I,1}\left(z,\xi\right)\right)+R_{2,0}\left(z,\xi\right), \hspace{0.1 cm}\mbox{such that:}\\&\quad\quad\quad\quad   \quad \quad\quad\quad R_{2,0}\left(z,\xi\right)\in \left(\displaystyle\bigcap_{I\in\mathbb{N}^{N},\hspace{0.1 cm} I\not\in \mathcal{S} \atop {\left|I\right|=p-2}}\left( \ker  C^{\star}_{I}  \cap  \ker  \tilde{C}^{\star}_{I}\right) \right) \bigcap \left(\displaystyle\bigcap_{I\in\mathbb{N}^{N},\hspace{0.1 cm} I \in \mathcal{S} \atop {\left|I\right|=p-2}}\left(\ker  \left(z^{I}\right)^{\star}   \cap  \ker  \left(\xi^{I}\right)^{\star}\right)\right),\end{split}\label{ha2lala}
\end{equation}
where  all  occurring polynomials   are obtained iteratively according to (\ref{new1lala1}). We obtain
\begin{equation}G_{r-2,1}(z)=\displaystyle\sum_{I\in\mathbb{N}^{N},\hspace{0.1 cm} I\not\in \mathcal{S} \atop {\left|I\right|=r-2}}a_{I,1}z^{I}+\displaystyle\sum_{I\in\mathbb{N}^{N},\hspace{0.1 cm} I \in \mathcal{S} \atop {\left|I\right|=r-2}}g_{I,1}z^{I} ,\quad \tilde{G}_{r-2,1}\left(\xi\right)=-\displaystyle\sum_{I\in\mathbb{N}^{N},\hspace{0.1 cm} I\not\in \mathcal{S} \atop {\left|I\right|=r-2}}\tilde{a}_{I,1}\xi^{I}+\displaystyle\sum_{I\in\mathbb{N}^{N},\hspace{0.1 cm} I \in \mathcal{S} \atop {\left|I\right|=r-2}}\tilde{g}_{I,1}\xi^{I},\label{lele2lala}
\end{equation} 
where the terms defining the above second sum are directly computed by identifying corresponding coefficients  in (\ref{ecuatieXYY11lala}).

These iterative computations determine inductively the $\left(G,\tilde{G}\right)$-component of the formal transformation using   the uniqueness of the Fischer Decomposition\cite{sh}, because  we can apply  an induction process with respect to $k=1,\dots,\left[\frac{r}{2}\right]$. We have

\begin{equation}\begin{split}&\displaystyle\sum_{\left|I\right|=r-2k\atop I\in\mathbb{N}^{N}}\frac{g_{I,k}z^{I}-\tilde{g}_{I,k}\xi^{I}}{2\sqrt{-1}}=     P_{k}\left(z,\xi\right)+\displaystyle\sum_{\left|I\right|=r-2k+2\atop I\in\mathbb{N}^{N}}\frac{A_{I,k-1}(z)-\tilde{A}_{I,k-1}(z)}{2\sqrt{-1}}\\& \quad \quad \quad\quad\quad\quad\quad\quad \quad\quad\quad=P_{k+1}\left(z,\xi\right)Q\left(z,\xi\right)+R_{k+1}\left(z,\xi\right),\quad\mbox{where $\tr\left(R_{k+1}\left(z,\xi\right)\right)=0$ and:}\\&  \quad\quad \quad\quad\quad\quad\quad\quad \quad \hspace{0.1 cm} R_{k+1}\left(z,\xi\right)=\displaystyle\sum_{I\in\mathbb{N}^{N}\atop {\left|I\right|=r-2k }} \left(a_{I,k}C_{I,k}\left(z,\xi\right)+\tilde{a}_{I,k}\tilde{C}_{I,k}\left(z,\xi\right)\right)+R_{k+1,0}\left(z,\xi\right),\hspace{0.1 cm}\mbox{where:}\\&\quad\quad\quad\quad  \quad\quad \quad\quad\quad R_{k+1,0}\left(z,\xi\right)\in \left(\displaystyle\bigcap_{I\in\mathbb{N}^{N}, \hspace{0.1 cm}I\not\in \mathcal{S} \atop {\left|I\right|=p-2k}} \left(\ker  C^{\star}_{I}  \cap  \ker  \tilde{C}^{\star}_{I}\right)\right)\bigcap\left(   \displaystyle\bigcap_{I\in\mathbb{N}^{N},\hspace{0.1 cm} I \in \mathcal{S} \atop {\left|I\right|=p-2k}}\left(\ker  \left(z^{I}\right)^{\star}   \cap  \ker \left(\xi^{I}\right)^{\star}\right) \right), \end{split}\end{equation}
where   all  occurring polynomials   are obtained iteratively according to (\ref{new1lala1}). We obtain 
\begin{equation}G_{r-2k }(z)=\displaystyle\sum_{I\in\mathbb{N}^{N},\hspace{0.1 cm} I\not\in \mathcal{S} \atop {\left|I\right|=r-2k }}a_{I,k}z^{I}+\displaystyle\sum_{I\in\mathbb{N}^{N},\hspace{0.1 cm} I \in \mathcal{S} \atop {\left|I\right|=r-2k }}g_{I,k}z^{I} ,\quad \tilde{G}_{r-2k }\left(\xi\right)=-\displaystyle\sum_{I\in\mathbb{N}^{N},\hspace{0.1 cm} I\not\in \mathcal{S} \atop {\left|I\right|=r-2k }}\tilde{a}_{I,k}\xi^{I}+\displaystyle\sum_{I\in\mathbb{N}^{N},\hspace{0.1 cm} I \in \mathcal{S} \atop {\left|I\right|=r-2k }}\tilde{g}_{I,k}\xi^{I},\label{lele3lala}
\end{equation} 
where the terms defining the above second sum are directly computed by identifying corresponding coefficients  in (\ref{ecuatieXYY11lala}).

 Now, according to (\ref{lele1lala}), (\ref{lele2lala}) and (\ref{lele3lala}), we obtain
\begin{equation}\displaystyle\sum_{m+2n=r}G_{m,n}(z)w^{n}=G_{r,0}(z)+G_{r-2,1}(z)w+\dots+G_{1,\frac{r-1}{2}}(z)w^{\frac{r-1}{2}}+  g_{0,\frac{r}{2} }(z)w^{\frac{r}{2}},\quad\mbox{for $r$ even, }\label{gee1la}
\end{equation} 
and respectively, we obtain
\begin{equation}\displaystyle\sum_{m+2n=r}G_{m,n}(z)w^{n}=G_{r,0}(z)+G_{r-2,1}(z)w+\dots+G_{1,\frac{r-1}{2}}(z)w^{\frac{r-1}{2}},\quad\mbox{for $r$ odd}.\label{gee2la}
\end{equation}

Analogously, according to (\ref{lele1lala}), (\ref{lele2lala}) and (\ref{lele3lala}), we obtain   
\begin{equation}\displaystyle\sum_{m+2n=r}G_{m,n}\left(\xi\right)\nu^{n}=G_{r,0}\left(\xi\right)+G_{r-2,1}\left(\xi\right)\nu+\dots+G_{1,\frac{r-1}{2}}\left(\xi\right)\nu^{\frac{r-1}{2}}-  g_{0,\frac{r}{2} }(z)\nu^{\frac{r}{2}},\quad\mbox{for $r$ even, }\label{gee1lala}
\end{equation} 
and respectively, we obtain  
\begin{equation}\displaystyle\sum_{m+2n=r}G_{m,n}(z)w^{n}=G_{r,0}\left(\xi\right)+G_{r-2,1}\left(\xi\right)\nu+\dots+G_{1,\frac{r-1}{2}}\left(\xi\right)\nu{\frac{r-1}{2}},\quad\mbox{for $r$ odd}.\label{gee2lala}
\end{equation}
 
Now, we apply the Fischer Decompositions (\ref{new1})   choosing
\begin{equation} \begin{split}       P\left(z,\xi\right)=&\frac{\varphi_{r}\left(z,\xi\right)+\overline{\varphi_{r}}\left( \xi,z\right) }{2}+V\left(z,\xi\right)\\&+\displaystyle\sum_{m+2n=r}\displaystyle\sum_{\left|I\right|=m\atop I\in\mathbb{N}^{N}}\frac{ g_{I,n}z^{I}+\overline{ g_{I,n}z^{I}}}{2 }\left(z_{1}\overline{z}_{1}+\dots+z_{N}\overline{z}_{N}+\lambda_{1}\left(z_{1}^{2}+\overline{z}_{1}^{2}\right)+\dots+\lambda_{N}\left(z_{N}^{2}+\overline{z}_{N}^{2}\right)\right)^{n},\end{split} \end{equation}

We determine now the following terms
$$ f_{J}^{\left(l\right)},\quad \tilde{f}_{J}^{\left(l\right)},\quad\mbox{for all $J=\left(j_{1},j_{2},\dots,j_{N}\right)\in\mathbb{N}^{N}$  such that $j_{1}+\dots+j_{N}=n$ with $m+2n=r-1$, for all $l=1,\dots,N$}.
$$

We consider the following Fischer Decompositions 
\begin{equation}\begin{split}&   \displaystyle\sum_{l=1}^{N}\displaystyle\sum_{\left|J\right|=r-1\atop J\in\mathbb{N}^{N}} \left(\xi_{l}f_{J}^{\left(l\right)}z^{J}+z_{l}\tilde{f}_{J}^{\left(l\right)\xi^{I}}\right) = P_{1}\left(z,\xi\right)Q\left(z,\xi\right)+R_{1}\left(z,\xi\right),\quad\mbox{where $\tr\left(R_{1}\left(z,\xi\right)\right)=0$ such that:}\\& \quad\quad\quad\quad\quad\quad\quad  R_{1}\left(z,\xi\right)=\displaystyle\sum_{l=1}^{N}\displaystyle\sum_{J\in\mathbb{N}^{N}\atop {\left|J\right|=r-1}} \left(a_{l,J,0}C_{l,J,0}\left(z,\xi\right)+\tilde{a}_{l,J,0}\tilde{C}_{l,J,0}\left(z,\xi\right) \right)+R_{1,0}\left(z,\xi\right),\quad\mbox{where:}\\&   \quad\quad\quad \quad\quad\quad\quad R_{1,0}\left(z,\xi\right)\in \displaystyle\bigcap_{l=1}^{N}\left(\left(\displaystyle\bigcap_{J\in\mathbb{N}^{N},\hspace{0.1 cm}J\not\in \mathcal{T}_{l}\atop {\left|J\right|=p-1}} \left(\ker  C^{\star}_{l,J}  \bigcap  \ker  \tilde{C}^{\star}_{l,J} \right)\right)\bigcap   \left(\displaystyle\bigcap_{J\in\mathbb{N}^{N},\hspace{0.1 cm}J\in \mathcal{T}_{l}\atop {\left|J\right|=p-1}} \left(\ker  \left( z_{l}\xi^{J} \right)^{\star} \bigcap  \ker  \left( \xi_{l}z^{J} \right)^{\star} \right)\right)\right),\end{split}
\end{equation}
where  all  occurring polynomials   are obtained iteratively according to (\ref{new2}).  We obtain \begin{equation}F_{r-1,0}^{\left(l\right)}(z)=\displaystyle\sum_{J\in\mathbb{N}^{N},\hspace{0.1 cm} J\not\in \mathcal{T}_{l} \atop {\left|J\right|=r-1}}a_{l,J,0}z^{J}+\displaystyle\sum_{J\in\mathbb{N}^{N},\hspace{0.1 cm} J \in \mathcal{T}_{l} \atop {\left|J\right|=r-1}}f_{J,0}^{\left(l\right)}z^{J},\quad \tilde{F}_{r-1,0}^{\left(l\right)}(z)=\displaystyle\sum_{J\in\mathbb{N}^{N},\hspace{0.1 cm} J\not\in \mathcal{T}_{l} \atop {\left|J\right|=r-1}}\tilde{a}_{l,J,0}z^{J}+\displaystyle\sum_{J\in\mathbb{N}^{N},\hspace{0.1 cm} J \in \mathcal{T}_{l} \atop {\left|J\right|=r-1}}\tilde{f}_{J,0}^{\left(l\right)}z^{J}, \label{lele1selala}
\end{equation}  
for all $l=1,\dots,N$, where all terms  are directly computed by identifying corresponding coefficients in (\ref{ecuatieXYY11lala}).  
 
Then,  (\ref{lele1selala}) provides contributions to the other terms, because we further deal with the following Fischer   Decompositions 
 \begin{equation}\begin{split}& \displaystyle\sum_{l=1}^{N}\displaystyle\sum_{\left|J\right|=r-3\atop J\in\mathbb{N}^{N}} \left(\xi_{l}f_{J}^{\left(l\right)}z^{J}+z_{l}\tilde{f}_{J}^{\left(l\right)}\xi^{J}\right)=P_{1}\left(z,\xi\right)+\displaystyle\sum_{l=1}^{N}\displaystyle\sum_{\left|J\right|=r-1\atop J\in\mathbb{N}^{N}}\left(A_{J,l}(z)+\tilde{A}_{J,l}\left(\xi\right)\right)\\& \quad\quad\quad\quad\quad\quad\quad\quad\quad\quad\quad\quad\quad\quad\quad   =P_{2}\left(z,\xi\right)Q\left(z,\xi\right)+R_{2}\left(z,\xi\right),\quad\mbox{where $\tr\left(R_{2}\left(z,\xi\right)\right)=0$ and:}\\& \quad\quad \quad \quad \quad \quad     R_{2}\left(z,\xi\right)=\displaystyle\sum_{l=1}^{N}\displaystyle\sum_{J\in\mathbb{N}^{N}\atop {\left|J\right|=r-3}} \left(a_{l,J,1}C_{l,J,1}\left(z,\xi\right)+\tilde{a}_{l,J,1}\tilde{C}_{l,J,1}\left(z,\xi\right)\right)+R_{2,0}\left(z,\xi\right),\quad\mbox{such that:} \\&    \quad\quad \quad\quad\quad\quad  R_{2,0}\left(z,\xi\right)\in \displaystyle\bigcap_{l=1}^{N}\left(\left(\displaystyle\bigcap_{J\in\mathbb{N}^{N},\hspace{0.1 cm}J\not\in \mathcal{T}_{l}\atop {\left|J\right|=p-3}}\left( \ker  C^{\star}_{l,J}  \bigcap   \ker  \tilde{C}^{\star}_{l,J}\right)\right)\bigcap\left(   \displaystyle\bigcap_{J\in\mathbb{N}^{N},\hspace{0.1 cm}J\in \mathcal{T}_{l}\atop {\left|J\right|=p-3}} \ker  \left( z_{l}\xi^{J} \right)^{\star} \bigcap   \ker  \left( \xi_{l}z^{J} \right)^{\star} \right)\right), \end{split} 
 \end{equation}
where all  occurring polynomials   are obtained iteratively using  generalized versions of the Fischer Decomposition\cite{sh}. We obtain
\begin{equation}F_{r-3,0}^{\left(l\right)}(z)=\displaystyle\sum_{J\in\mathbb{N}^{N},\hspace{0.1 cm} J\not\in \mathcal{T}_{l} \atop {\left|J\right|=r-3}}a_{l,J,1}z^{J}+\displaystyle\sum_{J\in\mathbb{N}^{N},\hspace{0.1 cm} J \in \mathcal{T}_{l} \atop {\left|J\right|=r-3}}f_{J,1}^{\left(l\right)}z^{J} ,\quad \tilde{F}_{r-3,0}^{\left(l\right)}\left(\xi\right)=\displaystyle\sum_{J\in\mathbb{N}^{N},\hspace{0.1 cm} J\not\in \mathcal{T}_{l} \atop {\left|J\right|=r-3}}\tilde{a}_{l,J,1}\xi^{J}+\displaystyle\sum_{J\in\mathbb{N}^{N},\hspace{0.1 cm} J \in \mathcal{T}_{l} \atop {\left|J\right|=r-3}}\tilde{f}_{J,1}^{\left(l\right)}\xi^{J},\label{lele2selala}
\end{equation}    
for all $l=1,\dots,N$, where all terms   are directly computed by identifying corresponding coefficients in (\ref{ecuatieXYY11lala}). 

These iterative computations determine inductively the $F$-component of the formal transformation using   the uniqueness of the Fischer Decomposition\cite{sh}, because  we can apply  an induction process with respect to $k=1,\dots,\left[\frac{r-1}{2}\right]$. We have 
\begin{equation}\begin{split}& \displaystyle\sum_{l=1}^{N}\displaystyle\sum_{\left|J\right|=r-1-2k\atop J\in\mathbb{N}^{N}} \left(\xi_{l}f_{J}^{\left(l\right)}z^{J}+z_{l}\tilde{f}_{J}^{\left(l\right)\xi^{J}}\right)=P_{k}\left(z,\xi\right)+ \displaystyle\sum_{l=1}^{N}\displaystyle\sum_{\left|J\right|=r+1-2k\atop J\in\mathbb{N}^{N}}\left(A_{J,l}(z)+\tilde{A}_{J,l}\left(\xi\right)\right)\\&\quad\quad   \quad\quad\quad \quad \quad \quad \quad\quad\quad\quad\quad\quad\quad\quad\quad =P_{k+1}\left(z,\xi\right)Q\left(z,\xi\right)+R_{k+1}\left(z,\xi\right),\quad\mbox{where $\tr\left(R_{k+1}\left(z,\xi\right)\right)=0$ such that:}\\&   \quad\quad\quad   R_{k+1}\left(z,\xi\right)=\displaystyle\sum_{l=1}^{N}\displaystyle\sum_{J\in\mathbb{N}^{N}\atop {\left|J\right|=r-2k-1}} \left(a_{l,J,k}C_{l,J,k}\left(z,\xi\right)+\tilde{a}_{l,J,k}\tilde{C}_{l,J,k}\left(z,\xi\right) \right)+R_{k+1,0}\left(z,\xi\right),\quad\mbox{where:}\\&     \quad\quad\quad R_{k+1,0}\left(z,\xi\right)\in \left(\displaystyle\bigcap_{I\in\mathbb{N}^{N}, \hspace{0.1 cm}I\not\in \mathcal{S} \atop {\left|I\right|=p-2k}} \left(\ker  C^{\star}_{I}  \cap  \ker  \tilde{C}^{\star}_{I}\right)\right)\bigcap\left(   \displaystyle\bigcap_{I\in\mathbb{N}^{N},\hspace{0.1 cm} I \in \mathcal{S} \atop {\left|I\right|=p-2k}}\left(\ker  \left(z^{I}\right)^{\star}   \cap  \ker \left(\xi^{I}\right)^{\star}\right) \right), \end{split} 
 \end{equation}
where  all  occurring polynomials   are obtained iteratively using  generalized versions of the Fischer Decomposition\cite{sh}. We obtain \begin{equation}F_{r-1-2k,0}^{\left(l\right)}(z)=\displaystyle\sum_{J\in\mathbb{N}^{N},\hspace{0.1 cm} J\not\in \mathcal{T}_{l} \atop {\left|J\right|=r-1-2k}}a_{l,J,k}z^{J}+\displaystyle\sum_{J\in\mathbb{N}^{N},\hspace{0.1 cm} J \in \mathcal{T}_{l} \atop {\left|J\right|=r-1-2k}}f_{J,k}^{\left(l\right)}z^{J} ,\quad \tilde{F}_{r-1-2k,0}^{\left(l\right)}\left(\xi\right)=\displaystyle\sum_{J\in\mathbb{N}^{N},\hspace{0.1 cm} J\not\in \mathcal{T}_{l} \atop {\left|J\right|=r-1-2k}}\tilde{a}_{l,J,k}\xi^{J}+\displaystyle\sum_{J\in\mathbb{N}^{N},\hspace{0.1 cm} J \in \mathcal{T}_{l} \atop {\left|J\right|=r-1-2k}}\tilde{f}_{J,k}^{\left(l\right)}\xi^{J},\label{lele3selala}
\end{equation}    
for all $l=1,\dots,N$, where all terms   are directly computed by identifying corresponding coefficients in (\ref{ecuatieXYY11lala}). 

Assume  $r$   even. We write then as follows
\begin{equation}F_{1,n}(z):=F_{1,n}\left(z_{1},z_{2},\dots,z_{N}\right)=\begin{pmatrix}
v_{11} & v_{12} & \dots & v_{1N} \\ v_{21} & v_{22} & \dots & v_{2N} \\  \vdots  & \vdots  & \ddots & \vdots  \\ v_{N1} & v_{N2} & \dots & v_{NN}  
\end{pmatrix}\begin{pmatrix}
z_{1}\\ z_{2}\\ \vdots \\ z_{N}\end{pmatrix},
\end{equation}
and respectively, as follows
 \begin{equation}\tilde{F}_{1,n}\left(\xi\right):=F_{1,n}\left(\xi_{1},\xi_{2},\dots,\xi_{N}\right)=\begin{pmatrix}
\tilde{v}_{11} & \tilde{v}_{12} & \dots & \tilde{v}_{1N} \\ \tilde{v}_{21} & \tilde{v}_{22} & \dots & \tilde{v}_{2N} \\  \vdots  & \vdots  & \ddots & \vdots  \\ \tilde{v}_{N1} & \tilde{v}_{N2} & \dots & \tilde{v}_{NN}  
\end{pmatrix}\begin{pmatrix}
\xi_{1}\\ \xi_{2}\\ \vdots \\ \xi_{N}\end{pmatrix}.
\end{equation}

 Then, it remains to study the following sum of terms
\begin{equation}S=\displaystyle\sum_{j=1}^{N}\displaystyle\sum_{i=1}^{k_{0}}2\lambda_{i}\xi_{i} \tilde{v}_{ij}\xi_{j}+\displaystyle\sum_{i,j=1}^{N}\xi_{i} v_{ij}z_{j} +\displaystyle\sum_{i,j=1}^{N}z_{i}\tilde{v}_{ij}\xi_{j}+\displaystyle\sum_{j=1}^{N}\displaystyle\sum_{i=1}^{k_{0}}2\lambda_{i}z_{i} v_{ij}z_{j}. 
\end{equation}

Now, the corresponding normalization follows. Then, $\Im V$ is determined by the first normalization condition from (\ref{cnlala1}). We move forward:

 Assume  $r$ is odd. Then, we write as follows
\begin{equation}\begin{split}&F_{2,r-1}(z)=\left(F_{2,r-1}^{\left(1\right)}(z),F_{2,r-1}^{\left(2\right)}(z),\dots,F_{2,r-1}^{\left(N\right)}(z)\right),\quad F_{0,r}(z)=a:=\left(a_{1},a_{2},\dots,a_{N}\right),\\& \tilde{F}_{2,r-1}(z)=\left(\tilde{F}_{2,r-1}^{\left(1\right)}(z),\tilde{F}_{2,r-1}^{\left(2\right)}(z),\dots,\tilde{F}_{2,r-1}^{\left(N\right)}(z)\right),\quad \tilde{F}_{0,r}(z)=\tilde{a}:=\left(\tilde{a}_{1},\tilde{a}_{2},\dots,\tilde{a}_{N}\right).\end{split}\label{788alala}
\end{equation}  

Now, it remains to understand the interactions of these two polynomials  in 
the local defining equation, according to the corresponding Fischer Decompositions by using the standard  inner product $\left<\cdot,\cdot\right>$. We have
\begin{equation}\begin{split}&    \left<z,\tilde{F}_{2,r-1}\left(\xi\right) \right>+\left<a,\xi\right>\left<z,\xi\right>+\left<F_{2,r-1}(z),\xi\right>+\left<z,\tilde{a}\right>\left<z,\xi\right> +\left(\left<z,\tilde{a}\right>+\left<a,\xi\right>\right)\left(\displaystyle\sum_{i=1}^{k_{0}}\lambda_{i}z_{i}^{2}+\displaystyle\sum_{i=1}^{k_{0}}\lambda_{i}\xi_{i}^{2}\right)+\displaystyle\sum_{i=1}^{k_{0}}2\lambda_{i}z_{i}F_{2,r-1}^{(i)}(z)\\&\quad\quad\quad\quad+\displaystyle\sum_{i=1}^{k_{0}}2\lambda_{i}\xi_{i}\tilde{F}_{2,r-1}^{(i)}\left(\xi\right)   +\displaystyle\sum_{i=1}^{k_{0}}2\lambda_{i}z_{i}a_{i}\left(\left<z,\xi\right>+\displaystyle\sum_{i=1}^{k_{0}}\lambda_{i}z_{i}^{2}+\displaystyle\sum_{i=1}^{k_{0}}\lambda_{i}\xi_{i}^{2}\right)  +\displaystyle\sum_{i=1}^{k_{0}}2\lambda_{i}\xi_{i}\tilde{a}_{i}\left(\left<z,\xi\right>+\displaystyle\sum_{i=1}^{k_{0}}\lambda_{i}z_{i}^{2}+\displaystyle\sum_{i=1}^{k_{0}}\lambda_{i}\xi_{i}^{2}\right)=K\left(z,\xi\right),\end{split}
\end{equation} 
where $K\left(z,\xi\right)$ is a  determined polynomial of degree $2$ provided by corresponding Fischer Decompositions  from the local defining equation. 

Then, we write  $K\left(z,\xi\right)$ as a sum of homogeneous terms as follows
$$ K\left(z,\xi\right)=K_{3,0}\left(z,\xi\right)+K_{2,1}\left(z,\xi\right)+K_{1,2}\left(z,\xi\right)+K_{0,3}\left(z,\xi\right).
$$

  We obtain
$$
 F_{2,r-1}^{\left(j\right)}(z)+\left<z,\tilde{a}\right>z_{j}+a_{j}\displaystyle\sum_{i=1}^{k_{0}}\lambda_{i}z_{i}^{2}+\displaystyle\sum_{i=1}^{k_{0}}2\lambda_{i}z_{i}a_{i}  z_{j}+2\lambda_{j}a_{j}\displaystyle\sum_{i=1}^{k_{0}}\lambda_{i}z_{i}^{2}= \frac{\partial }{\partial \xi_{l}}\left(K_{2,1}\left(z,\xi\right)\right),\quad\mbox{for all $l=1,\dots,k_{0}$}.
$$ 

Returning in the first equation in the complexification of (\ref{011}), we have
$$ \displaystyle\sum_{j=1}^{k_{0}}2\lambda_{j}z_{j}F_{2,r-1}^{\left(j\right)}(z)+\left<z,\tilde{a}\right>\displaystyle\sum_{j=1}^{k_{0}}\lambda_{j}z_{j}^{2}+\displaystyle\sum_{i,j=1}^{k_{0}}2\lambda_{i}z_{i}a_{i} \lambda_{j}z_{j}^{2}=K_{3,0}\left(z,\xi\right), $$
or equivalently, we have 
\begin{equation*} \begin{split}&\displaystyle\sum_{j=1}^{k_{0}}2\lambda_{j}z_{j}\left(\frac{\partial }{\partial \xi_{j}}\left(K_{2,1}\left(z,\xi\right)\right)-\left<z,\tilde{a}\right>z_{j}-a_{j}\displaystyle\sum_{i=1}^{k_{0}}\lambda_{i}z_{i}^{2}-\displaystyle\sum_{i=1}^{k_{0}}2\lambda_{i}z_{i}a_{i}  z_{j}-2\lambda_{j}a_{j}\displaystyle\sum_{i=1}^{k_{0}}\lambda_{i}z_{i}^{2}\right)\\&\quad\quad\quad\quad\quad\quad\quad\quad\quad\quad\quad\quad\quad\quad\quad\quad \quad\quad\quad\quad \quad\quad    +\left<z,\tilde{a}\right>\displaystyle\sum_{j=1}^{k_{0}}\lambda_{j}z_{j}^{2}+\displaystyle\sum_{i,j=1}^{k_{0}}2\lambda_{i}z_{i}a_{i} \lambda_{j}z_{j}^{2}=K_{3,0}\left(z,\xi\right). \end{split}
\end{equation*}  
 
It follows that
\begin{equation} \begin{split}&    0=\displaystyle\sum_{j=1}^{k_{0}}2\lambda_{j}z_{j}\frac{\partial }{\partial \xi_{j}}\left(K_{2,1}\left(z,\xi\right)\right)-K_{3,0}\left(z,\xi\right)-\displaystyle\sum_{j=1}^{k_{0}}2\lambda_{j}z_{j} \left( \left<z,\tilde{a}\right>z_{j}+a_{j}\displaystyle\sum_{i=1}^{k_{0}}\lambda_{i}z_{i}^{2}+\displaystyle\sum_{i=1}^{k_{0}}2\lambda_{i}z_{i}a_{i}  z_{j}+2\lambda_{j}a_{j}\displaystyle\sum_{i=1}^{k_{0}}\lambda_{i}z_{i}^{2}\right)\\& \quad\quad\quad\quad\quad\quad\quad\quad\quad\quad\quad\quad\quad\quad\quad\quad\quad\quad\quad\quad\quad\quad\quad\quad\quad\quad\quad\quad\quad\quad\quad\quad\quad   +\left<z,\tilde{a}\right>\displaystyle\sum_{j=1}^{k_{0}}\lambda_{j}z_{j}^{2}+\displaystyle\sum_{i,j=1}^{k_{0}}2\lambda_{i}z_{i}a_{i} \lambda_{j}z_{j}^{2}, \end{split}
\end{equation}
which, after few simplifications, confirms the following
\begin{equation} \begin{split}&    0=\displaystyle\sum_{j=1}^{k_{0}}2\lambda_{j}z_{j}\frac{\partial }{\partial \xi_{j}}\left(K_{2,1}\left(z,\xi\right)\right)-K_{3,0}\left(z,\xi\right)-\left<z,\tilde{a}\right>\displaystyle\sum_{j=1}^{k_{0}}\lambda_{j}z_{j}^{2}  -\displaystyle\sum_{j=1}^{k_{0}}2\lambda_{j}z_{j} \left(  \displaystyle\sum_{i=1}^{k_{0}}2\lambda_{i}z_{i}a_{i}  z_{j}+2\lambda_{j}a_{j}\displaystyle\sum_{i=1}^{k_{0}}\lambda_{i}z_{i}^{2}\right). \end{split}
\end{equation}

Now, we focus only on the coefficients of the above terms $z_{1}^{3},z_{2}^{3},\dots, z_{k_{0}}^{3}$ and $\xi_{1}^{3},\xi_{2}^{3},\dots, \xi_{k_{0}}^{3}$. Then, we compute easily $\tilde{a}_{1},\tilde{a}_{2},\dots, \tilde{a}_{k_{0}}$ and $a_{1},a_{1},\dots, a_{k_{0}}$  according to the assumption (\ref{03}), by computing the real parts of the complex numbers $2\lambda_{1}a_{1}+\tilde{a}_{1},\dots, 2\lambda_{k_{0}}a_{k_{0}}+\tilde{a}_{k_{0}}$ and $2\lambda_{1}\tilde{a}_{1}+a_{1},\dots, 2\lambda_{k_{0}}\tilde{a}_{k_{0}}+a_{k_{0}}$, because 
$$\left|\begin{matrix} 2\lambda_{k} & 1 \\ 1 & 2\lambda_{k}
\end{matrix}\right|\neq 0, \quad \mbox{for all $k=1,\dots, k_{0}$.}$$

 Then, we compute $F_{2,r-1}(z)$ and $\tilde{F}_{2,r-1}\left(\xi\right)$, for all $i=1,\dots, k_{0}$. We obtain 
\begin{equation}\left\{\begin{split}&\displaystyle\sum_{m+2n=r-1}F_{m,n}^{\left(l\right)}(z)w^{n}=F_{r-1,0}^{\left(l\right)}(z)+F_{r-3,1}^{\left(l\right)}(z)w+\dots +a_{k}w^{\frac{r}{2}} ,\quad\mbox{for all $l=1,\dots,k_{0}$,}\\& \displaystyle\sum_{m+2n=r-1}F_{m,n}^{\left(l\right)}(z)w^{n}=F_{r-1,0}^{\left(l\right)}(z)+F_{r-3,1}^{\left(l\right)}(z)w+\dots+a_{k}w^{\frac{r}{2}} ,\quad\mbox{for all $l=k_{0}+1,\dots,N$,} \end{split}\right.\label{londra1la} \end{equation}  
when $r$ is even, and respectively
\begin{equation}\left\{\begin{split}&\displaystyle\sum_{m+2n=r-1}F_{m,n}^{\left(l\right)}(z)w^{n}=F_{r-1,0}^{\left(l\right)}(z)+F_{r-3,1}^{\left(l\right)}(z)w+\dots+F_{1,\frac{r-1}{2}}\left(\xi\right)w^{\frac{r-1}{2}}   ,\quad\mbox{for all $l=1,\dots,k_{0}$,}\\& \displaystyle\sum_{m+2n=r-1}F_{m,n}^{\left(l\right)}(z)w^{n}=F_{r-1,0}^{\left(l\right)}(z)+F_{r-3,1}^{\left(l\right)}(z)w+\dots+F_{1,\frac{r-1}{2}}(z)w^{\frac{r-1}{2}}  ,\quad\mbox{for all $l=k_{0}+1,\dots,N$,} \end{split}\right. \label{londra2la}\end{equation} 
when $r$ is odd. 

Analogously, we have
\begin{equation}\left\{\begin{split}&\displaystyle\sum_{m+2n=r-1}\tilde{F}_{m,n}^{\left(l\right)}\left(\xi\right)\nu^{n}=\tilde{F}_{r-1,0}^{\left(l\right)}\left(\xi\right)+\tilde{F}_{r-3,1}^{\left(l\right)}\left(\xi\right)\nu+\dots +\tilde{a}_{k}\nu^{\frac{r}{2}} ,\quad\mbox{for all $l=1,\dots,k_{0}$,}\\& \displaystyle\sum_{m+2n=r-1}\tilde{F}_{m,n}^{\left(l\right)}\left(\xi\right)\nu^{n}=\tilde{F}_{r-1,0}^{\left(l\right)}\left(\xi\right)+\tilde{F}_{r-3,1}^{\left(l\right)}\left(\xi\right)\nu+\dots+\tilde{a}_{k}\nu^{\frac{r}{2}} ,\quad\mbox{for all $l=k_{0}+1,\dots,N$,} \end{split}\right.\label{londra1lala} \end{equation}  
when $r$ is even, and respectively
\begin{equation}\left\{\begin{split}&\displaystyle\sum_{m+2n=r-1}\tilde{F}_{m,n}^{\left(l\right)}\left(\xi\right)\nu^{n}=\tilde{F}_{r-1,0}^{\left(l\right)}(z)+\tilde{F}_{r-3,1}^{\left(l\right)}\left(\xi\right)\nu+\dots+\tilde{F}_{1,\frac{r-1}{2}}\left(\xi\right)\nu^{\frac{r-1}{2}}   ,\quad\mbox{for all $l=1,\dots,k_{0}$,}\\& \displaystyle\sum_{m+2n=r-1}\tilde{F}_{m,n}^{\left(l\right)}\left(\xi\right)\nu^{n}=\tilde{F}_{r-1,0}^{\left(l\right)}\left(\xi\right)\nu+\tilde{F}_{r-3,1}^{\left(l\right)}\left(\xi\right)\nu+\dots+\tilde{F}_{1,\frac{r-1}{2}}\left(\xi\right)\nu^{\frac{r-1}{2}}  ,\quad\mbox{for all $l=k_{0}+1,\dots,N$,} \end{split}\right. \label{londra2lala}\end{equation} 
when $r$ is odd.
\end{proof}
  
\section{Proofs of Theorems \ref{tA1sec} and \ref{tA2sec}}
We proceed in the light of  (\ref{var1A}),(\ref{var1B}),(\ref{var1Ase}),(\ref{var1Bse}),(\ref{var1AX}),(\ref{var1BX}). Let's   compute: 
\subsection{Transforming Equations and Models}We consider the following non-constant formal mapping
\begin{equation*}\begin{split}&w'=G(w,z),\quad z_{1}'=F_{1}(w,z),z_{2}'=F_{2}(w,z),\dots,z_{N}'=F_{N}(w,z), \\& \quad\quad \quad\quad\quad\quad\quad\hspace{0.1 cm} z_{N+1}'=F_{N+1}(w,z),z_{N+2}'=F_{N+2}(w,z),\dots,z_{N'}'=F_{N'}(w,z),\\&\nu'=\tilde{G}\left(\nu,\xi\right),\quad z_{1}'=\tilde{F}_{1}\left(\nu,\xi\right),z_{2}'=\tilde{F}_{2}\left(\nu,\xi\right),\dots,z_{N}'=\tilde{F}_{N}\left(\nu,\xi\right), \\& \quad\quad \quad\quad\quad\quad\quad\hspace{0.1 cm} z_{N+1}'=\tilde{F}_{N+1}\left(\nu,\xi\right),z_{N+2}'=\tilde{F}_{N+2}\left(\nu,\xi\right),\dots,z_{N'}'=\tilde{F}_{N'}\left(\nu,\xi\right),\end{split}
\end{equation*}
which sends the complexification of (\ref{var1A}) into the complexification of (\ref{var1B})  according to the following formal expansions
\begin{equation} G(w,z)=\displaystyle\sum_{m,n\geq 0}G_{m,n}(z)w^{n},\quad \left(F^{\left(1\right)}(w,z),\dots,F^{\left(N'\right)}(w,z)\right)=\left(\displaystyle\sum_{m,n\geq 0}F_{m,n}^{\left(1\right)}(z)w^{n},\dots,\displaystyle\sum_{m,n\geq 0}F_{m,n}^{\left(N'\right)}(z)w^{n}\right),\label{440la1} 
\end{equation}
where $G_{m,n}(z)$, $F_{m,n}^{\left(1\right)}(z),F_{m,n}^{\left(2\right)}(z),\dots,F_{m,n}^{\left(N\right)}(z),F_{m,n}^{\left(N+1\right)}(z),\dots,F_{m,n}^{\left(N'\right)}(z)$ are homogeneous polynomials of degree $m$ in $z$, and respectively  
\begin{equation} \tilde{G}\left(\nu,\xi\right)=\displaystyle\sum_{m,n\geq 0}\tilde{G}_{m,n}\left(\xi\right)\nu^{n},\quad \left(\tilde{F}^{\left(1\right)}\left(\nu,\xi\right),\dots,\tilde{F}^{\left(N'\right)}\left(\nu,\xi\right)\right)=\left(\displaystyle\sum_{m,n\geq 0}\tilde{F}_{m,n}^{\left(1\right)}\left(\xi\right)\nu^{n},\dots,\displaystyle\sum_{m,n\geq 0}\tilde{F}_{m,n}^{\left(N'\right)}\left(\xi\right)\nu^{n}\right),\label{440la2} 
\end{equation}
where $\tilde{G}_{m,n}\left(\xi\right)$, $\tilde{F}_{m,n}^{\left(1\right)}\left(\xi\right),\tilde{F}_{m,n}^{\left(2\right)}\left(\xi\right),\dots,\tilde{F}_{m,n}^{\left(N\right)}\left(\xi\right),\tilde{F}_{m,n}^{\left(N+1\right)}\left(\xi\right),\dots,\tilde{F}_{m,n}^{\left(N'\right)}\left(\xi\right)$ are homogeneous polynomials of degree $m$ in $\xi$. We obtain  
\begin{equation}
\displaystyle\sum_{m,n\geq 0}G_{m,n}(z)w^{n}= \displaystyle\sum_{k=1}^{N'} \left(\displaystyle\sum_{m,n\geq 0}\left(F_{m,n}^{\left(k\right)}(z)w^{n}\right) \left( \displaystyle\sum_{m,n\geq 0}\tilde{F}_{m,n}^{\left(k\right)}\left(\xi\right)\nu^{n}\right) +\lambda_{k}\left(\displaystyle\sum_{m,n\geq 0}F_{m,n}^{\left(k\right)}(z)w^{n}\right)^{2}+\lambda_{k}  \left(\displaystyle\sum_{m,n\geq 0}\tilde{F}_{m,n}^{\left(k\right)}\left(\xi\right)\nu^{n}\right)^{2} \right) , \label{ecuatieXYYla1}
\end{equation}
and respectively, we obtain 
\begin{equation}
\displaystyle\sum_{m,n\geq 0}\tilde{G}_{m,n}\left(\xi\right)\nu^{n}=\displaystyle\sum_{m,n\geq 0}G_{m,n}(z)w^{n}. \label{ecuatieXYYla2}
\end{equation}

We also observe the following important aspect
\begin{equation}\mbox{rank}\begin{pmatrix}\frac{\partial F_{1}}{\partial z_{1}}(0) &\dots & \frac{\partial F_{N}}{\partial z_{1}}(0)& \frac{\partial F_{N+1}}{\partial z_{1}}(0)&\dots&\frac{\partial F_{N'}}{\partial z_{1}}(0) \\ \vdots &\ddots & \vdots & \vdots &\ddots&\vdots \\ \frac{\partial F_{1}}{\partial z_{N}}(0) &\dots & \frac{\partial F_{N}}{\partial z_{N}}(0)& \frac{\partial F_{N+1}}{\partial z_{N}}(0)&\dots&\frac{\partial F_{N'}}{\partial z_{N}}(0) \\
\end{pmatrix}=N,
\end{equation}
and respectively, we also observe that
\begin{equation}\mbox{rank}\begin{pmatrix}\frac{\partial \tilde{F}_{1}}{\partial \xi_{1}}(0) &\dots & \frac{\partial \tilde{F}_{N}}{\partial \xi_{1}}(0)& \frac{\partial \tilde{F}_{N+1}}{\partial \xi_{1}}(0)&\dots&\frac{\partial \tilde{F}_{N'}}{\partial \xi_{1}}(0) \\ \vdots &\ddots & \vdots & \vdots &\ddots&\vdots \\ \frac{\partial \tilde{F}_{1}}{\partial \xi_{N}}(0) &\dots & \frac{\partial \tilde{F}_{N}}{\partial z_{N}}(0)& \frac{\partial \tilde{F}_{N+1}}{\partial \xi_{N}}(0)&\dots&\frac{\partial \tilde{F}_{N'}}{\partial \xi_{N}}(0) \\
\end{pmatrix}=N,
\end{equation}
because otherwise we would obtain a contraction with (\ref{transverse}). 

Indeed, (\ref{transverse}) holds, because contrary we would have 
\begin{equation*}G(z,w)=\mbox{O}(2),\quad
\end{equation*}
which implies that $G(z,w)=0$ and $\tilde{G}\left(\xi,\nu\right)=0$, in the light of the Fischer Decompositions defining (\ref{spartiuGlala1}), 
and also that  $F(z,w)=0$ and $\tilde{F}\left(\xi,\nu\right)=0$, because of (\ref{ecuatieXYYla2}), obtaining again a contradiction because the formal embedding was assumed to be not constant.

Imposing the normalization procedure of Baouendi-Huang\cite{BH} similarly as in \cite{bu5}, we obtain 
\begin{equation}\begin{split}&G(z,w)=w+\mbox{O}(2),\quad F^{\left(1\right)}(w,z)=z_{1}+\mbox{O}(2),\dots, F^{\left(N\right)}(w,z)=z_{N}+\mbox{O}(2),\\&\quad\quad\quad\quad\quad\quad\quad\quad\quad\quad F^{\left(N+1\right)}(w,z)=\mbox{O}(2),\dots, F^{\left(N'\right)}(w,z)= \mbox{O}(2),\\& \tilde{G}\left(\xi,\nu\right)=\nu+\mbox{O}(2),\quad F^{\left(1\right)}\left(\xi,\nu\right)=\xi_{1}+\mbox{O}(2),\dots, F^{\left(N\right)}\left(\xi,\nu\right)=\xi_{N}+\mbox{O}(2),\\&\quad\quad\quad\quad\quad\quad\quad\quad\quad\quad F^{\left(N+1\right)}\left(\xi,\nu\right)=\mbox{O}(2),\dots, F^{\left(N'\right)}\left(\xi,\nu\right)= \mbox{O}(2).\end{split}
 \end{equation} 

It follows that 
\begin{equation}G(z,w)=w,\quad \tilde{G}\left(\xi,\nu\right)=\nu,
\end{equation}
in the light of the Fischer Decompositions defining (\ref{spartiuGlala1}).

Now, we use the writing (\ref{882}) according to (\ref{lalala1}) and  (\ref{lalala}). Moreover, we write as follows\begin{equation} \left(F^{\left(1\right)},F^{\left(2\right)},\dots,F^{\left(N\right)}\right)\left(\xi,\nu\right)=\left(\xi_{1},\xi_{2},\dots,\xi_{N}\right)+A\left(\xi_{1},\xi_{2},\dots,\xi_{N}\right)\nu+B_{0}\left(\xi_{1},\xi_{2},\dots,\xi_{N}\right)+B_{1}\left(\xi_{1},\xi_{2},\dots,\xi_{N}\right)+\mbox{O}(4),\label{882la}
\end{equation}
where we deal with:
\begin{itemize}
\item $B_{0}\left(\xi_{1},\xi_{2},\dots,\xi_{N}\right)$ is a vector   polynomial of degree $2$ in  $\left(\xi_{1},\xi_{2},\dots,\xi_{N}\right)$,
\item $B_{1}\left(\xi_{1},\xi_{2},\dots,\xi_{N}\right)$ is a vector   polynomial of degree $3$ in  $\left(\xi_{1},\xi_{2},\dots,\xi_{N}\right)$,
\item   $A\left(\xi_{1},\xi_{2},\dots,\xi_{N}\right)$ is a linear form in $\left(\xi_{1},\xi_{2},\dots,\xi_{N}\right)$,   written   as follows
 \begin{equation} \tilde{A}\left(\xi_{1},\xi_{1},\dots,\xi_{N}\right)=\begin{pmatrix}
\tilde{a}_{11} & \tilde{a}_{12} &\dots& \tilde{a}_{1N} \\ \tilde{a}_{21} & \tilde{a}_{22} &\dots& \tilde{a}_{2N} \\  \vdots & \vdots &\ddots& \vdots \\ \tilde{a}_{N1} & \tilde{a}_{N2} &\dots& \tilde{a}_{NN} \\ 
\end{pmatrix}\begin{pmatrix}
\xi_{1} \\ \xi_{2}\\ \vdots \\ \xi_{N}
\end{pmatrix},\label{lalalap}
\end{equation}\end{itemize}
having in mind (\ref{bev}), and writing respectively as  s
\begin{equation}\begin{split}&\tilde{B}_{0}\left(\xi_{1},\xi_{2},\dots,\xi_{N}\right)=\left(\tilde{B}_{0}^{(1)},\tilde{B}_{0}^{(2)},\dots,\tilde{B}_{0}^{\left(N\right)}\right)\left(\xi_{1},\xi_{1},\dots,\xi_{N}\right),\\& \tilde{B}_{1}\left(\xi_{1},\xi_{2},\dots,\xi_{N}\right)=\left(\tilde{B}_{1}^{(2)},\tilde{B}_{2}^{(2)},\dots,\tilde{B}_{1}^{\left(N\right)}\right)\left(\xi_{1},\xi_{2},\dots,\xi_{N}\right).\end{split}\label{lalala1lap}
\end{equation}

Focusing on  (\ref{ecuatieXYYla1}), we obtain 
\begin{equation}\displaystyle\sum_{l=1}^{N}B_{0}^{\left(l\right)} \xi_{l}+\displaystyle\sum_{l=1}^{N}\tilde{B}_{0}^{\left(l\right)} z_{l}=0,
\end{equation}
because there are no homogeneous terms of degree $3$ in the right hand side of (\ref{ecuatieXYYla1}), concluding that
$$ B_{0}\left(z_{1},z_{2},\dots,z_{N}\right)=0,\quad \tilde{B}_{0}\left(\xi_{1},\xi_{2}\dots,\xi_{N}\right)=0.
$$

Moreover, we have
\begin{equation}\begin{split}&   \left(\displaystyle\sum_{k,l=1}^{N}\tilde{a}_{kl}\xi_{l}z_{k} +\displaystyle\sum_{l =1}^{N} 2\lambda_{l}\tilde{a}_{lk}\xi_{l}\xi_{k}     +\displaystyle\sum_{k,l=1}^{N}a_{kl}z_{l}\xi_{k} +\displaystyle\sum_{l =1}^{N} 2\lambda_{l}a_{lk}z_{l}z_{k}\right)\left(z_{1}\xi_{1}+\dots+z_{N}\xi_{N}+\lambda_{1}\left(z_{1}^{2}+\xi_{1}^{2}\right)+\dots+\lambda_{N}\left(z_{N}^{2}+\xi_{N}^{2}\right)\right)\\& +\\&\left(\displaystyle\sum_{l=1}^{N}B^{\left(l\right)}_{1} \xi_{l}+\sum_{l=1}^{N}B^{\left(l\right)}_{1} 2\lambda_{l} z_{l}\right)+   \left(\displaystyle\sum_{l=1}^{N}\tilde{B}^{\left(l\right)}_{1} z_{l}+\sum_{l=1}^{N}\tilde{B}^{\left(l\right)}_{1} 2\lambda_{l} \xi_{l}\right)\\&=\\&-\left(\displaystyle\sum_{l=N+1}^{N'}F_{2,0}^{\left(l\right)}(z)\right)\left(\displaystyle\sum_{l=N+1}^{N'}\tilde{F}_{2,0}^{\left(l\right)}\left(\xi\right)\right)                
-\displaystyle\sum_{l=N+1}^{N'}\lambda'_{l}\left( F_{2,0}^{\left(l\right)}(z)\right)^{2}
 -\displaystyle\sum_{l=N+1}^{N'}\lambda'_{l}  \left(\tilde{F}_{2,0}^{\left(l\right)}\left(\xi\right) \right)^{2} .\end{split}\label{ectlala}
\end{equation}

It follows that
  $$F_{2,0}^{\left(l\right)}(z)=0,\quad \tilde{F}_{2,0}^{\left(l\right)}\left(\xi\right)=0,\hspace{0.1 cm}\mbox{for all $l=N+1,\dots, N'$.}
$$

According to the Fischer Decompositions defined by (\ref{spartiuFlala}), we obtain
$$ A\left(z_{1},z_{2},\dots,z_{N}\right)=0,\quad B_{1}\left(z_{1},z_{2},\dots,z_{N}\right)=0,\quad \tilde{A}\left(\xi_{1},\xi_{2},\dots,\xi_{N}\right)=0,\quad \tilde{B}_{1}\left(\xi_{1},\xi_{2},\dots,\xi_{N}\right)=0.$$

We move forward:

\subsection{Proof of Theorem \ref{tA1sec}}  It is continued the previous analysis in 
(\ref{ecuatieXYYla1}) and (\ref{ecuatieXYYla2})  according to  (\ref{801}) and the following expansion
 \begin{equation} \left(\tilde{F}^{\left(1\right)},\tilde{F}^{\left(2\right)},\dots,\tilde{F}^{\left(N\right)}\right)\left(\xi,\nu\right)=\left(\xi_{1},\xi_{2},\dots,\xi_{N}\right)+\tilde{A}_{0}\left(\xi_{1},\xi_{2},\dots,\xi_{N}\right)\nu+\tilde{B}_{0}\left(\xi_{1},\xi_{2},\dots,\xi_{N}\right) +v\nu^{2}+\mbox{O}(5),
\end{equation} 
where we deal with
\begin{itemize}
\item $\tilde{B}_{0}\left(\xi_{1},\xi_{2},\dots,\xi_{N}\right)$ is a vector homogeneous polynomial of degree $4$ in  $\left(\xi_{1},\xi_{2},\dots,\xi_{N}\right)$,
\item  $\tilde{A}_{0}\left(\xi_{1},\xi_{2},\dots,\xi_{N}\right)$ is a vector homogeneous polynomial of degree $2$ in  $\left(\xi_{1},\xi_{2},\dots,\xi_{N}\right)$,
\item $\tilde{v}\in\mathbb{C}^{N}$, 
\end{itemize}
in respect to the following notations
$$\tilde{B}_{0}\left(\xi_{1},\xi_{2},\dots,\xi_{N}\right)=\left(\tilde{B}_{0}^{(1)},\tilde{B}_{0}^{(2)},\dots,\tilde{B}_{0}^{\left(N\right)}\right)\left(\xi_{1},\xi_{2},\dots,\xi_{N}\right),$$  $$ \tilde{A}_{0}\left(\xi_{1},\xi_{2},\dots,\xi_{N}\right)=\left(\tilde{A}_{0}^{(1)},\tilde{A}_{0}^{(2)},\dots,\tilde{A}_{0}^{\left(N\right)}\right)\left(\xi_{1},\xi_{2},\dots,\xi_{N}\right),$$ $$ \tilde{v}=\left(\tilde{v}_{1},\tilde{v}_{2},\dots,\tilde{v}_{N}\right)\in\mathbb{C}^{N}.$$

Extracting  terms of degree $5$ in $\left(z,\xi\right)$ from (\ref{ecuatieXYY}), we obtain
  \begin{equation}\begin{split}&  \left(\displaystyle\sum_{l=1}^{N} z_{l}\tilde{A}_{0}^{\left(l\right)} +  \displaystyle\sum_{l=1}^{N} 2\lambda_{l}\xi_{l}\tilde{A}_{0}^{\left(l\right)}     + \displaystyle\sum_{l=1}^{N} \xi_{l} A_{0}^{\left(l\right)} +  \displaystyle\sum_{l=1}^{N} 2\lambda_{l}z_{l}A_{0}^{\left(l\right)}\right) \left(z_{1}\xi_{1}+\dots+z_{N}\xi_{N}+\lambda_{1}\left(z_{1}^{2}+\xi_{1}^{2}\right)+\dots+\lambda_{N}\left(z_{N}^{2}+\xi_{N}^{2}\right)  \right)\\&\quad\quad\quad\quad\quad\quad \quad\quad\quad\quad\quad\quad\quad \quad\quad\quad\quad\quad\quad\quad \quad\quad\quad\quad\quad\quad\quad \quad\quad\quad\quad\quad\quad\quad \quad\quad\quad        +\\&    \left(\displaystyle\sum_{l=1}^{N} z_{l}\tilde{v}_{l}  +\displaystyle\sum_{l=1}^{N}2\lambda_{l}\xi_{l} \tilde{v}_{l} +        \displaystyle\sum_{l=1}^{N} \xi_{l} v_{l}  +\displaystyle\sum_{l=1}^{N}2\lambda_{l}z_{l} v_{l}\right)  \left(z_{1}\xi_{1}+\dots+z_{N}\xi_{N}+\lambda_{1}\left(z_{1}^{2}+\xi_{1}^{2}\right)+\dots+\lambda_{N}\left(z_{N}^{2}+\xi_{N}^{2}\right) \right)^{2}
\\&\quad\quad\quad\quad\quad\quad \quad\quad\quad\quad\quad\quad\quad \quad\quad\quad\quad\quad\quad\quad \quad\quad\quad\quad\quad\quad\quad \quad\quad\quad\quad\quad\quad\quad \quad\quad\quad  +\\& \quad\quad\quad\quad\quad\quad\quad\quad\quad\quad\quad\quad\quad\quad\quad\quad\quad\quad\quad\quad\quad\quad\quad\quad   \left( \displaystyle\sum_{l=1}^{N} \xi_{l}B_{0}^{\left(l\right)} +\displaystyle\sum_{l=1}^{N} 2\lambda_{l}z_{l}B_{0}^{\left(l\right)}\right)+\left( \displaystyle\sum_{l=1}^{N} z_{l}\tilde{B}_{0}^{\left(l\right)} +\displaystyle\sum_{l=1}^{N} 2\lambda_{l}\xi_{l}\tilde{B}_{0}^{\left(l\right)}\right)  =0.
\end{split}\end{equation}

According to the Fischer Decompositions defined by (\ref{spartiuFlala}), it follows that
$$B_{0}\left(z_{1},\dots,z_{N}\right)=0,\hspace{0.1 cm} \tilde{A}_{0}\left(z_{1},\dots,z_{N}\right)=0,\hspace{0.1 cm}  v=0,\quad \tilde{B}_{0}\left(\xi_{1},\xi_{2},\dots,\xi_{N}\right)=0,\hspace{0.1 cm} \tilde{A}_{0}\left(\xi_{1},\xi_{2},\dots,\xi_{N}\right)=0,\hspace{0.1 cm}  v=0.$$ 

We move forward to the next degree in  (\ref{ecuatieXYY}) according to the writing (\ref{901}). We have
\begin{equation}\left(\tilde{F}^{\left(1\right)},\tilde{F}^{\left(2\right)},\dots,
\tilde{F}^{\left(N\right)}\right)\left(\xi,\nu\right)=\left(\xi_{1},\xi_{2},\dots,\xi_{N}\right)+\tilde{A}\left(\xi_{1},\xi_{2},\dots,\xi_{N}\right) \nu^{2}+B_{1}\left(\xi_{1},\xi_{2},\dots,\xi_{N}\right)\nu+\tilde{B}_{0}\left(\xi_{1},\xi_{2},\dots,\xi_{N}\right)+ \mbox{O}(5),
\end{equation}
where we deal with
\begin{itemize}
\item $\tilde{B}_{1}\left(\xi_{1},\xi_{2},\dots,\xi_{N}\right)$ is a vector homogeneous polynomial of degree $5$ in  $\left(\xi_{1},\xi_{2},\dots,\xi_{N}\right)$,
\item  $\tilde{B}_{2}\left(\xi_{1},\xi_{2},\dots,\xi_{N}\right)$ is a vector homogeneous polynomial of degree $3$ in  $\left(\xi_{1},\xi_{2},\dots,\xi_{N}\right)$,
\item $\tilde{A}_{0}\left(\xi_{1},\xi_{2},\dots,\xi_{N}\right)$ is a linear form in  $\left(\xi_{1},\xi_{2},\dots,\xi_{N}\right)$  defined as in (\ref{lalala}),
\end{itemize}
according to the following notations
\begin{equation*}\begin{split}&\tilde{B}_{0}\left(\xi_{1},\xi_{2},\dots,\xi_{N}\right)=\left(\tilde{B}_{0}^{(1)},\tilde{B}_{0}^{(2)},\dots,\tilde{B}_{0}^{\left(N\right)}\right)\left(\xi_{1},\xi_{2},\dots,\xi_{N}\right),\\& \tilde{B}_{1}\left(\xi_{1},\xi_{2},\dots,\xi_{N}\right)=\left(\tilde{B}_{1}^{(1)},\tilde{B}_{1}^{(2)},\dots,\tilde{B}_{1}^{\left(N\right)}\right)\left(\xi_{1},\xi_{2},\dots,\xi_{N}\right) .\end{split}
\end{equation*}

We obtain
 \begin{equation}\begin{split}& \left(\displaystyle\sum_{k,l=1}^{N}\tilde{a}_{kl}\xi_{l}z_{k} +\displaystyle\sum_{k =1}^{N} 2\lambda_{k}\tilde{a}_{lk}\xi_{l}\xi_{k}         +\displaystyle\sum_{k,l=1}^{N}\tilde{a}_{kl}z_{l}\xi_{k} +\displaystyle\sum_{k =1}^{N} 2\lambda_{k}\tilde{a}_{lk}z_{l}z_{k}\right) \left(z_{1}\xi_{1}+\dots+z_{N}\xi_{N}+\lambda_{1}\left(z_{1}^{2}+\xi_{1}^{2}\right)+\dots+\lambda_{N}\left(z_{N}^{2}+\xi_{N}^{2}\right) \right)^{2}\\& \quad\quad\quad\quad\quad\quad\quad\quad\quad\quad\quad\quad\quad\quad\quad\quad+\\&  \quad\quad\hspace{0.1 cm} \left(\displaystyle\sum_{l=1}^{N}\tilde{B}_{0}^{\left(l\right)}z_{k} +  \displaystyle\sum_{l=1}^{N}2\lambda_{l} \xi_{l}\tilde{B}_{0}^{\left(l\right)}+      \displaystyle\sum_{l=1}^{N}B_{0}^{\left(l\right)}\xi_{k} +  \displaystyle\sum_{l=1}^{N}2\lambda_{l} z_{l}B_{0}^{\left(l\right)} \right) \left(z_{1}\xi_{1}+\dots+z_{N}\xi_{N}+\lambda_{1}\left(z_{1}^{2}+\xi_{1}^{2}\right)+\dots+\lambda_{N}\left(z_{N}^{2}+\xi_{N}^{2}\right)\right)\\&\quad\quad\quad\quad\quad\quad\quad\quad\quad\quad\quad\quad\quad\quad\quad\quad + \\& \quad\quad\quad\quad\quad\quad\quad\quad \left(\displaystyle\sum_{l=1}^{N}\tilde{B}_{1}^{\left(l\right)}z_{l} +  \displaystyle\sum_{l=1}^{N}2\lambda_{l} \xi_{l}\tilde{B}_{1}^{\left(l\right)}+        \displaystyle\sum_{l=1}^{N}B_{1}^{\left(l\right)}\xi_{l} +  \displaystyle\sum_{l=1}^{N}2\lambda_{l} z_{l}B_{1}^{\left(l\right)} \right)+\left(\displaystyle\sum_{l=N+1}^{N'}F_{3,0}^{\left(l\right)}(z)\right) \left(\displaystyle\sum_{l=N+1}^{N'}\tilde{F}_{3,0}^{\left(l\right)}\left(\xi\right)\right)         
\\&\quad\quad\quad\quad\quad\quad\quad\quad\quad\quad\quad\quad\quad\quad\quad\quad+\\&\quad\quad\quad\quad\quad\quad\quad\quad\quad\quad\quad\quad\displaystyle\sum_{l=N+1}^{N'}\lambda'_{l}\left( F_{3,0}^{\left(l\right)}(z)\right)^{2}+\displaystyle\sum_{l=N+1}^{N'}\lambda'_{l}  \left(  \tilde{F}_{3,0}^{\left(l\right)} \left(\xi\right)\right)^{2}=0.\end{split}\label{ect1lala}
\end{equation}

It follows that
 $$F_{3,0}^{\left(l\right)}(z)=0,\quad \tilde{F}_{3,0}^{\left(l\right)}\left(\xi\right)=0,\hspace{0.1 cm}\mbox{for all $l=N+1,\dots, N'$.}
$$

According to the Fischer Decompositions defined by (\ref{spartiuFlala}), it follows that
\begin{equation*}
\begin{split}& B_{0}\left(z_{1},\dots,z_{N}\right)=0,\quad B_{1}\left(z_{1},\dots,z_{N}\right)=0,\\& \tilde{B}_{0}\left(\xi_{1},\xi_{2},\dots,\xi_{N}\right)=0,\quad\tilde{B}_{1}\left(\xi_{1},\xi_{2},\dots,\xi_{N}\right)=0.\end{split}
\end{equation*}

It becomes clear now the induction process considered  
in (\ref{ecuatieXYYla1}) and (\ref{ecuatieXYYla2}) providing obviously the equivalence class from (\ref{bebe1se}).  
\subsection{Proof of Theorem \ref{tA2sec}}Let $\mathcal{M}\subset\mathbb{C}^{2N+2}$ be a real-formal
submanifold defined near  $p=0$ as in (\ref{ecuatie}), and respectively, let 
 another real-formal submanifold $M'\subset\mathbb{C}^{N+1}$ defined near $p=0$ as in (\ref{ecuatie1}) such that (\ref{diag1}) holds. It is important to observe that any formal holomorphic Segre preserving embedding of $\mathcal{M}$ into the corresponding model, can be written as follows
\begin{equation}
(w,z;\nu,\xi)\rightarrow \left(w+\mbox{O}(2),z+\mbox{O}(2),\mbox{O}(2);\nu+\mbox{O}(2),\xi+\mbox{O}(2),\mbox{O}(2)\right).
\end{equation}

More precisely, in the light of the previous computations, we have
\begin{equation}
(w,z;\nu,\xi)\rightarrow \left(w+\mbox{O}(2),z+\mbox{O}(2),0;\nu+\mbox{O}(2),\xi+\mbox{O}(2),0\right).
\end{equation}

We take an embedding of similar type for $\mathcal{M}'$. It follows that
\begin{equation}\left.\varphi_{k}\left(z,\xi\right)=\varphi' _{k}\left(z',\overline{\xi'}\right)\right\vert_{z'=(z,0)\atop{ x'=\left(\xi,0\right)}},\quad\mbox{for all
 $k\geq 3$. } \label{yuklala}\end{equation}

\section{Appendix }
 
 \end{document}